\documentclass[11pt]{article}
\newcounter{ourcount}
\usepackage{amsmath,amsfonts,amssymb,graphicx,epsfig,pdflscape,multirow,ntheorem,gensymb,enumerate}
\usepackage[matrix,arrow,curve]{xy} 
\numberwithin{equation}{section}
\graphicspath{{figpdf/}}
\usepackage[nosort]{cite}
\usepackage[usenames]{color}
\usepackage{pstricks}
\usepackage{pst-node,pst-blur,pstricks-add} 
\usepackage[backref=false]{hyperref}
\hypersetup{
colorlinks=true,
citecolor=red,
linkcolor=darkblue,
urlcolor=darkblue,
}

\usepackage{thmtools} 
\usepackage[capitalise,noabbrev]{cleveref}

\usepackage{array}
\usepackage{ragged2e}

\definecolor{darkblue}{rgb}{0,0,.8}
\definecolor{red}{rgb}{1,0,0}

\theorembodyfont{\itshape} 
\theoremheaderfont{\scshape}
\theoremstyle{plain}  

\newtheorem{Environment}{Environment}[section]

\newtheorem{Lemma}[Environment]{Lemma}
\newtheorem{Corollary}[Environment]{Corollary}
\newtheorem{Proposition}[Environment]{Proposition}
\newtheorem{Definition}[Environment]{Definition}
\newtheorem{Remark}[Environment]{Remark}
\numberwithin{equation}{section}

\newcommand{\nc}{\newcommand}
\def\arxiv#1#2{\href{http://arxiv.org/abs/#1}{\textsf{arXiv:#1\,#2}}}

\crefname{Lemma}{Lemma}{Lemmas}
\crefname{Proposition}{Proposition}{Propositions}
\crefname{Conjecture}{Conjecture}{Conjectures}
\nc{\iI}{\mathrm{i}}
\nc{\dd}{\mathrm{d}}   
\nc{\eE}{\mathsf{e}}
\renewcommand{\ge}{\geqslant}
\renewcommand{\le}{\leqslant}
\renewcommand{\geq}{\geqslant}

\nc{\pa}{\mathsf{a}}
\nc{\pb}{\mathsf{b}}
\nc{\pc}{\mathsf{c}}
\nc{\pab}{\mathsf{\bar a}}
\nc{\pbb}{\mathsf{\bar b}}

\nc{\proof}{{\scshape Proof.\ }} 
\nc{\eproof}{{\hfill \rule{0.5em}{0.5em}\medskip}}
\nc{\bib}{\bibitem}
\nc{\be}{\begin{equation}}
\nc{\ee}{\end{equation}}
\nc{\chit}{\raisebox{0.25ex}{$\chi$}}
\nc{\llangle}{\langle\!\langle}
\nc{\tl}{\mathsf{TL}}
\nc{\ptl}{\mathsf{pTL}}
\nc{\atl}{\mathsf{aTL}}
\nc{\qptl}{\mathsf{upTL}}
\nc{\qptla}{\mathsf{upTL}^{\tinyx{1}}}
\nc{\qptlb}{\mathsf{upTL}^{\tinyx{2}}}
\nc{\qatl}{\mathsf{uaTL}}
\nc{\qatla}{\mathsf{uaTL}^{\tinyx{1}}}
\nc{\qatlb}{\mathsf{uaTL}^{\tinyx{2}}}
\nc{\repA}{\mathsf{A}}
\nc{\repB}{\mathsf{B}}
\nc{\repE}{\mathsf{E}}
\nc{\repM}{\mathsf{M}}
\nc{\repS}{\mathsf{S}}
\nc{\repV}{\mathsf{V}}
\nc{\repW}{\mathsf{W}}
\nc{\repI}{\mathsf{I}}
\nc{\repJ}{\mathsf{J}}
\nc{\repX}{\mathsf{X}}
\nc{\repN}{\mathsf{N}}
\nc{\setB}{\mathcal{B}}
\nc{\setT}{\mathcal{T}}
\nc{\setE}{\mathcal{E}}
\nc{\setS}{\mathcal{S}}
\nc{\Fb}{\mbox{\boldmath $F$}}
\nc{\Fbb}{\mbox{\boldmath $\bar F$}}
\nc{\Gb}{\mbox{\boldmath $G$}}
\nc{\Gbb}{\mbox{\boldmath $\bar G$}}
\nc{\Tb}{\mbox{\boldmath $T$}}
\nc{\Hb}{\mbox{\boldmath $H$}}
\nc{\Hbb}{\mbox{\boldmath $\bar H$}}
\nc{\id}{\mathbf{1}}
\nc{\wh}{\widehat}
\nc{\tinyL}{\textrm{\tiny$(\ell)$}}
\nc{\tinyx}[1]{\textrm{\tiny$(#1)$}}

\definecolor{lightblue}{rgb}{.7,.7,1}
\definecolor{lightestblue}{rgb}{.95,.95,1}
\definecolor{lightlightblue}{rgb}{.9,.9,1}
\definecolor{midblue}{rgb}{.7,.7,1}
\definecolor{purple}{rgb}{0.5,0,0.5}
\definecolor{peach}{rgb}{1, 0.854902, 0.72549}
\definecolor{weirdblue}{rgb}{0.567, 0.767, 0.867}
\definecolor{aquamarine}{rgb}{0.498039, 1., 0.83137}
\newrgbcolor{darkgreen}{0., 0.733333, 0.0621042}
\newrgbcolor{mygreen}{0., 0.87, 0.0621042}
\newrgbcolor{myred}{0.9 0.085 0.031}

\def\facegrid#1#2{
\psframe[fillstyle=solid,fillcolor=lightlightblue,linewidth=0pt]#1#2
\psgrid[gridlabels=0pt,subgriddiv=1]#1#2}

\def\loopa{
\psframe[linewidth=.25pt](0,0)(1,1)
\psarc[linewidth=1.5pt,linecolor=blue](1,0){.5}{90}{180}
\psarc[linewidth=1.5pt,linecolor=blue](0,1){.5}{-90}{0}
}

\def\loopa{
\psframe[linewidth=.25pt](0,0)(1,1)
\psarc[linewidth=1.25pt,linecolor=blue](1,0){.5}{90}{180}
\psarc[linewidth=1.25pt,linecolor=blue](0,1){.5}{-90}{0}
}

\def\loopb{
\psframe[linewidth=.25pt](0,0)(1,1)
\psarc[linewidth=1.25pt,linecolor=blue](0,0){.5}{0}{90}
\psarc[linewidth=1.25pt,linecolor=blue](1,1){.5}{180}{270}
}

\def\loopab{
\psline[linewidth=1pt,linestyle=dashed,linecolor=lightgray,dash=2pt 1.5pt](1,0)(0,1)
}
\def\loopar{
\psline[linewidth=2pt,linecolor=myred](0,0)(1,1)
}

\def\loopbb{
\psline[linewidth=1pt,linestyle=dashed,linecolor=lightgray,dash=2pt 1.5pt](0,0)(1,1)
}
\def\loopbr{
\psline[linewidth=2pt,linecolor=myred](1,0)(0,1)
}

\def\Isinga{
\pspolygon[fillstyle=solid,fillcolor=myred,linewidth=0cm,linecolor=myred](-0.2,-0.05)(0.2,-0.05)(0.2,0.05)(-0.2,0.05)
\pspolygon[fillstyle=solid,fillcolor=myred,linewidth=0cm,linecolor=myred](-0.05,-0.2)(-0.05,0.2)(0.05,0.2)(0.05,-0.2)
}
\def\Isingb{
\pspolygon[fillstyle=solid,fillcolor=mygreen,linewidth=0cm,linecolor=mygreen](-0.2,-0.05)(0.2,-0.05)(0.2,0.05)(-0.2,0.05)
}

\nc{\elegant}{1.5pt}
\nc{\moyen}{1.0pt}
\nc{\mince}{0.5pt}

\nc{\alg}[1]{\mathsf{#1}}
\nc{\mc}[1]{\mathcal{#1}}
\nc{\field}[1]{\mathbb{#1}}

\begin{document}

\topmargin -5mm
\oddsidemargin 5mm

\vspace*{-2cm}

\setcounter{page}{1}

\vspace{22mm}
\begin{center}
{\huge {\bf 
Uncoiled affine Temperley--Lieb algebras\\[0.2cm] and their Wenzl--Jones projectors}}
\end{center}

\vspace{1cm}
\begin{center}
{\vspace{-5mm}\Large Alexis Langlois-R\'emillard$^{ab}$, Alexi Morin-Duchesne$^{ac}$}\\[.5cm]
{\em { }$^a$Department of Applied Mathematics, Computer Science and Statistics \\ Ghent University, 9000 Ghent, Belgium}
 \\[.2cm] 
{\em { }$^b$Hausdorff Center for Mathematics, 
University of Bonn, 53115 Bonn, Germany}
\\[.2cm] 
{\em { }$^c$Department of Mathematics, Royal Military Academy, 1000 Brussels, Belgium}
\\[.2cm] 
{\tt alexis.langlois-remillard\,@\,tutanota.com \qquad
alexi.morin.duchesne\,@\,gmail.com}
\end{center}

\vspace{0.6cm}
\centerline{\bf{Abstract}}
\vspace{0.3cm}
Affine and periodic Temperley--Lieb algebras are families of diagrammatic algebras that find diverse applications in mathematics and physics. These algebras are infinite dimensional, yet most of their interesting modules are finite. In this paper, we introduce finite quotients for these algebras, which we term {\it uncoiled affine Temperley--Lieb algebras} and {\it uncoiled periodic Temperley--Lieb algebras}. We study some of their properties, including their defining relations, their descriptions with diagrams, their dimensions, and their relations with  affine and skew sandwich cellular algebras. The uncoiled algebras all have finitely many one-dimensional modules. We construct a family of Wenzl--Jones idempotents, each of which projects onto one of these one-dimensional modules. Our construction is explicit and uses the similar projectors for the ordinary Temperley--Lieb algebras, as well as the diagrammatic description of the uncoiled algebras in terms of sandwich diagrams. We also discuss the Markov traces for the uncoiled algebras and their evaluations on the newly defined projectors, and find expressions involving Chebyshev polynomials of the first kind.

%
%

\tableofcontents

%
\section{Introduction}
%

The Temperley--Lieb algebra $\tl_n(\beta)$ is an algebra that describes the product of connectivity diagrams made of curves drawn inside a rectangle, in a way that keeps track of the non-local connectivity properties of those curves. It has a long history. This algebra of diagrams was first discovered by physicists in the context of percolation problems~\cite{TL71}, and subsequently used to describe ice-type models~\cite{B82}, XXZ spin chains \cite{PS90} and the $Q$-state Potts model~\cite{M91}. It then reappeared in the work of Jones on knot invariants~\cite{J83}, and its representation theory was worked out over the next years~\cite{GW93,W95,RSA14}. The algebra $\tl_n(\beta)$ is also in Schur--Weyl duality with the quantum group $U_q(\mathfrak{sl}_2)$, an idea that has continued to develop as it was later applied to study the representation theory of other diagrammatic algebras like the Brauer algebra~\cite{B37}, the BMW algebra~\cite{M87,BW89,CVZ21,CPAVZ23} and the partition algebra~\cite{J94,MS94,M96,HR05}. The Temperley--Lieb algebra is, in fact, a subalgebra of the partition algebra. Many results about the representation theory of this larger algebra were obtained using combinatorial methods which are also applicable to its subalgebras~\cite{COSSZ20}. Another approach is to view $\tl_n(\beta)$ as a quotient of the Hecke algebra of type $A$. From the latter, one can generalise the Temperley--Lieb algebra in various ways, either by studying other quotients of this Hecke algebra, or quotients of other types of Hecke algebras; see for example~\cite{GL04, HMR09}.\medskip

One generalisation of the Temperley--Lieb algebra is obtained by studying connectivity diagrams endowed with periodic boundary conditions. The precise names and the choices of relations used to define these algebras vary in the literature. In this paper, the {\it periodic Temperley--Lieb algebra} $\ptl_n(\beta)$ refers to the algebra defined with only the generators $e_0, e_1, \dots, e_{n-1}$, whereas the {\it affine Temperley--Lieb algebra} $\atl_n(\beta)$ also involves the translation generators $\Omega$ and $\Omega^{-1}$. These algebras were initially uncovered by physicists in their analysis of XXZ spin chains \cite{PS90,L91}. They also arise as certain quotients of the affine Hecke algebra of type $\hat{A}$~\cite{J94b,FG99}. In contrast to $\tl_n(\beta)$, both $\atl_n(\beta)$ and $\ptl_n(\beta)$ are infinite-dimensional algebras. Nonetheless, many of the methods employed to study $\tl_n(\beta)$ go through for $\atl_n(\beta)$ and $\ptl_n(\beta)$ with suitable modifications. The investigation of the representation theory of $\atl_n(\beta)$~\cite{MS93,G98,GL98,EG99,PSA22} allowed for a full characterisation of its irreducible modules, but also revealed a much richer panorama of indecomposable yet reducible modules than for the ordinary Temperley--Lieb algebra. The classification of these indecomposable modules remains an open problem.
\medskip

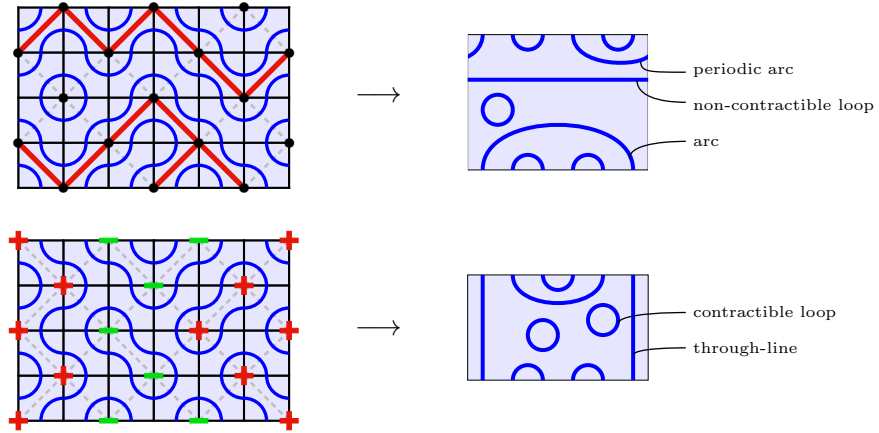
\begin{figure}
\begin{equation*}
\psset{unit=.6cm}
\begin{pspicture}[shift=-1.9](0,0)(6,4)
\facegrid{(0,0)}{(6,4)}
\rput(0,3){\loopar}\rput(1,3){\loopbr}\rput(2,3){\loopar}\rput(3,3){\loopbr}\rput(4,3){\loopbb}\rput(5,3){\loopab}
\rput(0,2){\loopab}\rput(1,2){\loopbb}\rput(2,2){\loopab}\rput(3,2){\loopbb}\rput(4,2){\loopbr}\rput(5,2){\loopar}
\rput(0,1){\loopbb}\rput(1,1){\loopab}\rput(2,1){\loopar}\rput(3,1){\loopbr}\rput(4,1){\loopbb}\rput(5,1){\loopab}
\rput(0,0){\loopbr}\rput(1,0){\loopar}\rput(2,0){\loopab}\rput(3,0){\loopar}\rput(4,0){\loopbr}\rput(5,0){\loopbb}
\rput(0,3){\loopa}\rput(1,3){\loopb}\rput(2,3){\loopa}\rput(3,3){\loopb}\rput(4,3){\loopb}\rput(5,3){\loopa}
\rput(0,2){\loopa}\rput(1,2){\loopb}\rput(2,2){\loopa}\rput(3,2){\loopb}\rput(4,2){\loopb}\rput(5,2){\loopa}
\rput(0,1){\loopb}\rput(1,1){\loopa}\rput(2,1){\loopa}\rput(3,1){\loopb}\rput(4,1){\loopb}\rput(5,1){\loopa}
\rput(0,0){\loopb}\rput(1,0){\loopa}\rput(2,0){\loopa}\rput(3,0){\loopa}\rput(4,0){\loopb}\rput(5,0){\loopb}
\multiput(0,0)(2,0){3}{\psdots(1,0)(1,2)(1,4)}
\multiput(0,0)(2,0){4}{\psdots(0,1)(0,3)}
\end{pspicture}
\qquad
\longrightarrow
\qquad
\psset{unit=1cm}
\begin{pspicture}[shift=-0.8](0,-0.8)(2.4,0.8)
\pspolygon[fillstyle=solid,fillcolor=lightlightblue,linewidth=0pt](0,-0.9)(2.4,-0.9)(2.4,0.9)(0,0.9)
\psarc[linecolor=blue,linewidth=1.5pt]{-}(2.0,0.9){0.2}{180}{0}
\psarc[linecolor=blue,linewidth=1.5pt]{-}(0.8,0.9){0.2}{180}{0}
\psarc[linecolor=blue,linewidth=1.5pt]{-}(0.8,-0.9){0.2}{0}{180}
\psarc[linecolor=blue,linewidth=1.5pt]{-}(1.6,-0.9){0.2}{0}{180}
\psbezier[linecolor=blue,linewidth=1.5pt]{-}(0.2,-0.9)(0.2,-0.1)(2.2,-0.1)(2.2,-0.9)
\psbezier[linecolor=blue,linewidth=1.5pt]{-}(2.42,0.63)(2.3,0.46)(1.4,0.46)(1.4,0.9)
\psbezier[linecolor=blue,linewidth=1.5pt]{-}(0.2,0.9)(0.2,0.72)(0.05,0.62)(-0.02,0.6)
\psframe[fillstyle=solid,linecolor=white,linewidth=0pt](-0.05,-0.91)(-0.005,0.91)
\psframe[fillstyle=solid,linecolor=white,linewidth=0pt](2.405,-0.91)(2.45,0.91)
\psarc[linecolor=blue,linewidth=1.5pt]{-}(0.4,-0.1){0.2}{0}{360}
\psline[linecolor=blue,linewidth=1.5pt]{-}(0,0.3)(2.4,0.3)
\rput(4.2,0){$\begin{array}{l}  \textrm{\tiny periodic arc} \\ \textrm{\tiny non-contractible loop} \\ \textrm{\tiny arc}\end{array}$}
\psbezier[linecolor=black,linewidth=0.5pt]{-}(2.3,0.55)(2.3,0.4)(2.4,0.4)(2.9,0.4)
\psbezier[linecolor=black,linewidth=0.5pt]{-}(2.25,0.3)(2.25,0)(2.8,0)(2.9,0)
\psbezier[linecolor=black,linewidth=0.5pt]{-}(2.15,-.75)(2.15,-.52)(2.8,-.52)(2.9,-.52)
\end{pspicture}
\end{equation*}
\vspace{0.4cm}
\begin{equation*}
\psset{unit=.6cm}
\begin{pspicture}[shift=-1.9](0,0)(6,4)
\facegrid{(0,0)}{(6,4)}
\rput(0,3){\loopab}\rput(1,3){\loopbb}\rput(2,3){\loopab}\rput(3,3){\loopbb}\rput(4,3){\loopab}\rput(5,3){\loopbb}
\rput(0,2){\loopbb}\rput(1,2){\loopab}\rput(2,2){\loopbb}\rput(3,2){\loopab}\rput(4,2){\loopbb}\rput(5,2){\loopab}
\rput(0,1){\loopab}\rput(1,1){\loopbb}\rput(2,1){\loopab}\rput(3,1){\loopbb}\rput(4,1){\loopab}\rput(5,1){\loopbb}
\rput(0,0){\loopbb}\rput(1,0){\loopab}\rput(2,0){\loopbb}\rput(3,0){\loopab}\rput(4,0){\loopbb}\rput(5,0){\loopab}
\rput(0,3){\loopb}\rput(1,3){\loopb}\rput(2,3){\loopb}\rput(3,3){\loopa}\rput(4,3){\loopa}\rput(5,3){\loopb}
\rput(0,2){\loopa}\rput(1,2){\loopa}\rput(2,2){\loopb}\rput(3,2){\loopa}\rput(4,2){\loopa}\rput(5,2){\loopa}
\rput(0,1){\loopb}\rput(1,1){\loopb}\rput(2,1){\loopb}\rput(3,1){\loopb}\rput(4,1){\loopa}\rput(5,1){\loopb}
\rput(0,0){\loopa}\rput(1,0){\loopa}\rput(2,0){\loopb}\rput(3,0){\loopa}\rput(4,0){\loopb}\rput(5,0){\loopb}
\rput(0,4){\Isinga}\rput(2,4){\Isingb}\rput(4,4){\Isingb}\rput(6,4){\Isinga}
\rput(1,3){\Isinga}\rput(3,3){\Isingb}\rput(5,3){\Isinga}
\rput(0,2){\Isinga}\rput(2,2){\Isingb}\rput(4,2){\Isinga}\rput(6,2){\Isinga}
\rput(1,1){\Isinga}\rput(3,1){\Isingb}\rput(5,1){\Isinga}
\rput(0,0){\Isinga}\rput(2,0){\Isingb}\rput(4,0){\Isingb}\rput(6,0){\Isinga}
\end{pspicture}
\qquad
\longrightarrow
\qquad
\psset{unit=1cm}
\begin{pspicture}[shift=-0.6](0,-0.7)(2.4,0.7)
\pspolygon[fillstyle=solid,fillcolor=lightlightblue,linewidth=0pt](0,-0.7)(2.4,-0.7)(2.4,0.7)(0,0.7)
\psarc[linecolor=blue,linewidth=1.5pt]{-}(1.2,0.7){0.2}{180}{0}
\psarc[linecolor=blue,linewidth=1.5pt]{-}(1.6,-0.7){0.2}{0}{180}
\psarc[linecolor=blue,linewidth=1.5pt]{-}(0.8,-0.7){0.2}{0}{180}
\psline[linecolor=blue,linewidth=1.5pt]{-}(0.2,-0.7)(0.2,0.7)
\psline[linecolor=blue,linewidth=1.5pt]{-}(2.2,-0.7)(2.2,0.7)
\psbezier[linecolor=blue,linewidth=1.5pt]{-}(0.6,0.7)(0.6,0.2)(1.8,0.2)(1.8,0.7)
\psarc[linecolor=blue,linewidth=1.5pt]{-}(1.0,-0.1){0.2}{0}{360}
\psarc[linecolor=blue,linewidth=1.5pt]{-}(1.8,0.1){0.2}{0}{360}
\rput(3.95,0.){$\begin{array}{l}  \textrm{\tiny contractible loop} \\ \textrm{\tiny through-line} \end{array}$}
\psbezier[linecolor=black,linewidth=0.5pt]{-}(2.0,0.05)(2.0,0.2)(2.8,0.2)(2.9,0.2)
\psbezier[linecolor=black,linewidth=0.5pt]{-}(2.2,-0.37)(2.2,-0.27)(2.8,-0.27)(2.9,-0.27)
\end{pspicture}
\end{equation*}
\caption{Top left panel: a configuration of the model of bond percolation on the square lattice tilted by $45^\circ$. Bottom left panel: a configuration of the Fortuin--Kasteleyen expansion of the Ising model on the square lattice. Right panel: the affine diagrams that these configurations are mapped to, along with the names we use for the various features in these diagrams.}
\label{fig:maps.to.ATL}
\end{figure}

The Temperley--Lieb algebra and its affine and periodic versions find further applications in statistical mechanics to describe special families of models defined in terms of non-local observables~\cite{PRZ06,J09}, with the Temperley--Lieb diagrams conveniently encoding the connectivity properties of these objects. \cref{fig:maps.to.ATL} gives examples of configurations of two statistical models on the square lattice, namely the model of bond percolation and the Ising model, and their maps to connectivity diagrams of the algebra $\atl_n(\beta)$. Interestingly, the representations of $\atl_n(\beta)$ and $\ptl_n(\beta)$ relevant in physics are finite-dimensional, marking a clear contrast with the algebras' infiniteness. The reasons behind this are clear. First, on a periodic geometry, the curves that form the Temperley--Lieb diagrams may wind arbitrarily many times along this axis. Second, such curves may close after winding along this periodicity axis and thus form {\it non-contractible loops}, of which there can also be arbitrarily many. It thus makes sense to define and study quotients of $\atl_n(\beta)$ and $\ptl_n(\beta)$, with extra relations that render the resulting algebras finite-dimensional. A first objective of this paper is to study and define such algebras, where (i) the non-contractible loops are removed with weights $\alpha$, and (ii) the winding of the curves around the full $n$ nodes is equivalent to no winding at all, times an {\it unwinding factor} $\gamma$. We refer to these quotient algebras as {\it uncoiled affine Temperley--Lieb algebras} and {\it uncoiled periodic Temperley--Lieb algebras}. These algebras turn out to depend crucially on the parity of the system size~$n$. In particular, it is natural to define two types of uncoiled algebras for $n$ even compared to just one type for $n$ odd. This will thus lead to a total of six algebras, half of them associated with $\atl_n(\beta)$ and the other half with $\ptl_n(\beta)$. In a sense, by studying these quotients, we will be fragmenting the infinite-dimensional algebras $\atl_n(\beta)$ and $\ptl_n(\beta)$ into families of finite-dimensional ones, whose representation theory can be expected to be simpler. In particular, we hope that this fragmentation will help make the study of the indecomposable modules more amenable.
\medskip

The Temperley--Lieb algebra is also a prime example of a unital associative algebra falling under the umbrella of cellular algebras~\cite{GL96,GL04}, and whose representation theory can be completely unravelled using this framework. The notion of cellularity was originally introduced for finite-dimensional algebras, but was already extended in~\cite{GL98} to an infinite-dimensional case: the affine Temperley--Lieb algebra. The extension to other infinite-dimensional cases was formalised by the notion of affine cellularity given by K\"onig and Xi~\cite{KX12}, which adjoins commutative algebras to the cells. The notion of affine cellularity can still be weakened, notably by removing the requirement of the commutativity of the associated algebras and by removing the involution from the cellular datum, which gives rise to sandwich cellularity~\cite{MT22,TV22,T22}. We refer to \cite{T22} for more on sandwich cell theory. Notably, the algebra $\atl_n(\beta)$ is affine cellular in all cases. To the best of our knowledge, the same question for the algebra $\ptl_n(\beta)$ has not been previously addressed. We will show in this paper that it is affine cellular for $n$ odd. Our construction, however, does not carry over to the case with $n$ even. We will also show that the uncoiled algebras are affine cellular, except in one special case. For  this special case and for $\ptl_n(\beta)$ with $n$ even, the first three axioms defining sandwich cellularity hold, but not the fourth axiom that dictates the action of the anti-involution. This motivate us to modify the fourth axiom and introduce the notion of \emph{skew sandwich cellular algebras}, building on Hu, Mathas and Rostam's notion of skew cellularity~\cite{HMR23}. The sandwich cellular and skew affine cellular algebras are, in fact, special subcases of skew sandwich cellular algebras. We will then show that the two problematic algebras described above are skew sandwich cellular.
\medskip

There is another reason, from a physics standpoint, to launch the study of these quotient algebras. The Temperley--Lieb algebra is particularly useful to study the continuum scaling limits ($n \to \infty$) of critical statistical models and their description with conformal field theories (CFTs). The connection between the fusion of modules over $\tl_n(\beta)$ and the fusion of fields for boundary CFTs is by now well established \cite{PRZ06,RS07,RP07,GV13}, and has proven very useful to compute conformal correlation functions. The search for a similar connection tying the fusion of modules over $\atl_n(\beta)$ and the fusion of bulk fields is less immediate. Despite some recent progress \cite{GS16,GJS18,BSA18,IMD22}, many questions remain unanswered, again leading us to wonder whether $\atl_n(\beta)$ is really the ideal algebraic structure to probe the scaling limit. We hope that the uncoiled algebras defined in this paper will provide a new tool to tackle these issues.\medskip

One of the main results that we provide for the uncoiled algebras is the definition of a family of idempotents mimicking the famed Wenzl--Jones projectors of $\tl_n(\beta)$.
The Wenzl--Jones projectors $P_1, P_2, \dots, P_n$~\cite{J83,W88,KL94} form an important family of elements of $\tl_n(\beta)$. They are defined recursively and play a crucial role for the representation theory. In particular, $P_n$ projects on the unique one-dimensional module of $\tl_n(\beta)$. In the study of Kazhdan--Lusztig theory~\cite{KL79}, Wenzl--Jones projectors encode particularly important idempotent endomorphisms projecting to indecomposable summands, allowing for explicit expressions for the bases dual to the canonical bases in the tensor products of finite-dimensional simple modules of $U_q(\mathfrak{sl}_2)$~\cite{FK97}. They have recently resurfaced in the trail of Soergel bimodule theory~\cite{EMTW20}, as slight generalisations of these projectors were used to give a categorical link between the Temperley--Lieb algebra and Bott--Samelson bimodules, specifically for the infinite dihedral group~\cite{E16} and for the universal Coxeter group~\cite{EL17}. This categorification should also hold for the $p$-Kazhdan--Lusztig basis~\cite{JW17}. Recently, Martin and Spencer~\cite{MS21} studied this correspondence and provided a construction, building on~\cite{BLS19}, of the Wenzl--Jones projectors in positive characteristic  intrinsic to the diagrammatic construction and to the Temperley--Lieb category without resorting to the tilting theory of $U_q(\mathfrak{sl}_2)$~\cite{STWZ21}. The Wenzl--Jones projectors also find multiple applications in physics. For instance, they arise in the construction of fused $R$-matrices for integrable models \cite{S88,MDPR14}. They are also used to define loop models with thick boundary conditions and thus to gain insight on the boundary conformal field theory describing the continuum scaling limit of these models \cite{PRZ06,MDRR15,LRSA19}.\medskip

Differently from $\tl_n(\beta)$, the algebra $\atl_n(\beta)$ has a one-parameter family of one-dimensional representations, labelled by the {\it twist parameter} $z \in \mathbb C^\times$. Defining a projector on each representation in this continuous family is not possible. This can, however, be achieved for the uncoiled algebras, as they turn out to have finitely many one-dimensional representations. We will construct these projectors by expressing them as linear combinations of connectivities sandwiched between the Wenzl--Jones projectors of $\tl_n(\beta)$, seen as a subalgebra of $\atl_n(\beta)$, and we will give closed-form formulas for the coefficients arising in these sums.\medskip

One of the tools used by Jones in his construction of the Jones invariant of knots~\cite{J85} was the Markov trace. It is originally a trace on braid groups amenable to Markov moves~\cite{B75}. It has multiple interpretations, and can be extended to other types of quotients of the braid group~\cite{WW11}. The Wenzl--Jones projectors for $\tl_n(\beta)$ behave nicely with respect to this trace: they evaluate to Chebyshev polynomials of the second kind. This fact shows up in the categorification of the integer polynomial ring by the Temperley--Lieb category where the Wenzl--Jones projectors and these polynomials are in direct correspondence~\cite{KS21}. In contrast, for the uncoiled algebras, the Markov traces of the Wenzl--Jones projectors evaluate to  Chebyshev polynomials of the first kind and sums thereof.\medskip 

We mention in passing that the quotients that we study are not the only possible quotients of $\atl_n(\beta)$ and $\ptl_n(\beta)$. For instance, the quotient considered in \cite{QW18a,BGJST20} involves an extra relation that equates the {\it braid transfer matrix} (or {\it hoop operator}) to a constant. These quotients then fragment $\atl_n(\beta)$ into finite-dimensional algebras that are isomorphic to certain one-boundary Temperley--Lieb algebras \cite{BGJST20}. The corresponding dimensions, projectors and representations are then different from the ones that we consider here. The interest of~\cite{QW18a} was also to categorify the integer polynomial ring by a diagrammatic category. The projectors that they define also correspond via their Markov trace to Chebyshev polynomials of the first kind. It will be interesting to interpret the uncoiled algebras in a similar framework, to study the resulting categorification, and compare it with the one studied in~\cite{QW18a}.

\paragraph{A note on the different definitions and naming conventions.}

There is vast literature on periodic/affine Temperley--Lieb algebras and their quotients. With no claim to exhaustiveness, we provide here some insight and comments on the definitions and terminology used, to help the readers compare our results with those obtained previously by other authors.\medskip

\begin{table}[h!]
\centering
\begin{tabular}{c | l | l}
Papers  & $\ptl_n(\beta)$ & $\atl_n(\beta)$\\\hline
This paper & periodic & affine\\\hline
\cite{BGJST20,GL98,PSA22} & periodic & affine\\\hline
\cite{EG99,G98} & affine & algebra of affine diagrams\\\hline
\cite{MS93,MS94} & periodic & algebra of affine diagrams\\\hline
\cite{IMD22,IMD24,MDSA13} & periodic & enlarged periodic\\\hline
\cite{RJ07} & periodic & NA\\\hline
\cite{BSA18,GJS18,GS16} & NA & affine or periodic\\\hline
\cite{GRSV15} & periodic & periodic with translation or 
\\&&algebra of affine diagrams\\\hline
 \cite{GL98,KX12,GL04,HMR09,QW18a} & Temperley--Lieb quotient & affine (type $GL_n$)
\\&
of the (unextended) affine 
\\&Hecke algebra  of type $\tilde{A}$ & 
\end{tabular}
\caption{
Different naming conventions in the articles cited here for the periodic and affine Temperley--Lieb algebras $\ptl_n(\beta)$ and $\atl_n(\beta)$ defined in \cref{sec:def.aptl.algebras}.}\label{table:naming_convention}
\end{table}

First, the two infinite extensions of the Temperley--Lieb algebras that we denote here by $\ptl_n(\beta)$ and $\atl_n(\beta)$ are assigned various names in the literature, which are often specific to certain communities and can change over time. To help the readers, we summarise in \cref{table:naming_convention} the terminology used for these algebras in the papers that we cite later on.\medskip

Second, finite-dimensional quotients of the affine and periodic Temperley--Lieb algebras also appeared early on in the literature, either because of their relevance in physics or as tools in the study of the representation theory of the infinite-dimensional diagrammatic algebras. We review the various naming conventions for these quotients in \cref{table:different_quotients} and comment on their relations to the uncoiled algebras that we study here.

\begin{table}[h]
\centering
\begin{tabular}{>{\RaggedRight}p{2.5cm}|c||>{\RaggedRight}p{5cm}|>{\RaggedRight}p{5cm}}
Name & Papers cited & Informal definition & Remarks \\\hline
Jones annular algebra & \cite{J94b} & $\Omega^n=\id$ and there are no periodic arcs on diagrams with no through-lines. & This corresponds to $\qatl_n(\beta,1)$ for $n$ odd, but is another algebra for $n$ even.\\\hline 
$q$-Jones algebra & \cite{G98,EG99} & $\Omega^n=\gamma\,\id$ 
and there are no periodic arcs on diagrams with no through-lines for $n$ even. & This corresponds to $\qatl_n(\beta,\gamma)$ for $n$ odd (their $q$ is our $\gamma$), but is another algebra for $n$ even.\\\hline
augmented Jones--Temperley--Lieb algebra & \cite{GRSV15} & This is the subalgebra $\mathsf{JTL}_n^{au}(\beta)\subset \qatla_n(\beta,\beta)$ with $n$ even
generated by $e_0,\dots, e_{n-1}, \Omega^2$ (so only even powers of $\Omega$ are allowed).
& This algebra is not defined for $n$ odd. The word {\it augmented} refers to the fact that the algebra contains all periodic diagrams with no through-lines, whereas the Jones--Temperley--Lieb algebra (or Jones annular algebra) contains only the planar ones.\\\hline
Cyl-algebra & \cite{CHMPTRSS24} & $\Omega^n =1$ (called there the Dehn twist) is the identity, so this is $\qatl_n(1,1)$ for $n$ odd and the quotient $\atl_n(1)/
\langle \Omega^n =1\rangle$ for $n$ even.& 
This work defines, in fact, the stripped cylinder cobordism category $\mathrm{Cyl}$. The identification with the uncoiled algebras requires additional enrichments; see~\cite[Sect.~4]{CHMPTRSS24}.
\\\hline
$\mathbf{ATL}^{\text{ess}}(n)$ & \cite{QW18a} & Quotient of the affine Temperley--Lieb algebra where the braid transfer matrix is $0$: $\atl_n(2)/\langle\mathbf{F}=0\rangle$ where $\mathbf F$ is defined in~\eqref{eq:Fb.def}.  & It is defined as a category; the algebra is understood as  the endomorphisms of the object $n$.
\end{tabular}
\caption{
Various quotients of the periodic or affine Temperley--Lieb algebras and their naming conventions, with remarks and comparisons with the definitions used in our paper.}\label{table:different_quotients}
\end{table}

\paragraph{Main results.} 

We list here the main results of our paper: 
(i) \cref{def:SandwichCell}
of skew sandwich cellular algebras, (ii) the proof that the algebra $\ptl_n(\beta)$ is affine cellular for $n$ odd and skew sandwich cellular for $n$ even in \cref{prop:pTL.skew.sandwich}, 
(iii)~the definition of the uncoiled affine and periodic algebras in~\eqref{eq:odd.quotients}, \eqref{eq:even.quotients.1} and~\eqref{eq:even.quotients.2};  (iv)~the proof that these algebras are sandwich cellular 
in \cref{prop:uaTLupTLSandwichCellular};  (v)~the construction of the Wenzl--Jones projectors in \eqref{eq:WJ.upTL} and \eqref{eq:WJ.uaTL};  (vi)~the linear relations satisfied by the coefficients arising in the definition of these projectors in \cref{Prop.const.rel}; (vii)~the closed forms of the solutions to these relations in \cref{Prop.const.uptla.uptlb,Prop.const.uptl,Prop.const.qatl.qatla.qatlb}; and (viii)~the Markov trace of the Wenzl--Jones projectors in \cref{prop:MarkovTrace.FQnr.et.FQn}.

\paragraph{Outline.} 

The paper is structured as follows. In \cref{sec:def.aptl.algebras}, we review the definitions of $\atl_n(\beta)$ and $\ptl_n(\beta)$, describe their standard modules $\repW_{n,d,z}$ and present their description with sandwich diagrams. A short discussion on central elements of these two algebras is presented in \cref{sec:non.isomorphic}. In \cref{sec:def.uaptl}, we define the uncoiled affine and periodic Temperley--Lieb algebras, both in terms of generators and connectivity diagrams. We study the basic properties of these algebras, namely their dimensions, their standard representations, and their description as affine cellular algebras and skew sandwich cellular algebras. The proofs of certain binomial identities, used to compute the dimensions of the uncoiled algebras, are relegated to \cref{app:binom.identities}. In \cref{sec:WJ.projs}, we define the Wenzl--Jones projectors, first for the uncoiled periodic Temperley--Lieb algebras and then for the uncoiled affine ones. We discuss the relation between the periodic and affine projectors and present their explicit construction. The construction involves some constants $\Gamma_{k,\ell}$, whose computation is non-trivial. In \cref{sec:Gammas}, we derive linear relations satisfied by these constants that determine them uniquely and use them to find
closed forms for the constants $\Gamma_{k,\ell}$, for each of the six uncoiled algebras. The proof uses contour integral expressions for these constants and is given in \cref{Sec:proof.closed.form}. 
In \cref{sec:Markov}, we discuss the Markov trace and its evaluation on the Wenzl--Jones projectors. We conclude with some remarks in \cref{sec:conclusion}.

%
\section{Affine and periodic Temperley--Lieb algebras}\label{sec:def.aptl.algebras}
%

\subsection{Definitions with generators}\label{ssec:def}

Temperley--Lieb algebras and their generalisations are families of associative unital algebras that have equivalent definitions in terms of generators and diagrams. Here we define two algebras, the affine Temperley--Lieb algebra $\atl_n(\beta)$ and the periodic Temperley--Lieb algebra $\ptl_n(\beta)$, for $n\in\mathbb N$ and $\beta\in \mathbb C$.
For $\atl_n(\beta)$, the definition with generators involves a unit $\id$, elements $e_0, e_1, \dots, e_{n-1}$ (with $n > 2$), and two translation generators $\Omega$ and $\Omega^{-1}$. These are subject to the relations
\begin{subequations}
\begin{alignat}{4}
\label{eq:PTL.relations}
e_j^2 &= \beta \, e_j, \qquad &&e_j  e_{j \pm 1} e_j = e_j, \qquad &&e_i e_j = e_j e_i \qquad \text{for } |i-j|>1,
\\
\label{eq:ATL.relations}
\Omega \, e_j \, \Omega^{-1} &= e_{j-1}, \qquad &&\Omega \, \Omega^{-1} = \Omega^{-1} \, \Omega = \id, \qquad  &&\Omega^2 e_1 = e_{n-1}e_{n-2} \cdots e_2e_1,
\end{alignat}
\end{subequations}
where $i,j$ are taken modulo $n$, and $\beta \in \mathbb C$ is a free parameter. In contrast, the algebra $\ptl_n(\beta)$ involves only the generators $e_0, e_1, \dots, e_{n-1}$ and the unit $\id$, and these satisfy the relations~\eqref{eq:PTL.relations}. We therefore have
\begin{subequations}
\label{eq:atl.ptl.def}
\begin{alignat}{2}
\label{eq:ptl.def}
\atl_n(\beta) &= \big\langle \Omega,\Omega^{-1}, e_0, e_1,\dots, e_{n-1} \big\rangle \Big/ \big\{\eqref{eq:PTL.relations},\eqref{eq:ATL.relations}\big\}, \\[0.1cm]
\label{eq:atl.def}
\ptl_n(\beta) &= \big\langle \id, e_0, e_1,\dots, e_{n-1} \big\rangle \Big/ \big\{\eqref{eq:PTL.relations}\big\}.
\end{alignat}
\end{subequations}
Clearly, $\ptl_n(\beta)$ is a subalgebra of $\atl_n(\beta)$. Because the words $\Omega^k$ and $(e_0 e_1 \cdots e_{n-1})^k$ with $k \in \mathbb N$ can never be reduced using the above relations, it is clear that $\atl_n(\beta)$ and $\ptl_n(\beta)$ are infinite-dimensional algebras. Finally, we note that the ordinary Temperley--Lieb algebra $\tl_n(\beta)$ is the finite subalgebra of $\ptl_n(\beta)$ defined as
\be
\tl_n(\beta) = \big\langle \id, e_1, e_2,\dots, e_{n-1} \big\rangle \Big/ \big\{\eqref{eq:PTL.relations}\big\}.
\ee
We parameterise the variable $\beta$ as
\be
\beta = -q-q^{-1}, \qquad q \in \mathbb C^\times.
\ee

\begin{Proposition}
\label{prop:3.forms}
Any word $c$ in the generators of $\atl_n(\beta)$ can be rewritten into one of the three following forms: 
\begin{enumerate}
\item[$\mathrm{(i)}$] $c = \Omega^k$, for some $k \in \mathbb Z$ (with $\Omega^0 = \id$);
\item[$\mathrm{(ii)}$] $c = w$, where $w$ is a word in the generators $e_j$;
\item[$\mathrm{(iii)}$] $c = \Omega w$, where $w$ is a word in the generators $e_j$.
\end{enumerate}
\end{Proposition}
\proof
If $c=\Omega^k$ or $c = \id$, there is nothing to prove, and likewise if $c$ is already a word in only the generators $e_j$. The only pending case is when $c$ involves both a non-zero power of $\Omega$ and at least one generator $e_j$.
In this case, one can use \eqref{eq:ATL.relations} to commute all powers of $\Omega$ to the left, which yields $c = \Omega^k w$ for some $k \in \mathbb Z$ and $w$ a word in the generators $e_j$. Finally, the powers of $\Omega$ are reduced in steps of two using the relation
\be
\Omega^2 e_i = e_{n-2+i}e_{n-3+i} \cdots e_{i+1}e_{i},
\ee
where the labels are again taken modulo $n$. This last relation is obtained by combining the first and last relations of \eqref{eq:ATL.relations}. This ends the proof.
\eproof

\subsection{Definitions with diagrams}

In its diagrammatic presentation, $\atl_n(\beta)$ is constructed on the vector space of objects called {\it connectivities}. These are diagrams drawn inside a horizontal box with $n$ marked nodes on the top edge of the box and $n$ more nodes on its bottom edge. A connectivity is then a collection of non-intersecting curves that we call {\it loop segments} that connect the nodes pairwise. The box has periodic boundary conditions in the horizontal direction, so that a loop segment connecting two nodes may travel to the left edge of the box and reappear on the right edge at the same height. Moreover, a connectivity may have an arbitrary number of {\it non-contractible loops}, namely loop segments that directly connect the left and right edges of the box at the same height. Here are examples of connectivities in $\atl_6(\beta)$:
\be
\label{eq:c1c2}
c_1 = \ 
\begin{pspicture}[shift=-0.6](0,-0.7)(2.4,0.7)
\pspolygon[fillstyle=solid,fillcolor=lightlightblue,linewidth=0pt](0,-0.7)(2.4,-0.7)(2.4,0.7)(0,0.7)
\psarc[linecolor=blue,linewidth=1.5pt]{-}(0.4,0.7){0.2}{180}{0}
\psarc[linecolor=blue,linewidth=1.5pt]{-}(1.6,0.7){0.2}{180}{0}
\psarc[linecolor=blue,linewidth=1.5pt]{-}(0.8,-0.7){0.2}{0}{180}
\psarc[linecolor=blue,linewidth=1.5pt]{-}(1.6,-0.7){0.2}{0}{180}
\psbezier[linecolor=blue,linewidth=1.5pt]{-}(0.2,-0.7)(0.2,-0.1)(2.2,-0.1)(2.2,-0.7)
\psbezier[linecolor=blue,linewidth=1.5pt]{-}(-0.02,0.43)(0.1,0.26)(1.0,0.26)(1.0,0.7)
\psbezier[linecolor=blue,linewidth=1.5pt]{-}(2.2,0.7)(2.2,0.52)(2.35,0.42)(2.42,0.4)
\psframe[fillstyle=solid,linecolor=white,linewidth=0pt](-0.05,-0.71)(-0.005,0.71)
\psframe[fillstyle=solid,linecolor=white,linewidth=0pt](2.405,-0.71)(2.45,0.71)
\psline[linecolor=blue,linewidth=1.5pt]{-}(0,0)(2.4,0)
\end{pspicture}\ ,
\qquad
c_2 = \ 
\begin{pspicture}[shift=-0.6](0,-0.7)(2.4,0.7)
\pspolygon[fillstyle=solid,fillcolor=lightlightblue,linewidth=0pt](0,-0.7)(2.4,-0.7)(2.4,0.7)(0,0.7)
\psarc[linecolor=blue,linewidth=1.5pt]{-}(0,0.7){0.2}{-90}{0}
\psarc[linecolor=blue,linewidth=1.5pt]{-}(2.4,0.7){0.2}{180}{270}
\psarc[linecolor=blue,linewidth=1.5pt]{-}(1.2,0.7){0.2}{180}{0}
\psarc[linecolor=blue,linewidth=1.5pt]{-}(1.6,-0.7){0.2}{0}{180}
\psbezier[linecolor=blue,linewidth=1.5pt]{-}(0.2,-0.7)(0.2,0)(0.6,0)(0.6,0.7)
\psbezier[linecolor=blue,linewidth=1.5pt]{-}(0.6,-0.7)(0.6,0)(1.8,0)(1.8,0.7)
\psbezier[linecolor=blue,linewidth=1.5pt]{-}(1.0,-0.7)(1.0,-0.2)(2.2,-0.2)(2.2,-0.7)
\end{pspicture}\ .
\ee

The product $c_1 c_2$ of two connectivities is defined by stacking $c_2$ over $c_1$. The resulting connectivity is then read off from the connection of the nodes at the top and bottom of the new bigger box. If there are {\it contractible loops}, they are erased and the resulting connectivity is multiplied by a factor of $\beta$ for each such loop. The parameter $\beta \in \mathbb C$ is the {\it loop fugacity} for the contractible loops. If new non-contractible loops are created, they are not erased and are instead added to the stack of such loops, if there were some in the first place. Here is an example of this product:
\be
c_1 c_2 = \ 
\begin{pspicture}[shift=-1.3](0,-0.7)(2.4,2.1)
\pspolygon[fillstyle=solid,fillcolor=lightlightblue,linewidth=0pt](0,-0.7)(2.4,-0.7)(2.4,0.7)(0,0.7)
\psarc[linecolor=blue,linewidth=1.5pt]{-}(0.4,0.7){0.2}{180}{0}
\psarc[linecolor=blue,linewidth=1.5pt]{-}(1.6,0.7){0.2}{180}{0}
\psarc[linecolor=blue,linewidth=1.5pt]{-}(0.8,-0.7){0.2}{0}{180}
\psarc[linecolor=blue,linewidth=1.5pt]{-}(1.6,-0.7){0.2}{0}{180}
\psbezier[linecolor=blue,linewidth=1.5pt]{-}(0.2,-0.7)(0.2,-0.1)(2.2,-0.1)(2.2,-0.7)
\psbezier[linecolor=blue,linewidth=1.5pt]{-}(-0.02,0.43)(0.1,0.26)(1.0,0.26)(1.0,0.7)
\psbezier[linecolor=blue,linewidth=1.5pt]{-}(2.2,0.7)(2.2,0.52)(2.35,0.42)(2.42,0.4)
\psframe[fillstyle=solid,linecolor=white,linewidth=0pt](-0.05,-0.71)(-0.005,0.71)
\psframe[fillstyle=solid,linecolor=white,linewidth=0pt](2.405,-0.71)(2.45,0.71)
\psline[linecolor=blue,linewidth=1.5pt]{-}(0,0)(2.4,0)
\rput(0,1.4){\pspolygon[fillstyle=solid,fillcolor=lightlightblue,linewidth=0pt](0,-0.7)(2.4,-0.7)(2.4,0.7)(0,0.7)
\psarc[linecolor=blue,linewidth=1.5pt]{-}(0,0.7){0.2}{-90}{0}
\psarc[linecolor=blue,linewidth=1.5pt]{-}(2.4,0.7){0.2}{180}{270}
\psarc[linecolor=blue,linewidth=1.5pt]{-}(1.2,0.7){0.2}{180}{0}
\psarc[linecolor=blue,linewidth=1.5pt]{-}(1.6,-0.7){0.2}{0}{180}
\psbezier[linecolor=blue,linewidth=1.5pt]{-}(0.2,-0.7)(0.2,0)(0.6,0)(0.6,0.7)
\psbezier[linecolor=blue,linewidth=1.5pt]{-}(0.6,-0.7)(0.6,0)(1.8,0)(1.8,0.7)
\psbezier[linecolor=blue,linewidth=1.5pt]{-}(1.0,-0.7)(1.0,-0.2)(2.2,-0.2)(2.2,-0.7)}
\end{pspicture}
\ = \beta \ 
\begin{pspicture}[shift=-0.6](0,-0.7)(2.4,0.7)
\pspolygon[fillstyle=solid,fillcolor=lightlightblue,linewidth=0pt](0,-0.7)(2.4,-0.7)(2.4,0.7)(0,0.7)
\psarc[linecolor=blue,linewidth=1.5pt]{-}(0,0.7){0.2}{-90}{0}
\psarc[linecolor=blue,linewidth=1.5pt]{-}(1.2,0.7){0.2}{180}{0}
\psarc[linecolor=blue,linewidth=1.5pt]{-}(2.4,0.7){0.2}{180}{270}
\psarc[linecolor=blue,linewidth=1.5pt]{-}(0.8,-0.7){0.2}{0}{180}
\psarc[linecolor=blue,linewidth=1.5pt]{-}(1.6,-0.7){0.2}{0}{180}
\psbezier[linecolor=blue,linewidth=1.5pt]{-}(0.6,0.7)(0.6,0.2)(1.8,0.2)(1.8,0.7)
\psframe[fillstyle=solid,linecolor=white,linewidth=0pt](-0.05,-0.71)(-0.005,0.71)
\psframe[fillstyle=solid,linecolor=white,linewidth=0pt](2.4005,-0.71)(2.45,0.71)
\psline[linecolor=blue,linewidth=1.5pt]{-}(0,-0.1)(2.4,-0.1)
\psline[linecolor=blue,linewidth=1.5pt]{-}(0,0.1)(2.4,0.1)
\psbezier[linecolor=blue,linewidth=1.5pt]{-}(0.2,-0.7)(0.2,-0.1)(2.2,-0.1)(2.2,-0.7)
\end{pspicture}\ .
\ee
The product is then extended bilinearly to any linear combinations of connectivities. The algebra $\atl_n(\beta)$ is the free vector space spanned by the set of connectivities on $n$ nodes, endowed with this product rule. There can be an arbitrary number of non-contractible loops, and likewise the through-lines may cross the periodic boundary condition an arbitrary number of times. It is then clear that this algebra is infinite dimensional.\medskip

The two presentations of $\atl_n(\beta)$, the first with generators and their relations and the second in terms of the diagrams, are known~\cite{MS93,FG99} to be equivalent for $n\geq 2$. In particular, the generators are represented by the connectivities
\begin{subequations}
\label{eq:generator.diagrams}
\begin{alignat}{3}
 e_0 &= \
\begin{pspicture}[shift=-0.45](0,-0.55)(2.0,0.35)
\pspolygon[fillstyle=solid,fillcolor=lightlightblue,linewidth=0pt](0,-0.35)(2.0,-0.35)(2.0,0.35)(0,0.35)
\rput(0.2,-0.55){$_1$}\rput(0.6,-0.55){$_2$}\rput(1.0,-0.55){\small$...$}\rput(1.8,-0.55){$_n$}
\rput(1.0,0.0){\small$...$}
\psarc[linecolor=blue,linewidth=1.5pt]{-}(0.0,0.35){0.2}{-90}{0}
\psarc[linecolor=blue,linewidth=1.5pt]{-}(0.0,-0.35){0.2}{0}{90}
\psline[linecolor=blue,linewidth=1.5pt]{-}(0.6,0.35)(0.6,-0.35)
\psline[linecolor=blue,linewidth=1.5pt]{-}(1.4,0.35)(1.4,-0.35)
\psarc[linecolor=blue,linewidth=1.5pt]{-}(2.0,-0.35){0.2}{90}{180}
\psarc[linecolor=blue,linewidth=1.5pt]{-}(2.0,0.35){0.2}{180}{-90}
\end{pspicture}\ \ ,
\qquad
&&e_j =  \
\begin{pspicture}[shift=-0.525](-0.0,-0.55)(3.2,0.35)
\pspolygon[fillstyle=solid,fillcolor=lightlightblue,linecolor=black,linewidth=0pt](0,-0.35)(0,0.35)(3.2,0.35)(3.2,-0.35)(0,-0.35)
\psline[linecolor=blue,linewidth=1.5pt]{-}(0.2,-0.35)(0.2,0.35)
\rput(0.6,0){$...$}
\psline[linecolor=blue,linewidth=1.5pt]{-}(1.0,-0.35)(1.0,0.35)
\psarc[linecolor=blue,linewidth=1.5pt]{-}(1.6,0.35){0.2}{180}{360}
\psarc[linecolor=blue,linewidth=1.5pt]{-}(1.6,-0.35){0.2}{0}{180}
\psline[linecolor=blue,linewidth=1.5pt]{-}(2.2,-0.35)(2.2,0.35)
\rput(2.6,0){$...$}
\psline[linecolor=blue,linewidth=1.5pt]{-}(3.0,-0.35)(3.0,0.35)
\rput(0.2,-0.6){$_1$}
\rput(1.38,-0.6){$_{j}$}
\rput(1.9,-0.6){$_{j+1}$}
\rput(3.0,-0.6){$_n$}
\end{pspicture}
 \qquad (1\le j \le n-1),  
\\[0.5cm] 
\Omega &= \ 
\begin{pspicture}[shift=-0.45](0,-0.55)(2.0,0.35)
\rput(0.2,-0.55){$_1$}\rput(0.6,-0.55){$_2$}\rput(1.0,-0.55){\small$...$}\rput(1.8,-0.55){$_n$}
\pspolygon[fillstyle=solid,fillcolor=lightlightblue,linewidth=0pt](0,-0.35)(2.0,-0.35)(2.0,0.35)(0,0.35)
\psbezier[linecolor=blue,linewidth=1.5pt]{-}(0.2,0.35)(0.2,0.1)(-0.02,0)(-0.05,-0.05)
\rput(2,0){\psbezier[linecolor=blue,linewidth=1.5pt]{-}(-0.2,-0.35)(-0.2,-0.1)(0.02,0)(0.05,0.05)}
\multiput(0.4,0)(0.4,0){4}{\psbezier[linecolor=blue,linewidth=1.5pt]{-}(-0.2,-0.35)(-0.2,-0.0)(0.2,0.0)(0.2,0.35)}
\psframe[fillstyle=solid,linecolor=white,linewidth=0pt](-0.1,-0.4)(-0.005,0.4)
\psframe[fillstyle=solid,linecolor=white,linewidth=0pt](2.005,-0.4)(2.1,0.4)
\end{pspicture}
\ \ ,
\qquad
&&\Omega^{-1}\,= \
 \begin{pspicture}[shift=-0.45](0,-0.55)(2.0,0.35)
\rput(0.2,-0.55){$_1$}\rput(0.6,-0.55){$_2$}\rput(1.0,-0.55){\small$...$}\rput(1.8,-0.55){$_n$}
\pspolygon[fillstyle=solid,fillcolor=lightlightblue,linewidth=0pt](0,-0.35)(2.0,-0.35)(2.0,0.35)(0,0.35)
\psbezier[linecolor=blue,linewidth=1.5pt]{-}(0.2,-0.35)(0.2,-0.1)(-0.02,0)(-0.05,0.05)
\rput(2,0){\psbezier[linecolor=blue,linewidth=1.5pt]{-}(-0.2,0.35)(-0.2,0.1)(0.02,0)(0.05,-0.05)}
\multiput(0.4,0)(0.4,0){4}{\psbezier[linecolor=blue,linewidth=1.5pt]{-}(-0.2,0.35)(-0.2,-0.0)(0.2,0.0)(0.2,-0.35)}
\psframe[fillstyle=solid,linecolor=white,linewidth=0pt](-0.1,-0.4)(-0.005,0.4)
\psframe[fillstyle=solid,linecolor=white,linewidth=0pt](2.005,-0.4)(2.1,0.4)
\end{pspicture}
\ \ ,
\qquad 
\id = \ 
\begin{pspicture}[shift=-0.525](0,-0.25)(2.4,0.8)
\pspolygon[fillstyle=solid,fillcolor=lightlightblue,linecolor=black,linewidth=0pt](0,0)(0,0.8)(2.4,0.8)(2.4,0)(0,0)
\psline[linecolor=blue,linewidth=1.5pt]{-}(0.2,0)(0.2,0.8)
\psline[linecolor=blue,linewidth=1.5pt]{-}(0.6,0)(0.6,0.8)
\psline[linecolor=blue,linewidth=1.5pt]{-}(1.0,0)(1.0,0.8)
\rput(1.4,0.4){$...$}
\psline[linecolor=blue,linewidth=1.5pt]{-}(1.8,0)(1.8,0.8)
\psline[linecolor=blue,linewidth=1.5pt]{-}(2.2,0)(2.2,0.8)
\rput(0.2,-0.25){$_1$}
\rput(0.6,-0.25){$_2$}
\rput(1.0,-0.25){\small$...$}
\rput(2.2,-0.25){$_n$}
\end{pspicture}
\ \ .
\end{alignat}
\end{subequations}
With the above rule for the product of connectivities, it is straightforward to show that these connectivities satisfy the relations \eqref{eq:PTL.relations} and \eqref{eq:ATL.relations}.\medskip

We note that $\atl_n(\beta)$ is well defined in terms of connectivities for $n=1$ and $n=2$. This is therefore the definition that we choose for $\atl_1(\beta)$ and $\atl_2(\beta)$. Their equivalent algebraic definitions are as follows. For $n=1$, the algebra is $\atl_1(\beta) = \langle\id,\Omega,\Omega^{-1}\rangle\big/\{\Omega \, \Omega^{-1} = \Omega^{-1} \, \Omega = \id\}$. For $n=2$, the generators of the algebra are $\Omega$, $\Omega^{-1}$, $e_0$ and $e_1$, and these satisfy all the relations in \eqref{eq:PTL.relations} and \eqref{eq:ATL.relations} except for the relation $e_j e_{j\pm1}e_j = e_j$ which is removed in this case. Finally, it also makes sense to define $\atl_n(\beta)$ for $n = 0$, for which the connectivities are stacks of non-contractible loops. We denote by $\id$ and $f$ the connectivities of $\atl_0(\beta)$ with zero and one non-contractible loop, respectively, and write $\atl_0(\beta) = \langle\id,f\rangle$ with no extra relations other than $\id f = f \id = f$. For uniformity of notation, we leave the dependence on $\beta$ in writing $\atl_0(\beta)$ and $\atl_1(\beta)$, even though these two algebras do not depend on $\beta$.\medskip

Because it is a subalgebra of $\atl_n(\beta)$, it is clear that $\ptl_n(\beta)$ also has a presentation in terms of diagrams, with the corresponding vector space only a subspace of the full space of connectivities. In preparation for the proposition below, let us define the parity of a connectivity as the parity of the number of loop segments that cross the periodic boundary condition. For instance, the connectivities $c_1$ and $c_2$ in \eqref{eq:c1c2} are respectively even and odd, and their product $c_1c_2$ is odd. Likewise $\Omega$ and $\Omega^{-1}$ are odd, whereas $\id$ and each $e_j$ are even.

\begin{Proposition}
The algebra $\ptl_n(\beta)$ is the subalgebra of $\atl_n(\beta)$ generated by the set of even connectivities excluding the strictly positive even powers of $\Omega$ and $\Omega^{-1}$.
\end{Proposition}
\proof This proposition defines a diagrammatic algebra and claims that it is equivalent to the algebraic definition \eqref{eq:ptl.def} of $\ptl_n(\beta)$. To prove this claim, we show the inclusion both ways. The first part consists in showing that any word in the generators $e_j$ is represented by an even connectivity. This follows readily by noting that the product of even connectivities is always even.\medskip

For the reverse inclusion, we must show that if $c$ is an even connectivity different from $\Omega^{2k}$ with $k \in \mathbb Z\setminus \{0\}$, then it is either equal to $\id$ or can be written as a product of $e_j$ generators. First, we note that such a connectivity is also an element of $\atl_n(\beta)$. It can thus be written as a word in the generators $e_j$, $\Omega$ and $\Omega^{-1}$. Then either $c = \id$ or else $c$ contains at least one generator $e_j$. In the latter case, by \cref{prop:3.forms}, it can be written either as $c = w$ or $c = \Omega\, w$, with $w$ a word in the generators~$e_j$. Because $c$ is even, only $c=w$ is possible, thus ending the proof.
\eproof

Finally, we note that the algebras $\atl_n(\beta)$  and $\ptl_n(\beta)$ have central elements that can be useful in the study of their representation theory. They are included in \cref{sec:non.isomorphic} for the benefit of interested readers.

\subsection{Standard modules}\label{sec:standard.modules}

The standard modules $\repW_{n,d,z}$ of $\atl_n(\beta)$ are constructed on vector spaces spanned by {\it link states}. A link state is a diagram drawn above a horizontal segment with $n$ marked nodes. These are either connected pairwise by non-intersecting loop segments, or occupied by straight loop segments called {\it defects} that extend to infinity and cannot be overarched. The boundary conditions are again periodic in the horizontal direction, meaning that loop segments can travel to the left and reappear to the right at the same height. We denote by $\setB_{n,d}$ the set of link states with $n$ nodes and $d$ defects, with $0\le d\le n$ and $d \equiv n \textrm{ mod } 2$. For
example, the sets of link states for $n=4$ and $n=5$ are
\begin{subequations}
\begin{alignat}{2}
\setB_{4,0}: \quad
\Big \{\ 
&
\psset{unit=0.8}
\begin{pspicture}[shift=-0.05](0,0)(1.6,0.6)
\psline[linewidth=0.5pt](0,0)(1.6,0)
\psarc[linecolor=blue,linewidth=1.5pt]{-}(0.4,0){0.2}{0}{180}
\psarc[linecolor=blue,linewidth=1.5pt]{-}(1.2,0){0.2}{0}{180}
\end{pspicture}\ ,
\quad
\begin{pspicture}[shift=-0.05](0,0)(1.6,0.6)
\psline[linewidth=0.5pt](0,0)(1.6,0)
\psarc[linecolor=blue,linewidth=1.5pt]{-}(0.8,0){0.2}{0}{180}
\psbezier[linecolor=blue,linewidth=1.5pt]{-}(0.2,0)(0.2,0.5)(1.4,0.5)(1.4,0)
\end{pspicture}\ ,
\quad
\begin{pspicture}[shift=-0.05](0,0)(1.6,0.6)
\psline[linewidth=0.5pt](0,0)(1.6,0)
\psbezier[linecolor=blue,linewidth=1.5pt]{-}(1.62,0.27)(1.5,0.44)(0.6,0.44)(0.6,0)
\psbezier[linecolor=blue,linewidth=1.5pt]{-}(0.2,0)(0.2,0.18)(0.05,0.28)(-0.02,0.3)
\psarc[linecolor=blue,linewidth=1.5pt]{-}(1.2,0){0.2}{0}{180}
\psframe[fillstyle=solid,linecolor=white,linewidth=0pt](-0.05,0)(-0.005,0.4)
\psframe[fillstyle=solid,linecolor=white,linewidth=0pt](1.605,0)(1.65,0.4)
\end{pspicture}\ ,
\quad
\begin{pspicture}[shift=-0.05](0,0)(1.6,0.6)
\psline[linewidth=0.5pt](0,0)(1.6,0)
\psarc[linecolor=blue,linewidth=1.5pt]{-}(0,0){0.2}{0}{90}
\psarc[linecolor=blue,linewidth=1.5pt]{-}(1.6,0){0.2}{90}{180}
\psarc[linecolor=blue,linewidth=1.5pt]{-}(0.8,0){0.2}{0}{180}
\end{pspicture}\ ,
\quad
\begin{pspicture}[shift=-0.05](0,0)(1.6,0.6)
\psline[linewidth=0.5pt](0,0)(1.6,0)
\psbezier[linecolor=blue,linewidth=1.5pt]{-}(-0.02,0.27)(0.1,0.44)(1.0,0.44)(1.0,0)
\psbezier[linecolor=blue,linewidth=1.5pt]{-}(1.4,0)(1.4,0.18)(1.55,0.28)(1.62,0.3)
\psarc[linecolor=blue,linewidth=1.5pt]{-}(0.4,0){0.2}{0}{180}
\psframe[fillstyle=solid,linecolor=white,linewidth=0pt](-0.05,0)(-0.005,0.4)
\psframe[fillstyle=solid,linecolor=white,linewidth=0pt](1.605,0)(1.65,0.4)
\end{pspicture}\ ,
\quad
\begin{pspicture}[shift=-0.05](0,0)(1.6,0.6)
\psline[linewidth=0.5pt](0,0)(1.6,0)
\psarc[linecolor=blue,linewidth=1.5pt]{-}(0,0){0.2}{0}{90}
\psarc[linecolor=blue,linewidth=1.5pt]{-}(1.6,0){0.2}{90}{180}
\psbezier[linecolor=blue,linewidth=1.5pt]{-}(0.6,0)(0.6,0.36)(0.05,0.4)(0,0.4)
\psbezier[linecolor=blue,linewidth=1.5pt]{-}(1.0,0)(1.0,0.36)(1.45,0.4)(1.6,0.4)
\end{pspicture}
\ \Big\},
\\[0.2cm]
\setB_{4,2}: \quad
\Big\{\ 
&
\psset{unit=0.8}
\begin{pspicture}[shift=-0.05](0,0)(1.6,0.6)
\psline[linewidth=0.5pt](0,0)(1.6,0)
\psline[linecolor=blue,linewidth=1.5pt]{-}(0.2,0)(0.2,0.5)
\psline[linecolor=blue,linewidth=1.5pt]{-}(0.6,0)(0.6,0.5)
\psarc[linecolor=blue,linewidth=1.5pt]{-}(1.2,0){0.2}{0}{180}
\end{pspicture}\ , 
\quad
\begin{pspicture}[shift=-0.05](0,0)(1.6,0.6)
\psline[linewidth=0.5pt](0,0)(1.6,0)
\psline[linecolor=blue,linewidth=1.5pt]{-}(0.2,0)(0.2,0.5)
\psline[linecolor=blue,linewidth=1.5pt]{-}(1.4,0)(1.4,0.5)
\psarc[linecolor=blue,linewidth=1.5pt]{-}(0.8,0){0.2}{0}{180}
\end{pspicture}\ , 
\quad
\begin{pspicture}[shift=-0.05](0,0)(1.6,0.6)
\psline[linewidth=0.5pt](0,0)(1.6,0)
\psline[linecolor=blue,linewidth=1.5pt]{-}(1.0,0)(1.0,0.5)
\psline[linecolor=blue,linewidth=1.5pt]{-}(1.4,0)(1.4,0.5)
\psarc[linecolor=blue,linewidth=1.5pt]{-}(0.4,0){0.2}{0}{180}
\end{pspicture}\ ,
\quad
\begin{pspicture}[shift=-0.05](0,0)(1.6,0.6)
\psline[linewidth=0.5pt](0,0)(1.6,0)
\psline[linecolor=blue,linewidth=1.5pt]{-}(0.6,0)(0.6,0.5)
\psline[linecolor=blue,linewidth=1.5pt]{-}(1.0,0)(1.0,0.5)
\psarc[linecolor=blue,linewidth=1.5pt]{-}(0.0,0){0.2}{0}{90}
\psarc[linecolor=blue,linewidth=1.5pt]{-}(1.6,0){0.2}{90}{180}
\end{pspicture}\ \Big\},
\\[0.2cm]
\setB_{4,4}: \quad
\Big\{\ 
&
\psset{unit=0.8}
\begin{pspicture}[shift=-0.05](0,0)(1.6,0.6)
\psline[linewidth=0.5pt](0,0)(1.6,0)
\psline[linecolor=blue,linewidth=1.5pt]{-}(0.2,0)(0.2,0.5)
\psline[linecolor=blue,linewidth=1.5pt]{-}(0.6,0)(0.6,0.5)
\psline[linecolor=blue,linewidth=1.5pt]{-}(1.0,0)(1.0,0.5)
\psline[linecolor=blue,linewidth=1.5pt]{-}(1.4,0)(1.4,0.5)
\end{pspicture}\ \Big\} ,
\end{alignat}
\end{subequations}
and
\begin{subequations}
\begin{alignat}{2}
\setB_{5,1}: \quad
\Big \{\ 
&
\psset{unit=0.8}
\begin{pspicture}[shift=-0.05](0,0)(2.0,0.6)
\psline[linewidth=0.5pt](0,0)(2.0,0)
\psline[linecolor=blue,linewidth=1.5pt]{-}(0.2,0)(0.2,0.5)
\psarc[linecolor=blue,linewidth=1.5pt]{-}(0.8,0){0.2}{0}{180}
\psarc[linecolor=blue,linewidth=1.5pt]{-}(1.6,0){0.2}{0}{180}
\end{pspicture}\ ,
\quad
\begin{pspicture}[shift=-0.05](0,0)(2.0,0.6)
\psline[linewidth=0.5pt](0,0)(2.0,0)
\psline[linecolor=blue,linewidth=1.5pt]{-}(0.6,0)(0.6,0.5)
\psarc[linecolor=blue,linewidth=1.5pt]{-}(0,0){0.2}{0}{90}
\psarc[linecolor=blue,linewidth=1.5pt]{-}(1.2,0){0.2}{0}{180}
\psarc[linecolor=blue,linewidth=1.5pt]{-}(2,0){0.2}{90}{180}
\end{pspicture}\ ,
\quad
\begin{pspicture}[shift=-0.05](0,0)(2.0,0.6)
\psline[linewidth=0.5pt](0,0)(2.0,0)
\psline[linecolor=blue,linewidth=1.5pt]{-}(1.0,0)(1.0,0.5)
\psarc[linecolor=blue,linewidth=1.5pt]{-}(0.4,0){0.2}{0}{180}
\psarc[linecolor=blue,linewidth=1.5pt]{-}(1.6,0){0.2}{0}{180}
\end{pspicture}\ ,
\quad
\begin{pspicture}[shift=-0.05](0,0)(2.0,0.6)
\psline[linewidth=0.5pt](0,0)(2.0,0)
\psline[linecolor=blue,linewidth=1.5pt]{-}(1.4,0)(1.4,0.5)
\psarc[linecolor=blue,linewidth=1.5pt]{-}(0,0){0.2}{0}{90}
\psarc[linecolor=blue,linewidth=1.5pt]{-}(0.8,0){0.2}{0}{180}
\psarc[linecolor=blue,linewidth=1.5pt]{-}(2,0){0.2}{90}{180}
\end{pspicture}\ , 
\quad
\begin{pspicture}[shift=-0.05](0,0)(2.0,0.6)
\psline[linewidth=0.5pt](0,0)(2.0,0)
\psline[linecolor=blue,linewidth=1.5pt]{-}(1.8,0)(1.8,0.5)
\psarc[linecolor=blue,linewidth=1.5pt]{-}(0.4,0){0.2}{0}{180}
\psarc[linecolor=blue,linewidth=1.5pt]{-}(1.2,0){0.2}{0}{180}
\end{pspicture}\ , 
\nonumber\\[0.2cm]& 
\psset{unit=0.8}
\begin{pspicture}[shift=-0.05](0,0)(2.0,0.6)
\psline[linewidth=0.5pt](0,0)(2.0,0)
\psline[linecolor=blue,linewidth=1.5pt]{-}(0.2,0)(0.2,0.5)
\psarc[linecolor=blue,linewidth=1.5pt]{-}(1.2,0){0.2}{0}{180}
\psbezier[linecolor=blue,linewidth=1.5pt]{-}(0.6,0)(0.6,0.5)(1.8,0.5)(1.8,0)
\end{pspicture}\ ,
\quad
\begin{pspicture}[shift=-0.05](0,0)(2.0,0.6)
\psline[linewidth=0.5pt](0,0)(2.0,0)
\psbezier[linecolor=blue,linewidth=1.5pt]{-}(2.02,0.27)(1.9,0.44)(1.0,0.44)(1.0,0)
\psbezier[linecolor=blue,linewidth=1.5pt]{-}(0.2,0)(0.2,0.18)(0.05,0.28)(-0.02,0.3)
\psline[linecolor=blue,linewidth=1.5pt]{-}(0.6,0)(0.6,0.5)
\psarc[linecolor=blue,linewidth=1.5pt]{-}(1.6,0){0.2}{0}{180}
\end{pspicture}\ ,
\quad
\begin{pspicture}[shift=-0.05](0,0)(2.0,0.6)
\psline[linewidth=0.5pt](0,0)(2.0,0)
\psline[linecolor=blue,linewidth=1.5pt]{-}(1.0,0)(1.0,0.5)
\psarc[linecolor=blue,linewidth=1.5pt]{-}(0,0){0.2}{0}{90}
\psarc[linecolor=blue,linewidth=1.5pt]{-}(2,0){0.2}{90}{180}
\psbezier[linecolor=blue,linewidth=1.5pt]{-}(0.6,0)(0.6,0.36)(0.05,0.4)(0,0.4)
\psbezier[linecolor=blue,linewidth=1.5pt]{-}(1.4,0)(1.4,0.36)(1.85,0.4)(2,0.4)
\end{pspicture}\ ,
\quad
\begin{pspicture}[shift=-0.05](0,0)(2.0,0.6)
\psline[linewidth=0.5pt](0,0)(2.0,0)
\psbezier[linecolor=blue,linewidth=1.5pt]{-}(-0.02,0.27)(0.1,0.44)(1.0,0.44)(1.0,0)
\psbezier[linecolor=blue,linewidth=1.5pt]{-}(1.8,0)(1.8,0.18)(1.95,0.28)(2.02,0.3)
\psline[linecolor=blue,linewidth=1.5pt]{-}(1.4,0)(1.4,0.5)
\psarc[linecolor=blue,linewidth=1.5pt]{-}(0.4,0){0.2}{0}{180}
\psframe[fillstyle=solid,linecolor=white,linewidth=0pt](-0.05,0)(-0.005,0.4)
\psframe[fillstyle=solid,linecolor=white,linewidth=0pt](2.005,0)(2.05,0.4)
\end{pspicture}\ , 
\quad
\begin{pspicture}[shift=-0.05](0,0)(2.0,0.6)
\psline[linewidth=0.5pt](0,0)(2.0,0)
\psline[linecolor=blue,linewidth=1.5pt]{-}(1.8,0)(1.8,0.5)
\psarc[linecolor=blue,linewidth=1.5pt]{-}(0.8,0){0.2}{0}{180}
\psbezier[linecolor=blue,linewidth=1.5pt]{-}(0.2,0)(0.2,0.5)(1.4,0.5)(1.4,0)
\end{pspicture}
\ \Big\} ,
\\[0.2cm]
\setB_{5,3}: \quad
\Big\{\ 
&
\psset{unit=0.8}
\begin{pspicture}[shift=-0.05](0,0)(2.0,0.6)
\psline[linewidth=0.5pt](0,0)(2.0,0)
\psline[linecolor=blue,linewidth=1.5pt]{-}(0.2,0)(0.2,0.5)
\psline[linecolor=blue,linewidth=1.5pt]{-}(0.6,0)(0.6,0.5)
\psline[linecolor=blue,linewidth=1.5pt]{-}(1.0,0)(1.0,0.5)
\psarc[linecolor=blue,linewidth=1.5pt]{-}(1.6,0){0.2}{0}{180}
\end{pspicture}\ ,
\quad
\begin{pspicture}[shift=-0.05](0,0)(2.0,0.6)
\psline[linewidth=0.5pt](0,0)(2.0,0)
\psline[linecolor=blue,linewidth=1.5pt]{-}(0.2,0)(0.2,0.5)
\psline[linecolor=blue,linewidth=1.5pt]{-}(0.6,0)(0.6,0.5)
\psline[linecolor=blue,linewidth=1.5pt]{-}(1.8,0)(1.8,0.5)
\psarc[linecolor=blue,linewidth=1.5pt]{-}(1.2,0){0.2}{0}{180}
\end{pspicture}\ , 
\quad
\begin{pspicture}[shift=-0.05](0,0)(2.0,0.6)
\psline[linewidth=0.5pt](0,0)(2.0,0)
\psline[linecolor=blue,linewidth=1.5pt]{-}(0.2,0)(0.2,0.5)
\psline[linecolor=blue,linewidth=1.5pt]{-}(1.4,0)(1.4,0.5)
\psline[linecolor=blue,linewidth=1.5pt]{-}(1.8,0)(1.8,0.5)
\psarc[linecolor=blue,linewidth=1.5pt]{-}(0.8,0){0.2}{0}{180}
\end{pspicture}\ , 
\quad
\begin{pspicture}[shift=-0.05](0,0)(2.0,0.6)
\psline[linewidth=0.5pt](0,0)(2.0,0)
\psline[linecolor=blue,linewidth=1.5pt]{-}(1.0,0)(1.0,0.5)
\psline[linecolor=blue,linewidth=1.5pt]{-}(1.4,0)(1.4,0.5)
\psline[linecolor=blue,linewidth=1.5pt]{-}(1.8,0)(1.8,0.5)
\psarc[linecolor=blue,linewidth=1.5pt]{-}(0.4,0){0.2}{0}{180}
\end{pspicture}\ ,
\quad
\begin{pspicture}[shift=-0.05](0,0)(2.0,0.6)
\psline[linewidth=0.5pt](0,0)(2.0,0)
\psline[linecolor=blue,linewidth=1.5pt]{-}(0.6,0)(0.6,0.5)
\psline[linecolor=blue,linewidth=1.5pt]{-}(1.0,0)(1.0,0.5)
\psline[linecolor=blue,linewidth=1.5pt]{-}(1.4,0)(1.4,0.5)
\psarc[linecolor=blue,linewidth=1.5pt]{-}(0.0,0){0.2}{0}{90}
\psarc[linecolor=blue,linewidth=1.5pt]{-}(2.0,0){0.2}{90}{180}
\end{pspicture}\ \Big\} ,
\\[0.2cm] 
\setB_{5,5}: \quad
\Big\{\ 
&
\psset{unit=0.8}
\begin{pspicture}[shift=-0.05](0,0)(2.0,0.6)
\psline[linewidth=0.5pt](0,0)(2.0,0)
\psline[linecolor=blue,linewidth=1.5pt]{-}(0.2,0)(0.2,0.5)
\psline[linecolor=blue,linewidth=1.5pt]{-}(0.6,0)(0.6,0.5)
\psline[linecolor=blue,linewidth=1.5pt]{-}(1.0,0)(1.0,0.5)
\psline[linecolor=blue,linewidth=1.5pt]{-}(1.4,0)(1.4,0.5)
\psline[linecolor=blue,linewidth=1.5pt]{-}(1.8,0)(1.8,0.5)
\end{pspicture}\ \Big\} .
\end{alignat}
\end{subequations}
The cardinality of $\setB_{n,d}$ is
\be
|\setB_{n,d}| = \binom{n}{\frac{n-d}2}.
\ee

The module $\repW_{n,d,z}$ is the free vector space over $\setB_{n,d}$ endowed with an action of $\atl_n(\beta)$. To define this standard action, we first define $c\cdot w$ for $c$ a connectivity of $\atl_n(\beta)$ and $w$ a link state $w \in \setB_{n,d}$, and then extend it linearly on both entries. To compute $c\cdot w$, we draw $w$ above $c$. If two defects are attached together in the resulting diagram, the result is set to zero. Otherwise, a new link state in $\setB_{n,d}$ is read off from the connections of the nodes at the bottom of the diagram. The result of $c\cdot w$ is this link state times a number of prefactors. First, each contractible loop is erased and replaced by a factor of~$\beta$. Second, each non-contractible loop is also erased and replaced by a factor of $\alpha \in \mathbb C$. This may occur only for $d=0$. The parameter $\alpha$ is the fugacity of non-contractible loops. Third, for $d \ge 0$, each defect that crosses the periodic boundary condition is unwound, at the cost of a twist factor $z$ or $z^{-1}$, with $z \in \mathbb C^\times$. This factor is equal to $z$ if the defect crossed the boundary condition towards the left as it propagated downwards, and $z^{-1}$ if it travelled towards the right. Here are examples for $n=4$ and $n=5$:
\be
\begin{array}{lll}
\psset{unit=0.8}
e_1 \cdot \,
\begin{pspicture}[shift=-0.05](0,0)(1.6,0.5)
\psline[linewidth=0.5pt](0,0)(1.6,0)
\psline[linecolor=blue,linewidth=1.5pt]{-}(1.0,0)(1.0,0.5)
\psline[linecolor=blue,linewidth=1.5pt]{-}(1.4,0)(1.4,0.5)
\psarc[linecolor=blue,linewidth=1.5pt]{-}(0.4,0){0.2}{0}{180}
\end{pspicture} 
\ = \beta\ 
\begin{pspicture}[shift=-0.05](0,0)(1.6,0.5)
\psline[linewidth=0.5pt](0,0)(1.6,0)
\psline[linecolor=blue,linewidth=1.5pt]{-}(1.0,0)(1.0,0.5)
\psline[linecolor=blue,linewidth=1.5pt]{-}(1.4,0)(1.4,0.5)
\psarc[linecolor=blue,linewidth=1.5pt]{-}(0.4,0){0.2}{0}{180}
\end{pspicture}\ ,
&\quad&
\psset{unit=0.8}
e_4 \cdot \,
\begin{pspicture}[shift=-0.05](0,0)(2.0,0.5)
\psline[linewidth=0.5pt](0,0)(2.0,0)
\psline[linecolor=blue,linewidth=1.5pt]{-}(1.8,0)(1.8,0.5)
\psarc[linecolor=blue,linewidth=1.5pt]{-}(0.4,0){0.2}{0}{180}
\psarc[linecolor=blue,linewidth=1.5pt]{-}(1.2,0){0.2}{0}{180}
\end{pspicture} 
\ = \ 
\begin{pspicture}[shift=-0.05](0,0)(2.0,0.5)
\psline[linewidth=0.5pt](0,0)(2.0,0)
\psline[linecolor=blue,linewidth=1.5pt]{-}(1.0,0)(1.0,0.5)
\psarc[linecolor=blue,linewidth=1.5pt]{-}(0.4,0){0.2}{0}{180}
\psarc[linecolor=blue,linewidth=1.5pt]{-}(1.6,0){0.2}{0}{180}
\end{pspicture} \ ,
\\[0.3cm]
\psset{unit=0.8}
e_0 \cdot \,
\begin{pspicture}[shift=-0.05](0,0)(1.6,0.5)
\psline[linewidth=0.5pt](0,0)(1.6,0)
\psarc[linecolor=blue,linewidth=1.5pt]{-}(0.8,0){0.2}{0}{180}
\psbezier[linecolor=blue,linewidth=1.5pt]{-}(0.2,0)(0.2,0.5)(1.4,0.5)(1.4,0)
\end{pspicture} 
\ = \alpha \ \,
\begin{pspicture}[shift=-0.05](0,0)(1.6,0.5)
\psline[linewidth=0.5pt](0,0)(1.6,0)
\psarc[linecolor=blue,linewidth=1.5pt]{-}(0.8,0){0.2}{0}{180}
\psarc[linecolor=blue,linewidth=1.5pt]{-}(0.0,0){0.2}{0}{90}
\psarc[linecolor=blue,linewidth=1.5pt]{-}(1.6,0){0.2}{90}{180}
\end{pspicture}\ ,
&\quad&
\psset{unit=0.8}
\Omega \cdot \,\,
\begin{pspicture}[shift=-0.05](0,0)(2.0,0.5)
\psline[linewidth=0.5pt](0,0)(2.0,0)
\psline[linewidth=0.5pt](0,0)(2.0,0)
\psline[linecolor=blue,linewidth=1.5pt]{-}(0.2,0)(0.2,0.5)
\psarc[linecolor=blue,linewidth=1.5pt]{-}(1.2,0){0.2}{0}{180}
\psbezier[linecolor=blue,linewidth=1.5pt]{-}(0.6,0)(0.6,0.5)(1.8,0.5)(1.8,0)
\end{pspicture} 
\ = z\ 
\begin{pspicture}[shift=-0.05](0,0)(2.0,0.5)
\psline[linewidth=0.5pt](0,0)(2.0,0)
\psline[linecolor=blue,linewidth=1.5pt]{-}(1.8,0)(1.8,0.5)
\psarc[linecolor=blue,linewidth=1.5pt]{-}(0.8,0){0.2}{0}{180}
\psbezier[linecolor=blue,linewidth=1.5pt]{-}(0.2,0)(0.2,0.5)(1.4,0.5)(1.4,0)
\end{pspicture}\ ,
\\[0.3cm]
\psset{unit=0.8}
e_3 \cdot \,
\begin{pspicture}[shift=-0.05](0,0)(1.6,0.5)
\psline[linewidth=0.5pt](0,0)(1.6,0)
\psline[linecolor=blue,linewidth=1.5pt]{-}(0.6,0)(0.6,0.5)
\psline[linecolor=blue,linewidth=1.5pt]{-}(1.0,0)(1.0,0.5)
\psarc[linecolor=blue,linewidth=1.5pt]{-}(0.0,0){0.2}{0}{90}
\psarc[linecolor=blue,linewidth=1.5pt]{-}(1.6,0){0.2}{90}{180}
\end{pspicture}
\ = z^{-1}\ 
\begin{pspicture}[shift=-0.05](0,0)(1.6,0.5)
\psline[linewidth=0.5pt](0,0)(1.6,0)
\psline[linecolor=blue,linewidth=1.5pt]{-}(0.2,0)(0.2,0.5)
\psline[linecolor=blue,linewidth=1.5pt]{-}(0.6,0)(0.6,0.5)
\psarc[linecolor=blue,linewidth=1.5pt]{-}(1.2,0){0.2}{0}{180}
\end{pspicture} \ ,
&\quad&
\psset{unit=0.8}
e_3 \cdot \,
\begin{pspicture}[shift=-0.05](0,0)(2.0,0.5)
\psline[linewidth=0.5pt](0,0)(2.0,0)
\psline[linecolor=blue,linewidth=1.5pt]{-}(1.0,0)(1.0,0.5)
\psline[linecolor=blue,linewidth=1.5pt]{-}(1.4,0)(1.4,0.5)
\psline[linecolor=blue,linewidth=1.5pt]{-}(1.8,0)(1.8,0.5)
\psarc[linecolor=blue,linewidth=1.5pt]{-}(0.4,0){0.2}{0}{180}
\end{pspicture}
\ = 0\, .
\end{array}
\ee
Each standard module thus depends on one complex parameter: $\alpha$ for $d = 0$ and $z$ for $d\ge 1$. It is then customary to parameterize the weight of non-contractible loops as $\alpha = z+z^{-1}$, and the corresponding module as $\repW_{n,0,z}$, thus providing a unified notation for all the standard modules. The set of standard modules over $\atl_n(\beta)$ is
\be
\label{eq:atl.standard}
\atl_n(\beta): 
\left\{\begin{array}{ll}
\big\{\repW_{n,d,z} \ | \ z \in \mathbb C^\times, \ d\in \{1,3,\dots, n\}\big\} & n \textrm{ odd,} \\[0.15cm]
\{\repW_{n,0,r \eE^{\iI \theta}}\ | \ r>0, \ 0 \le \theta < \pi\}\ \cup \  \big\{\repW_{n,d,z} \ | \ z \in \mathbb C^\times, \ d \in \{2,4,\dots, n\}\big\} & n \textrm{ even.}
\end{array}\right.
\ee 
For $n$ even, we used the property $\repW_{n,0,z} = \repW_{n,0,z^{-1}}$ to restrict the parameter $z$ to half of the complex plane. These modules satisfy the following property.
\begin{Proposition}
\label{prop:atl.iso}
Over $\atl_n(\beta)$, the standard modules in \eqref{eq:atl.standard} are all pairwise non-isomorphic.
\end{Proposition}
\proof
The dimension of $\repW_{n,d,z}$ is strictly decreasing with $d$, so two standard modules $\repW_{n,d,z}$ and $\repW_{n,d',z}$ with $d\neq d'$ are never isomorphic. Moreover, Graham and Lehrer \cite[Theorem 2.8]{GL98} have shown that the only possible homomorphisms are between standard modules with different defect numbers.
\eproof

Over $\ptl_n(\beta)$, the modules $\repW_{n,d,z}$ are not all non-isomorphic. First, because each $e_j$ is represented by the zero matrix of size one in $\repW_{n,n,z}$, it follows that $\repW_{n,n,y}\simeq \repW_{n,n,z}$ over $\ptl_n(\beta)$, for all $y,z \in \mathbb C^\times$. The next proposition shows that some pairs of standard modules for $d<n$ are also isomorphic.

\begin{Lemma} 
\label{lem:iso.minus.z}
Over $\ptl_n(\beta)$, the modules $\repW_{n,d,z}$ and $\repW_{n,d,-z}$ are isomorphic.
\end{Lemma}
\proof Let us define the parity $\sigma_v$ of a link state $v$ to equal $0$ if an even number of its arcs cross the periodic boundary condition, and $1$ if this number is odd. For $d >0$, if the action of $e_j$ on a link state changes the parity, then it produces a factor of $z$ or $z^{-1}$. One can then perform a diagonal change of basis in $\repW_{n,d,z}$ that multiplies each link state $v$ by its parity $(-1)^{\sigma_v}$. The resulting action of the algebra is then precisely the one for $\repW_{n,d,-z}$, thus showing that $\repW_{n,d,z}$ and $\repW_{n,d,-z}$ are indeed isomorphic. The proof for $d = 0$ works in the same way, with the new action after the change of basis obtained by changing $\alpha \to -\alpha$.\eproof

\noindent We also note that the same proof does not hold for $\atl_n(\beta)$, because $\Omega$ and $\Omega^{-1}$ may produce factors of $z$ or $z^{-1}$ without changing the parity of the link states.

\begin{Proposition} Over $\ptl_n(\beta)$, the modules $\repW_{n,d,y}$ and $\repW_{n,d,z}$ with $0<d \le n-2$ are isomorphic if and only if $y^2 = z^2$. For $d=0$, the modules $\repW_{n,0,y}$ and $\repW_{n,0,z}$ are isomorphic over $\ptl_n(\beta)$ if and only if $y^2 \in \{z^2,z^{-2}\}$.
\end{Proposition}
\proof
\cref{lem:iso.minus.z} showed that $\repW_{n,d,y} \simeq \repW_{n,d,z}$ for $0<d \le n-2$ with $y^2 = z^2$. For $d=0$, because $\repW_{n,0,z}=\repW_{n,0,z^{-1}}$, it readily follows from the same lemma that $\repW_{n,0,y} \simeq \repW_{n,0,z}$ for $y^2 \in \{z^2,z^{-2}\}$. It thus only remains to show that $\repW_{n,d,y}$ is not isomorphic to $\repW_{n,d,z}$ for the other values of $y$. Let us define the connectivities
\be
a_d = 
\left\{
\begin{array}{ll}
\ \begin{pspicture}[shift=-0.7](0,-0.8)(5.6,1.3) 
\pspolygon[fillstyle=solid,fillcolor=lightlightblue,linewidth=0pt](0,-0.7)(5.6,-0.7)(5.6,0.7)(0,0.7)
\psarc[linecolor=blue,linewidth=1.5pt]{-}(0.4,0.7){0.2}{180}{0}
\rput(1.2,0.6){$...$}
\psarc[linecolor=blue,linewidth=1.5pt]{-}(2.0,0.7){0.2}{180}{0}
\psarc[linecolor=blue,linewidth=1.5pt]{-}(0.8,-0.7){0.2}{0}{180}
\rput(1.6,-0.6){$...$}
\psarc[linecolor=blue,linewidth=1.5pt]{-}(2.4,-0.7){0.2}{0}{180}
\rput(5.6,0){
\psbezier[linecolor=blue,linewidth=1.5pt]{-}(-0.6,-0.7)(-0.6,-0.29)(0.02,-0.2)(0.02,-0.13)
\psbezier[linecolor=blue,linewidth=1.5pt]{-}(-0.2,-0.7)(-0.2,-0.5)(0.02,-0.4)(0.02,-0.37)
}
\psbezier[linecolor=blue,linewidth=1.5pt]{-}(-0.02,-0.17)(0.5,0.09)(2.6,0.15)(2.6,0.7)
\psbezier[linecolor=blue,linewidth=1.5pt]{-}(-0.02,-0.42)(0.5,0.03)(3.0,0.03)(3.0,0.7)
\psbezier[linecolor=blue,linewidth=1.5pt]{-}(0.2,-0.7)(0.2,0)(3.4,-0.2)(3.4,0.7)
\psbezier[linecolor=blue,linewidth=1.5pt]{-}(3.0,-0.7)(3.0,0)(3.8,0)(3.8,0.7)
\psbezier[linecolor=blue,linewidth=1.5pt]{-}(3.4,-0.7)(3.4,0)(4.2,0)(4.2,0.7)
\rput(4.25,0){$...$}
\psbezier[linecolor=blue,linewidth=1.5pt]{-}(4.2,-0.7)(4.2,0)(5.0,0)(5.0,0.7)
\psbezier[linecolor=blue,linewidth=1.5pt]{-}(4.6,-0.7)(4.6,0)(5.4,0)(5.4,0.7)
\psframe[fillstyle=solid,linecolor=white,linewidth=0pt](-0.05,-0.71)(-0.005,0.71)
\psframe[fillstyle=solid,linecolor=white,linewidth=0pt](5.605,-0.71)(5.65,0.71)
\rput(4.0,1.1){$\overbrace{\ \hspace{2.75cm}\ }^{d}$} 
\end{pspicture}
\quad & 3 \le d \le n-2,\\[0.8cm]
\ \begin{pspicture}[shift=-0.7](0,-0.8)(5.6,0.8) 
\pspolygon[fillstyle=solid,fillcolor=lightlightblue,linewidth=0pt](0,-0.7)(5.6,-0.7)(5.6,0.7)(0,0.7)
\psarc[linecolor=blue,linewidth=1.5pt]{-}(0.4,0.7){0.2}{180}{0}
\psarc[linecolor=blue,linewidth=1.5pt]{-}(1.2,0.7){0.2}{180}{0}
\psarc[linecolor=blue,linewidth=1.5pt]{-}(2.0,0.7){0.2}{180}{0}
\rput(2.8,0.6){$...$}
\psarc[linecolor=blue,linewidth=1.5pt]{-}(3.6,0.7){0.2}{180}{0}
\psarc[linecolor=blue,linewidth=1.5pt]{-}(4.4,0.7){0.2}{180}{0}
\psarc[linecolor=blue,linewidth=1.5pt]{-}(0.8,-0.7){0.2}{0}{180}
\psarc[linecolor=blue,linewidth=1.5pt]{-}(1.6,-0.7){0.2}{0}{180}
\psarc[linecolor=blue,linewidth=1.5pt]{-}(2.4,-0.7){0.2}{0}{180}
\rput(3.2,-0.6){$...$}
\psarc[linecolor=blue,linewidth=1.5pt]{-}(4.0,-0.7){0.2}{0}{180}
\psarc[linecolor=blue,linewidth=1.5pt]{-}(4.8,-0.7){0.2}{0}{180}
\psarc[linecolor=blue,linewidth=1.5pt]{-}(2.4,-0.7){0.2}{0}{180}
\psbezier[linecolor=blue,linewidth=1.5pt]{-}(5.4,0.7)(5.4,-0.3)(-0.2,0.2)(-0.2,-0.7)
\psbezier[linecolor=blue,linewidth=1.5pt]{-}(5.0,0.7)(5.0,0)(-0.6,0)(0.03,0)
\psbezier[linecolor=blue,linewidth=1.5pt]{-}(5.6,0)(5.0,0)(0.2,-0.12)(0.2,-0.7)
\psbezier[linecolor=blue,linewidth=1.5pt]{-}(5.4,-0.7)(5.4,-0.55)(5.5,-0.45)(5.6,-0.4)
\rput(5.6,0){
}
\psframe[fillstyle=solid,linecolor=white,linewidth=0pt](-0.3,-0.71)(-0.005,0.71)
\psframe[fillstyle=solid,linecolor=white,linewidth=0pt](5.605,-0.71)(5.65,0.71)
\end{pspicture}
& d=2,
\\[0.8cm]
\ \begin{pspicture}[shift=-0.7](0,-0.8)(5.6,0.8) 
\pspolygon[fillstyle=solid,fillcolor=lightlightblue,linewidth=0pt](0,-0.7)(5.6,-0.7)(5.6,0.7)(0,0.7)
\psarc[linecolor=blue,linewidth=1.5pt]{-}(0.4,0.7){0.2}{180}{0}
\psarc[linecolor=blue,linewidth=1.5pt]{-}(1.2,0.7){0.2}{180}{0}
\psarc[linecolor=blue,linewidth=1.5pt]{-}(2.0,0.7){0.2}{180}{0}
\rput(3.0,0.6){$...$}
\psarc[linecolor=blue,linewidth=1.5pt]{-}(4.0,0.7){0.2}{180}{0}
\psarc[linecolor=blue,linewidth=1.5pt]{-}(4.8,0.7){0.2}{180}{0}
\psarc[linecolor=blue,linewidth=1.5pt]{-}(0.8,-0.7){0.2}{0}{180}
\psarc[linecolor=blue,linewidth=1.5pt]{-}(1.6,-0.7){0.2}{0}{180}
\psarc[linecolor=blue,linewidth=1.5pt]{-}(2.4,-0.7){0.2}{0}{180}
\rput(3.4,-0.6){$...$}
\psarc[linecolor=blue,linewidth=1.5pt]{-}(4.4,-0.7){0.2}{0}{180}
\psarc[linecolor=blue,linewidth=1.5pt]{-}(5.2,-0.7){0.2}{0}{180}
\psarc[linecolor=blue,linewidth=1.5pt]{-}(2.4,-0.7){0.2}{0}{180}
\psbezier[linecolor=blue,linewidth=1.5pt]{-}(5.4,0.7)(5.4,0)(-0.6,0.1)(0.03,0.1)
\psbezier[linecolor=blue,linewidth=1.5pt]{-}(5.6,0.1)(5.4,0.1)(0.2,-0.1)(0,-0.1)
\psbezier[linecolor=blue,linewidth=1.5pt]{-}(5.6,-0.1)(5.4,-0.1)(0.2,0)(0.2,-0.7)
\rput(5.6,0){
}
\psframe[fillstyle=solid,linecolor=white,linewidth=0pt](-0.3,-0.71)(-0.005,0.71)
\psframe[fillstyle=solid,linecolor=white,linewidth=0pt](5.605,-0.71)(5.65,0.71)
\end{pspicture}
& d=1,\\[.8cm]
\ \begin{pspicture}[shift=-0.7](0,-0.8)(5.6,0.8)
\pspolygon[fillstyle=solid,fillcolor=lightlightblue,linewidth=0pt](0,-0.7)(5.6,-0.7)(5.6,0.7)(0,0.7)
\multiput(0,0)(0.8,0){3}{\psarc[linecolor=blue,linewidth=1.5pt]{-}(0.4,0.7){0.2}{180}{0}}
\multiput(0,0)(0.8,0){3}{\psarc[linecolor=blue,linewidth=1.5pt]{-}(3.6,0.7){0.2}{180}{0}}
\rput(2.8,0.6){$...$}
\multiput(0,0)(0.8,0){3}{\psarc[linecolor=blue,linewidth=1.5pt]{-}(0.8,-0.7){0.2}{0}{180}}
\multiput(0,0)(0.8,0){2}{\psarc[linecolor=blue,linewidth=1.5pt]{-}(4.0,-0.7){0.2}{0}{180}}
\rput(3.25,-0.6){$...$}
\psline[linecolor=blue,linewidth=1.5pt]{-}(0,0)(5.6,0)
\psarc[linecolor=blue,linewidth=1.5pt]{-}(0,-0.7){0.2}{0}{90}
\psarc[linecolor=blue,linewidth=1.5pt]{-}(5.6,-0.7){0.2}{90}{180}
\end{pspicture}
\quad & d = 0.
\end{array}
\right.
\ee
The connectivity $a_d$ has $\frac{n-d}2$ simple arcs on both the bottom and top segments of the box, as well as $d$~through-lines. For $2 \le d \le n-2$, two of the through-lines cross the periodic boundary condition. For $d=1$, the unique through-line crosses the periodic boundary condition twice. For $d=0$, the connectivity contains one non-contractible loop. Clearly, $a_d$ is an even connectivity and is an element of $\ptl_n(\beta)$ for $0 \le d\le n-2$. It is not hard to see that, of all the link states in the basis $\setB_{n,d}$, 
the state 
\be
v = \ 
\left\{\begin{array}{ll}
\begin{pspicture}[shift=-0.05](0,-0.1)(5.6,0.7)
\psline[linewidth=0.5pt](0,0)(5.6,0)
\psline[linecolor=blue,linewidth=1.5pt]{-}(0.2,0)(0.2,0.5)
\psarc[linecolor=blue,linewidth=1.5pt]{-}(0.8,0){0.2}{0}{180}
\rput(1.6,0.2){$...$}
\psarc[linecolor=blue,linewidth=1.5pt]{-}(2.4,0){0.2}{0}{180}
\multiput(0,0)(0.4,0){3}{\psline[linecolor=blue,linewidth=1.5pt]{-}(3.0,0)(3.0,0.5)}
\multiput(0,0)(0.4,0){3}{\psline[linecolor=blue,linewidth=1.5pt]{-}(4.6,0)(4.6,0.5)}
\rput(4.2,.2){$...$}
\end{pspicture}
\quad & 0 < d \le n-2,
\\[0.2cm]
\begin{pspicture}[shift=-0.05](0,-0.1)(5.6,0.5)
\psline[linewidth=0.5pt](0,0)(5.6,0)
\psarc[linecolor=blue,linewidth=1.5pt]{-}(0,0){0.2}{0}{90}
\psarc[linecolor=blue,linewidth=1.5pt]{-}(5.6,0){0.2}{90}{180}
\psarc[linecolor=blue,linewidth=1.5pt]{-}(0.8,0){0.2}{0}{180}
\psarc[linecolor=blue,linewidth=1.5pt]{-}(1.6,0){0.2}{0}{180}
\psarc[linecolor=blue,linewidth=1.5pt]{-}(2.4,0){0.2}{0}{180}
\multiput(0,0)(0.8,0){3}{\psarc[linecolor=blue,linewidth=1.5pt]{-}(0.8,0){0.2}{0}{180}}
\multiput(0,0)(0.8,0){2}{\psarc[linecolor=blue,linewidth=1.5pt]{-}(4.0,0){0.2}{0}{180}}
\rput(3.2,.2){$...$}
\end{pspicture}
\quad & d=0,
\end{array}\right.
\ee
is the only one that is an eigenstate of $a_d$. Its eigenvalue is $\alpha^2$ for $d=0$ and $z^2$ for $d>0$. It then follows that 
\be
\textrm{tr}_{\raisebox{-0.05cm}{\tiny$\repW_{n,d,z}$}}a_d = 
\left\{\begin{array}{cl}
(z+z^{-1})^2 & d = 0,
\\[0.15cm]
z^2 & d>0.
\end{array}\right.
\ee
The trace over a module of any element of the algebra is invariant under isomorphisms. As a result, $\repW_{n,d,y} \simeq \repW_{n,d,z}$ is only possible if $y^2 \in \{z^2,z^{-2}\}$ for $d=0$, and $y^2 = z^2$ for $d>0$.\eproof

Thus it suffices for $\ptl_n(\beta)$ to consider standard modules $\repW_{n,d,z}$ with the parameter $z$ covering only half of the complex plane for $0 \le d \le n-2$, and one quadrant for $d=0$. The set of standard modules over $\ptl_n(\beta)$ is
\be 
\label{eq:standards.ptl}
\ptl_n(\beta): \quad 
\left\{\begin{array}{ll}
\big\{\repW_{n,n,1}\big\} \ \cup \ \big\{\repW_{n,d,r \eE^{\iI \theta}} \ | \ r>0, \ 0 \le \theta < \pi,  \ d \in \{1,3, \dots, n-2\}\big\}
  & n \textrm{ odd}, \\[0.25cm]
\hspace{-0.2cm}\begin{array}{l}
\big\{\repW_{n,n,1}\big\} \ \cup \ \big\{ \repW_{n,0,r \eE^{\iI \theta}} \ | \ r>0, \ 0 \le \theta < \frac{\pi}2\big\} \\[0.15cm] \ \cup \ 
\big\{\repW_{n,d,r \eE^{\iI \theta}} \ | \ r>0, \ 0 \le \theta < \pi,  \ d \in \{2,4, \dots, n-2\}\big\}
\end{array}
& n \textrm{ even}.
\end{array}\right.
\ee

\paragraph{The modules $\boldsymbol{\wh\repW_{n,d,z}}$ for $\boldsymbol{\ptl_n(\beta)}$.}

For the description in \cref{sec:sandwich.cell} of the cellular theory of the uncoiled algebras, it is useful to define another set of modules over $\ptl_n(\beta)$, denoted $\wh\repW_{n,d,z}$. This module is spanned by link states in a basis $\wh\setB_{n,d}$ that consists only of even diagrams. This basis includes all the even link states of $\setB_{n,d}$. For $d=0$, it also includes each odd link state of $\setB_{n,0}$, to which we add one non-contractible loop above the arches, thus making it an even diagram. For $d>0$, $\wh\setB_{n,d}$ also includes each odd link state of $\setB_{n,d}$, transformed into an even diagram by adding a winding to the defects, in such a way that exactly one defect crosses the periodic boundary condition following the connectivity of $\Omega$ on $d$ strands. For example, for $n=4$, the link states in the bases $\wh\setB_{4,d}$ are
\begin{subequations}
\begin{alignat}{2}
\wh\setB_{4,0}: \quad
\Big \{\ 
&
\psset{unit=0.8}
\begin{pspicture}[shift=-0.05](0,0)(1.6,0.6)
\psline[linewidth=0.5pt](0,0)(1.6,0)
\psarc[linecolor=blue,linewidth=1.5pt]{-}(0.4,0){0.2}{0}{180}
\psarc[linecolor=blue,linewidth=1.5pt]{-}(1.2,0){0.2}{0}{180}
\end{pspicture}\ ,
\quad
\begin{pspicture}[shift=-0.05](0,0)(1.6,0.6)
\psline[linewidth=0.5pt](0,0)(1.6,0)
\psarc[linecolor=blue,linewidth=1.5pt]{-}(0.8,0){0.2}{0}{180}
\psbezier[linecolor=blue,linewidth=1.5pt]{-}(0.2,0)(0.2,0.5)(1.4,0.5)(1.4,0)
\end{pspicture}\ ,
\quad
\begin{pspicture}[shift=-0.05](0,0)(1.6,0.6)
\psline[linewidth=0.5pt](0,0)(1.6,0)
\psbezier[linecolor=blue,linewidth=1.5pt]{-}(1.62,0.27)(1.5,0.44)(0.6,0.44)(0.6,0)
\psbezier[linecolor=blue,linewidth=1.5pt]{-}(0.2,0)(0.2,0.18)(0.05,0.28)(-0.02,0.3)
\psarc[linecolor=blue,linewidth=1.5pt]{-}(1.2,0){0.2}{0}{180}
\psline[linecolor=blue,linewidth=1.5pt]{-}(0,0.6)(1.6,0.6)
\psframe[fillstyle=solid,linecolor=white,linewidth=0pt](-0.05,0)(-0.005,0.4)
\psframe[fillstyle=solid,linecolor=white,linewidth=0pt](1.605,0)(1.65,0.4)
\end{pspicture}\ ,
\quad
\begin{pspicture}[shift=-0.05](0,0)(1.6,0.6)
\psline[linewidth=0.5pt](0,0)(1.6,0)
\psarc[linecolor=blue,linewidth=1.5pt]{-}(0,0){0.2}{0}{90}
\psarc[linecolor=blue,linewidth=1.5pt]{-}(1.6,0){0.2}{90}{180}
\psarc[linecolor=blue,linewidth=1.5pt]{-}(0.8,0){0.2}{0}{180}
\psline[linecolor=blue,linewidth=1.5pt]{-}(0,0.6)(1.6,0.6)
\end{pspicture}\ ,
\quad
\begin{pspicture}[shift=-0.05](0,0)(1.6,0.6)
\psline[linewidth=0.5pt](0,0)(1.6,0)
\psbezier[linecolor=blue,linewidth=1.5pt]{-}(-0.02,0.27)(0.1,0.44)(1.0,0.44)(1.0,0)
\psbezier[linecolor=blue,linewidth=1.5pt]{-}(1.4,0)(1.4,0.18)(1.55,0.28)(1.62,0.3)
\psarc[linecolor=blue,linewidth=1.5pt]{-}(0.4,0){0.2}{0}{180}
\psline[linecolor=blue,linewidth=1.5pt]{-}(0,0.6)(1.6,0.6)
\psframe[fillstyle=solid,linecolor=white,linewidth=0pt](-0.05,0)(-0.005,0.4)
\psframe[fillstyle=solid,linecolor=white,linewidth=0pt](1.605,0)(1.65,0.4)
\end{pspicture}\ ,
\quad
\begin{pspicture}[shift=-0.05](0,0)(1.6,0.6)
\psline[linewidth=0.5pt](0,0)(1.6,0)
\psarc[linecolor=blue,linewidth=1.5pt]{-}(0,0){0.2}{0}{90}
\psarc[linecolor=blue,linewidth=1.5pt]{-}(1.6,0){0.2}{90}{180}
\psbezier[linecolor=blue,linewidth=1.5pt]{-}(0.6,0)(0.6,0.36)(0.05,0.4)(0,0.4)
\psbezier[linecolor=blue,linewidth=1.5pt]{-}(1.0,0)(1.0,0.36)(1.45,0.4)(1.6,0.4)
\end{pspicture}
\ \Big\},
\\[0.2cm]
\wh\setB_{4,2}: \quad
\Big\{\ 
&
\psset{unit=0.8}
\begin{pspicture}[shift=-0.05](0,0)(1.6,0.6)
\psline[linewidth=0.5pt](0,0)(1.6,0)
\psline[linecolor=blue,linewidth=1.5pt]{-}(0.2,0)(0.2,0.5)
\psline[linecolor=blue,linewidth=1.5pt]{-}(0.6,0)(0.6,0.5)
\psarc[linecolor=blue,linewidth=1.5pt]{-}(1.2,0){0.2}{0}{180}
\end{pspicture}\ , 
\quad
\begin{pspicture}[shift=-0.05](0,0)(1.6,0.6)
\psline[linewidth=0.5pt](0,0)(1.6,0)
\psline[linecolor=blue,linewidth=1.5pt]{-}(0.2,0)(0.2,0.5)
\psline[linecolor=blue,linewidth=1.5pt]{-}(1.4,0)(1.4,0.5)
\psarc[linecolor=blue,linewidth=1.5pt]{-}(0.8,0){0.2}{0}{180}
\end{pspicture}\ , 
\quad
\begin{pspicture}[shift=-0.05](0,0)(1.6,0.6)
\psline[linewidth=0.5pt](0,0)(1.6,0)
\psline[linecolor=blue,linewidth=1.5pt]{-}(1.0,0)(1.0,0.5)
\psline[linecolor=blue,linewidth=1.5pt]{-}(1.4,0)(1.4,0.5)
\psarc[linecolor=blue,linewidth=1.5pt]{-}(0.4,0){0.2}{0}{180}
\end{pspicture}\ ,
\quad 
\begin{pspicture}[shift=-0.05](0,0)(1.6,0.6)
\psline[linewidth=0.5pt](0,0)(1.6,0)
\psbezier[linecolor=blue,linewidth=1.5pt]{-}(0.6,0)(0.6,0.25)(1.0,0.25)(1.0,0.55)
\psbezier[linecolor=blue,linewidth=1.5pt]{-}(1.0,0)(1.0,0.35)(1.5,0.35)(1.6,0.35)
\psbezier[linecolor=blue,linewidth=1.5pt]{-}(0,0.35)(0.5,0.35)(0.6,0.45)(0.6,0.55)
\psarc[linecolor=blue,linewidth=1.5pt]{-}(0.0,0){0.2}{0}{90}
\psarc[linecolor=blue,linewidth=1.5pt]{-}(1.6,0){0.2}{90}{180}
\end{pspicture}\ \Big\},
\\[0.2cm]
\wh\setB_{4,4}: \quad
\Big\{\ 
&
\psset{unit=0.8}
\begin{pspicture}[shift=-0.05](0,0)(1.6,0.6)
\psline[linewidth=0.5pt](0,0)(1.6,0)
\psline[linecolor=blue,linewidth=1.5pt]{-}(0.2,0)(0.2,0.5)
\psline[linecolor=blue,linewidth=1.5pt]{-}(0.6,0)(0.6,0.5)
\psline[linecolor=blue,linewidth=1.5pt]{-}(1.0,0)(1.0,0.5)
\psline[linecolor=blue,linewidth=1.5pt]{-}(1.4,0)(1.4,0.5)
\end{pspicture}\ \Big\}.
\end{alignat}
\end{subequations}

The action of $\ptl_n(\beta)$ on $\wh\repW_{n,d,z}$ is similar to the action on $\repW_{n,d,z}$, namely one draws the connectivity below the link state, reads the new link state at the bottom, and multiplies by factors to account for the non-contractible loops and the winding of the defects. The only differences are that (i) contractible loops are only removed in pairs, and (ii) defects are only unwound across the periodic boundary condition in pairs. It is thus clear that the action of $\ptl_n(\beta)$ on $\wh\repW_{n,d,z}$ depends on the parameter $\alpha$ only through factors of $\alpha^2$ for $d=0$, and on the parameter $z$ only through factors of $z^2$ for $d>0$. Moreover, we note that $\wh\repW_{n,d,z}$ is not a module over $\atl_n(\beta)$, as the action of $\Omega$ is not defined in this module.\medskip

There exists a diagonal change of basis from $\repW_{n,d,z}$ to $\wh\repW_{n,d,z}$, except in the case $d = 0$ with $\alpha = 0$. The associated weight is $1$ for the even link states of $\repW_{n,d,z}$, whereas it is either $\alpha^{-1}$ or $z^{-1}$ for the odd link states. For instance, for $(n,d) = (4,0)$ and $(4,2)$, we respectively have
\be
\psset{unit=0.8}
\begin{pspicture}[shift=-0.05](0,0)(1.6,0.6)
\psline[linewidth=0.5pt](0,0)(1.6,0)
\psbezier[linecolor=blue,linewidth=1.5pt]{-}(1.62,0.27)(1.5,0.44)(0.6,0.44)(0.6,0)
\psbezier[linecolor=blue,linewidth=1.5pt]{-}(0.2,0)(0.2,0.18)(0.05,0.28)(-0.02,0.3)
\psarc[linecolor=blue,linewidth=1.5pt]{-}(1.2,0){0.2}{0}{180}
\psframe[fillstyle=solid,linecolor=white,linewidth=0pt](-0.05,0)(-0.005,0.4)
\psframe[fillstyle=solid,linecolor=white,linewidth=0pt](1.605,0)(1.65,0.4)
\end{pspicture}
\ \mapsto \frac 1\alpha \
\begin{pspicture}[shift=-0.05](0,0)(1.6,0.6)
\psline[linewidth=0.5pt](0,0)(1.6,0)
\psbezier[linecolor=blue,linewidth=1.5pt]{-}(1.62,0.27)(1.5,0.44)(0.6,0.44)(0.6,0)
\psbezier[linecolor=blue,linewidth=1.5pt]{-}(0.2,0)(0.2,0.18)(0.05,0.28)(-0.02,0.3)
\psarc[linecolor=blue,linewidth=1.5pt]{-}(1.2,0){0.2}{0}{180}
\psline[linecolor=blue,linewidth=1.5pt]{-}(0,0.6)(1.6,0.6)
\psframe[fillstyle=solid,linecolor=white,linewidth=0pt](-0.05,0)(-0.005,0.4)
\psframe[fillstyle=solid,linecolor=white,linewidth=0pt](1.605,0)(1.65,0.4)
\end{pspicture}\ ,
\qquad
\begin{pspicture}[shift=-0.05](0,0)(1.6,0.6)
\psline[linewidth=0.5pt](0,0)(1.6,0)
\psline[linecolor=blue,linewidth=1.5pt]{-}(0.6,0)(0.6,0.5)
\psline[linecolor=blue,linewidth=1.5pt]{-}(1.0,0)(1.0,0.5)
\psarc[linecolor=blue,linewidth=1.5pt]{-}(0.0,0){0.2}{0}{90}
\psarc[linecolor=blue,linewidth=1.5pt]{-}(1.6,0){0.2}{90}{180}
\end{pspicture}
\ \mapsto \frac 1z \
\begin{pspicture}[shift=-0.05](0,0)(1.6,0.6)
\psline[linewidth=0.5pt](0,0)(1.6,0)
\psbezier[linecolor=blue,linewidth=1.5pt]{-}(0.6,0)(0.6,0.25)(1.0,0.25)(1.0,0.55)
\psbezier[linecolor=blue,linewidth=1.5pt]{-}(1.0,0)(1.0,0.35)(1.5,0.35)(1.6,0.35)
\psbezier[linecolor=blue,linewidth=1.5pt]{-}(0,0.35)(0.5,0.35)(0.6,0.45)(0.6,0.55)
\psarc[linecolor=blue,linewidth=1.5pt]{-}(0.0,0){0.2}{0}{90}
\psarc[linecolor=blue,linewidth=1.5pt]{-}(1.6,0){0.2}{90}{180}
\end{pspicture}\ .
\ee
It is easy to check that this transformation is an isomorphism of $\ptl_n(\beta)$-modules, leading to the following result.

\begin{Proposition}
Let $d>0$ with $z \in \mathbb C^\times$, or $d=0$ with $\alpha \in \mathbb C^\times$. As modules over $\ptl_n(\beta)$, we have $\repW_{n,d,z} \simeq \wh\repW_{n,d,z}$.
\end{Proposition}

\paragraph{The special case $\boldsymbol{\alpha = 0}$.}

The case $\alpha = z+z^{-1}= 0$ must be treated separately. First, we note that the modules $\repW_{n,0,z=\iI}$ and $\wh\repW_{n,0,z=\iI}$ are both well defined, but the isomorphism constructed above breaks down. We now show that they are, in general, 
not isomorphic.\medskip

The standard module $\repW_{n,0,\iI}$ splits as the direct sum
\be
\repW_{n,0,\iI} = \repW^{\tinyx 0}_{n,0,\iI} \oplus \repW^{\tinyx 1}_{n,0,\iI}\,,
\ee 
where $\repW^{\tinyx 0}_{n,0,\iI}$ and $\repW^{\tinyx 1}_{n,0,\iI}$ are respectively spanned by the even and odd link states of $\setB_{n,0}$. Indeed, any connectivity that acts on a link state of a given parity and produces a link state of the opposite parity does so by creating a non-contractible loop, which gets a weight $\alpha = 0$ in $\repW_{n,0,\iI}$. We denote by $\setB^{\tinyx 0}_{n,0}$ and $\setB^{\tinyx 1}_{n,0}$ the bases of odd and even link states in $\setB_{n,0}$, respectively. Each odd link state is obtained from an even one by translation by one node, so we have
\be
|\setB^{\tinyx 0}_{n,0}| = |\setB^{\tinyx 1}_{n,0}| = \tfrac12|\setB_{n,0}|.
\ee

To investigate the module $\wh\repW_{n,0,\iI}$, let us denote by $\wh\setB^{\tinyx 0}_{n,0}$ and $\wh\setB^{\tinyx 1}_{n,0}$ the subsets of $\wh\setB_{n,0}$ made of link states with zero and one non-contractible loops, respectively. The module $\wh\repW_{n,0,\iI}$ has a non-trivial submodule denoted by $\wh\repW^{\tinyx 1}_{n,0,\iI}$, spanned by the link states in $\wh\setB^{\tinyx 1}_{n,0}$. It is, moreover, easy to see that this submodule is isomorphic to $\repW^{\tinyx 1}_{n,0,\iI}$. The following proposition shows that the vector space spanned by the link states in $\wh\setB^{\tinyx 0}_{n,0}$ is not closed under the action of $\ptl_n(\beta)$. Its proof uses the elements 
\be
\label{eq:E.and.F}
E=e_0e_2e_4\cdots e_{n-2} 
\qquad \textrm{and} \qquad
F=e_1e_3e_5\cdots e_{n-1}
\ee
in $\ptl_n(\beta)$ and $\atl_n(\beta)$ for $n$ even, whose diagrammatic representations are
\be 
E=\  \begin{pspicture}[shift=-0.25](0,-0.35)(4.4,0.35)
\pspolygon[fillstyle=solid,fillcolor=lightlightblue,linewidth=0pt](0,-0.35)(4.4,-0.35)(4.4,0.35)(0,0.35)
\multiput(0,0)(0.8,0){2}{\psarc[linecolor=blue,linewidth=1.5pt]{-}(0.8,0.35){0.2}{180}{0}}
\multiput(0,0)(0.8,0){2}{\psarc[linecolor=blue,linewidth=1.5pt]{-}(2.8,0.35){0.2}{180}{0}}
\rput(2.2,0){$...$}
\multiput(0,0)(0.8,0){2}{\psarc[linecolor=blue,linewidth=1.5pt]{-}(0.8,-0.35){0.2}{0}{180}}
\multiput(0,0)(0.8,0){2}{\psarc[linecolor=blue,linewidth=1.5pt]{-}(2.8,-0.35){0.2}{0}{180}}
\psarc[linecolor=blue,linewidth=1.5pt]{-}(0,-0.35){0.2}{0}{90}
\psarc[linecolor=blue,linewidth=1.5pt]{-}(4.4,-0.35){0.2}{90}{180}
\psarc[linecolor=blue,linewidth=1.5pt]{-}(0,0.35){0.2}{270}{00}
\psarc[linecolor=blue,linewidth=1.5pt]{-}(4.4,0.35){0.2}{180}{270}
\end{pspicture}\ , 
\qquad F=\ \begin{pspicture}[shift=-0.25](0,-0.35)(4.4,0.35)
\pspolygon[fillstyle=solid,fillcolor=lightlightblue,linewidth=0pt](0,-0.35)(4.4,-0.35)(4.4,0.35)(0,0.35)
\multiput(0,0)(0.8,0){3}{\psarc[linecolor=blue,linewidth=1.5pt]{-}(0.4,0.35){0.2}{180}{0}}
\multiput(0,0)(0.8,0){2}{\psarc[linecolor=blue,linewidth=1.5pt]{-}(3.2,0.35){0.2}{180}{0}}
\rput(2.6,0){$...$}
\multiput(0,0)(0.8,0){3}{\psarc[linecolor=blue,linewidth=1.5pt]{-}(0.4,-0.35){0.2}{0}{180}}
\multiput(0,0)(0.8,0){2}{\psarc[linecolor=blue,linewidth=1.5pt]{-}(3.2,-0.35){0.2}{0}{180}}
\end{pspicture}\ .
\ee 

\begin{Proposition}\label{prop:alpha0.submodule}
For $\beta \neq 0$, the submodule $\wh\repW^{\tinyx 1}_{n,0,\iI}$ does not appear as a direct summand in $\wh\repW_{n,0,\iI}$.
\end{Proposition}
\proof
Let $v$ be any link state in $\wh\setB^{\tinyx 0}_{n,0}$. The state $EF\cdot v$ is equal to a link state in $\wh\setB^{\tinyx 1}_{n,0}$ times a non-zero power of $\beta$. For $\beta \neq 0$, this state is non-zero and lies in the submodule $\wh\repW^{\tinyx 1}_{n,0,\iI}$. Moreover, we have $(EF)^2\cdot v=0$. We have thus shown that, for $\beta \neq 0$, $EF$ has a rank-two Jordan cell of eigenvalue $0$ that couples the submodule with its complement.\eproof 

\subsection{Sandwich bases}\label{ssec:SandwichBases}

One way to produce bases of connectivities for $\ptl_n(\beta)$ or $\atl_n(\beta)$ is to consider {\it sandwich diagrams}. Let $\setE$ be a subset of connectivities of $\atl_d(\beta)$ that have $d$ through-lines. In other words, for $d>0$, the set $\setE$ may only contain the unit $\id$ and powers of $\Omega$ and $\Omega^{-1}$, and for $d=0$, the set $\setE$ may only contain the unit $\id$ and powers of $f$, namely the spanning elements of $\atl_0(\beta)$. A sandwich diagram is therefore a way to draw a connectivity by dividing it in three parts: its bottom, middle and top. We define the sets $\setS_d(\setE)$ of connectivities with $d$ through-lines as
\be
\label{eq:sandwich.diag}
\setS_d(\setE) = \Big\{
c_{b,m,t}^d = \
\begin{pspicture}[shift=-1.2](0,-1.3)(2.0,1.3)
\psline[linewidth=0.5pt](0,-1.3)(2,-1.3)
\pspolygon[linewidth=0.5pt,fillstyle=solid,fillcolor=lightlightblue](0.2,-1.3)(1.8,-1.3)(1.8,-0.8)(0.2,-0.8)\rput(1,-1.05){$b$}
\psline[linewidth=0.5pt](0,1.3)(2,1.3)
\pspolygon[linewidth=0.5pt,fillstyle=solid,fillcolor=lightlightblue](0.2,1.3)(1.8,1.3)(1.8,0.8)(0.2,0.8)\rput(1,1.05){${\iota(t)}$}
\pspolygon[fillstyle=solid,fillcolor=lightlightblue,linewidth=0.5pt](0.55,-0.4)(1.45,-0.4)(1.45,0.4)(0.55,0.4)
\rput(1,0){$m$}
\psbezier[linecolor=blue,linewidth=1.5pt]{-}(0.3,-0.8)(0.3,-0.6)(0.65,-0.6)(0.65,-0.4)
\psbezier[linecolor=blue,linewidth=1.5pt]{-}(0.6,-0.8)(0.6,-0.6)(0.85,-0.6)(0.85,-0.4)
\psbezier[linecolor=blue,linewidth=1.5pt]{-}(1.4,-0.8)(1.4,-0.6)(1.15,-0.6)(1.15,-0.4)
\psbezier[linecolor=blue,linewidth=1.5pt]{-}(1.7,-0.8)(1.7,-0.6)(1.35,-0.6)(1.35,-0.4)
\rput(1.0,-0.6){...}
\psbezier[linecolor=blue,linewidth=1.5pt]{-}(0.3,0.8)(0.3,0.6)(0.65,0.6)(0.65,0.4)
\psbezier[linecolor=blue,linewidth=1.5pt]{-}(0.6,0.8)(0.6,0.6)(0.85,0.6)(0.85,0.4)
\psbezier[linecolor=blue,linewidth=1.5pt]{-}(1.4,0.8)(1.4,0.6)(1.15,0.6)(1.15,0.4)
\psbezier[linecolor=blue,linewidth=1.5pt]{-}(1.7,0.8)(1.7,0.6)(1.35,0.6)(1.35,0.4)
\rput(1.0,0.6){...}
\end{pspicture} \ 
\Big | \ b,t \in \setB_{n,d}\,, 
\ m \in \setE
\Big\}.
\ee
In this construction, the operator $\iota$ applies a reflection to the link states $t$ with respect to a horizontal axis. The $k$ bottom nodes of the connectivity $m$ are attached to the $k$ defects of $b$, and likewise the top nodes of $m$ are attached to the $k$ defects of $t$ drawn upside down. \medskip

The sandwich construction allows us to write the set of connectivities of $\atl_n(\beta)$ as
\be
\label{eq:SandwichBaseaTL}
\atl_n(\beta): \quad \left\{ 
\begin{array}{cl}
\displaystyle\bigcup_{d=1,3,\dots,n}\setS_d(\{\Omega^r\,|\, r \in \mathbb Z\}) & n \textrm{ odd,}
\\[0.6cm]
\displaystyle\setS_0(\{f^r\,|\, r \in \mathbb Z_{\ge 0}\}) \ \cup \bigcup_{d=2,4,\dots,n}\setS_d(\{\Omega^r\,|\, r \in \mathbb Z\}) & n \textrm{ even.}
\end{array}\right.
\ee
Indeed, the connectivities of $\atl_n(\beta)$ with $d\ge 1$ allow for an arbitrary winding of the through-lines, accounted for by the arbitrary power of $\Omega$ (with the convention $\Omega^0=\id$). For $n$ even and $d=0$, the arbitrary number of non-contractible loops is accounted for by the powers of $f$ (with the convention $f^0=\id$). Because the indices $r$ in~\eqref{eq:SandwichBaseaTL} are unbounded, this algebra is infinite dimensional.\medskip

For $\ptl_n(\beta)$, the same construction requires that we retain only the even sandwich diagrams and also that we exclude the diagrams $\Omega^{2j}$ and $\Omega^{-2j}$ for $j\in \mathbb{Z}_{>0}$. 
For $b,t\in \setB_{n,d}$, let us write $\sigma_{bt}=0$ if $b$ and $t$ have the same parity, and $\sigma_{bt}=1$ if they have different parities. Then the set of connectivities of $\ptl_n(\beta)$ is
\be \label{eq:sandwich.basis.ptl}
\ptl_n(\beta): \left\{
\begin{array}{cl}
\setS_n(\{\id\}) \ \cup \displaystyle\bigcup_{d=1,3,\dots,n-2}\setS_d(\{\Omega^r\,|\, r \in 2\mathbb Z + \sigma_{bt}\})& n \textrm{ odd,}
\\[0.6cm]
\displaystyle\setS_0(\{f^r\,|\, r \in 2\mathbb Z_{\ge 0}+\sigma_{bt}\})  \cup  \setS_n(\{\id\})  \cup\!\!\!\!\!\bigcup_{d=2,4,\dots,n-2}\!\!\!\!\setS_d(\{\Omega^r\,|\, r \in 2\mathbb Z+\sigma_{bt}\}) & n \textrm{ even.}
\end{array}\right.
\ee
Thus, by restricting $r$ to take only odd and even integers according to $\sigma_{bt}$ in this way, we indeed produce only the even diagrams.
In fact, we can give an equivalent construction for the sandwich bases for $d>0$ by defining the sets
\be 
\label{eq:sandwich.whS}
\wh \setS_d(\setE) = \Big\{
 c_{b,m,t}^{\,d} \ 
\Big | \ b,t \in \wh \setB_{n,d}\,, 
\ m \in \setE
\Big\}.
\ee
The set of connectivities of $\ptl_n(\beta)$ can then be equivalently written as 
\be \label{eq:sandwich.basis.ptl.v2}
\ptl_n(\beta): \left\{
\begin{array}{cl}
\wh\setS_n(\{\id\}) \ \cup \displaystyle\bigcup_{d=1,3,\dots,n-2}\wh\setS_d(\{\Omega^r\,|\, r \in 2\mathbb Z\})& n \textrm{ odd,}
\\[0.6cm]
\displaystyle\setS_0(\{f^r\,|\, r \in 2\mathbb Z_{\ge 0}+\sigma_{bt}\})  \cup \wh\setS_n(\{\id\})  \cup\bigcup_{d=2,4,\dots,n-2}\wh\setS_d(\{\Omega^r\,|\, r \in 2\mathbb Z\}) & n \textrm{ even,}
\end{array}\right.
\ee
with the added benefit that the sets $\mathcal E$ for $d>0$ do not depend on the parities $\sigma_{bt}$. Crucially, for $n$ even and $d=0$, $\setS_0(\{f^r\,|\, r \in 2\mathbb Z_{\ge 0}+\sigma_{bt}\})$ and $\wh \setS_0(\{f^r\,|\, r \in 2\mathbb Z_{\ge 0}\})$ are not identical, as for instance the connectivity $E$ is in the first set but not in the second. This remark plays an important role in \cref{sec:sandwich.cell} in the description of $\ptl_n(\beta)$ as an affine cellular or a skew sandwich algebra.\medskip

The set of connectivities of $\tl_n(\beta)$ can also be obtained from sandwich diagrams. In this case however, the link states $b$ and $t$ are taken from the subset $\mathcal B_{n,d}'\subset \setB_{n,d}$ of link states without periodic arcs, and the corresponding sandwiches have no filling, namely $\mathcal E$ contains only the identity. We denote the corresponding subset of sandwich diagrams as $\setS_d'(\{\id\})$, and find
\be
\tl_n(\beta): \quad \bigcup_{d}\ \setS_d'(\{\id\}),
\ee
where $d$ runs over the integers in the range $0 \le d \le n$ with $d \equiv n \textrm{ mod }2$. The dimension of $\tl_n(\beta)$ is then given by
\be
\tl_n(\beta) = \sum_d |\setB'_{n,d}|^2 = \sum_{d} \left( \binom{n}{\frac{n-d}{2}} - \binom{n}{\frac{n-d-2}{2}}\right)^2 = \frac1{n+1}\binom{2n}{n} = C_n,
\ee
the Catalan number.\medskip

%
\section{Uncoiled affine and periodic Temperley--Lieb algebras}\label{sec:def.uaptl}
%

\subsection{Definition of the algebras}

In this section, we write down extra relations for $\atl_n(\beta)$ and $\ptl_n(\beta)$ that ultimately allow us to define finite-dimensional quotients of these algebras, which we respectively call {\it uncoiled affine Temperley--Lieb algebras} and {\it uncoiled periodic Temperley--Lieb algebras}. As we shall see, the cases with $n$ odd and even are quite different, as we shall define two types of uncoiled algebras for $n$ even and one kind for $n$ odd. We thus treat the two parities of $n$ separately.

\paragraph{$\boldsymbol n$ odd.}
For $n$ odd, we consider the quotient relations
\begin{subequations}
\label{eq:unwinding.ptl.atl.odd}
\begin{alignat}{2}
&\atl_n(\beta): \qquad \Omega^n = \gamma\,\id, \label{eq:unwinding.atl.odd}
\\[0.1cm]
&\ptl_n(\beta): \qquad e_0(e_{n-1}e_{n-2} \cdots e_1 e_0)^{n-2} = \gamma^2 e_0.
\label{eq:unwinding.ptl.odd}
\end{alignat}
\end{subequations}
In the diagram presentation, these relations allow for the unwinding of the through-lines, by $n$ nodes for $\atl_n(\beta)$ and by $2n$ nodes for $\ptl_n(\beta)$, and involve a free parameter $\gamma \in \mathbb C^\times$. The uncoiled algebras $\qatl_n(\beta,\gamma)$ and $\qptl_n(\beta,\gamma)$ are then defined as
\be
\label{eq:odd.quotients}
\hspace{-0.1cm}n \textrm{ odd}: \quad
\left\{\begin{array}{l}
\qatl_n(\beta,\gamma) = 
\big\langle \Omega, \Omega^{-1}, e_0, e_1,\dots, e_{n-1} \big\rangle \Big/ \big\{\eqref{eq:PTL.relations},\eqref{eq:ATL.relations},\eqref{eq:unwinding.atl.odd}\big\},
\\[0.1cm]
\qptl_n(\beta,\gamma) = \big\langle \id, e_0, e_1,\dots, e_{n-1} \big\rangle \Big/ \big\{\eqref{eq:PTL.relations},\eqref{eq:unwinding.ptl.odd}\big\}.
\end{array}\right.
\ee
The relation \eqref{eq:unwinding.ptl.odd} in fact also holds in $\qatl_n(\beta,\gamma)$. It is proven as follows:
\begin{alignat}{2}
e_0(e_{n-1}e_{n-2} \cdots e_1 e_0)^{n-2} &= e_0(\Omega^2 e_1 e_0)^{n-2} = \Omega^{2n-4} e_{n-4}e_{n-5} \cdots e_1e_0e_{n-1}e_{n-2}\cdots e_1e_0 
\nonumber\\
&= \Omega^{2n-2} e_{n-2}e_{n-3} \cdots e_1e_0 = \Omega^{2n} e_0 = \gamma^2 e_0. 
\end{alignat}
This implies that $\qptl_n(\beta,\gamma)$ is a subalgebra of $\qatl_n(\beta,\gamma)$.

\paragraph{$\boldsymbol n$ even.}

For $n$ even, we define two distinct pairs of uncoiled algebras. For the first pair of algebras, we respectively consider the extra relations 
\begin{subequations}
\label{eq:quotients.ptl.and.atl.1}
\begin{alignat}{3}
\label{eq:extra.relations.atl.even.1} 
&\atl_n(\beta): \qquad\qquad\quad E\, \Omega\, E = \alpha\, E, \qquad &&\Omega^n = \id, 
\\[0.1cm]
\label{eq:extra.relations.ptl.even.1}
&\ptl_n(\beta): \qquad E F E = \alpha^2 E, \ FEF = \alpha^2 F, \qquad &&e_0(e_{n-1}e_{n-2} \cdots e_1 e_0)^{(n-2)/2} = e_0,
\end{alignat}
\end{subequations}
where we recall that $E$ and $F$ are defined in~\eqref{eq:E.and.F}. In the diagram presentation, the first relations in \eqref{eq:extra.relations.atl.even.1} and \eqref{eq:extra.relations.ptl.even.1} allow one to erase non-contractible loops and replace them by factors of $\alpha$. In $\qatl_n(\beta,\alpha)$, each such loop may be removed individually, whereas these loops are removed in pairs for $\qptl_n(\beta,\alpha)$. The second relations allow for the unwinding of the through-lines, by $n$ nodes in both cases. The algebras $\qatla_n(\beta,\alpha)$ and $\qptla_n(\beta,\alpha)$ are then defined as
\be
\label{eq:even.quotients.1}
n \textrm{ even}: \quad 
\left\{\begin{array}{l}
\qatla_n(\beta,\alpha) = \big\langle \Omega, \Omega^{-1}, e_0, e_1,\dots, e_{n-1} \big\rangle \Big/ \big\{\eqref{eq:PTL.relations},\eqref{eq:ATL.relations},\eqref{eq:extra.relations.atl.even.1}\big\},
\\[0.2cm]
\qptla_n(\beta,\alpha) = \big\langle \id, e_0, e_1,\dots, e_{n-1} \big\rangle \Big/ \big\{\eqref{eq:PTL.relations},\eqref{eq:extra.relations.ptl.even.1}\big\}.
\end{array}\right.
\ee
Seen as relations of $\atl_n(\beta)$, the identities in \eqref{eq:extra.relations.ptl.even.1} can be proven using those of \eqref{eq:extra.relations.atl.even.1}, therefore implying that $\qptla_n(\beta,\alpha)$ is a subalgebra of $\qatla_n(\beta,\alpha)$.\medskip

The pair of quotient algebras~\eqref{eq:even.quotients.1} is very different from those defined in~\eqref{eq:odd.quotients} for $n$ odd. First, in the relation~\eqref{eq:unwinding.ptl.odd} that unwinds the through-lines in $\qptla_n(\beta,\alpha)$ for $n$ odd, the exponent of $(e_{n-1}e_{n-2}\dots e_1 e_0)$ is double the one for $n$ even in relation~\eqref{eq:extra.relations.ptl.even.1}. It is indeed not possible in $\qptl_n(\beta, \gamma)$ with $n$ odd to unwind the through-lines by $n$ nodes, as this would create connectivities with odd parities.  Secondly, for $n$ odd, the unwinding of the through-lines allows for the introduction of the parameter~$\gamma$. For $n$ even, it is not possible to include this parameter in~\eqref{eq:quotients.ptl.and.atl.1}. Indeed, the relation $\Omega^n c = c$, valid in $\atl_n(\beta)$ for $n$ even and any connectivity $c$ without through-lines, shows that only $\gamma = 1$ is possible for $n$ even.\medskip

It is, however, possible to consider a second kind of quotient for $n$ even, with the relations
\begin{subequations}
\label{eq:quotients.ptl.and.atl.2}
\begin{alignat}{3}
\label{eq:extra.relations.atl.even.2}
&\atl_n(\beta): \qquad \qquad E = 0, \qquad &&\Omega^n = \gamma\,\id,
\\[0.1cm]
\label{eq:extra.relations.ptl.even.2}
&\ptl_n(\beta): \qquad E=0,\ F=0, \qquad &&e_0(e_{n-1}e_{n-2} \cdots e_1 e_0)^{(n-2)/2} = \gamma e_0. 
\end{alignat}
\end{subequations}
In this case, it is not hard to see that the conditions $E=0$ for $\atl_n(\beta)$ and $E=0=F$ for $\ptl_n(\beta)$ imply that any connectivity $c$ with no through-lines is set to zero, as all such connectivities can be obtained by acting with the algebra on~$E$ or $F$ on the left and right. The identity $\Omega^n c = c$ for these connectivities is then satisfied. Moreover, the unwinding of the defects now involves the free parameter~$\gamma$. The second pair of uncoiled algebras for $n$ even, denoted as $\qatlb_n(\beta,\gamma)$ and $\qptlb_n(\beta,\gamma)$, are then defined as
\be
\label{eq:even.quotients.2}
\hspace{-0.1cm}n \textrm{ even}: \quad 
\left\{\begin{array}{l}
\qatlb_n(\beta,\gamma) = \big\langle \Omega, \Omega^{-1}, e_0, e_1,\dots, e_{n-1} \big\rangle \Big/ \big\{\eqref{eq:PTL.relations},\eqref{eq:ATL.relations},\eqref{eq:extra.relations.atl.even.2}\big\},
\\[0.15cm]
\qptlb_n(\beta,\gamma) = \big\langle \id, e_0, e_1,\dots, e_{n-1} \big\rangle \Big/ \big\{\eqref{eq:PTL.relations},\eqref{eq:extra.relations.ptl.even.2}\big\}.
\end{array}\right.
\ee

We end this section by illustrating in examples the differences between the three uncoiled periodic algebras in how the defects are unwound. For $\qptl_5(\beta,\gamma)$, $\qptla_6(\beta,\alpha)$ and $\qptlb_6(\beta,\gamma)$, we have
\begin{subequations}
\begin{alignat}{2}
\label{eq:quotient.eptl.n56}
\begin{pspicture*}[shift=-0.9](0,-1)(2.0,1)
\pspolygon[fillstyle=solid,fillcolor=lightlightblue,linewidth=0pt](0,-1)(2.0,-1)(2.0,1)(0,1)
\psarc[linecolor=blue,linewidth=1.5pt]{-}(0,1){0.2}{-90}{0}
\psarc[linecolor=blue,linewidth=1.5pt]{-}(0,-1){0.2}{0}{90}
\psarc[linecolor=blue,linewidth=1.5pt]{-}(2.0,1){0.2}{-180}{-90}
\psarc[linecolor=blue,linewidth=1.5pt]{-}(2.0,-1){0.2}{90}{180}
\multiput(0,0)(-2,0){3}{\multiput(0,0)(0.4,0){3}{\psbezier[linecolor=blue,linewidth=1.5pt]{-}(0.6,-1)(0.6,-0.2)(4.6,0.2)(4.6,1)}}
\end{pspicture*} 
\ &= \gamma^2  \ 
\begin{pspicture}[shift=-0.25](0,-0.35)(2.0,0.35)
\pspolygon[fillstyle=solid,fillcolor=lightlightblue,linewidth=0pt](0,-0.35)(2.0,-0.35)(2.0,0.35)(0,0.35)
\psarc[linecolor=blue,linewidth=1.5pt]{-}(0,0.35){0.2}{-90}{0}
\psarc[linecolor=blue,linewidth=1.5pt]{-}(0,-0.35){0.2}{0}{90}
\psarc[linecolor=blue,linewidth=1.5pt]{-}(2.0,0.35){0.2}{-180}{-90}
\psarc[linecolor=blue,linewidth=1.5pt]{-}(2.0,-0.35){0.2}{90}{180}
\psline[linecolor=blue,linewidth=1.5pt]{-}(0.6,-0.35)(0.6,0.35)
\psline[linecolor=blue,linewidth=1.5pt]{-}(1.0,-0.35)(1.0,0.35)
\psline[linecolor=blue,linewidth=1.5pt]{-}(1.4,-0.35)(1.4,0.35)
\qquad 
\end{pspicture}
\qquad \qptl_5(\beta,\gamma), 
\\
\begin{pspicture*}[shift=-0.6](0,-0.7)(2.4,0.7)
\pspolygon[fillstyle=solid,fillcolor=lightlightblue,linewidth=0pt](0,-0.7)(2.4,-0.7)(2.4,0.7)(0,0.7)
\psarc[linecolor=blue,linewidth=1.5pt]{-}(0,0.7){0.2}{-90}{0}
\psarc[linecolor=blue,linewidth=1.5pt]{-}(0,-0.7){0.2}{0}{90}
\psarc[linecolor=blue,linewidth=1.5pt]{-}(2.4,0.7){0.2}{-180}{-90}
\psarc[linecolor=blue,linewidth=1.5pt]{-}(2.4,-0.7){0.2}{90}{180}
\multiput(0,0)(0.4,0){4}{\psbezier[linecolor=blue,linewidth=1.5pt]{-}(0.6,-0.7)(0.6,-0.2)(3.0,0.2)(3.0,0.7)}
\multiput(-2.4,0)(0.4,0){4}{\psbezier[linecolor=blue,linewidth=1.5pt]{-}(0.6,-0.7)(0.6,-0.2)(3.0,0.2)(3.0,0.7)}
\end{pspicture*} 
\ &=  \ \left\{\begin{array}{rl}
\begin{pspicture}[shift=-0.25](0,-0.35)(2.4,0.35)
\pspolygon[fillstyle=solid,fillcolor=lightlightblue,linewidth=0pt](0,-0.35)(2.4,-0.35)(2.4,0.35)(0,0.35)
\psarc[linecolor=blue,linewidth=1.5pt]{-}(0,0.35){0.2}{-90}{0}
\psarc[linecolor=blue,linewidth=1.5pt]{-}(0,-0.35){0.2}{0}{90}
\psarc[linecolor=blue,linewidth=1.5pt]{-}(2.4,0.35){0.2}{-180}{-90}
\psarc[linecolor=blue,linewidth=1.5pt]{-}(2.4,-0.35){0.2}{90}{180}
\psline[linecolor=blue,linewidth=1.5pt]{-}(0.6,-0.35)(0.6,0.35)
\psline[linecolor=blue,linewidth=1.5pt]{-}(1.0,-0.35)(1.0,0.35)
\psline[linecolor=blue,linewidth=1.5pt]{-}(1.4,-0.35)(1.4,0.35)
\psline[linecolor=blue,linewidth=1.5pt]{-}(1.8,-0.35)(1.8,0.35)
\end{pspicture}  & \qptla_6(\beta,\alpha), \\[.6cm]
\gamma \	\begin{pspicture}[shift=-0.25](0,-0.35)(2.4,0.35)
\pspolygon[fillstyle=solid,fillcolor=lightlightblue,linewidth=0pt](0,-0.35)(2.4,-0.35)(2.4,0.35)(0,0.35)
\psarc[linecolor=blue,linewidth=1.5pt]{-}(0,0.35){0.2}{-90}{0}
\psarc[linecolor=blue,linewidth=1.5pt]{-}(0,-0.35){0.2}{0}{90}
\psarc[linecolor=blue,linewidth=1.5pt]{-}(2.4,0.35){0.2}{-180}{-90}
\psarc[linecolor=blue,linewidth=1.5pt]{-}(2.4,-0.35){0.2}{90}{180}
\psline[linecolor=blue,linewidth=1.5pt]{-}(0.6,-0.35)(0.6,0.35)
\psline[linecolor=blue,linewidth=1.5pt]{-}(1.0,-0.35)(1.0,0.35)
\psline[linecolor=blue,linewidth=1.5pt]{-}(1.4,-0.35)(1.4,0.35)
\psline[linecolor=blue,linewidth=1.5pt]{-}(1.8,-0.35)(1.8,0.35)
\end{pspicture} & \qptlb_6(\beta,\gamma).
\end{array}\right.
\end{alignat} 
\end{subequations}

\subsection{Properties of the algebras}\label{sec:properties.quotients}

This section discusses various properties of the uncoiled algebras. First, it is clear that all six algebras are finite-dimensional. Indeed, the quotient relations allow us to write their connectivities in terms of diagrams with at most one non-contractible loop, and where each through-line crosses the periodic boundary condition at most twice. For each quotient algebra, we give a basis for its connectivities in terms of sandwich diagrams and use this to compute their dimensions. The resulting formulas involve sums with squares of binomial coefficients, which are evaluated using the identities derived in \cref{app:binom.identities}. We then give the finite set of standard modules for each of these algebras. These are, of course, a subset of the standard modules of $\atl_n(\beta)$ or $\ptl_n(\beta)$. We describe separately the cases of odd and even~$n$.

\paragraph{$\boldsymbol n$ odd.} 

In this case, there are two algebras to consider. 
\begin{Proposition}\label{Prop:sandwich.bases.odd}
Let $n\in \mathbb{N}$ be odd. The full sets of sandwich
diagrams of $\qatl_n(\beta,\gamma)$ and $\qptl_n(\beta,\gamma)$ are given by
\begin{subequations}
\begin{alignat}{2}
\qatl_n(\beta,\gamma)&: \quad
\bigcup_{d=1,3,\dots,n}\setS_d(\{\Omega^r \, | \, 0\le r < d\}),\label{eq:Sandwich.Basis.uatl}
\\[0.3cm]
\qptl_n(\beta,\gamma)&: \quad
\wh\setS_n(\id)\, \cup \, \bigcup_{d=1,3,\dots,n-2}\wh\setS_d(\{\Omega^{2r} \, | \, 0\le r < d\}).\label{eq:Sandwich.Basis.uptl}
\end{alignat}
\end{subequations}
\end{Proposition}
\proof
For $n$ odd, by using the relations \eqref{eq:unwinding.atl.odd}, we can write any connectivity in $\qatl_n(\beta,\gamma)$ in terms of diagrams where (i)~at least one through-line does not cross the periodic boundary, and (ii)~the through-lines that cross the periodic boundary, if there are any, do so by travelling to the left as they move from the top to the bottom edges of the box. This also implies that each through-line can travel by at most $n$ nodes. For $\qptl_n(\beta,\gamma)$, the relation \eqref{eq:unwinding.ptl.odd} only allows for the unwinding of the through-lines by $2n$ nodes. Thus in this case, the connectivities
have the following properties: (i) there is at least one through-line that crosses the periodic boundary at most once, (ii) all crossing through-lines go leftwards as they proceed from the top to the bottom of the box, (iii) the only connectivity in $\qptl_n(\beta,\gamma)$ with $n$ through-lines is the unit $\id$, and (iv) in the sandwich basis, the powers of $\Omega$ must be even as the parities of the link states in 
the bases $\wh \setB_{n,d}$ are all even. \eproof

The dimensions of these algebras follow directly.
\begin{Corollary}\label{coro:dim.uatl.odd}
The dimensions of $\qatl_n(\beta,\gamma)$ and $\qptl_n(\beta,\gamma)$ for $n$ odd are
\begin{subequations}
\begin{alignat}{2}
\dim \qatl_n(\beta,\gamma) &= \sum_{d=1,3,\dots,n} d\, |\setB_{n,d}|^2 = n\, \binom{n-1}{\frac{n-1}2}^2,
\\[0.3cm]
\dim \qptl_n(\beta,\gamma) &= 1+\sum_{d=1,3,\dots,n-2} d\, |\setB_{n,d}|^2 = n\, \binom{n-1}{\frac{n-1}2}^2 - (n-1).
\end{alignat}
\end{subequations}
\end{Corollary}
\proof
The two leftmost equalities follow directly from the previous proposition, and the evaluation of the sums follows from \cref{app:binom.identities}.
\eproof

\begin{Proposition}
\label{prop:standards.uncoiled.odd}
For $n$ odd, the subsets of standard modules of $\atl_n(\beta)$ that are also modules over $\qatl_n(\beta, \gamma)$ or $\qptl_n(\beta,\gamma)$ are
\begin{subequations}
\begin{alignat}{2}
\qatl_n(\beta,\gamma)&: \quad \big\{\repW_{n,d,\gamma^{1/d}\eE^{2 \pi \iI r/d}} \ | \ 0\le r < d, \ d = 1,3,\dots,n\big\}, \label{eq:Std.Mod.uaTL.odd}
\\[0.3cm]
\qptl_n(\beta,\gamma)&: \quad \big\{\repW_{n,n,1} \big\} \cup \big\{\repW_{n,d,\gamma^{1/d}\eE^{\pi \iI r/d}} \ | \ 0\le r < d, \ d = 1,3,\dots,n-2\big\}.
\end{alignat}
\end{subequations}
\end{Proposition}
\proof
To identify the standard modules for these algebras, we evaluate the quotient relations \eqref{eq:unwinding.ptl.atl.odd} in $\repW_{n,d,z}$. For $\qatl_n(\beta,\gamma)$, we read from \eqref{eq:FO.eigenvalues} that $\Omega^n = z^d \id$ on this module, which constrains~$z$ to the values satisfying $z^d = \gamma$. For $\qptl_n(\beta,\gamma)$, we instead have $e_0(e_{n-1}e_{n-2} \cdots e_1 e_0)^{n-2} = z^{2d} e_0$ in $\repW_{n,d,z}$. This is easily seen from the diagram representation of $e_0(e_{n-1}e_{n-2} \cdots e_1 e_0)^{n-2}$ (given in the first connectivity of~\eqref{eq:quotient.eptl.n56} for $n=5$), where it is clear that each defect crosses the boundary condition twice. The constraint on $z$ in this case reads $z^{d} = \pm\gamma$. For $\qptl_n(\beta,\gamma)$, we additionally use \cref{lem:iso.minus.z}, which states that $\repW_{n,d,z}\simeq\repW_{n,d,-z}$. 
\eproof

\paragraph{$\boldsymbol n$ even.} In this case, there are four algebras to consider.

\begin{Proposition}\label{Prop:sandwich.bases.even}
Let $n\in \mathbb{N}$ be even. The full sets of sandwich diagrams for $\qatla_n(\beta,\alpha)$ and $\qptla_n(\beta,\alpha)$ with $\alpha \neq 0$ are given by
\begin{subequations}
\begin{alignat}{2}
\qatla_n(\beta,\alpha)&: \quad
\displaystyle\setS_0(\{\id\}) \ \cup \bigcup_{d=2,4,\dots,n}\setS_d(\{\Omega^r \, | \, 0\le r <d\}), \label{eq:Sandwich.Basis.uatla}
\\[0.3cm]
\qptla_n(\beta,\alpha)&: \quad
\displaystyle \wh\setS_0(\{\id\}) \ \cup\ \wh \setS_n(\{\id\}) \ \cup\bigcup_{d=2,4,\dots,n-2}\wh \setS_d(\{\Omega^{2r}\,|\, 0\le r < \tfrac d 2\}), \qquad \alpha \neq 0,\label{eq:Sandwich.Basis.uptla}
\end{alignat}
\end{subequations}
and for $\qatlb_n(\beta,\gamma)$ and $\qptlb_n(\beta,\gamma)$, they are given by
\begin{subequations}
\begin{alignat}{2}
\qatlb_n(\beta,\gamma)&: \quad
\bigcup_{d=2,4,\dots,n}\setS_d(\{\Omega^r \, | \, 0\le r <d\}),\label{eq:Sandwich.Basis.uatlb}
\\[0.3cm]
\qptlb_n(\beta,\gamma)&: \quad
\wh \setS_n(\{\id\}) \ \cup\bigcup_{d=2,4,\dots,n-2} \wh \setS_d(\{\Omega^{2r}\,|\, 0\le r < \tfrac d 2\}).\label{eq:Sandwich.Basis.uptlb}
\end{alignat}
\end{subequations}
\end{Proposition}
\proof
For $n$ even, the relations \eqref{eq:quotients.ptl.and.atl.1} imply that any element in $\qatla_n(\beta,\alpha)$ can be written in terms of a connectivity where (i)~at least one through-line does not cross the periodic boundary, (ii)~all crossing through-lines go to the left as they move from the top to the bottom of the box, and (iii)~there are no non-contractible loops. For $\qptla_n(\beta,\alpha)$, the connectivities have the following properties: (i)~at~least one through-line does not cross the periodic boundary, (ii)~all crossing through-lines go to the left, (iii)~in the case of zero through-lines, the connectivities have zero and one non-contractible loops for $\sigma_{bt} = 0$ and $\sigma_{bt} = 1$ respectively, which is correctly implemented in the sets $\wh\setS_{0}(\{\id\})$,
(iv)~the only connectivity with $n$ through-lines is the unit $\id$, and (v)~in the sandwich bases with $0 \le d \le n-2$, 
the powers of $\Omega$ must be even as the link states of $\wh \setB_{n,d}$ are even. We note that, for $d=0$, certain elements in~$\wh\setS_0$ include prefactors of $\alpha^2$, which can be divided out since $\alpha \neq 0$ by assumption.\medskip

Finally, for $\qatlb_n(\beta,\gamma)$, the set of connectivities is the same as for $\qatla_n(\beta,\alpha)$, except that those with zero through-lines are absent.
\eproof

For the special case of $\qatla_n(\beta,\alpha)$ with $\alpha = 0$, some elements in $\wh\setS_0(\{\id\})$ are zero, and this set does not give a full set of connectivities without through-lines. The full set of connectivities can instead be obtained from \eqref{eq:Sandwich.Basis.uptla} by replacing $\wh\setS_0(\{\id\})$ by $\setS_0(\{f^{\sigma_{bt}}\})$. In particular, this does not play a role in the following corollary and proposition. 
However, as discussed in \cref{sec:sandwich.cell}, it prevents us from describing this algebra as an affine cellular algebra.

\begin{Corollary}\label{coro:dim.uatl.even}
The dimensions of $\qatla_n(\beta,\alpha)$ and $\qptla_n(\beta,\alpha)$ are
\begin{subequations}
\begin{alignat}{2}
\hspace{-0.3cm}\dim \qatla_n(\beta,\alpha) &= |\setB_{n,0}|^2 + \sum_{d=2,4,\dots,n} d\, |\setB_{n,d}|^2 = (n+4)\, \binom{n-1}{\frac{n}2}^2,
\\[0.3cm]
\hspace{-0.3cm}\dim \qptla_n(\beta,\alpha) &=  |\setB_{n,0}|^2+1+\sum_{d=2,4,\dots,n-2} \tfrac{d}2\, |\setB_{n,d}|^2 = (\tfrac n 2 + 4)\, \binom{n-1}{\frac{n}2}^2 - (\tfrac n2-1),
\end{alignat}
\end{subequations}
and those of $\qatlb_n(\beta,\gamma)$ and $\qptlb_n(\beta,\gamma)$ are
\begin{subequations}
\begin{alignat}{2}
\dim \qatlb_n(\beta,\gamma) &= \sum_{d=2,4,\dots,n} d\, |\setB_{n,d}|^2 = n\, \binom{n-1}{\frac{n}2}^2,
\\[0.3cm]
\dim \qptlb_n(\beta,\gamma) &= 1+\sum_{d=2,4,\dots,n-2} \tfrac d2\, |\setB_{n,d}|^2=\tfrac n 2 \, \binom{n-1}{\frac{n}2}^2 - (\tfrac n2-1).
\end{alignat}
\end{subequations}
\end{Corollary}
\proof
The leftmost equalities follow from the previous proposition, and the sums are evaluated using the results of \cref{app:binom.identities}. 
\eproof

\begin{Proposition}
\label{prop:standards.uncoiled.even}
For $n$ even, the subsets of standard modules of $\atl_n(\beta)$ that are also modules over $\qatla_n(\beta, \alpha)$ or $\qptla_n(\beta,\alpha)$ are
\begin{subequations}
\begin{alignat}{2}
\qatla_n(\beta,\alpha)&: \quad \big\{\repW_{n,0,z}\big\}\cup \big\{\repW_{n,d,\eE^{2 \pi \iI r/d}} \ | \ 0\le r < d, \ d = 2,4,\dots,n\big\},\label{eq:Std.Mod.uaTLa.even}
\\[0.3cm]
\qptla_n(\beta,\alpha)&: \quad \big\{\repW_{n,0,z}\big\}\cup\big\{\repW_{n,n,1} \big\} \cup \big\{\repW_{n,d,\eE^{2 \pi \iI r/d}} \ | \ 0\le r < \tfrac d2, \ d = 2,4,\dots,n-2\big\}.
\end{alignat}
\end{subequations}
Likewise, for $\qatlb_n(\beta,\gamma)$ and $\qptlb_n(\beta,\gamma)$, these modules are those with $d>0$ and $z^d = \gamma$, namely
 \begin{subequations}
\begin{alignat}{2}
\qatlb_n(\beta,\gamma)&: \quad \big\{\repW_{n,d,\gamma^{1/d}\eE^{2 \pi \iI r/d}} \ | \ 0\le r < d, \ d = 2,4,\dots,n\big\},
\\[0.3cm]
\qptlb_n(\beta,\gamma)&: \quad\big\{\repW_{n,n,1} \big\} \cup \big\{\repW_{n,d,\gamma^{1/d}\eE^{2 \pi \iI r/d}} \ | \ 0\le r < \tfrac d2, \ d = 2,4,\dots,n-2\big\}.
\end{alignat}
\end{subequations}
\end{Proposition}
\proof
 In the standard module $\repW_{n,d,z}$ with $d>0$, the quotient relations that allow for the unwinding of the through-lines in \eqref{eq:quotients.ptl.and.atl.1} read $\Omega^n = z^d \id$ and $e_0(e_{n-1}e_{n-2} \cdots e_1 e_0)^{(n-2)/2} = z^d e_0$. This imposes that $z^d = 1$  for $\qatla_n(\beta,\alpha)$ and $\qptla_n(\beta,\alpha)$, and $z^d = \gamma$ for $\qatlb_n(\beta,\gamma)$ and $\qptlb_n(\beta,\gamma)$.
\eproof

For $n$ odd or even, we say that an uncoiled affine or periodic algebra is \emph{generic} if all of its standard modules are simple and pairwise inequivalent. 
\begin{Corollary}\label{Coro:uncoiled_semisimple}
An uncoiled affine or periodic algebra that is generic is also semisimple.
\end{Corollary}
\proof
Let $\alg{A}$ be one of the uncoiled affine or periodic algebras and $\{\repW\}$ be its set of standard modules, as in \cref{prop:standards.uncoiled.odd,prop:standards.uncoiled.even}. In all cases, the dimension of the algebra is given by
\be 
\dim \alg{A} = \sum_{\{\repW\}} (\dim \repW)^2.
\ee
If all these standard modules are simple and pairwise inequivalent, by Wedderburn's Theorem, 
the finite-dimensional $\mathbb C$-algebra~$\repA$ is semisimple.
\eproof

\subsection{Skew sandwich cellularity}\label{sec:sandwich.cell}

In this section, we first introduce the notion of skew sandwich cellular algebras by adapting the ideas of \cite[Def.~2A.12]{T22} to include a skew involution as in \cite{HMR23}. We then describe how the affine, periodic and uncoiled Temperley--Lieb algebras fit within this framework.
 
\begin{Definition}[Skew sandwich cellular algebra]\label{def:SandwichCell}
Let $\alg{A}$ be an associative unital $\field{C}$-algebra, and let $\big((\mc{P},\sigma), (\mc{B},\mc{T}, \sigma_{\lambda},\sigma^{\lambda}),(\alg{M}_{\lambda}, \mc{M}(\lambda)),C,(\iota,\iota_{\lambda})\big)$ be a tuple of elements defined as follows:
\begin{itemize}
\item $\mc{P}$ is a poset with a partial order $>_{\mc{P}}$
and $\sigma$ is a poset involution of $\mc P$ that preserves the order $>_{\mc P}$;
\item $\mc{B} = \bigcup_{\lambda\in\mc{P}} \mc{B}(\lambda)$ and $\mc{T} = \bigcup_{\lambda\in\mc{P}} \mc{T}(\lambda)$ are collections of finite sets, and $\sigma_{\lambda}: \mc B(\lambda) \to \mc T(\sigma(\lambda))$ and $\sigma^{\lambda}: \mc T(\lambda) \to \mc B(\sigma(\lambda))$ are bijections such that $\sigma_{\sigma(\lambda)}\sigma^{\lambda} = \id_{\mc T(\lambda)}$ and $\sigma^{\sigma(\lambda)}\sigma_{\lambda} = \id_{\mc B(\lambda)}$;
\item for each $\lambda \in \mc P$, $\alg{M}_{\lambda}$ is a $\field{C}$-algebra, and  $\mc M(\lambda)$ is a basis of $\alg{M}_{\lambda}$; 
\item $C: \bigcup_{\lambda\in \mc P} \mc B(\lambda) \times \mc{M}(\lambda) \times  \mc T(\lambda) \to \alg{A}$ given by the assignment $(b,m,t)\mapsto c_{b,m,t}^{\lambda}$ 
 is an injective map; 
 \item $\iota:\alg{A}\to \alg{A}$ is an anti-involution of $\alg{A}$, and for each $\lambda\in \mc P$, $\iota_{\lambda}:\alg{M}_{\lambda}\to \alg{M}_{\sigma(\lambda)}$ is an anti-morphism of order 2.
\end{itemize}
We call the tuple 
$\big((\mc{P},\sigma), (\mc{B},\mc{T}, \sigma_{\lambda},\sigma^{\lambda}),(\alg{M}_{\lambda}, \mc{M}(\lambda)),C,(\iota,\iota_{\lambda})\big)$
a {\bf skew sandwich cell datum} 
for $\alg{A}$ if the following axioms hold:
\begin{enumerate}
\item[A1)]
the  set $\setS_{\alg{A}}= \{ c_{b,m,t}^{\lambda} \mid \lambda \in \mc P,\, b\in \mc B (\lambda), \, m\in \mathcal M(\lambda),\, t\in \mc T(\lambda)\}$ is a basis of $\alg{A}$ (the {\bf sandwich cellular basis});
\item[A2)]
for all $a\in \alg{A}$, there exist scalars $r_{b,b',m'}(a)\in \field{C}$
independent of $t$
such that
\begin{equation}\label{eq:CellAxiom action}
a\, c_{b,m,t}^{\lambda} \equiv \sum_{b'\in \mc B(\lambda),\, m'\in \mathcal M(\lambda)}  r_{b,b',m'}(a) c_{b',m',t}^{\lambda} \mod \alg{A}^{>_{\mc P} \lambda},
\end{equation}
where $\alg{A}^{>_{\mc P}\lambda}$ is the ideal generated by elements of the form $c_{u,n,v}^{\mu}$ for $\mu >_{\mc P} \lambda$, $u\in \mc B(\mu)$, $n\in \mc M(\mu)$ and $v\in \mc T(\mu)$; 
\item[A3)]
there exist a free $\alg{A}$-$\alg{M}_{\lambda}$--bimodule $\Delta(\lambda)$, a free $\alg{M}_{\lambda}$-$\alg{A}$--bimodule $\nabla (\lambda)$ and an $\alg{A}$-$\alg{A}$--bimodule isomorphism
\begin{equation}
\alg{A}_{\lambda} = \alg{A}^{\geq_{\mc P} \lambda}/\alg{A}^{>_{\mc P}\lambda} \simeq \Delta(\lambda)\otimes_{\alg{M}_{\lambda}} \nabla(\lambda);
\end{equation}
\item[A4)]
for all $\lambda\in \mc P$ and for all $t\in \mc T(\lambda)$, $m\in \mc{M}(\lambda)$, $b\in \mc B(\lambda)$,
\begin{equation}\label{eq:CellInvo}
\iota(c_{b,m,t}^{\lambda}) \equiv c_{\sigma^{\lambda}(t),\iota_{\lambda}(m),\sigma_{\lambda}(b)}^{\sigma(\lambda)} \mod \alg{A}^{>_{\mc P}\lambda}.
\end{equation}
\end{enumerate}
We call an algebra $\alg{A}$ equipped with a skew sandwich cellular datum, a {\bf skew sandwich cellular algebra}.
\end{Definition}
\noindent We now compare this definition to the other similar concepts already present in the literature.
\begin{Remark}
Let $\alg{A}$ be an associative unital $\mathbb C$-algebra.
\begin{itemize}
\item A tuple $\big((\mc{P}, \sigma), (\mc{B},\mc{T}, \sigma_\lambda,\sigma^\lambda),(\alg{M}_{\lambda}, \mc{M}(\lambda)),C\big)$ respecting Axioms~A1--A3) with $\sigma$, $\sigma_\lambda$ and $\sigma^\lambda$ all given by the identity on their respective spaces is called a {\bf sandwich cellular datum}, and $\alg{A}$ is called {\bf sandwich cellular} in the sense of~\cite{T22}. 
\item A skew cellular datum $\big((\mc{P},\sigma), (\mc{B},\mc{T}),(\alg{M}_{\lambda}, \mc{M}(\lambda)),C,(\iota,\iota_{\lambda})\big)$ with $\sigma$, $\sigma_\lambda$ and $\sigma^\lambda$ all given by the identity on their respective spaces and with $\mc B = \mc T$ is called an {\bf involutive sandwich cellular datum}, and $\alg{A}$ is called {\bf involutive sandwich cellular} in the sense of~\cite{T22}.
\item An involutive sandwich cellular $\alg{A}$ where $\alg{M}_{\lambda}$ is commutative for all $\lambda \in \mc P$ is called {\bf affine cellular} in the sense of~\cite{KX12}.
\item A skew sandwich cellular algebra $\alg{A}$ with $\alg{M}_{\lambda} = \mathbb C$ and $\setB=\setT$ for all $\lambda \in \mc P$ is a {\bf skew cellular algebra} in the (ungraded) sense of~{\cite[Def.~2.2]{HMR23}}.
\item A skew cellular algebra $\alg{A}$ with $\sigma = \id$ is a {\bf cellular algebra} in the sense of~\cite{GL96}.
\end{itemize}
\end{Remark}
\noindent The definition of affine cellularity was stated in~\cite{KX12} using the language of cell ideals. It is, however, easier for computational purposes to have the cellular datum as in \cref{def:SandwichCell}.
\begin{Proposition}[Prop.~2.5 \cite{KX12}]
\label{Prop:atl_affine_cell}
The affine Temperley--Lieb algebra $\atl_n(\beta)$ is affine cellular for all $n\in \mathbb N$ and all $\beta\in \mathbb C$.
\end{Proposition}
We do not repeat the proof of \cite{KX12}, but nonetheless describe the structure of the ideals of $\atl_n(\beta)$ and give its cellular datum. The algebra $\atl_n(\beta)$ has the filtration of two-sided ideals
\be
\label{eq:filtration.of.ideals}
\repJ_{n_0} \subset \repJ_{n_0+2} \subset \dots \subset \repJ_n = \atl_n(\beta), 
\qquad
n_0 = \left\{\begin{array}{ll}
0 & n \textrm{ even,}\\[0.1cm]
1 & n \textrm{ odd,}
\end{array}\right.
\ee
where
\be
\label{def:Jd}
\repJ_d = \{c_1 E_d c_2\,|\, c_1,c_2 \in \atl_n(\beta) \}, 
\qquad
E_d = e_0 e_2 \cdots e_{n-d-2},
\qquad
d=n_0,n_0+2,\dots, n.
\ee
The associated cellular datum is 
\begin{subequations}
\begin{alignat}{2}
\mc P &=\{n_0>_{\mc P}n_0+2>_{\mc P}\dots >_{\mc P} n\},\\[0.15cm]
\qquad
\setB(d) &= \setT(d) = \setB_{n,d}\,,\\[0.15cm]
\qquad
\alg{M}_d &= \left\{\begin{array}{ll}
\mathbb C[f] & d=0,
\\[0.15cm]
\mathbb C [\Omega_d,\Omega_d^{-1}] / \langle \Omega_d\Omega_d^{-1} = \Omega_d^{-1}\Omega_d = \id_d \rangle & d>0,
\end{array}\right.
\end{alignat} 
\end{subequations}
where the partial order is the reverse order on the natural numbers: $d>_{\mc P} d'$ if and only if $d<_{\mathbb N} d'$. Moreover, we denote by $\id_d$, $\Omega_d$ and $\Omega^{-1}_d$ the identity and the two translation generators on $d$ strands. The involution $\iota$ is the vertical reflection of the diagrams, which on $\alg{A}_d$ sends $\Omega_d^{\pm 1} \mapsto \Omega_d^{\mp 1}$ and $f\mapsto f$. The cell bases are those given in~\eqref{eq:SandwichBaseaTL}.\medskip

We illustrate the assignment with an example in $\atl_{6}(\beta)$:
\be
\label{eq:ExMultiaTL}
\begin{array}{r}
\psset{unit=0.8cm}
(b,m,t) = \Big( \
\begin{pspicture}[shift=-0.05](0,0)(2.4,0.6)
\psline[linewidth=0.5pt](0,0)(2.4,0)
\psline[linecolor=blue,linewidth=1.5pt]{-}(0.2,0)(0.2,0.5)
\psline[linecolor=blue,linewidth=1.5pt]{-}(0.6,0)(0.6,0.5)
\psarc[linecolor=blue,linewidth=1.5pt]{-}(1.6,0){0.2}{0}{180}
\psbezier[linecolor=blue,linewidth=1.5pt](1.0,0)(1.0,0.7)(2.2,0.7)(2.2,0)
\end{pspicture}
\ ,\id_2,\  
\begin{pspicture}[shift=-0.05](0,0)(2.4,0.6)
\psline[linewidth=0.5pt](0,0)(2.4,0)
\psline[linecolor=blue,linewidth=1.5pt]{-}(0.6,0)(0.6,0.5)
\psline[linecolor=blue,linewidth=1.5pt]{-}(1.8,0)(1.8,0.5)
\psarc[linecolor=blue,linewidth=1.5pt]{-}(0,0){0.2}{0}{90}
\psarc[linecolor=blue,linewidth=1.5pt]{-}(2.4,0){0.2}{90}{180}
\psarc[linecolor=blue,linewidth=1.5pt]{-}(1.2,0){0.2}{0}{180}
\end{pspicture}\ \Big)\ , \qquad
a = \
\begin{pspicture}[shift=-0.4](0,-0.5)(2.4,0.5)
\pspolygon[fillstyle=solid,fillcolor=lightlightblue,linewidth=0pt](0,-0.5)(2.4,-0.5)(2.4,0.5)(0,0.5)
\psbezier[linecolor=blue,linewidth=1.5pt](0.2,-0.5)(0.2,0.1)(1.4,0.1)(1.4,-0.5)
\psarc[linecolor=blue,linewidth=1.5pt]{-}(0.8,-0.5){0.2}{0}{180}
\psarc[linecolor=blue,linewidth=1.5pt]{-}(0.8,0.5){0.2}{180}{0}
\psarc[linecolor=blue,linewidth=1.5pt]{-}(1.6,0.5){0.2}{180}{0}
\psbezier[linecolor=blue,linewidth=1.5pt]{-}(1.8,-0.5)(1.8,0)(2.2,0)(2.2,0.5)
\psbezier[linecolor=blue,linewidth=1.5pt]{-}(2.2,-0.5)(2.2,-0.3)(2.35,-0.05)(2.41,0)
\psbezier[linecolor=blue,linewidth=1.5pt]{-}(0.2,0.5)(0.2,0.3)(0.05,0.05)(-0.01,0)
\psframe[fillstyle=solid,linecolor=white,linewidth=0pt](-0.05,-0.51)(-0.005,0.51)
\psframe[fillstyle=solid,linecolor=white,linewidth=0pt](2.405,-0.51)(2.45,0.51)
\end{pspicture} 
\ + \
\begin{pspicture}[shift=-0.4](0,-0.5)(2.4,0.5)
\pspolygon[fillstyle=solid,fillcolor=lightlightblue,linewidth=0pt](0,-0.5)(2.4,-0.5)(2.4,0.5)(0,0.5)
\psarc[linecolor=blue,linewidth=1.5pt]{-}(0.4,0.5){0.2}{180}{0}
\psarc[linecolor=blue,linewidth=1.5pt]{-}(1.6,0.5){0.2}{180}{0}
\psarc[linecolor=blue,linewidth=1.5pt]{-}(0.8,-0.5){0.2}{0}{180}
\psarc[linecolor=blue,linewidth=1.5pt]{-}(1.6,-0.5){0.2}{0}{180}
\psbezier[linecolor=blue,linewidth=1.5pt]{-}(0.2,-0.5)(0.2,0.1)(2.2,0.1)(2.2,-0.5)
\psbezier[linecolor=blue,linewidth=1.5pt]{-}(-0.02,0.17)(0.1,0)(1.0,-0.06)(1.0,0.5)
\psbezier[linecolor=blue,linewidth=1.5pt]{-}(2.2,0.5)(2.2,0.32)(2.35,0.22)(2.42,0.2)
\psframe[fillstyle=solid,linecolor=white,linewidth=0pt](-0.05,-0.51)(-0.005,0.51)
\psframe[fillstyle=solid,linecolor=white,linewidth=0pt](2.405,-0.51)(2.45,0.51)
\end{pspicture}\ ,
\\[0.6cm]
\Longrightarrow \quad
\psset{unit=0.8cm}
c_{b,m,t}^2 = \
\begin{pspicture}[shift=-0.4](0,-0.5)(2.4,0.5)
\pspolygon[fillstyle=solid,fillcolor=lightlightblue,linewidth=0pt](0,-0.5)(2.4,-0.5)(2.4,0.5)(0,0.5)
\psarc[linecolor=blue,linewidth=1.5pt]{-}(1.6,-0.5){0.2}{0}{180}
\psbezier[linecolor=blue,linewidth=1.5pt](1.0,-0.5)(1.0,0.1)(2.2,0.1)(2.2,-0.5)
\psarc[linecolor=blue,linewidth=1.5pt]{-}(0,0.5){0.2}{-90}{0}
\psarc[linecolor=blue,linewidth=1.5pt]{-}(1.2,0.5){0.2}{180}{0}
\psarc[linecolor=blue,linewidth=1.5pt]{-}(2.4,0.5){0.2}{180}{270}
\psbezier[linecolor=blue,linewidth=1.5pt](0.2,-0.5)(0.2,0)(0.6,0)(0.6,0.5)
\psbezier[linecolor=blue,linewidth=1.5pt](0.6,-0.5)(0.6,0)(1.8,0)(1.8,0.5)
\end{pspicture}
\ , \qquad
a\, c_{b,m,t}^2 = \beta \
\begin{pspicture}[shift=-0.4](0,-0.65)(2.4,0.65)
\pspolygon[fillstyle=solid,fillcolor=lightlightblue,linewidth=0pt](0,-0.65)(2.4,-0.65)(2.4,0.65)(0,0.65)
\psarc[linecolor=blue,linewidth=1.5pt]{-}(0.8,-0.65){0.2}{0}{180}
\psbezier[linecolor=blue,linewidth=1.5pt](0.2,-0.65)(0.2,-0.05)(1.4,-0.05)(1.4,-0.65)
\psarc[linecolor=blue,linewidth=1.5pt]{-}(0,0.65){0.2}{-90}{0}
\psarc[linecolor=blue,linewidth=1.5pt]{-}(1.2,0.65){0.2}{180}{0}
\psarc[linecolor=blue,linewidth=1.5pt]{-}(2.4,0.65){0.2}{180}{270}
\psline[linecolor=blue,linewidth=1.5pt](1.8,-0.65)(1.8,0.65)
\psbezier[linecolor=blue,linewidth=1.5pt](0.6,0.65)(0.6,0.1)(0.1,0.03)(-0.02,-0.17)
\psbezier[linecolor=blue,linewidth=1.5pt](2.2,-0.65)(2.2,-0.3)(2.37,-0.16)(2.42,-0.13)
\psframe[fillstyle=solid,linecolor=white,linewidth=0pt](-0.05,-0.651)(-0.005,0.651)
\psframe[fillstyle=solid,linecolor=white,linewidth=0pt](2.405,-0.651)(2.45,0.651)
\end{pspicture}
\ + \beta \
\begin{pspicture}[shift=-0.4](0,-0.65)(2.4,0.65)
\pspolygon[fillstyle=solid,fillcolor=lightlightblue,linewidth=0pt](0,-0.65)(2.4,-0.65)(2.4,0.65)(0,0.65)
\psarc[linecolor=blue,linewidth=1.5pt]{-}(0,0.65){0.2}{-90}{0}
\psarc[linecolor=blue,linewidth=1.5pt]{-}(2.4,0.65){0.2}{180}{270}
\psarc[linecolor=blue,linewidth=1.5pt]{-}(1.2,0.65){0.2}{180}{0}
\psbezier[linecolor=blue,linewidth=1.5pt]{-}(0.6,0.65)(0.6,0.05)(1.8,0.05)(1.8,0.65)
\psarc[linecolor=blue,linewidth=1.5pt]{-}(0.8,-0.65){0.2}{0}{180}
\psarc[linecolor=blue,linewidth=1.5pt]{-}(1.6,-0.65){0.2}{0}{180}
\psbezier[linecolor=blue,linewidth=1.5pt]{-}(0.2,-0.65)(0.2,-0.05)(2.2,-0.05)(2.2,-0.65)
\psframe[fillstyle=solid,linecolor=white,linewidth=0pt](-0.05,-0.651)(-0.005,0.651)
\psframe[fillstyle=solid,linecolor=white,linewidth=0pt](2.405,-0.651)(2.45,0.651)
\psline[linecolor=blue, linewidth=1.5pt](0,0)(2.4,0)
\end{pspicture}\ .
\end{array}
\ee
The result of the product $a\, c_{b,m,t}^2$ has a first term involving a connectivity $c^2_{b',m',t}$ with $b' = \,
\psset{unit=0.6cm}
\begin{pspicture}[shift=-0.05](0,0)(2.4,0.6)
\psline[linewidth=0.5pt](0,0)(2.4,0)
\psline[linecolor=blue,linewidth=1.5pt]{-}(1.8,0)(1.8,0.5)
\psline[linecolor=blue,linewidth=1.5pt]{-}(2.2,0)(2.2,0.5)
\psarc[linecolor=blue,linewidth=1.5pt]{-}(0.8,0){0.2}{0}{180}
\psbezier[linecolor=blue,linewidth=1.5pt](0.2,0)(0.2,0.7)(1.4,0.7)(1.4,0)
\end{pspicture}\ $ and $m' = \Omega_2$, and the coefficient is $r_{b,b',m'}(a) = \beta$. The second term has no through-lines, and is therefore absent if we consider the action of the algebra modulo the elements with fewer than two through-lines, as in~\eqref{eq:CellAxiom action}.\medskip

An important feature of affine cellularity is that the cellular datum of an affine cellular algebra may not be unique. In particular, for $\atl_n(\beta)$ with $n$ even, one can refine the cellular order by decomposing the $d=0$ cell further in terms of the two-sided ideals
\be
\repJ_{0_r} = \{c_1 E(\Omega E)^{r} c_2\,|\, c_1,c_2 \in \atl_n(\beta) \},\qquad r \in \mathbb Z_{\ge 0},
\ee
which satisfy the filtration
\be
\dots \subset \repJ_{0_2} \subset \repJ_{0_1} \subset \repJ_{0_0} = \repJ_0.
\ee
The resulting poset is 
\be
\mc P = \{\dots >_{\mc P} 0_2 >_{\mc P} 0_1 >_{\mc P} 0 >_{\mc P} 2 >_{\mc P} \dots >_{\mc P} n\}
\ee 
and is infinite. For each $r$, the bottom and top sets are $\setB(0_r) = \setT(0_r) = \setB_{n,0} $ and the middle algebras are all $\repM(0_r)=\mathbb C$. The assignment $C_{0_r}$ sandwiches $r$ non-contractible loops between two link states $b,t\in \setB_{n,0}$. The action of a connectivity cannot reduce the number of non-contractible loops, so the order is respected. Moreover, because the non-contractible loops separate the top and bottom, it is easy to check that Axiom~A2) is satisfied. In the proof of the next proposition, we will apply a similar argument to show that the algebra $\ptl_n(\beta)$ for $n$ even is skew sandwich cellular.

\begin{Proposition}
\label{prop:pTL.skew.sandwich}
Let $n\in \mathbb Z_{\ge 0}$ and $\beta\in \mathbb C$. The periodic Temperley--Lieb algebra $\ptl_n(\beta)$ is affine cellular for $n$ odd, whereas it is skew sandwich cellular for $n$ even.
\end{Proposition}
\proof
For $n$ odd, the algebra $\ptl_n(\beta)$ has the filtration of ideals 
\be
\label{eq:Jdptl}
\repJ_{1} \subset \repJ_{3} \subset \dots \subset \repJ_n = \ptl_n(\beta), 
\qquad
\repJ_d = \{c_1 E_d c_2\,|\, c_1,c_2 \in \ptl_n(\beta) \}.
\ee
The cellular datum is
\begin{subequations}
\label{eq:cell.data.ptl.odd}
\begin{alignat}{2}
\mc P &=\{1>_{\mc P}3 >_{\mc P}\dots >_{\mc P} n\},\\[0.15cm]
\qquad
\setB(d) &= \setT(d) = \wh\setB_{n,d} ,\\[0.15cm]
\qquad
\alg{M}_d &= \left\{\begin{array}{ll}
\mathbb C[\Omega_d^2, \Omega_d^{-2}]/\langle \Omega_d^2\Omega_d^{-2} = \Omega_d^{-2}\Omega_d^2= \id_d \rangle  & d<n.
\\[0.15cm]
\mathbb C & d=n.
\end{array}\right.
\end{alignat} 
\end{subequations}
The assignment $C_d$ produces the diagrams $\wh c_{b,m,t}^{\,d}$ in the sets $\wh\setS_d$ defined in \eqref{eq:sandwich.whS} and produces the cell bases given in \eqref{eq:sandwich.basis.ptl.v2}. The involution $\iota$ is the vertical reflection of the diagrams. \medskip
 
Axiom A2) is satisfied because multiplication can only decrease the number of through-lines, showing that the equivalence is respected. Moreover, one can check that the coefficients in the multiplication $a\wh c_{b,m,t}^{\, d}$ depend on the bottom state $b$ and on the connectivity $a$ that is being multiplied, but do not depend on the top state $t$, since 
these coefficients depend only on the number of contractible loops.\medskip 

We now turn to Axiom A3). Let $\alg{A} = \ptl_n(\beta)$ with $n$ odd. We need to construct an isomorphism of $\alg{A}$-$\alg{A}$--bimodules $\Delta(d)\otimes_{\alg{M}_{d}}\nabla(d) \simeq \alg{A}^{\geq_{\mc P} d}/\alg{A}^{>_{\mc P}d}$. For $d\in \mc P$, the module $\Delta(d)$ is the module given as vector space by $\repW_{n,d,z}\otimes_{\field{C}} \alg{M}_{d}$. It inherits the regular right $\alg{M}_{d}$-module structure of $\alg{M}_d$, and is a left $\alg{A}$-module since $\repW_{n,d,z}$ is itself such a module. The module $\nabla(d)$ is the dual of this module, namely $\alg{M}_{d} \otimes_{\field{C}} \repW_{n,d,z^{-1}}$. The isomorphism is given by a construction from~\cite[Prop.~2.2]{KX12}. We note that $\Delta(d)\otimes_{\alg{M}_{d}}\nabla(d) \simeq \repW_{n,d,z} \otimes_{\field{C}} \alg{M}_{d} \otimes_{\field{C}} \repW_{n,d,z^{-1}}$ as vector spaces. Given a bilinear form 
$\psi_d:\repW_{n,d,z^{-1}}\otimes_{\field{C}}\repW_{n,d,z} \to \alg{M}_{d}$,
there exists a multiplication on $\repW_{n,d,z} \otimes_{\field{C}} \alg{M}_{d} \otimes_{\field{C}} \repW_{n,d,z^{-1}}$, given by the assignment 
\be
\label{eq:BiliAssignmentAffine}
(b,m,t)\times(b',m',t') \mapsto (b,m\psi_d(t,b')m',t'),
\ee 
for $b,b',t,t'\in \setB_{n,d}$ and $m,m'\in \mc M(\lambda)$.
\medskip

In the present case, the bilinear form takes two link states $t,b$ in $\wh \setB_{n,d}$ and computes their diagrammatic product $\iota(t)b$ where $\iota(t)$ is the vertical reflection of $t$. The product $\psi_d(t,b)$ is non-zero only if each defect of $b$ is connected to a defect of $t$. In this case, $\psi_d(t,b)$ is equal to $\beta^m \Omega_d^\ell$, where $m \in \mathbb Z_{\ge 0}$ is the number of contractible loops formed in the process, and $\ell \in \mathbb Z$ is the total winding number of the defects. This number $\ell$ measures the number of times the defects cross the periodic boundary condition, and is a sum of contributions $+1$ and $-1$ for each defect crossing it towards the left and right, respectively, as it goes downwards. For example in $\ptl_5(\beta)$, we have
\begin{subequations}
\begin{alignat}{2}
v &= \ 
\begin{pspicture}[shift=-0.05](0,0)(2.0,0.6)
\psline[linewidth=0.5pt](0,0)(2.0,0)
\psline[linecolor=blue,linewidth=1.5pt]{-}(1.0,0)(1.0,0.5)
\psarc[linecolor=blue,linewidth=1.5pt]{-}(0,0){0.2}{0}{90}
\psarc[linecolor=blue,linewidth=1.5pt]{-}(2.0,0){0.2}{90}{180}
\psbezier[linecolor=blue,linewidth=1.5pt]{-}(0.6,0)(0.6,0.36)(0.05,0.4)(0,0.4)
\psbezier[linecolor=blue,linewidth=1.5pt]{-}(1.4,0)(1.4,0.36)(1.85,0.4)(2.0,0.4)
\end{pspicture}
\ ,
\qquad&
\psi_1(v,v) &= \ 
\begin{pspicture}[shift=-0.55](0,-0.6)(2.0,0.6)
\psline[linewidth=0.5pt](0,0)(2.0,0)
\psline[linecolor=blue,linewidth=1.5pt]{-}(1.0,0)(1.0,0.5)
\psarc[linecolor=blue,linewidth=1.5pt]{-}(0,0){0.2}{0}{90}
\psarc[linecolor=blue,linewidth=1.5pt]{-}(2.0,0){0.2}{90}{180}
\psbezier[linecolor=blue,linewidth=1.5pt]{-}(0.6,0)(0.6,0.36)(0.05,0.4)(0,0.4)
\psbezier[linecolor=blue,linewidth=1.5pt]{-}(1.4,0)(1.4,0.36)(1.85,0.4)(2.0,0.4)
\psline[linecolor=blue,linewidth=1.5pt]{-}(1.0,0)(1.0,-0.5)
\psarc[linecolor=blue,linewidth=1.5pt]{-}(0,0){0.2}{270}{0}
\psarc[linecolor=blue,linewidth=1.5pt]{-}(2.0,0){0.2}{180}{270}
\psbezier[linecolor=blue,linewidth=1.5pt]{-}(0.6,0)(0.6,-0.36)(0.05,-0.4)(0,-0.4)
\psbezier[linecolor=blue,linewidth=1.5pt]{-}(1.4,0)(1.4,-0.36)(1.85,-0.4)(2.0,-0.4)
\end{pspicture} \ = \beta^2 \id_1,
\\
w &= \
\begin{pspicture}[shift=-0.05](0,0)(2.0,0.6)
\psline[linewidth=0.5pt](0,0)(2.0,0)
\psbezier[linecolor=blue,linewidth=1.5pt]{-}(0.6,0)(0.6,0.4)(1.8,0.4)(2.0,0.4)
\psarc[linecolor=blue,linewidth=1.5pt]{-}(1.2,0){0.2}{0}{180}
\psbezier[linecolor=blue,linewidth=1.5pt]{-}(0,0.35)(0.5,0.35)(0.6,0.45)(0.6,0.55)
\psarc[linecolor=blue,linewidth=1.5pt]{-}(0.0,0){0.2}{0}{90}
\psarc[linecolor=blue,linewidth=1.5pt]{-}(2.0,0){0.2}{90}{180}
\end{pspicture}
\ , 
\qquad &
\psi_1(v,w) & =\  
\begin{pspicture}[shift=-0.55](0,-0.6)(2.0,0.6)
\psline[linewidth=0.5pt](0,0)(2.0,0)
\psbezier[linecolor=blue,linewidth=1.5pt]{-}(0.6,0)(0.6,0.4)(1.8,0.4)(2.0,0.4)
\psarc[linecolor=blue,linewidth=1.5pt]{-}(1.2,0){0.2}{0}{180}
\psbezier[linecolor=blue,linewidth=1.5pt]{-}(0,0.35)(0.5,0.35)(0.6,0.45)(0.6,0.55)
\psarc[linecolor=blue,linewidth=1.5pt]{-}(0.0,0){0.2}{0}{90}
\psarc[linecolor=blue,linewidth=1.5pt]{-}(2.0,0){0.2}{90}{180}
\psline[linecolor=blue,linewidth=1.5pt]{-}(1.0,0)(1.0,-0.5)
\psarc[linecolor=blue,linewidth=1.5pt]{-}(0,0){0.2}{270}{0}
\psarc[linecolor=blue,linewidth=1.5pt]{-}(2.0,0){0.2}{180}{270}
\psbezier[linecolor=blue,linewidth=1.5pt]{-}(0.6,0)(0.6,-0.36)(0.05,-0.4)(0,-0.4)
\psbezier[linecolor=blue,linewidth=1.5pt]{-}(1.4,0)(1.4,-0.36)(1.85,-0.4)(2.0,-0.4)
\end{pspicture}\ = \beta\, \Omega_1^2 .
\end{alignat}
\end{subequations}
Then $\repW_{n,d,z}\otimes_{\field{C}}\alg{M}_d\otimes_{\field{C}}\repW_{n,d,z^{-1}}$ is sent to $\alg{A}^{\geq_{\mc P}d}/\alg{A}^{>_{\mc P}d}$ by the map~\eqref{eq:BiliAssignmentAffine}, since it contains only the connectivities of $\alg{A}$ with exactly $d$ through-lines and the multiplication sends to zero any element with less than $d$ through-lines. This constructs the isomorphism of $\alg{A}$-$\alg{A}$--bimodules required by Axiom~A3).\medskip

Lastly, the anti-involution $\iota_d$ is defined on the middle algebras $\alg{M}_d$ by $\Omega\mapsto \Omega^{-1}$ and $f\mapsto f$. Then indeed, we see that $\iota(c_{b,m,t}^d) = c_{t,\iota_d(m),b}^d$, so Axiom~A4) holds. This ends the proof that $\ptl_n(\beta)$ is affine cellular for $n$ odd.\medskip

For $n$ even, we have the filtration 
\be
\repJ_{0} \subset \repJ_{2} \subset \dots \subset \repJ_n = \ptl_n(\beta),
\ee
with $\repJ_d$ defined as in \eqref{eq:Jdptl}. The main difference in the construction arises because of our inability in~\eqref{eq:sandwich.basis.ptl.v2} to write the set of diagrams spanning $\repJ_0$ in terms of the set $\wh\setS_0(\{f^r\,|\,r \in 2 \mathbb Z_{\ge 0}\})$. To circumvent this issue, we note that the description of $\repJ_0$ can be refined by defining the ideals
\begin{subequations}
\begin{alignat}{2}
&r \in 2\mathbb Z_{\ge 0}: \qquad
&&\left\{\begin{array}{l}
\repJ^{\tinyx{E}}_{0_r} = \{c_1 E(FE)^{r/2} c_2\,|\, c_1,c_2 \in \ptl_n(\beta) \},
\\[0.25cm]
\repJ^{\tinyx{F}}_{0_r} = \{c_1 F(EF)^{r/2} c_2\,|\, c_1,c_2 \in \ptl_n(\beta) \},
\end{array}\right.
\\[0.3cm]
&r \in 2\mathbb Z_{\ge 0}+1: \qquad
&&\left\{\begin{array}{l}
\repJ^{\tinyx{E}}_{0_r} = \{c_1 (EF)^{(r+1)/2} c_2\,|\, c_1,c_2 \in \ptl_n(\beta) \},
\\[0.25cm]
\repJ^{\tinyx{F}}_{0_r} = \{c_1 (FE)^{(r+1)/2} c_2\,|\, c_1,c_2 \in \ptl_n(\beta) \}.
\end{array}\right.
\end{alignat}
\end{subequations}
With these definitions, the diagrams in $\repJ^{\tinyx{E}}_{0_r}$ and $\repJ^{\tinyx{F}}_{0_r}$ all have at least $r$ non-contractible loops. These ideals satisfy
\be
\dots \subset (\repJ^{\tinyx E}_{0_2} \oplus \repJ^{\tinyx F}_{0_2}) \subset (\repJ^{\tinyx E}_{0_1}\oplus \repJ^{\tinyx F}_{0_1}) \subset (\repJ^{\tinyx E}_{0_0} \oplus \repJ^{\tinyx F}_{0_0}) \subset \repJ_2 \subset \repJ_4 \subset \dots \subset \repJ_n = \ptl_n(\beta).
\ee
Crucially, the right and left action of $\ptl_n(\beta)$ on $\repJ^{\tinyx E}_{0_r}$ can produce diagrams in $\repJ^{\tinyx E}_{0_{r'}}$ or $\repJ^{\tinyx F}_{0_{r'}}$ for some $r'<r$, but not diagrams in $\repJ^{\tinyx F}_{0_r}$. The same statement holds with $E$ and $F$ interchanged. The resulting poset is therefore
\be
\mc P = \big\{\dots  >_{\mc P} 0^{\tinyx{E}}_3, 0^{\tinyx{F}}_3 >_{\mc P} 0^{\tinyx{E}}_2, 0^{\tinyx{F}}_2 >_{\mc P} 0^{\tinyx{E}}_1, 0^{\tinyx{F}}_1 >_{\mc P} 0^{\tinyx{E}}_0, 0^{\tinyx{F}}_0 >_{\mc P} 2 >_{\mc P} 4 >_{\mc P} \dots >_{\mc P} n\big\}.
\ee 
For $d\ge 2$, $\setB(d)$, $\setT(d)$, and $\alg{M}_d$ are as in \eqref{eq:cell.data.ptl.odd}, and the assignment~$C$ is constructed in the same way as in the case with $n$ odd. For the cells with labels $d=0_r^{\tinyx E}$ and $d=0_r^{\tinyx F}$, we instead have
\begin{subequations}
\label{eq:sets.ptla.skew}
\begin{alignat}{2}
&r \in 2 \mathbb Z_{\ge 0}: \qquad&&
\left\{\begin{array}{l}
\mc B = \setB^{\tinyx 1}_{n,0}\,,\quad
\mc T = \setB^{\tinyx 1}_{n,0}\,, \qquad d=0_r^{\tinyx E},
\\[0.15cm]
\mc B = \setB^{\tinyx 0}_{n,0}\,,\quad
\mc T = \setB^{\tinyx 0}_{n,0}\,, \qquad d=0_r^{\tinyx F},
\end{array}\right.
\\[0.15cm]
&r \in 2 \mathbb Z_{\ge 0}+1: \qquad&&
\left\{\begin{array}{l}
\mc B = \wh\setB^{\tinyx 1}_{n,0}\,,\quad
\mc T = \wh\setB^{\tinyx 0}_{n,0}\,, \qquad d=0_r^{\tinyx E},
\\[0.15cm]
\mc B = \wh\setB^{\tinyx 0}_{n,0}\,,\quad
\mc T = \wh\setB^{\tinyx 1}_{n,0}\,, \qquad d=0_r^{\tinyx F},
\end{array}\right.
\end{alignat}
\end{subequations}
and $\alg{M}_{0^{\tinyx E}_r} = \alg{M}_{0^{\tinyx F}_r} = \mathbb C$ in all cases. In these cases, the assignment $C$ sandwiches the link states $b$ and $t$ on $r$ non-contractible loops for $r$ even, and $r-1$ non-contractible loops for $r$ odd. In the latter case, either $b$ or $t$ contains an extra non-contractible loop, as one of them is in 
$\wh \setB^{\tinyx 1}_{n,0}$, so the total number of such loops is also $r$. The resulting diagram is even in all cases.\medskip

Axiom~A2) holds from the above remark that the action of $\ptl_n(\beta)$ cannot decrease the non-contractible loops, so the equivalence is respected. Furthermore, each time the action of $\ptl_n(\beta)$ on a link state produces a state with the opposite parity, a non-contractible loop is produced, and has the effect of increasing $r$ and thus changing the order. It is also easy to see that the coefficients in the multiplication do not depend on the top state $t$. \medskip

For Axiom~A3), the bilinear form $\psi_d(t,b)$ for $d\geq 2$ is constructed in the same way as in the case with $n$~odd. For $d= 0_r^{\tinyx{E}}$ and $0_r^{\tinyx{F}}$, we construct the bilinear form as before from the diagram $\iota(t)b$, where $t$ and $b$ are link states in the top and bottom bases, respectively. The resulting diagram is even and contains an even number of non-contractible loops. If this number is non-zero, the result is $\psi_d(t,b) = 0$. Otherwise, $\psi_d(t,b) = \beta^m$, where $m$ is the number of contractible loops. We note in particular that $\psi_d(t,b) = 0$ for $d= 0_r^{\tinyx{E}}$ and $d=0_r^{\tinyx{F}}$ with $r$ odd, for all states $b$ and $t$. This is consistent with the fact that the ideals $\repJ_{0_r}^{\tinyx{E}}$ and $\repJ_{0_r}^{\tinyx{F}}$  are nilpotent in these cases.\medskip

The involution is the vertical reflection of the diagrams. It sends $c_{b,m,t}$ to $c_{t,\iota(m),b}$ for $d\geq 2$, and also for $d=0_r^{\tinyx{E}}$ and $d=0_r^{\tinyx{F}}$ with $r$ even. For $r$ odd, it instead maps the diagrams indexed by $0_r^{\tinyx{E}}$ to those indexed by $0_r^{\tinyx{F}}$, and vice versa. The permutation of cells thus acts as the identity in all cases except for $0_r^{\tinyx{E}}$ and $0_r^{\tinyx{F}}$, whereby it permutes these two cells, and both $\sigma_\lambda$ and $\sigma^\lambda$ act as the identity on their respective sets. This shows that Axiom~4) holds.\medskip

This ends the proof that $\ptl_n(\beta)$ is skew sandwich cellular for $n$ even.
\eproof

We deduce from this proof that the algebra $\ptl_n(\beta)/\{E=F=0\}$ without connectivities with $d=0$ through-lines is a proper affine cellular algebra. The following proposition addresses the same questions for the uncoiled algebras.

\begin{Proposition}\label{prop:uaTLupTLSandwichCellular}
The uncoiled algebras $\qatl_n(\beta,\gamma)$, $\qptl_n(\beta,\gamma)$, $\qatla_n(\beta, \alpha)$, $\qptla_n(\beta, \alpha)$ with $\alpha\neq 0$, $\qatlb_n(\beta,\gamma)$ and $\qptlb_n(\beta,\gamma)$ are affine cellular.
\end{Proposition}
\proof
The posets for the uncoiled algebras are
\be
\mc P = \left\{\begin{array}{ll}
\mc P^{\tinyx{0}} = \{1,3,\dots,n-2,n\} 
& \qatl_{n}(\beta,\gamma) \textrm{ and } \qptl_{n}(\beta,\gamma),
\\[0.1cm]
\mc P^{\tinyx{1}} = \{ 0,2,\dots, n-2,n\}
& \qatla_{n}(\beta,\alpha), \textrm{and } \qptla_{n}(\beta,\alpha) \textrm{ with } \alpha \neq 0,
\\[0.1cm]
\mc P^{\tinyx{2}} = \{2,4,\dots, n-2,n\}
& \qatlb_{n}(\beta,\gamma) \textrm{ and } \qptlb_{n}(\beta,\gamma).
\end{array}\right.
\ee
For the three uncoiled affine algebras, the top and bottom bases $\mc T(d)$ and $\mc B(d)$ are the sets of link states $\setB_{n,d}$, for all $d$. For the three uncoiled periodic algebras, the sets $\mc T(d)$ and $\mc B(d)$ are the sets of link states $\wh \setB_{n,d}$ for all $d$. The six uncoiled algebras also differ in their middle algebras $\alg{M}_{d}$, which we denote using upper labels indicating the algebra that they correspond to. They are given by
\begingroup
\allowdisplaybreaks
\begin{subequations}
\begin{alignat}{3}
\qatl_n(\beta,\gamma):\quad && \alg{M}_d^{{\sf a}} &= \field{C}[\Omega_d]/\langle \Omega_d^d = \gamma\,\id_d\rangle \quad d\in \mc P^{\tinyx{0}};\\[0.1cm]
\qatla_n(\beta,\alpha):\quad && \alg{M}_d^{{\sf a},\tinyx{1}}&= 
\left\{\begin{array}{ll}
\field{C}[f]/\langle f = \alpha\, \id_d\rangle & d=0,
\\[0.1cm]
\field{C}[\Omega_d]/\langle \Omega_d^d = \id_d\rangle & d\in \mc P^{\tinyx{1}}\setminus
\{0\};
\end{array}\right.
\\[0.1cm]
\qatlb_n(\beta,\gamma):\quad &&  \alg{M}_d^{{\sf a},\tinyx{2}} &= \field{C}[\Omega_d]/\langle \Omega_d^d = \gamma\,\id_d\rangle \quad d\in \mc P^{\tinyx{2}};\\[0.1cm]
\qptl_n(\beta,\gamma):\quad && 
\alg{M}_d^{{\sf p}} &= \left\{\begin{array}{ll}
\field{C}[\Omega_d^2]/\langle \Omega_d^{2d} = \gamma^2\id_d\rangle& d\in \mc P^{\tinyx{0}} \setminus\{n\},
\\[0.1cm]
\field{C}&d=n;
\end{array}\right.
\\[0.1cm] 
\qptla_n(\beta,\alpha) \textrm{ with } \alpha \neq 0 :\quad && \alg{M}_d^{{\sf p},\tinyx{1}} &= 
\left\{\begin{array}{ll}
\field{C}[f]/\langle f=\alpha \id_d\rangle & d=0, n,
\label{eq:Malgp1}
\\[0.1cm]
\field{C}[\Omega_d^2]/\langle 
\Omega_d^{d} 
= 
\id_d\rangle & d\in \mc P^{\tinyx{1}} \setminus\{0,n\};
\end{array}\right. \\[0.1cm]
\qptlb_n(\beta,\gamma):\quad && \alg{M}_d^{{\sf p},\tinyx{2}}&= 
\left\{\begin{array}{ll}
\field{C}[\Omega_d^2]/\langle 
\Omega_d^{d} = \gamma \, \id_d
\rangle & d\in \mc P^{\tinyx{2}}\setminus\{n\},
\\[0.1cm] 
\field{C}& d=n.
\end{array}\right.
\end{alignat}
\end{subequations}
\endgroup
The corresponding bases for the middle algebras are
\begingroup
\allowdisplaybreaks
\begin{subequations}
\begin{alignat}{3}
 \mathcal M^{{\sf a}}(d) &= \{ \Omega_d^j\mid j=0,1,\dots, d-1\} 
 \qquad d\in \mc P^{\tinyx{0}};\\[0.1cm]
 \mathcal M^{{\sf a},\tinyx{1}}(d)&= 
  \left\{\begin{array}{ll}
 \{\id_d \} & d = 0, \\[0.1cm]
 \{ \Omega_d^j\mid j=0,1,\dots, d-1\} & d \in \mc P^{\tinyx{1}}\setminus\{0\};
 \end{array}\right.
 \\[0.1cm]
  \mathcal M^{{\sf a},\tinyx{2}}(d) &= \{ \Omega_d^j\mid j=0,1,\dots, d-1\} \qquad d\in \mc P^{\tinyx{2}};\\[0.1cm]
\mathcal M^{{\sf p}}(d) &= 
\left\{\begin{array}{ll}
\{\Omega_d^{2j}\mid j=0,1,\dots, d-1\} & d\in \mc P^{\tinyx{0}} \setminus\{n\},\\[0.1cm]
\{\id_d\} & d = n; 
\end{array}\right. 
\\[0.1cm]
\mathcal M^{{\sf p},\tinyx{1}}(d)&= 
\left\{\begin{array}{ll}
\{\id_d\} & d = 0, n, \\[0.1cm] 
\{ \Omega_d^{2j}\mid j=0,1,\dots, \frac d2-1\} & d\in \mc P^{\tinyx{1}}\setminus\{0,n\};
\end{array}\right.
\label{eq:Mbasisp1}
\\[0.1cm]
\mathcal M^{{\sf p},\tinyx{2}}(d) &= 
\left\{\begin{array}{ll}
\{ \Omega_d^{2j}\mid j=0,1,\dots, \frac d2-1\} & d\in \mc P^{\tinyx{2}}\setminus\{n\}, \\[0.1cm]
\{\id_d\} & d=n. 
\end{array}\right.
\end{alignat}
\end{subequations}
\endgroup

The assignment $C$ for the uncoiled algebras sandwiches $m \in \mathcal M(d)$ between the states $b \in \mc B(d)$ and $t\in \mc T(d)$, producing the diagram $c_{b,m,t}^d$ appearing in~\eqref{eq:sandwich.diag}. The sandwich cellular bases for $\qatl_n(\beta,\gamma)$, $\qatla_n(\beta, \alpha)$ and $\qatlb_n(\beta,\gamma)$ are then given by~\eqref{eq:Sandwich.Basis.uatl}, \eqref{eq:Sandwich.Basis.uatla} and~\eqref{eq:Sandwich.Basis.uatlb} respectively. As usual $\iota$ is the vertical reflection of the diagrams. To illustrate, we return to the example~\eqref{eq:ExMultiaTL}. For $\qatla_6(\beta,\alpha)$, 
it instead reads
\be
\begin{array}{r}
 \psset{unit=0.8cm}
 (b,m,t) = \Big( \
 \begin{pspicture}[shift=-0.05](0,0)(2.4,0.6)
 \psline[linewidth=0.5pt](0,0)(2.4,0)
 \psline[linecolor=blue,linewidth=1.5pt]{-}(0.2,0)(0.2,0.5)
 \psline[linecolor=blue,linewidth=1.5pt]{-}(0.6,0)(0.6,0.5)
 \psarc[linecolor=blue,linewidth=1.5pt]{-}(1.6,0){0.2}{0}{180}
 \psbezier[linecolor=blue,linewidth=1.5pt](1.0,0)(1.0,0.7)(2.2,0.7)(2.2,0)
 \end{pspicture}
 \ ,\id_2,\  
 \begin{pspicture}[shift=-0.05](0,0)(2.4,0.6)
 \psline[linewidth=0.5pt](0,0)(2.4,0)
 \psline[linecolor=blue,linewidth=1.5pt]{-}(0.6,0)(0.6,0.5)
 \psline[linecolor=blue,linewidth=1.5pt]{-}(1.8,0)(1.8,0.5)
 \psarc[linecolor=blue,linewidth=1.5pt]{-}(0,0){0.2}{0}{90}
 \psarc[linecolor=blue,linewidth=1.5pt]{-}(2.4,0){0.2}{90}{180}
 \psarc[linecolor=blue,linewidth=1.5pt]{-}(1.2,0){0.2}{0}{180}
 \end{pspicture}\ \Big)\ , \qquad
 a = \
 \begin{pspicture}[shift=-0.4](0,-0.5)(2.4,0.5)
 \pspolygon[fillstyle=solid,fillcolor=lightlightblue,linewidth=0pt](0,-0.5)(2.4,-0.5)(2.4,0.5)(0,0.5)
 \psbezier[linecolor=blue,linewidth=1.5pt](0.2,-0.5)(0.2,0.1)(1.4,0.1)(1.4,-0.5)
 \psarc[linecolor=blue,linewidth=1.5pt]{-}(0.8,-0.5){0.2}{0}{180}
 \psarc[linecolor=blue,linewidth=1.5pt]{-}(0.8,0.5){0.2}{180}{0}
 \psarc[linecolor=blue,linewidth=1.5pt]{-}(1.6,0.5){0.2}{180}{0}
 \psbezier[linecolor=blue,linewidth=1.5pt]{-}(1.8,-0.5)(1.8,0)(2.2,0)(2.2,0.5)
 \psbezier[linecolor=blue,linewidth=1.5pt]{-}(2.2,-0.5)(2.2,-0.3)(2.35,-0.05)(2.41,0)
 \psbezier[linecolor=blue,linewidth=1.5pt]{-}(0.2,0.5)(0.2,0.3)(0.05,0.05)(-0.01,0)
 \psframe[fillstyle=solid,linecolor=white,linewidth=0pt](-0.05,-0.51)(-0.005,0.51)
 \psframe[fillstyle=solid,linecolor=white,linewidth=0pt](2.405,-0.51)(2.45,0.51)
 \end{pspicture} 
 \ + \
 \begin{pspicture}[shift=-0.4](0,-0.5)(2.4,0.5)
 \pspolygon[fillstyle=solid,fillcolor=lightlightblue,linewidth=0pt](0,-0.5)(2.4,-0.5)(2.4,0.5)(0,0.5)
 \psarc[linecolor=blue,linewidth=1.5pt]{-}(0.4,0.5){0.2}{180}{0}
 \psarc[linecolor=blue,linewidth=1.5pt]{-}(1.6,0.5){0.2}{180}{0}
 \psarc[linecolor=blue,linewidth=1.5pt]{-}(0.8,-0.5){0.2}{0}{180}
 \psarc[linecolor=blue,linewidth=1.5pt]{-}(1.6,-0.5){0.2}{0}{180}
 \psbezier[linecolor=blue,linewidth=1.5pt]{-}(0.2,-0.5)(0.2,0.1)(2.2,0.1)(2.2,-0.5)
 \psbezier[linecolor=blue,linewidth=1.5pt]{-}(-0.02,0.17)(0.1,0)(1.0,-0.06)(1.0,0.5)
 \psbezier[linecolor=blue,linewidth=1.5pt]{-}(2.2,0.5)(2.2,0.32)(2.35,0.22)(2.42,0.2)
 \psframe[fillstyle=solid,linecolor=white,linewidth=0pt](-0.05,-0.51)(-0.005,0.51)
 \psframe[fillstyle=solid,linecolor=white,linewidth=0pt](2.405,-0.51)(2.45,0.51)
 \end{pspicture}\ ,
 \\[0.6cm]
 \Longrightarrow \quad
 \psset{unit=0.8cm}
 c_{b,m,t}^2 = \
 \begin{pspicture}[shift=-0.4](0,-0.5)(2.4,0.5)
 \pspolygon[fillstyle=solid,fillcolor=lightlightblue,linewidth=0pt](0,-0.5)(2.4,-0.5)(2.4,0.5)(0,0.5)
 \psarc[linecolor=blue,linewidth=1.5pt]{-}(1.6,-0.5){0.2}{0}{180}
 \psbezier[linecolor=blue,linewidth=1.5pt](1.0,-0.5)(1.0,0.1)(2.2,0.1)(2.2,-0.5)
 \psarc[linecolor=blue,linewidth=1.5pt]{-}(0,0.5){0.2}{-90}{0}
 \psarc[linecolor=blue,linewidth=1.5pt]{-}(1.2,0.5){0.2}{180}{0}
 \psarc[linecolor=blue,linewidth=1.5pt]{-}(2.4,0.5){0.2}{180}{270}
 \psbezier[linecolor=blue,linewidth=1.5pt](0.2,-0.5)(0.2,0)(0.6,0)(0.6,0.5)
 \psbezier[linecolor=blue,linewidth=1.5pt](0.6,-0.5)(0.6,0)(1.8,0)(1.8,0.5)
 \end{pspicture}
 \ , \qquad
 a\, c_{b,m,t}^2 = \beta \
 \begin{pspicture}[shift=-0.4](0,-0.5)(2.4,0.5)
 \pspolygon[fillstyle=solid,fillcolor=lightlightblue,linewidth=0pt](0,-0.5)(2.4,-0.5)(2.4,0.5)(0,0.5)
 \psarc[linecolor=blue,linewidth=1.5pt]{-}(0.8,-0.5){0.2}{0}{180}
 \psbezier[linecolor=blue,linewidth=1.5pt](0.2,-0.5)(0.2,0.1)(1.4,0.1)(1.4,-0.5)
 \psarc[linecolor=blue,linewidth=1.5pt]{-}(0,0.5){0.2}{-90}{0}
 \psarc[linecolor=blue,linewidth=1.5pt]{-}(1.2,0.5){0.2}{180}{0}
 \psarc[linecolor=blue,linewidth=1.5pt]{-}(2.4,0.5){0.2}{180}{270}
 \psline[linecolor=blue,linewidth=1.5pt](1.8,-0.5)(1.8,0.5)
 \psbezier[linecolor=blue,linewidth=1.5pt](0.6,0.5)(0.6,0.1)(0.1,0.03)(-0.02,-0.17)
 \psbezier[linecolor=blue,linewidth=1.5pt](2.2,-0.5)(2.2,-0.3)(2.37,-0.16)(2.42,-0.13)
 \psframe[fillstyle=solid,linecolor=white,linewidth=0pt](-0.05,-0.51)(-0.005,0.51)
 \psframe[fillstyle=solid,linecolor=white,linewidth=0pt](2.405,-0.51)(2.45,0.51)
 \end{pspicture}
 \ + \alpha \beta \
 \begin{pspicture}[shift=-0.4](0,-0.5)(2.4,0.5)
 \pspolygon[fillstyle=solid,fillcolor=lightlightblue,linewidth=0pt](0,-0.5)(2.4,-0.5)(2.4,0.5)(0,0.5)
 \psarc[linecolor=blue,linewidth=1.5pt]{-}(0,0.5){0.2}{-90}{0}
 \psarc[linecolor=blue,linewidth=1.5pt]{-}(2.4,0.5){0.2}{180}{270}
 \psarc[linecolor=blue,linewidth=1.5pt]{-}(1.2,0.5){0.2}{180}{0}
 \psbezier[linecolor=blue,linewidth=1.5pt]{-}(0.6,0.5)(0.6,-0.05)(1.8,-0.05)(1.8,0.5)
 \psarc[linecolor=blue,linewidth=1.5pt]{-}(0.8,-0.5){0.2}{0}{180}
 \psarc[linecolor=blue,linewidth=1.5pt]{-}(1.6,-0.5){0.2}{0}{180}
 \psbezier[linecolor=blue,linewidth=1.5pt]{-}(0.2,-0.5)(0.2,0.05)(2.2,0.05)(2.2,-0.5)
 \psframe[fillstyle=solid,linecolor=white,linewidth=0pt](-0.05,-0.51)(-0.005,0.51)
 \psframe[fillstyle=solid,linecolor=white,linewidth=0pt](2.405,-0.51)(2.45,0.51)
 \end{pspicture}\ .
\end{array}
\ee
The only difference is that one non-contractible loop is removed and replaced by a factor of $\alpha$.\medskip

One needs to be careful with the assignment $C$ for the uncoiled periodic Temperley--Lieb algebras. We first note that all link states in the bases $\wh \setB_{n,d}$ are even, as are the diagrams in the middle bases~$\mc M^{\sf p}$. Combining $(b,m,t)\in \mc B(d)\times \mathcal M^{\sf p}(d) \times \mc T(d)$ in a sandwich diagram of the form~\eqref{eq:sandwich.diag} also produces an even diagram. In certain cases, the sandwich diagram thus created reduces further because of the defining relations of the algebra. In particular, this occurs for $\qptla_n(\beta,\alpha)$ with $\alpha\neq 0$ in the cases where $b$ and $t$ are both link states with $d=0$ defects and with a non-contractible loop. The resulting connectivity has two non-contractible loops which are removed in the quotient and replaced by a factor~$\alpha^2$. The multiplicative factor of $\alpha^2$ is kept for these states as elements of the basis. Here are two examples of the assignment for $\qptla_{4}(\beta,\alpha)$:
\begin{subequations}
\begin{alignat}{3}
(b,m,t) &= \Big( \
\psset{unit=0.8cm}
\begin{pspicture}[shift=-0.05](0,0)(1.6,0.6)
\psline[linewidth=0.5pt](0,0)(1.6,0)
\psbezier[linecolor=blue,linewidth=1.5pt]{-}(1.62,0.27)(1.5,0.44)(0.6,0.44)(0.6,0)
\psbezier[linecolor=blue,linewidth=1.5pt]{-}(0.2,0)(0.2,0.18)(0.05,0.28)(-0.02,0.3)
\psarc[linecolor=blue,linewidth=1.5pt]{-}(1.2,0){0.2}{0}{180}
\psline[linecolor=blue,linewidth=1.5pt]{-}(0,0.6)(1.6,0.6)
\psframe[fillstyle=solid,linecolor=white,linewidth=0pt](-0.05,0)(-0.005,0.4)
\psframe[fillstyle=solid,linecolor=white,linewidth=0pt](1.605,0)(1.65,0.4)
\end{pspicture}
\ , \id_0, \ 
\begin{pspicture}[shift=-0.05](0,0)(1.6,0.6)
\psline[linewidth=0.5pt](0,0)(1.6,0)
\psarc[linecolor=blue,linewidth=1.5pt]{-}(0.4,0){0.2}{0}{180}
\psarc[linecolor=blue,linewidth=1.5pt]{-}(1.2,0){0.2}{0}{180}
\end{pspicture}
 \ \Big)
\qquad \Longrightarrow \qquad
c_{b,m,t}^0 &&= \
\begin{pspicture}[shift=-0.4](0,-0.5)(1.6,0.5)
\pspolygon[fillstyle=solid,fillcolor=lightlightblue,linewidth=0pt](0,-0.5)(1.6,-0.5)(1.6,0.5)(0,0.5)
\psbezier[linecolor=blue,linewidth=1.5pt]{-}(1.62,-0.23)(1.5,-0.06)(0.6,-0.06)(0.6,-0.5)
\psbezier[linecolor=blue,linewidth=1.5pt]{-}(0.2,-0.5)(0.2,-0.32)(0.05,-0.22)(-0.02,-0.2)
\psarc[linecolor=blue,linewidth=1.5pt]{-}(1.2,-0.5){0.2}{0}{180}
\psline[linecolor=blue,linewidth=1.5pt]{-}(0,0)(1.6,0)
\psarc[linecolor=blue,linewidth=1.5pt]{-}(0.4,0.5){0.2}{180}{0}
\psarc[linecolor=blue,linewidth=1.5pt]{-}(1.2,0.5){0.2}{180}{0}
\psframe[fillstyle=solid,linecolor=white,linewidth=0pt](-0.05,0)(-0.005,-0.5)
\psframe[fillstyle=solid,linecolor=white,linewidth=0pt](1.605,0)(1.65,-0.5)
\end{pspicture} \ ,
\\[0.15cm]
(b,m,t') &= \Big( \
\psset{unit=0.8cm}
\begin{pspicture}[shift=-0.05](0,0)(1.6,0.6)
\psline[linewidth=0.5pt](0,0)(1.6,0)
\psbezier[linecolor=blue,linewidth=1.5pt]{-}(1.62,0.27)(1.5,0.44)(0.6,0.44)(0.6,0)
\psbezier[linecolor=blue,linewidth=1.5pt]{-}(0.2,0)(0.2,0.18)(0.05,0.28)(-0.02,0.3)
\psarc[linecolor=blue,linewidth=1.5pt]{-}(1.2,0){0.2}{0}{180}
\psline[linecolor=blue,linewidth=1.5pt]{-}(0,0.6)(1.6,0.6)
\psframe[fillstyle=solid,linecolor=white,linewidth=0pt](-0.05,0)(-0.005,0.4)
\psframe[fillstyle=solid,linecolor=white,linewidth=0pt](1.605,0)(1.65,0.4)
\end{pspicture}
\ , \id_0, \ 
\begin{pspicture}[shift=-0.05](0,0)(1.6,0.6)
\psline[linewidth=0.5pt](0,0)(1.6,0)
\psarc[linecolor=blue,linewidth=1.5pt]{-}(0,0){0.2}{0}{90}
\psarc[linecolor=blue,linewidth=1.5pt]{-}(1.6,0){0.2}{90}{180}
\psarc[linecolor=blue,linewidth=1.5pt]{-}(0.8,0){0.2}{0}{180}
\psline[linecolor=blue,linewidth=1.5pt]{-}(0,0.6)(1.6,0.6)
\end{pspicture} \ \Big)
\qquad \Longrightarrow \qquad
c_{b,m,t'}^0 &&= \ 
\begin{pspicture}[shift=-0.4](0,-0.5)(1.6,0.5)
\pspolygon[fillstyle=solid,fillcolor=lightlightblue,linewidth=0pt](0,-0.5)(1.6,-0.5)(1.6,0.5)(0,0.5)
\psbezier[linecolor=blue,linewidth=1.5pt]{-}(1.62,-0.23)(1.5,-0.06)(0.6,-0.06)(0.6,-0.5)
\psbezier[linecolor=blue,linewidth=1.5pt]{-}(0.2,-0.5)(0.2,-0.32)(0.05,-0.22)(-0.02,-0.2)
\psarc[linecolor=blue,linewidth=1.5pt]{-}(1.2,-0.5){0.2}{0}{180}
\psline[linecolor=blue,linewidth=1.5pt]{-}(0,0)(1.6,0)
\psline[linecolor=blue,linewidth=1.5pt]{-}(0,0.15)(1.6,0.15)
\psarc[linecolor=blue,linewidth=1.5pt]{-}(0,0.5){0.2}{-90}{0}
\psarc[linecolor=blue,linewidth=1.5pt]{-}(1.6,0.5){0.2}{180}{270}
\psarc[linecolor=blue,linewidth=1.5pt]{-}(0.8,0.5){0.2}{180}{0}
\psframe[fillstyle=solid,linecolor=white,linewidth=0pt](-0.05,0)(-0.005,-0.5)
\psframe[fillstyle=solid,linecolor=white,linewidth=0pt](1.605,0)(1.65,-0.5)
\end{pspicture} 
\ = \alpha^2 \ 
\begin{pspicture}[shift=-0.4](0,-0.5)(1.6,0.5)
\pspolygon[fillstyle=solid,fillcolor=lightlightblue,linewidth=0pt](0,-0.5)(1.6,-0.5)(1.6,0.5)(0,0.5)
\psbezier[linecolor=blue,linewidth=1.5pt]{-}(1.62,-0.23)(1.5,-0.06)(0.6,-0.06)(0.6,-0.5)
\psbezier[linecolor=blue,linewidth=1.5pt]{-}(0.2,-0.5)(0.2,-0.32)(0.05,-0.22)(-0.02,-0.2)
\psarc[linecolor=blue,linewidth=1.5pt]{-}(1.2,-0.5){0.2}{0}{180}
\psarc[linecolor=blue,linewidth=1.5pt]{-}(0,0.5){0.2}{-90}{0}
\psarc[linecolor=blue,linewidth=1.5pt]{-}(1.6,0.5){0.2}{180}{270}
\psarc[linecolor=blue,linewidth=1.5pt]{-}(0.8,0.5){0.2}{180}{0}
\psframe[fillstyle=solid,linecolor=white,linewidth=0pt](-0.05,0)(-0.005,-0.5)
\psframe[fillstyle=solid,linecolor=white,linewidth=0pt](1.605,0)(1.65,-0.5)
\end{pspicture} \ . \label{eq:ex_cell_C_uptla}
\end{alignat}
\end{subequations}
The presence of the coefficient $\alpha^2$ is important for Axiom~A2) to be satisfied. Indeed, setting $a=e_1$, we find in the above example
\be
\label{eq:example.coeff}
e_1 c^0_{b,m,t} = \alpha^2 c^0_{b',m,t}\,, 
\qquad
e_1 c^0_{b,m,t'} = \alpha^2 c^0_{b',m,t'}\,,
\qquad
\textrm{with}
\quad
b' = \ 
\begin{pspicture}[shift=-0.05](0,0)(1.6,0.6)
\psline[linewidth=0.5pt](0,0)(1.6,0)
\psarc[linecolor=blue,linewidth=1.5pt]{-}(0.4,0){0.2}{0}{180}
\psarc[linecolor=blue,linewidth=1.5pt]{-}(1.2,0){0.2}{0}{180}
\end{pspicture}\ .
\ee
Axiom~A2) requires that the coefficients arising in the right-hand side in these products be
identical in the two examples, as they share the same elements $b$ and $m$. If we had defined the assignment $C$ differently, with $c_{b,m,t'}^0$ given by the diagram on the right-hand side of~\eqref{eq:ex_cell_C_uptla} without the factor $\alpha^2$, the resulting coefficients  in \eqref{eq:example.coeff} would not agree.
\medskip

We also note that there are no multiplicative coefficients in the basis elements for $d>0$, even in the cases where both $b$ and $t$ are link states in $\wh{\mc B}_{n,d}$ with one defect winding across the boundary condition. For instance, we have in $\qptla_{4}(\beta,\gamma)$
\begin{subequations}
\begin{alignat}{3}
(b,m,t) &= \Big( \
\psset{unit=0.8cm}
\begin{pspicture}[shift=-0.05](0,0)(1.6,0.6)
\psline[linewidth=0.5pt](0,0)(1.6,0)
\psbezier[linecolor=blue,linewidth=1.5pt]{-}(0.6,0)(0.6,0.25)(1.0,0.25)(1.0,0.55)
\psbezier[linecolor=blue,linewidth=1.5pt]{-}(1.0,0)(1.0,0.35)(1.5,0.35)(1.6,0.35)
\psbezier[linecolor=blue,linewidth=1.5pt]{-}(0,0.35)(0.5,0.35)(0.6,0.45)(0.6,0.55)
\psarc[linecolor=blue,linewidth=1.5pt]{-}(0.0,0){0.2}{0}{90}
\psarc[linecolor=blue,linewidth=1.5pt]{-}(1.6,0){0.2}{90}{180}
\end{pspicture}
\ , \id_2, \ 
\begin{pspicture}[shift=-0.05](0,0)(1.6,0.6)
\psline[linewidth=0.5pt](0,0)(1.6,0)
\psline[linecolor=blue,linewidth=1.5pt]{-}(0.2,0)(0.2,0.5)
\psline[linecolor=blue,linewidth=1.5pt]{-}(0.6,0)(0.6,0.5)
\psarc[linecolor=blue,linewidth=1.5pt]{-}(1.2,0){0.2}{0}{180}
\end{pspicture}
 \ \Big)
\qquad \Longrightarrow \qquad 
c_{b,m,t}^0 &&= \
\begin{pspicture}[shift=-0.4](0,-0.5)(1.6,0.5)
\pspolygon[fillstyle=solid,fillcolor=lightlightblue,linewidth=0pt](0,-0.5)(1.6,-0.5)(1.6,0.5)(0,0.5)
\psline[linecolor=blue,linewidth=1.5pt]{-}(0.6,-0.5)(0.6,0.5)
\psbezier[linecolor=blue,linewidth=1.5pt]{-}(1.0,-0.5)(1.0,0)(1.8,0)(1.8,0.5)
\psbezier[linecolor=blue,linewidth=1.5pt]{-}(0.2,0.5)(0.2,0.32)(0.05,0.22)(-0.02,0.2)
\psarc[linecolor=blue,linewidth=1.5pt]{-}(1.2,0.5){0.2}{180}{360}
\psarc[linecolor=blue,linewidth=1.5pt]{-}(0,-0.5){0.2}{0}{90}
\psarc[linecolor=blue,linewidth=1.5pt]{-}(1.6,-0.5){0.2}{90}{180}
\psframe[fillstyle=solid,linecolor=white,linewidth=0pt](-0.05,0)(-0.005,0.5)
\psframe[fillstyle=solid,linecolor=white,linewidth=0pt](1.605,0)(1.82,0.5)
\end{pspicture} \ ,
\\[0.15cm]
(b,m,t') &= \Big( \
\psset{unit=0.8cm}
\begin{pspicture}[shift=-0.05](0,0)(1.6,0.6)
\psline[linewidth=0.5pt](0,0)(1.6,0)
\psbezier[linecolor=blue,linewidth=1.5pt]{-}(0.6,0)(0.6,0.25)(1.0,0.25)(1.0,0.55)
\psbezier[linecolor=blue,linewidth=1.5pt]{-}(1.0,0)(1.0,0.35)(1.5,0.35)(1.6,0.35)
\psbezier[linecolor=blue,linewidth=1.5pt]{-}(0,0.35)(0.5,0.35)(0.6,0.45)(0.6,0.55)
\psarc[linecolor=blue,linewidth=1.5pt]{-}(0.0,0){0.2}{0}{90}
\psarc[linecolor=blue,linewidth=1.5pt]{-}(1.6,0){0.2}{90}{180}
\end{pspicture}
\ , \id_2, \ 
\begin{pspicture}[shift=-0.05](0,0)(1.6,0.6)
\psline[linewidth=0.5pt](0,0)(1.6,0)
\psbezier[linecolor=blue,linewidth=1.5pt]{-}(0.6,0)(0.6,0.25)(1.0,0.25)(1.0,0.55)
\psbezier[linecolor=blue,linewidth=1.5pt]{-}(1.0,0)(1.0,0.35)(1.5,0.35)(1.6,0.35)
\psbezier[linecolor=blue,linewidth=1.5pt]{-}(0,0.35)(0.5,0.35)(0.6,0.45)(0.6,0.55)
\psarc[linecolor=blue,linewidth=1.5pt]{-}(0.0,0){0.2}{0}{90}
\psarc[linecolor=blue,linewidth=1.5pt]{-}(1.6,0){0.2}{90}{180}
\end{pspicture} \ \Big)
\qquad \Longrightarrow \qquad
c_{b,m,t'}^0 &&= \ 
\begin{pspicture}[shift=-0.4](0,-0.5)(1.6,0.5)
\pspolygon[fillstyle=solid,fillcolor=lightlightblue,linewidth=0pt](0,-0.5)(1.6,-0.5)(1.6,0.5)(0,0.5)
\psline[linecolor=blue,linewidth=1.5pt]{-}(0.6,-0.5)(0.6,0.5)
\psline[linecolor=blue,linewidth=1.5pt]{-}(1.0,-0.5)(1.0,0.5)
\psarc[linecolor=blue,linewidth=1.5pt]{-}(0,0.5){0.2}{-90}{0}
\psarc[linecolor=blue,linewidth=1.5pt]{-}(1.6,0.5){0.2}{180}{270}
\psarc[linecolor=blue,linewidth=1.5pt]{-}(0,-0.5){0.2}{0}{90}
\psarc[linecolor=blue,linewidth=1.5pt]{-}(1.6,-0.5){0.2}{90}{180}
\end{pspicture} \ .
\end{alignat}
\end{subequations}

Axiom~A2) is satisfied since, in all the algebras, the diagrammatic rules state that defects can only be removed, not created, so the action of the algebra can only increase the order. Finally, the scalars arising from the diagrammatic rules depend on the bottom link state and on the middle element, but not on the top link state. One can check that this also holds for the periodic cases with $d=0$, for which some basis elements contain $\alpha^2$ prefactors.
\medskip
 
To show that Axiom~A3) is satisfied, we construct the isomorphism using the same arguments as in the proof of \cref{prop:pTL.skew.sandwich}. We describe the bilinear form $\psi_d$ for $\qatla_n(\beta,\alpha)$ and $\qptla_n(\beta,\alpha)$, the others being similar. For $d=0$, we set $b,t\in \setB_{n,0}$ for $\qatla_n(\beta,\alpha)$ and $b,t\in \wh\setB_{n,0}$ for $\qptla_n(\beta,\alpha)$, and the corresponding bilinear forms are defined as $\psi_0(t,b) = \beta^{m}\alpha^{m'}\id_0 \in \alg M_0^{{\sf a},\tinyx{1}}$ and $\psi_0(t,b) = \beta^{m}\alpha^{m'}\id_0 \in \alg M_0^{{\sf p},\tinyx{1}}$ for $\qatla_n(\beta,\alpha)$ and $\qptla_n(\beta,\alpha)$, respectively. Here, $m$ and $m'$ are respectively the numbers of contractible and non-contractible loops in the diagram $\iota(t)b$. The construction of the basis $\widehat{\mc B}_{n,0}$ with only even link states implies that $m'$ is even for $\qptla_n(\beta,\alpha)$. 
For $d>0$, we set $b,t\in \setB_{n,d}$ or $b,t\in \wh\setB_{n,d}$ for the two algebras, and the bilinear form is defined in the same way as for the algebra $\ptl_n(\beta)$, namely $\psi_d(t,b) = \beta^{m} \Omega_d^{\ell}$ if all defects of $b$ connect to those of $t$, and $\psi_d(t,b)=0$ otherwise. This result is then simplified using the relations $\Omega_d^d = \id_d$ in both $\alg M_d^{{\sf a},\tinyx{1}}$ and $\alg M_d^{{\sf p},\tinyx{1}}$.
\medskip

Here are examples of applications of this bilinear form for $\qptlb_4(\beta,\gamma)$, $\qptla_{10}(\beta, \alpha)$, $\qatlb_{12}(\beta,\gamma)$ and $\qptl_{13}(\beta,\gamma)$:
\begingroup
\allowdisplaybreaks
\begin{subequations}
\begin{alignat}{2}
\psset{unit=0.8cm}
\begin{pspicture}[shift=-0.05](0,0)(1.6,0.6)
\psline[linewidth=0.5pt](0,0)(1.6,0)
\psline[linecolor=blue,linewidth=1.5pt]{-}(1,0)(1,0.5)
\psline[linecolor=blue,linewidth=1.5pt]{-}(1.4,0)(1.4,0.5)
\psarc[linecolor=blue,linewidth=1.5pt]{-}(.4,0){0.2}{0}{180}
\end{pspicture}
 \hspace{1cm}
&\otimes  \hspace{1cm}
\begin{pspicture}[shift=-0.05](0,0)(1.6,0.6)
\psline[linewidth=0.5pt](0,0)(1.6,0)
\psbezier[linecolor=blue,linewidth=1.5pt]{-}(0.6,0)(0.6,0.25)(1.0,0.25)(1.0,0.55)
\psbezier[linecolor=blue,linewidth=1.5pt]{-}(1.0,0)(1.0,0.35)(1.5,0.35)(1.6,0.35)
\psbezier[linecolor=blue,linewidth=1.5pt]{-}(0,0.35)(0.5,0.35)(0.6,0.45)(0.6,0.55)
\psarc[linecolor=blue,linewidth=1.5pt]{-}(0.0,0){0.2}{0}{90}
\psarc[linecolor=blue,linewidth=1.5pt]{-}(1.6,0){0.2}{90}{180}
\end{pspicture}\qquad \xrightarrow{\hspace{0.2cm} \psi_2 \hspace{0.2cm}} \qquad
\begin{pspicture}[shift=-0.5](0,-.6)(1.6,0.6)
\psline[linewidth=0.5pt](0,0)(1.6,0)
\psbezier[linecolor=blue,linewidth=1.5pt]{-}(0.6,0)(0.6,0.25)(1.0,0.25)(1.0,0.55)
\psbezier[linecolor=blue,linewidth=1.5pt]{-}(1.0,0)(1.0,0.35)(1.5,0.35)(1.6,0.35)
\psbezier[linecolor=blue,linewidth=1.5pt]{-}(0,0.35)(0.5,0.35)(0.6,0.45)(0.6,0.55)
\psarc[linecolor=blue,linewidth=1.5pt]{-}(0.0,0){0.2}{0}{90}
\psarc[linecolor=blue,linewidth=1.5pt]{-}(1.6,0){0.2}{90}{180}
\psline[linecolor=blue,linewidth=1.5pt]{-}(1,0)(1,-0.5)
\psline[linecolor=blue,linewidth=1.5pt]{-}(1.4,0)(1.4,-0.5)
\psarc[linecolor=blue,linewidth=1.5pt]{-}(.4,0){0.2}{180}{360}
\end{pspicture} \ 
= \Omega_2^2 = \gamma \id_2,\\[.2cm]
\begin{pspicture}[shift=-0.05](0,0)(4.0,0.8)
\psline[linewidth=0.5pt](0,0)(4.0,0)
\psarc[linecolor=blue,linewidth=1.5pt]{-}(0,0){0.2}{0}{90}
\psbezier[linecolor=blue,linewidth=1.5pt]{-}(0.6,0)(0.6,0.5)(0.1,0.5)(0,0.5)
\psbezier[linecolor=blue,linewidth=1.5pt]{-}(1.0,0)(1.0,0.8)(0.1,0.8)(0,0.8)
\psline[linecolor=blue,linewidth=1.5pt]{-}(0,1.1)(4,1.1)
\psarc[linecolor=blue,linewidth=1.5pt]{-}(1.6,0){0.2}{0}{180}
\psarc[linecolor=blue,linewidth=1.5pt]{-}(2.4,0){0.2}{0}{180}
\psbezier[linecolor=blue,linewidth=1.5pt]{-}(3.0,0)(3.0,0.8)(3.9,0.8)(4,0.8)
\psbezier[linecolor=blue,linewidth=1.5pt]{-}(3.4,0)(3.4,0.5)(3.9,0.5)(4,0.5)
\psarc[linecolor=blue,linewidth=1.5pt]{-}(4,0){0.2}{90}{180}
\end{pspicture} 
\ &\otimes \psset{unit=0.8cm} \ 
\begin{pspicture}[shift=-0.05](0,0)(4.0,1.4)
\psline[linewidth=0.5pt](0,0)(4.0,0)
\psbezier[linecolor=blue,linewidth=1.5pt]{-}(0.2,0)(0.2,1.6)(3.8,1.6)(3.8,0)
\psbezier[linecolor=blue,linewidth=1.5pt]{-}(0.6,0)(0.6,1.2)(3.4,1.2)(3.4,0)
\psbezier[linecolor=blue,linewidth=1.5pt]{-}(1.8,0)(1.8,0.7)(3,0.7)(3,0)
\psarc[linecolor=blue,linewidth=1.5pt]{-}(1.2,0){0.2}{0}{180}
\psarc[linecolor=blue,linewidth=1.5pt]{-}(2.4,0){0.2}{0}{180}
\end{pspicture} 
\ \xrightarrow{\hspace{0.2cm} \psi_0 \hspace{0.2cm}} \ 
\begin{pspicture}[shift=-1.35](0,-1.3)(4.0,0.8)
\psline[linewidth=0.5pt](0,0)(4.0,0)
\psarc[linecolor=blue,linewidth=1.5pt]{-}(0,0){0.2}{0}{90}
\psbezier[linecolor=blue,linewidth=1.5pt]{-}(0.6,0)(0.6,0.5)(0.1,0.5)(0,0.5)
\psbezier[linecolor=blue,linewidth=1.5pt]{-}(1.0,0)(1.0,0.8)(0.1,0.8)(0,0.8)
\psline[linecolor=blue,linewidth=1.5pt]{-}(0,1.1)(4,1.1)
\psarc[linecolor=blue,linewidth=1.5pt]{-}(1.6,0){0.2}{0}{180}
\psarc[linecolor=blue,linewidth=1.5pt]{-}(2.4,0){0.2}{0}{180}
\psbezier[linecolor=blue,linewidth=1.5pt]{-}(3.0,0)(3.0,0.8)(3.9,0.8)(4,0.8)
\psbezier[linecolor=blue,linewidth=1.5pt]{-}(3.4,0)(3.4,0.5)(3.9,0.5)(4,0.5)
\psarc[linecolor=blue,linewidth=1.5pt]{-}(4,0){0.2}{90}{180}
\psbezier[linecolor=blue,linewidth=1.5pt]{-}(0.2,0)(0.2,-1.6)(3.8,-1.6)(3.8,0)
\psbezier[linecolor=blue,linewidth=1.5pt]{-}(0.6,0)(0.6,-1.2)(3.4,-1.2)(3.4,0)
\psbezier[linecolor=blue,linewidth=1.5pt]{-}(1.8,0)(1.8,-0.7)(3,-0.7)(3,0)
\psarc[linecolor=blue,linewidth=1.5pt]{-}(1.2,0){0.2}{180}{0}
\psarc[linecolor=blue,linewidth=1.5pt]{-}(2.4,0){0.2}{180}{0}
\end{pspicture} \
= \beta \alpha^4 \id_0\,,
\\[0.2cm]
\psset{unit=0.8cm}
\begin{pspicture}[shift=-1.35](0,-1.3)(4.8,0.8)
\psline[linewidth=0.5pt](0,0)(4.8,0)
\psbezier[linecolor=blue,linewidth=1.5pt]{-}(1.0,0)(1.0,0.68)(0.1,0.58)(0,0.42)
\psarc[linecolor=blue,linewidth=1.5pt]{-}(2.4,0){0.2}{0}{180}
\psarc[linecolor=blue,linewidth=1.5pt]{-}(3.2,0){0.2}{0}{180}
\psline[linecolor=blue,linewidth=1.5pt]{-}(1.4,0)(1.4,0.8)
\psline[linecolor=blue,linewidth=1.5pt]{-}(4.2,0)(4.2,0.8)
\psbezier[linecolor=blue,linewidth=1.5pt]{-}(1.8,0)(1.8,1.0)(3.8,1.0)(3.8,0)
\psbezier[linecolor=blue,linewidth=1.5pt]{-}(4.6,0)(4.6,0.35)(4.75,0.375)(4.8,0.4)
\psarc[linecolor=blue,linewidth=1.5pt]{-}(0.4,0){0.2}{0}{180}
\end{pspicture} 
\ &\otimes \psset{unit=0.8cm} \ 
\begin{pspicture}[shift=-0.05](0,0)(4.8,1.4)
\psline[linewidth=0.5pt](0,0)(4.8,0)
\psarc[linecolor=blue,linewidth=1.5pt]{-}(0,0){0.2}{0}{90}
\psarc[linecolor=blue,linewidth=1.5pt]{-}(4.8,0){0.2}{90}{180}
\psarc[linecolor=blue,linewidth=1.5pt]{-}(1.6,0){0.2}{0}{180}
\psarc[linecolor=blue,linewidth=1.5pt]{-}(3.2,0){0.2}{0}{180}
\psarc[linecolor=blue,linewidth=1.5pt]{-}(4.0,0){0.2}{0}{180}
\psline[linecolor=blue,linewidth=1.5pt]{-}(0.6,0)(0.6,0.8)
\psline[linecolor=blue,linewidth=1.5pt]{-}(2.6,0)(2.6,0.8)
\psbezier[linecolor=blue,linewidth=1.5pt]{-}(1.0,0)(1.0,0.7)(2.2,0.7)(2.2,0)
\end{pspicture} 
\ \xrightarrow{\hspace{0.2cm} \psi_2 \hspace{0.2cm}} \ 
\begin{pspicture}[shift=-0.05](0,0)(4.8,0.8)
\psline[linewidth=0.5pt](0,0)(4.8,0)
\psbezier[linecolor=blue,linewidth=1.5pt]{-}(1.0,0)(1.0,0.68)(0.1,0.58)(0,0.42)
\psarc[linecolor=blue,linewidth=1.5pt]{-}(2.4,0){0.2}{0}{180}
\psarc[linecolor=blue,linewidth=1.5pt]{-}(3.2,0){0.2}{0}{180}
\psline[linecolor=blue,linewidth=1.5pt]{-}(1.4,0)(1.4,0.8)
\psline[linecolor=blue,linewidth=1.5pt]{-}(4.2,0)(4.2,0.8)
\psbezier[linecolor=blue,linewidth=1.5pt]{-}(1.8,0)(1.8,1.0)(3.8,1.0)(3.8,0)
\psbezier[linecolor=blue,linewidth=1.5pt]{-}(4.6,0)(4.6,0.35)(4.75,0.375)(4.8,0.4)
\psarc[linecolor=blue,linewidth=1.5pt]{-}(0.4,0){0.2}{0}{180}
\psarc[linecolor=blue,linewidth=1.5pt]{-}(0,0){0.2}{-90}{0}
\psarc[linecolor=blue,linewidth=1.5pt]{-}(4.8,0){0.2}{180}{270}
\psarc[linecolor=blue,linewidth=1.5pt]{-}(1.6,0){0.2}{180}{0}
\psarc[linecolor=blue,linewidth=1.5pt]{-}(3.2,0){0.2}{180}{0}
\psarc[linecolor=blue,linewidth=1.5pt]{-}(4.0,0){0.2}{180}{0}
\psline[linecolor=blue,linewidth=1.5pt]{-}(0.6,0)(0.6,-0.8)
\psline[linecolor=blue,linewidth=1.5pt]{-}(2.6,0)(2.6,-0.8)
\psbezier[linecolor=blue,linewidth=1.5pt]{-}(1.0,0)(1.0,-0.7)(2.2,-0.7)(2.2,0)
\end{pspicture} \ =0,
\\[0.2cm]
\psset{unit=0.8cm}
\begin{pspicture}[shift=-1.35](0,-1.3)(5.2,0.8)
\psline[linewidth=0.5pt](0,0)(5.2,0)
\psarc[linecolor=blue,linewidth=1.5pt]{-}(0,0){0.2}{0}{90}
\psbezier[linecolor=blue,linewidth=1.5pt]{-}(0.6,0)(0.6,0.5)(0.1,0.5)(0,0.5)
\psbezier[linecolor=blue,linewidth=1.5pt]{-}(1.0,0)(1.0,0.8)(0.1,0.8)(0,0.8)
\psbezier[linecolor=blue,linewidth=1.5pt]{-}(1.4,0)(1.4,1.1)(0.1,1.1)(0,1.1)
\psarc[linecolor=blue,linewidth=1.5pt]{-}(2.4,0){0.2}{0}{180}
\psline[linecolor=blue,linewidth=1.5pt]{-}(1.8,0)(1.8,0.8)
\psline[linecolor=blue,linewidth=1.5pt]{-}(3.0,0)(3.0,0.8)
\psline[linecolor=blue,linewidth=1.5pt]{-}(3.4,0)(3.4,0.8)
\psbezier[linecolor=blue,linewidth=1.5pt]{-}(3.8,0)(3.8,1.1)(5.1,1.1)(5.2,1.1)
\psbezier[linecolor=blue,linewidth=1.5pt]{-}(4.2,0)(4.2,0.8)(5.1,0.8)(5.2,0.8)
\psbezier[linecolor=blue,linewidth=1.5pt]{-}(4.6,0)(4.6,0.5)(5.1,0.5)(5.2,0.5)
\psarc[linecolor=blue,linewidth=1.5pt]{-}(5.2,0){0.2}{90}{180}
\end{pspicture} 
\ &\otimes \psset{unit=0.8cm} \ 
\begin{pspicture}[shift=-0.05](0,0)(5.2,1.4)
\psline[linewidth=0.5pt](0,0)(5.2,0)
\psarc[linecolor=blue,linewidth=1.5pt]{-}(0,0){0.2}{0}{90}
\psarc[linecolor=blue,linewidth=1.5pt]{-}(5.2,0){0.2}{90}{180}
\psarc[linecolor=blue,linewidth=1.5pt]{-}(1.6,0){0.2}{0}{180}
\psarc[linecolor=blue,linewidth=1.5pt]{-}(2.4,0){0.2}{0}{180}
\psbezier[linecolor=blue,linewidth=1.5pt]{-}(3.8,0)(3.8,0.9)(4.2,0.9)(4.2,1.8)
\psbezier[linecolor=blue,linewidth=1.5pt]{-}(4.2,0)(4.2,0.9)(4.6,0.9)(4.6,1.8)
\psbezier[linecolor=blue,linewidth=1.5pt]{-}(4.6,0)(4.6,0.9)(5.1,1.2)(5.2,1.2)
\psbezier[linecolor=blue,linewidth=1.5pt]{-}(0,1.2)(1.6,1.2)(3.8,1.2)(3.8,1.8)
\psbezier[linecolor=blue,linewidth=1.5pt]{-}(1.0,0)(1.0,1.0)(3.0,1.0)(3.0,0)
\psbezier[linecolor=blue,linewidth=1.5pt]{-}(0.6,0)(0.6,1.4)(3.4,1.4)(3.4,0)
\end{pspicture} 
\ \xrightarrow{\hspace{0.2cm} \psi_3 \hspace{0.2cm}} \ 
\begin{pspicture}[shift=-0.05](0,0)(5.2,0.8)
\psline[linewidth=0.5pt](0,0)(5.2,0)
\psarc[linecolor=blue,linewidth=1.5pt]{-}(0,0){0.2}{0}{90}
\psbezier[linecolor=blue,linewidth=1.5pt]{-}(0.6,0)(0.6,0.5)(0.1,0.5)(0,0.5)
\psbezier[linecolor=blue,linewidth=1.5pt]{-}(1.0,0)(1.0,0.8)(0.1,0.8)(0,0.8)
\psbezier[linecolor=blue,linewidth=1.5pt]{-}(1.4,0)(1.4,1.1)(0.1,1.1)(0,1.1)
\psarc[linecolor=blue,linewidth=1.5pt]{-}(2.4,0){0.2}{0}{180}
\psline[linecolor=blue,linewidth=1.5pt]{-}(1.8,0)(1.8,0.8)
\psline[linecolor=blue,linewidth=1.5pt]{-}(3.0,0)(3.0,0.8)
\psline[linecolor=blue,linewidth=1.5pt]{-}(3.4,0)(3.4,0.8)
\psbezier[linecolor=blue,linewidth=1.5pt]{-}(3.8,0)(3.8,1.1)(5.1,1.1)(5.2,1.1)
\psbezier[linecolor=blue,linewidth=1.5pt]{-}(4.2,0)(4.2,0.8)(5.1,0.8)(5.2,0.8)
\psbezier[linecolor=blue,linewidth=1.5pt]{-}(4.6,0)(4.6,0.5)(5.1,0.5)(5.2,0.5)
\psarc[linecolor=blue,linewidth=1.5pt]{-}(5.2,0){0.2}{90}{180}
\psarc[linecolor=blue,linewidth=1.5pt]{-}(0,0){0.2}{-90}{0}
\psarc[linecolor=blue,linewidth=1.5pt]{-}(5.2,0){0.2}{180}{270}
\psarc[linecolor=blue,linewidth=1.5pt]{-}(1.6,0){0.2}{180}{0}
\psarc[linecolor=blue,linewidth=1.5pt]{-}(2.4,0){0.2}{180}{0}
\psbezier[linecolor=blue,linewidth=1.5pt]{-}(3.8,0)(3.8,-0.9)(4.2,-0.9)(4.2,-1.8)
\psbezier[linecolor=blue,linewidth=1.5pt]{-}(4.2,0)(4.2,-0.9)(4.6,-0.9)(4.6,-1.8)
\psbezier[linecolor=blue,linewidth=1.5pt]{-}(4.6,0)(4.6,-0.9)(5.1,-1.2)(5.2,-1.2)
\psbezier[linecolor=blue,linewidth=1.5pt]{-}(0,-1.2)(1.6,-1.2)(3.8,-1.2)(3.8,-1.8)
\psbezier[linecolor=blue,linewidth=1.5pt]{-}(1.0,0)(1.0,-1.0)(3.0,-1.0)(3.0,0)
\psbezier[linecolor=blue,linewidth=1.5pt]{-}(0.6,0)(0.6,-1.4)(3.4,-1.4)(3.4,0)
\end{pspicture} \
= \beta^2 \Omega_3^2.
\end{alignat}
\end{subequations}
\endgroup

The anti-involution $\iota$ is again the vertical reflection of the diagrams, and one can check that Axiom~A4) indeed holds. We have thus given the sandwich cellular data for the uncoiled affine and periodic Temperley--Lieb algebras, and showed that the four axioms are satisfied, thus proving their sandwich cellularity. Furthermore, because all the middle algebras are commutative, this shows that the uncoiled algebras are also affine cellular.\eproof

\paragraph{The case $\boldsymbol{\qptla_n(\beta,\alpha)}$ with $\boldsymbol{\alpha = 0}$.}

For $\qptla_n(\beta,\alpha)$ with $\alpha=0$, we have not succeeded in constructing a basis as in Axiom A1) to describe the entire set of diagrams with no through-lines as a single cell. 

\begin{Proposition}
The uncoiled algebras $\qptla_n(\beta,\alpha)$ 
with $\alpha=0$ are skew sandwich cellular.
\end{Proposition}
\proof
We define the two-sided ideals
\begin{subequations}
\begin{alignat}{2}
&\repJ_{0_0}^{\tinyx{E}} = \{c_1 E c_2\,|\, c_1,c_2 \in \qptla_n(\beta,0) \},
\qquad
&&\repJ_{0_1}^{\tinyx{E}} = \{c_1 EF c_2\,|\, c_1,c_2 \in \qptla_n(\beta,0) \},
\\[0.15cm]
&\repJ_{0_0}^{\tinyx{F}} = \{c_1 F c_2\,|\, c_1,c_2 \in \qptla_n(\beta,0) \},
\qquad
&&\repJ_{0_1}^{\tinyx{F}} = \{c_1 FE c_2\,|\, c_1,c_2 \in \qptla_n(\beta,0) \},
\end{alignat}
\end{subequations}
satisfying
\be
(\repJ_{0_1}^{\tinyx{E}} \oplus \repJ_{0_1}^{\tinyx{F}}) \subset
(\repJ_{0_0}^{\tinyx{E}} \oplus \repJ_{0_0}^{\tinyx{F}}) \subset 
\repJ_2 \subset \dots \subset \repJ_n = \qptla_n(\beta,0).
\ee
The poset is
\be
\mc P = \{0_1^{\tinyx{E}}, 0_1^{\tinyx{F}} >_{\mc P} 0_0^{\tinyx{E}}, 0_0^{\tinyx{F}} >_{\mc P} 2 >_{\mc P} 4 >_{\mc P} \dots >_{\mc P} n\}.
\ee
For $d=2,4,\dots, n$, the sets $\mc B(d)$ and $\mc T(d)$ are the bases $\wh\setB_{n,d}$, and the middle algebras $\alg{M}_d^{{\sf p},\tinyx{1}}$ and their bases $\mathcal M^{{\sf p},\tinyx{1}}(d)$ are identical to those given in \eqref{eq:Malgp1} and \eqref{eq:Mbasisp1} for the case $\alpha \neq 0$. For $d=0_0^{\tinyx{E}}, 0_0^{\tinyx{F}}, 0_1^{\tinyx{E}}, 0_1^{\tinyx{F}}$
the sets $\mc B$ and $\mc T$ are
\begin{subequations}
\label{eq:sets.uptla.alpha0}
\begin{alignat}{2}
&d=0_0^{\tinyx{E}}: \qquad&&
\mc B = \setB^{\tinyx 1}_{n,0}\,,\quad
\mc T = \setB^{\tinyx 1}_{n,0}\,,
\\[0.15cm]
&d=0_{0}^{\tinyx{F}}: \qquad&&
\mc B = \setB^{\tinyx 0}_{n,0}\,,\quad
\mc T = \setB^{\tinyx 0}_{n,0}\,,
\\[0.15cm]
&d=0_1^{\tinyx{E}}: \qquad&&
\mc B = \wh\setB^{\tinyx 1}_{n,0}\,,\quad
\mc T = \wh\setB^{\tinyx 0}_{n,0}\,,
\\[0.15cm]
&d=0_1^{\tinyx{F}}: \qquad&&
\mc B = \wh\setB^{\tinyx 0}_{n,0}\,,\quad
\mc T = \wh\setB^{\tinyx 1}_{n,0}\,.
\end{alignat}
\end{subequations}
Moreover, we have $\alg{M} = \mathbb C$ and $\mc M = \{\id\}$ in all four cases. As usual, the involution is the vertical reflection of the diagrams. It acts as the identity on all cells with $d\geq 2$, and likewise for the cells $d=0_0^{\tinyx E}$ and $d=0_0^{\tinyx F}$. It, however, permutes the cells $d=0_1^{\tinyx{E}}$ and $d=0_1^{\tinyx{F}}$. The proof that the axioms hold follows the same arguments as those given in the proof of \cref{prop:pTL.skew.sandwich}.\eproof 

%
\section{Wenzl--Jones projectors}\label{sec:WJ.projs}
%

In this section, we first recall the definition of the Wenzl--Jones projectors for $\tl_n(\beta)$. We then construct the similar projectors for the uncoiled algebras, first for the periodic ones and subsequently for the affine ones. Throughout, we use the notation for parallel loop segments
\be
\begin{pspicture}[shift=-0.35](-0.2,-0.5)(0.2,0.4)
\psline[linecolor=blue,linewidth=1.5pt]{-}(0,-0.4)(0,0.4)\rput(0.2,0.0){$_k$}
\end{pspicture}
\ = \ 
\begin{pspicture}[shift=-0.725](0,-0.85)(1.1,0.40)
\multiput(0,0)(0.2,0){4}{\psline[linecolor=blue,linewidth=1.5pt]{-}(0,-0.4)(0,0.4)}
\rput(0.8,0){$_{...}$}
\psline[linecolor=blue,linewidth=1.5pt]{-}(1,-0.4)(1,0.4)
\rput(0.5,-0.75){$\underbrace{\phantom{\ \hspace{0.9cm}\ }}_k$}
\end{pspicture}
\ .
\ee

\subsection[Projectors for $\tl_n(\beta)$]{Projectors for $\boldsymbol{\tl_n(\beta)}$}

The Wenzl--Jones projectors $P_1, \dots, P_n$ \cite{J83,KL94,W88} are elements of the ordinary Temperley--Lieb algebra $\tl_n(\beta)$. They are defined recursively by $P_1 = \id$ and
\begin{equation} 
\label{eq:Pm.def}
P_{m} = P_{m-1} + \frac{[m-1]}{[m]} \, P_{m-1} e_{m-1} P_{m-1},\qquad 1 < m \le n.
\end{equation}
The $q$-numbers are
\begin{equation}
[k]=\frac{q^k-q^{-k}}{q-q^{-1}},
\end{equation}
whereby $\beta = -[2]$. The recursive definition can be used to show that
\be
(P_m)^2 = P_m, \qquad
e_j P_m = P_m e_j = 0, \qquad j = 1, 2,\dots, m-1,
\ee 
confirming that these indeed are projectors. The projector $P_n$ is particularly useful to study the representation theory of $\tl_n(\beta)$, as it projects onto the unique one-dimensional standard module of $\tl_n(\beta)$, spanned by the unique link state with $n$ defects.\medskip

The projectors satisfy the identities
\begin{subequations}\label{Pm.props}
\begin{alignat}{2}
e_m &P_m e_m = -\frac{[m+1]}{[m]}\,P_{m-1}e_m,\label{Pm.prop3}
\\
P_m &= (\id_1 \otimes P_{m-1})\bigg(\id + \sum_{j=1}^{m-1} \frac{[m-j]}{[m]} e_1e_2 \cdots e_j\bigg)\label{eq:other.rec.1}
\\&
= (P_{m-1} \otimes \id_1)\bigg(\id + \sum_{j=1}^{m-1} \frac{[j]}{[m]} e_{m-1}e_{m-2} \cdots e_j\bigg).\label{eq:other.rec.2}
\end{alignat}
\end{subequations}
Here, we denote as $a_1 \otimes a_2$ the element of $\tl_{n_1+n_2}(\beta)$ where $a_1 \in \tl_{n_1}(\beta)$ and $a_2 \in \tl_{n_2}(\beta)$ are drawn side-by-side.\medskip

In terms of the diagram presentation, we write the Wenzl--Jones projectors as
\be
P_m = \
\begin{pspicture}[shift=-0.05](-0.03,0.00)(1.47,0.3)
\pspolygon[fillstyle=solid,fillcolor=pink](0,0)(1.4,0)(1.4,0.3)(0,0.3)(0,0)\rput(0.7,0.15){$_m$}
\end{pspicture}\ .
\ee
For instance, the relations \eqref{eq:Pm.def} and \eqref{Pm.props} are written in terms of diagrams as
\begin{subequations}
\begin{alignat}{2}
\begin{pspicture}[shift=-0.05](0,0)(1.4,0.3)
\pspolygon[fillstyle=solid,fillcolor=pink](0,0)(1.4,0)(1.4,0.3)(0,0.3)(0,0)\rput(0.7,0.15){$_m$}
\end{pspicture} 
\ &= \  
\begin{pspicture}[shift=-0.05](0,0)(1.6,0.3)
\pspolygon[fillstyle=solid,fillcolor=pink](0,0)(1.4,0)(1.4,0.3)(0,0.3)(0,0)\rput(0.7,0.15){$_{m-1}$}
\psline[linecolor=blue,linewidth=1.5pt]{-}(1.6,0.3)(1.6,0.0)
\rput(1.6,0.5){$_1$}
\end{pspicture}
 \ + \frac{[m-1]}{[m]}\  
\begin{pspicture}[shift=-0.75](0,-0.2)(1.6,1.5)
\pspolygon[fillstyle=solid,fillcolor=pink](0,0)(1.4,0)(1.4,0.3)(0,0.3)(0,0)\rput(0.7,0.15){$_{m-1}$}
\rput(0,1){\pspolygon[fillstyle=solid,fillcolor=pink](0,0)(1.4,0)(1.4,0.3)(0,0.3)(0,0)\rput(0.7,0.15){$_{m-1}$}}
\psline[linecolor=blue,linewidth=1.5pt]{-}(0.4,0.3)(0.40,1.0)
\psarc[linecolor=blue,linewidth=1.5pt]{-}(1.4,0.3){0.2}{0}{180}
\psarc[linecolor=blue,linewidth=1.5pt]{-}(1.4,1){0.2}{180}{0}
\psline[linecolor=blue,linewidth=1.5pt]{-}(1.6,0.3)(1.6,0.0)
\psline[linecolor=blue,linewidth=1.5pt]{-}(1.6,1.3)(1.6,1.0)
\rput(1.6,-0.2){$_1$}
\rput(1.6,1.5){$_1$}
\rput(0.8,0.65){$_{m-2}$}
\end{pspicture}
\ \,, \qquad
\begin{pspicture}[shift=-0.05](0,0)(1.6,0.3)
\pspolygon[fillstyle=solid,fillcolor=pink](0,0)(1.4,0)(1.4,0.3)(0,0.3)(0,0)\rput(0.7,0.15){$_m$}
\psline[linecolor=blue,linewidth=1.5pt]{-}(1.6,0.3)(1.6,0.0)
\psarc[linecolor=blue,linewidth=1.5pt]{-}(1.4,0.3){0.2}{0}{180}
\psarc[linecolor=blue,linewidth=1.5pt]{-}(1.4,0){0.2}{180}{0}
\rput(1.6,0.65){$_1$}
\end{pspicture} 
\ = -\frac{[m+1]}{[m]}\ 
\begin{pspicture}[shift=-0.05](0,0)(1.4,0.3)
\pspolygon[fillstyle=solid,fillcolor=pink](0,0)(1.4,0)(1.4,0.3)(0,0.3)(0,0)\rput(0.7,0.15){$_{m-1}$}
\end{pspicture}\ ,
\\[0.5cm]
\begin{pspicture}[shift=-0.05](0,0)(1.4,0.3)
\pspolygon[fillstyle=solid,fillcolor=pink](0,0)(1.4,0)(1.4,0.3)(0,0.3)(0,0)\rput(0.7,0.15){$_m$}
\end{pspicture} 
\ &= \
\begin{pspicture}[shift=-0.05](0,0)(1.6,0.3)
\pspolygon[fillstyle=solid,fillcolor=pink](0.2,0)(1.6,0)(1.6,0.3)(0.2,0.3)(0.2,0)\rput(0.9,0.15){$_{m-1}$}
\psline[linecolor=blue,linewidth=1.5pt]{-}(0,0.3)(0,0.0)
\rput(0,0.5){$_1$}
\end{pspicture}
\ + \sum_{j=1}^{m-1} \frac{[m-j]}{[m]} \
\begin{pspicture}[shift=-0.4](-0.2,0)(1.9,1.0)
\pspolygon[fillstyle=solid,fillcolor=pink](0.2,0)(1.6,0)(1.6,0.3)(0.2,0.3)(0.2,0)\rput(0.9,0.15){$_{m-1}$}
\psline[linecolor=blue,linewidth=1.5pt]{-}(0,0.3)(0,0.0)
\psline[linecolor=blue,linewidth=1.5pt]{-}(1.4,0.3)(1.4,1.0)\rput(0.1,1.2){$_{j-1}$}
\psarc[linecolor=blue,linewidth=1.5pt]{-}(0.2,0.3){0.2}{0}{180}
\psarc[linecolor=blue,linewidth=1.5pt]{-}(0.8,1.0){0.2}{180}{360}\rput(0.6,1.2){$_{1}$}
\psbezier[linecolor=blue,linewidth=1.5pt]{-}(0.9,0.3)(0.9,0.60)(0.2,0.60)(0.2,1.0)
\rput(1.65,1.2){$_{m-j-1}$}
\rput(0,-0.2){$_1$}
\end{pspicture}
\ = \
\begin{pspicture}[shift=-0.05](0,0)(1.6,0.3)
\pspolygon[fillstyle=solid,fillcolor=pink](0,0)(1.4,0)(1.4,0.3)(0,0.3)(0,0)\rput(0.7,0.15){$_{m-1}$}
\psline[linecolor=blue,linewidth=1.5pt]{-}(1.6,0.3)(1.6,0.0)
\rput(1.6,0.5){$_1$}
\end{pspicture}
\ + \sum_{j=1}^{m-1} \frac{[j]}{[m]} \
\begin{pspicture}[shift=-0.4](-0.2,0)(2.1,1.0)
\pspolygon[fillstyle=solid,fillcolor=pink](0,0)(1.4,0)(1.4,0.3)(0,0.3)(0,0)\rput(0.7,0.15){$_{m-1}$}
\psline[linecolor=blue,linewidth=1.5pt]{-}(1.6,0.3)(1.6,0.0)
\psline[linecolor=blue,linewidth=1.5pt]{-}(0.2,0.3)(0.2,1.0)\rput(0.05,1.2){$_{j-1}$}
\psarc[linecolor=blue,linewidth=1.5pt]{-}(1.4,0.3){0.2}{0}{180}
\psarc[linecolor=blue,linewidth=1.5pt]{-}(0.8,1.0){0.2}{180}{360}\rput(0.6,1.2){$_{1}$}
\psbezier[linecolor=blue,linewidth=1.5pt]{-}(0.7,0.3)(0.7,0.60)(1.4,0.60)(1.4,1.0)
\rput(1.8,1.2){$_{m-j-1}$}
\rput(1.6,-0.2){$_1$}
\end{pspicture}\ .
\end{alignat}
\end{subequations}

\subsection[Projectors for $\qptl_n(\beta,\gamma)$, $\qptla_n(\beta,\alpha)$ and $\qptlb_n(\beta,\gamma)$]{Projectors for $\boldsymbol{\qptl_n(\beta,\gamma)}$, $\boldsymbol{\qptla_n(\beta, \alpha)}$ and $\boldsymbol{\qptlb_n(\beta,\gamma)}$}\label{sec:proj.qptl}

In this section, we construct the Wenzl--Jones projectors $Q_n$ for the uncoiled periodic Temperley--Lieb algebras, for $n$ odd and $n$ even. These elements satisfy
\be
(Q_n)^2 = Q_n ,\qquad e_j Q_n = Q_n e_j = 0, \qquad j=0, 1,\dots, n-1.
\ee
In fact, there can be only one element of the algebra satisfying these two properties. Indeed, the first property implies that $Q_n = \id + \sum_{j=0}^{n-1} a_j e_j$ for some elements $a_j$ of the algebra. Supposing that there exists a second element $\bar Q_n$ with the same properties, we compute $Q_n \bar Q_n$ and find
\be
Q_n \bar Q_n = \Big(\id + \sum_{j=0}^{n-1} a_j e_j \Big) \bar Q_n = \bar Q_n. 
\ee
The same calculation can be done by expanding $\bar Q_n$ as $\id + \sum_{j=0}^{n-1} \bar a_j e_j$ and then shows that $Q_n \bar Q_n = Q_n$, thus proving that $Q_n = \bar Q_n$. This proof was first given in \cite{KL94} to show the uniqueness of the projectors of $\tl_n(\beta)$.\medskip 

The projectors $Q_n$ will satisfy a number of properties. First, it is clear that $Q_n$ is also an element of $\atl_n(\beta)$. As such, it commutes with $\Omega$. Indeed, $\Omega\, Q_n\, \Omega^{-1}$ is an element of $\qptl_n(\beta)$ that squares to itself and is annihilated by each $e_j$. By the uniqueness, it follows that $\Omega\, Q_n\, \Omega^{-1} = Q_n$, which thus confirms that $\Omega Q_n = Q_n \Omega$. Secondly, for each quotient algebra, the projector $Q_n$ evaluates to $1$ on $\repW_{n,n,z}$, whereas it vanishes on the standard modules $\repW_{n,d,z}$ with $d<n$.\medskip

To construct these projectors, we write them as 
\be
\label{eq:Qn.idea}
Q_n = \begin{pspicture}[shift=-0.05](-0.03,0.00)(1.47,0.3)
 \pspolygon[fillstyle=solid,fillcolor=pink](0,0)(1.4,0)(1.4,0.3)(0,0.3)(0,0)\rput(0.7,0.15){$_n$}
 \end{pspicture}\ +   \sum_{c} \Gamma(c) \ 
\begin{pspicture}[shift=-0.75](0,-0.2)(1.4,1.5)
\pspolygon[fillstyle=solid,fillcolor=pink](0,0)(1.4,0)(1.4,0.3)(0,0.3)(0,0)\rput(0.7,0.15){$_{n}$}
\rput(0,1){\pspolygon[fillstyle=solid,fillcolor=pink](0,0)(1.4,0)(1.4,0.3)(0,0.3)(0,0)\rput(0.7,0.15){$_{n}$}}
\pspolygon[fillstyle=solid,fillcolor=lightlightblue](0,0.3)(1.4,0.3)(1.4,1)(0,1)(0,0.3)\rput(0.7,0.65){$c$}
\end{pspicture} \ .
\ee
There is in fact only a restricted subset of connectivities $c$ of the uncoiled periodic Temperley--Lieb algebras for which $P_n c P_n$ is non-zero. In terms of the sandwich diagrams of \eqref{eq:sandwich.diag}, the connectivities are those where the link states $b$ and $t$ have all of their arcs crossing the vertical line between the nodes $n$ and $1$, namely
\be
\label{eq:vwc}
b,t = \ \ 
\begin{pspicture}[shift=-0.05](-0.2,0)(1.5,0.7)
\psline[linewidth=0.5pt](0,0)(1.2,0)
\psline[linecolor=blue,linewidth=1.5pt]{-}(0.6,0)(0.6,0.5)
\psarc[linecolor=blue,linewidth=1.5pt]{-}(0,0){0.2}{0}{90}
\psarc[linecolor=blue,linewidth=1.5pt]{-}(1.2,0){0.2}{90}{180}
\rput(0.6,0.7){$_{n-2k}$}
\rput(-0.2,0.2){$_k$}
\rput(1.4,0.2){$_k$}
\end{pspicture}\ , 
\qquad
c = 
\begin{pspicture}[shift=-0.75](-0.2,0.3)(2.4,2)
\pspolygon[fillstyle=solid,fillcolor=lightlightblue](0,0.3)(2.2,0.3)(2.2,0.8)(0,0.8)(0,0.3)
\pspolygon[fillstyle=solid,fillcolor=lightlightblue](0,2)(2.2,2)(2.2,1.5)(0,1.5)(0,2)
\pspolygon[fillstyle=solid,fillcolor=lightlightblue](0.4,0.8)(1.8,0.8)(1.8,1.5)(0.4,1.5)(0.4,0.8)
\psarc[linecolor=blue,linewidth=1.5pt]{-}(2.2,0.3){0.2}{90}{180}
\psarc[linecolor=blue,linewidth=1.5pt]{-}(0,0.3){0.2}{0}{90}
\psarc[linecolor=blue,linewidth=1.5pt]{-}(2.2,2){0.2}{180}{270}
\psarc[linecolor=blue,linewidth=1.5pt]{-}(0,2){0.2}{270}{0}
\psline[linecolor=blue,linewidth=1.5pt]{-}(1.1,0.3)(1.1,0.8)
\psline[linecolor=blue,linewidth=1.5pt]{-}(1.1,1.5)(1.1,2)
\rput(-0.15,0.5){$_k$}\rput(2.35,0.5){$_k$}
\rput(-0.15,1.8){$_k$}\rput(2.35,1.8){$_k$}
\rput(1.55,0.55){$_{n-2k}$}
\rput(1.55,1.75){$_{n-2k}$}
\rput(1.1,1.15){$\Omega^{2\ell}$}
\end{pspicture} \ ,
\ee
where $\Omega^{2\ell}$ in the middle of the diagram is an element of $\atl_{n-2k}(\beta)$.\medskip

For $\qptl_n(\beta,\gamma)$, $\qptla_n(\beta, \alpha)$ and $\qptlb_n(\beta,\gamma)$, the sum over $c$ will run over even connectivities that have at least one generator $e_j$. These correspond to connectivities as in \eqref{eq:vwc} with $\ell \in \mathbb Z$ and $k \ge 1$. It is also clear from the proposed form \eqref{eq:Qn.idea} that $Q_n e_j = e_j Q_n =0$ for $j = 1, \dots, n-1$. Moreover, the first term is the only one with a component along the unit $\id$. It has coefficient $1$, which is the only choice consistent with $(Q_n)^2 = Q_n$. The constants $\Gamma(c)$ are then chosen to ensure that $Q_n e_0 = e_0 Q_n = 0$. To illustrate, here are some examples of projectors for small values of $n$:
\begin{subequations}
\begin{alignat}{2}
\qptla_2(\beta,\alpha)&: \quad Q_2 = \
\begin{pspicture}[shift=-0.05](0,0.00)(0.8,0.3)
\pspolygon[fillstyle=solid,fillcolor=pink](0,0)(0.8,0)(0.8,0.3)(0,0.3)(0,0)\rput(0.4,0.15){$_2$}
\end{pspicture} 
\ + \frac{\beta}{\alpha^2-\beta^2} \ 
\begin{pspicture}[shift=-0.5](0,0.00)(0.8,1.2)
\pspolygon[fillstyle=solid,fillcolor=pink](0,0)(0.8,0)(0.8,0.3)(0,0.3)(0,0)\rput(0.4,0.15){$_2$}
\rput(0,0.9){\pspolygon[fillstyle=solid,fillcolor=pink](0,0)(0.8,0)(0.8,0.3)(0,0.3)(0,0)\rput(0.4,0.15){$_2$}}
\pspolygon[fillstyle=solid,fillcolor=lightlightblue](0,0.3)(0.8,0.3)(0.8,0.9)(0,0.9)(0,0.3)
\psarc[linecolor=blue,linewidth=1.5pt]{-}(0,0.3){0.2}{0}{90}
\psarc[linecolor=blue,linewidth=1.5pt]{-}(0.8,0.3){0.2}{90}{180}
\psarc[linecolor=blue,linewidth=1.5pt]{-}(0,0.9){0.2}{-90}{0}
\psarc[linecolor=blue,linewidth=1.5pt]{-}(0.8,0.9){0.2}{180}{270}
\end{pspicture} \ ,
\\[0.2cm]
\qptl_3(\beta,\gamma)&: \quad Q_3 = \
\begin{pspicture}[shift=-0.05](0,0.00)(1.2,0.3)
\pspolygon[fillstyle=solid,fillcolor=pink](0,0)(1.2,0)(1.2,0.3)(0,0.3)(0,0)\rput(0.6,0.15){$_3$}
\end{pspicture} 
\ - \frac{\beta^2-1}{\gamma^2+\gamma^{-2}+\beta(\beta^2-3)} \ 
\begin{pspicture}[shift=-0.5](0,0.00)(1.2,1.2)
\pspolygon[fillstyle=solid,fillcolor=pink](0,0)(1.2,0)(1.2,0.3)(0,0.3)(0,0)\rput(0.6,0.15){$_3$}
\rput(0,0.9){\pspolygon[fillstyle=solid,fillcolor=pink](0,0)(1.2,0)(1.2,0.3)(0,0.3)(0,0)\rput(0.6,0.15){$_3$}}
\pspolygon[fillstyle=solid,fillcolor=lightlightblue](0,0.3)(1.2,0.3)(1.2,0.9)(0,0.9)(0,0.3)
\psarc[linecolor=blue,linewidth=1.5pt]{-}(0,0.3){0.2}{0}{90}
\psarc[linecolor=blue,linewidth=1.5pt]{-}(1.2,0.3){0.2}{90}{180}
\psarc[linecolor=blue,linewidth=1.5pt]{-}(0,0.9){0.2}{-90}{0}
\psarc[linecolor=blue,linewidth=1.5pt]{-}(1.2,0.9){0.2}{180}{270}
\psline[linecolor=blue,linewidth=1.5pt]{-}(0.6,0.3)(0.6,0.9)
\end{pspicture} \ ,
\\[0.2cm]
\qptla_4(\beta,\alpha) 
&: \quad Q_4 = \
\begin{pspicture}[shift=-0.05](0,0.00)(1.6,0.3)
\pspolygon[fillstyle=solid,fillcolor=pink](0,0)(1.6,0)(1.6,0.3)(0,0.3)(0,0)\rput(0.8,0.15){$_4$}
\end{pspicture} 
\ - \frac{\beta^2-2}{\beta(\beta^2-4)} \ 
\begin{pspicture}[shift=-0.5](0,0.00)(1.6,1.2)
\pspolygon[fillstyle=solid,fillcolor=pink](0,0)(1.6,0)(1.6,0.3)(0,0.3)(0,0)\rput(0.8,0.15){$_4$}
\rput(0,0.9){\pspolygon[fillstyle=solid,fillcolor=pink](0,0)(1.6,0)(1.6,0.3)(0,0.3)(0,0)\rput(0.8,0.15){$_4$}}
\pspolygon[fillstyle=solid,fillcolor=lightlightblue](0,0.3)(1.6,0.3)(1.6,0.9)(0,0.9)(0,0.3)
\psarc[linecolor=blue,linewidth=1.5pt]{-}(0,0.3){0.2}{0}{90}
\psarc[linecolor=blue,linewidth=1.5pt]{-}(1.6,0.3){0.2}{90}{180}
\psarc[linecolor=blue,linewidth=1.5pt]{-}(0,0.9){0.2}{-90}{0}
\psarc[linecolor=blue,linewidth=1.5pt]{-}(1.6,0.9){0.2}{180}{270}
\psline[linecolor=blue,linewidth=1.5pt]{-}(0.6,0.3)(0.6,0.9)
\psline[linecolor=blue,linewidth=1.5pt]{-}(1.0,0.3)(1.0,0.9)
\end{pspicture} 
\ - \frac{(\beta^2-2)(\beta^2-4)^{-1}}{\alpha^2-(\beta^2-2)^2} \ 
\begin{pspicture}[shift=-0.75](0,0.00)(1.6,1.7)
\pspolygon[fillstyle=solid,fillcolor=pink](0,0)(1.6,0)(1.6,0.3)(0,0.3)(0,0)\rput(0.8,0.15){$_4$}
\rput(0,1.4){\pspolygon[fillstyle=solid,fillcolor=pink](0,0)(1.6,0)(1.6,0.3)(0,0.3)(0,0)\rput(0.8,0.15){$_4$}}
\pspolygon[fillstyle=solid,fillcolor=lightlightblue](0,0.3)(1.6,0.3)(1.6,1.4)(0,1.4)(0,0.3)
\psarc[linecolor=blue,linewidth=1.5pt]{-}(0,0.3){0.2}{0}{90}
\psarc[linecolor=blue,linewidth=1.5pt]{-}(1.6,0.3){0.2}{90}{180}
\psarc[linecolor=blue,linewidth=1.5pt]{-}(0,1.4){0.2}{-90}{0}
\psarc[linecolor=blue,linewidth=1.5pt]{-}(1.6,1.4){0.2}{180}{270}
\psbezier[linecolor=blue,linewidth=1.5pt]{-}(0.6,0.3)(0.6,0.6)(0.1,0.7)(0,0.7)
\psbezier[linecolor=blue,linewidth=1.5pt]{-}(0.6,1.4)(0.6,1.1)(0.1,1.0)(0,1.0)
\psbezier[linecolor=blue,linewidth=1.5pt]{-}(1,0.3)(1,0.6)(1.5,0.7)(1.6,0.7)
\psbezier[linecolor=blue,linewidth=1.5pt]{-}(1,1.4)(1,1.1)(1.5,1.0)(1.6,1.0)
\end{pspicture} 
 \ ,
\\[0.2cm]
\qptlb_4(\beta,\gamma) 
&: \quad Q_4 = \
\begin{pspicture}[shift=-0.05](0,0.00)(1.6,0.3)
\pspolygon[fillstyle=solid,fillcolor=pink](0,0)(1.6,0)(1.6,0.3)(0,0.3)(0,0)\rput(0.8,0.15){$_4$}
\end{pspicture} 
\ + \frac{\beta(\beta^2-2)}{\gamma+\gamma^{-1}-\beta^4+4\beta^2-2} \ 
\begin{pspicture}[shift=-0.5](0,0.00)(1.6,1.2)
\pspolygon[fillstyle=solid,fillcolor=pink](0,0)(1.6,0)(1.6,0.3)(0,0.3)(0,0)\rput(0.8,0.15){$_4$}
\rput(0,0.9){\pspolygon[fillstyle=solid,fillcolor=pink](0,0)(1.6,0)(1.6,0.3)(0,0.3)(0,0)\rput(0.8,0.15){$_4$}}
\pspolygon[fillstyle=solid,fillcolor=lightlightblue](0,0.3)(1.6,0.3)(1.6,0.9)(0,0.9)(0,0.3)
\psarc[linecolor=blue,linewidth=1.5pt]{-}(0,0.3){0.2}{0}{90}
\psarc[linecolor=blue,linewidth=1.5pt]{-}(1.6,0.3){0.2}{90}{180}
\psarc[linecolor=blue,linewidth=1.5pt]{-}(0,0.9){0.2}{-90}{0}
\psarc[linecolor=blue,linewidth=1.5pt]{-}(1.6,0.9){0.2}{180}{270}
\psline[linecolor=blue,linewidth=1.5pt]{-}(0.6,0.3)(0.6,0.9)
\psline[linecolor=blue,linewidth=1.5pt]{-}(1.0,0.3)(1.0,0.9)
\end{pspicture} \ .
\end{alignat}
\end{subequations}
These examples involve only connectivities $c$ with $\ell = 0$, however contributions with $\ell>0$ do in fact arise for larger values of $n$. The details of the general constructions depend on the parity of $n$ and the uncoiled algebra considered.
They are given by
\begin{subequations}\label{eq:WJ.upTL}
\begin{alignat}{3}
&\qptl_n(\beta,\gamma): \quad 
&&Q_n = \
 \begin{pspicture}[shift=-0.05](-0.03,0.00)(1.62,0.3)
 \pspolygon[fillstyle=solid,fillcolor=pink](0,0)(1.4,0)(1.4,0.3)(0,0.3)(0,0)\rput(0.7,0.15){$_n$}
 \end{pspicture} + 
 \sum_{k=1}^{\frac{n-1}2} \sum_{\ell = 0}^{n-2k-1}\Gamma_{k,\ell} \
\begin{pspicture}[shift=-1.05](-0.2,0)(2.2,2.3)
\pspolygon[fillstyle=solid,fillcolor=pink](0,0)(2.2,0)(2.2,0.3)(0,0.3)(0,0)\rput(1.1,0.15){$_{n}$}
\rput(0,2){\pspolygon[fillstyle=solid,fillcolor=pink](0,0)(2.2,0)(2.2,0.3)(0,0.3)(0,0)\rput(1.1,0.15){$_{n}$}}
\pspolygon[fillstyle=solid,fillcolor=lightlightblue](0,0.3)(2.2,0.3)(2.2,0.8)(0,0.8)(0,0.3)
\pspolygon[fillstyle=solid,fillcolor=lightlightblue](0,2)(2.2,2)(2.2,1.5)(0,1.5)(0,2)
\pspolygon[fillstyle=solid,fillcolor=lightlightblue](0.4,0.8)(1.8,0.8)(1.8,1.5)(0.4,1.5)(0.4,0.8)
\psarc[linecolor=blue,linewidth=1.5pt]{-}(2.2,0.3){0.2}{90}{180}
\psarc[linecolor=blue,linewidth=1.5pt]{-}(0,0.3){0.2}{0}{90}
\psarc[linecolor=blue,linewidth=1.5pt]{-}(2.2,2){0.2}{180}{270}
\psarc[linecolor=blue,linewidth=1.5pt]{-}(0,2){0.2}{270}{0}
\psline[linecolor=blue,linewidth=1.5pt]{-}(1.1,0.3)(1.1,0.8)
\psline[linecolor=blue,linewidth=1.5pt]{-}(1.1,1.5)(1.1,2)
\rput(-0.15,0.5){$_k$}\rput(2.35,0.5){$_k$}
\rput(-0.15,1.8){$_k$}\rput(2.35,1.8){$_k$}
\rput(1.55,0.55){$_{n-2k}$}
\rput(1.55,1.75){$_{n-2k}$}
\rput(1.1,1.15){$\Omega^{2\ell}$}
\end{pspicture}\ ,
\\[0.4cm]
&\qptla_n(\beta,\alpha): \quad 
&&Q_n = \
\begin{pspicture}[shift=-0.05](-0.03,0.00)(1.47,0.3)
\pspolygon[fillstyle=solid,fillcolor=pink](0,0)(1.4,0)(1.4,0.3)(0,0.3)(0,0)\rput(0.7,0.15){$_n$}
\end{pspicture}\ + 
\sum_{k=1}^{\frac{n-2}2} \sum_{\ell = 0}^{\frac {n-2k}2-1}\Gamma_{k,\ell} \
\begin{pspicture}[shift=-1.05](-0.2,0)(2.2,2.3)
\pspolygon[fillstyle=solid,fillcolor=pink](0,0)(2.2,0)(2.2,0.3)(0,0.3)(0,0)\rput(1.1,0.15){$_{n}$}
\rput(0,2){\pspolygon[fillstyle=solid,fillcolor=pink](0,0)(2.2,0)(2.2,0.3)(0,0.3)(0,0)\rput(1.1,0.15){$_{n}$}}
\pspolygon[fillstyle=solid,fillcolor=lightlightblue](0,0.3)(2.2,0.3)(2.2,0.8)(0,0.8)(0,0.3)
\pspolygon[fillstyle=solid,fillcolor=lightlightblue](0,2)(2.2,2)(2.2,1.5)(0,1.5)(0,2)
\pspolygon[fillstyle=solid,fillcolor=lightlightblue](0.4,0.8)(1.8,0.8)(1.8,1.5)(0.4,1.5)(0.4,0.8)
\psarc[linecolor=blue,linewidth=1.5pt]{-}(2.2,0.3){0.2}{90}{180}
\psarc[linecolor=blue,linewidth=1.5pt]{-}(0,0.3){0.2}{0}{90}
\psarc[linecolor=blue,linewidth=1.5pt]{-}(2.2,2){0.2}{180}{270}
\psarc[linecolor=blue,linewidth=1.5pt]{-}(0,2){0.2}{270}{0}
\psline[linecolor=blue,linewidth=1.5pt]{-}(1.1,0.3)(1.1,0.8)
\psline[linecolor=blue,linewidth=1.5pt]{-}(1.1,1.5)(1.1,2)
\rput(-0.15,0.5){$_k$}\rput(2.35,0.5){$_k$}
\rput(-0.15,1.8){$_k$}\rput(2.35,1.8){$_k$}
\rput(1.55,0.55){$_{n-2k}$}
\rput(1.55,1.75){$_{n-2k}$}
\rput(1.1,1.15){$\Omega^{2\ell}$}
\end{pspicture}
\ \ + \Gamma_{n/2,0} \ \
\begin{pspicture}[shift=-0.7](-0.2,0)(1.6,1.6)
\pspolygon[fillstyle=solid,fillcolor=pink](0,0)(1.4,0)(1.4,0.3)(0,0.3)(0,0)\rput(0.7,0.15){$_{n}$}
\rput(0,1.3){\pspolygon[fillstyle=solid,fillcolor=pink](0,0)(1.4,0)(1.4,0.3)(0,0.3)(0,0)\rput(0.7,0.15){$_{n}$}}
\pspolygon[fillstyle=solid,fillcolor=lightlightblue](0,0.3)(1.4,0.3)(1.4,1.3)(0,1.3)(0,0.3)
\psarc[linecolor=blue,linewidth=1.5pt]{-}(1.4,0.3){0.2}{90}{180}
\psarc[linecolor=blue,linewidth=1.5pt]{-}(0,0.3){0.2}{0}{90}
\psarc[linecolor=blue,linewidth=1.5pt]{-}(1.4,1.3){0.2}{180}{270}
\psarc[linecolor=blue,linewidth=1.5pt]{-}(0,1.3){0.2}{270}{0}
\rput(-0.15,0.5){$_{\frac n2}$}\rput(1.55,0.5){$_{\frac n2}$}
\rput(-0.15,1.1){$_{\frac n2}$}\rput(1.55,1.1){$_{\frac n2}$}
\end{pspicture}\ ,
\\[0.4cm]
&\qptlb_n(\beta,\gamma): \quad 
&&Q_n = \
 \begin{pspicture}[shift=-0.05](-0.03,0.00)(1.47,0.3)
 \pspolygon[fillstyle=solid,fillcolor=pink](0,0)(1.4,0)(1.4,0.3)(0,0.3)(0,0)\rput(0.7,0.15){$_n$}
 \end{pspicture}\ + 
 \sum_{k=1}^{\frac{n-2}2} \sum_{\ell = 0}^{\frac {n-2k-2}2}\Gamma_{k,\ell} \
\begin{pspicture}[shift=-1.05](-0.2,0)(2.4,2.3)
\pspolygon[fillstyle=solid,fillcolor=pink](0,0)(2.2,0)(2.2,0.3)(0,0.3)(0,0)\rput(1.1,0.15){$_{n}$}
\rput(0,2){\pspolygon[fillstyle=solid,fillcolor=pink](0,0)(2.2,0)(2.2,0.3)(0,0.3)(0,0)\rput(1.1,0.15){$_{n}$}}
\pspolygon[fillstyle=solid,fillcolor=lightlightblue](0,0.3)(2.2,0.3)(2.2,0.8)(0,0.8)(0,0.3)
\pspolygon[fillstyle=solid,fillcolor=lightlightblue](0,2)(2.2,2)(2.2,1.5)(0,1.5)(0,2)
\pspolygon[fillstyle=solid,fillcolor=lightlightblue](0.4,0.8)(1.8,0.8)(1.8,1.5)(0.4,1.5)(0.4,0.8)
\psarc[linecolor=blue,linewidth=1.5pt]{-}(2.2,0.3){0.2}{90}{180}
\psarc[linecolor=blue,linewidth=1.5pt]{-}(0,0.3){0.2}{0}{90}
\psarc[linecolor=blue,linewidth=1.5pt]{-}(2.2,2){0.2}{180}{270}
\psarc[linecolor=blue,linewidth=1.5pt]{-}(0,2){0.2}{270}{0}
\psline[linecolor=blue,linewidth=1.5pt]{-}(1.1,0.3)(1.1,0.8)
\psline[linecolor=blue,linewidth=1.5pt]{-}(1.1,1.5)(1.1,2)
\rput(-0.15,0.5){$_k$}\rput(2.35,0.5){$_k$}
\rput(-0.15,1.8){$_k$}\rput(2.35,1.8){$_k$}
\rput(1.55,0.55){$_{n-2k}$}
\rput(1.55,1.75){$_{n-2k}$}
\rput(1.1,1.15){$\Omega^{2\ell}$}
\end{pspicture}\ .
\end{alignat}
\end{subequations}
In each case, all terms of the sum are even connectivities with at least one $e_j$, and are thus indeed elements of $\qptl_n(\beta,\gamma)$, $\qptla_n(\beta,\alpha)$ or $\qptlb_n(\beta,\gamma)$. Note that the constants $\Gamma_{k,\ell}$ of~\eqref{eq:WJ.upTL} are different in each case and share the same notation simply for convenience in the derivation of their explicit formulas in \cref{sec:Gammas}.

\subsection[Projectors for $\qatl_n(\beta,\gamma)$, $\qatla_n(\beta, \alpha)$ and $\qatlb_n(\beta,\gamma)$]{Projectors for $\boldsymbol{\qatl_n(\beta,\gamma)}$, $\boldsymbol{\qatla_n(\beta,\alpha)}$ and $\boldsymbol{\qatlb_n(\beta,\gamma)}$}\label{sec:proj.qatl}

In this section, we construct the Wenzl--Jones projectors for the quotient algebras $\qatl_n(\beta,\gamma)$, $\qatla_n(\beta,\alpha)$ and $\qatlb_n(\beta,\gamma)$. In these cases, as discussed in \cref{sec:properties.quotients}, there are $n$ one-dimensional standard modules, labelled by an integer $r=0, 1, \dots, n-1$. We construct a projector $Q_{n,r}$, for each~$r$, that projects on the one-dimensional standard modules labelled by $r$, as given in \cref{prop:standards.uncoiled.even,prop:standards.uncoiled.odd}. As a first remark, we note that these projectors are easily obtained from the projectors $Q_n$ constructed in \cref{sec:proj.qptl}. Indeed, let us write
\be
\omega_{n,r} = \left\{\begin{array}{cc}
\gamma^{1/n}\eE^{2 \pi \iI r/n}& \qatl_n(\beta,\gamma) \textrm{ and }\qatlb_{n}(\beta,\gamma),\\[0.15cm]
\eE^{2 \pi \iI r/n} &  \qatla_n(\beta,\alpha),
\end{array}\right.
\qquad 0 \le r < n.
\ee
These are the possible eigenvalues of $\Omega$, consistent with the quotient relations $\Omega^n = \gamma\, \id$ and $\Omega^n = \id$, respectively. Then the projectors $\Pi_{n,r}$ on the eigenspaces of $\Omega$ of eigenvalue $\omega_{n,r}$ are 
\be
\Pi_{n,r} = \frac1n\sum_{j=0}^{n-1} (\omega_{n,r})^{-j}\Omega^j.
\ee
One way to construct the Wenzl--Jones projectors $Q_{n,r}$ is then
\be
Q_{n,r} = \Pi_{n,r} Q_n.
\ee
The commutation relation $[\Omega,Q_n]=0$ then implies that $Q_{n,r} = Q_n \Pi_{n,r}$.
From the properties of $Q_n$, the projector $Q_{n,r}$ is then readily seen to satisfy the relations
\be
(Q_{n,r})^2 = Q_{n,r}, 
\qquad 
\Omega\, Q_{n,r}=Q_{n,r}\,\Omega = \omega_{n,r}Q_{n,r}, 
\qquad
e_j Q_{n,r}=Q_{n,r} e_j=0, \quad j = 0,1, \dots, n-1.
\ee
It is the unique element of the corresponding uncoiled affine Temperley--Lieb algebra satisfying these relations. The proof of this statement is analogous to the same proof given above for the projectors~$Q_n$. We also note that $Q_{n,r}$ evaluates to~$1$ on $\repW_{n,n,\omega_{n,r}}$, whereas it vanishes on all the other standard modules of the algebra. Finally, these projectors also satisfy the property 
\be
\label{eq:Qn=sumQnr}
Q_n = 
\sum_{r=0}^{n-1} Q_{n,r}\,.
\ee

A second equivalent construction of $Q_{n,r}$ is possible in terms of the sandwich diagrams:
\be
\label{eq:Qnr.idea}
Q_{n,r} = \
\begin{pspicture}[shift=-0.75](0,-0.2)(1.4,1.5)
\pspolygon[fillstyle=solid,fillcolor=pink](0,0)(1.4,0)(1.4,0.3)(0,0.3)(0,0)\rput(0.7,0.15){$_{n}$}
\rput(0,1){\pspolygon[fillstyle=solid,fillcolor=pink](0,0)(1.4,0)(1.4,0.3)(0,0.3)(0,0)\rput(0.7,0.15){$_{n}$}}
\pspolygon[fillstyle=solid,fillcolor=lightlightblue](0,0.3)(1.4,0.3)(1.4,1)(0,1)(0,0.3)\rput(0.7,0.65){$\Pi_{n,r}$}
\end{pspicture}
 \ + \sum_{c} \Gamma(c) \ 
\begin{pspicture}[shift=-0.75](0,-0.2)(1.4,1.5)
\pspolygon[fillstyle=solid,fillcolor=pink](0,0)(1.4,0)(1.4,0.3)(0,0.3)(0,0)\rput(0.7,0.15){$_{n}$}
\rput(0,1){\pspolygon[fillstyle=solid,fillcolor=pink](0,0)(1.4,0)(1.4,0.3)(0,0.3)(0,0)\rput(0.7,0.15){$_{n}$}}
\pspolygon[fillstyle=solid,fillcolor=lightlightblue](0,0.3)(1.4,0.3)(1.4,1)(0,1)(0,0.3)\rput(0.7,0.65){$c$}
\end{pspicture} \ .
\ee
Here, the connectivities $c$ in the sum are elements of the uncoiled affine Temperley--Lieb algebras that always have at least one generator~$e_j$. These can be both odd and even elements, namely the connectivities $c$ are of the form \eqref{eq:vwc} with $\ell \in \frac12\mathbb Z$. The prefactor in front of the first term is again equal to $1$, which is the only value ensuring that $(Q_{n,r})^2 = Q_{n,r}$. The details of this construction again depend on the algebra considered. For example, the projectors for $n=3$ are
\be
\qatl_3(\beta,\gamma): \quad
Q_{3,r} = \frac13 \left[\ \
\begin{pspicture}[shift=-0.5](0,0.00)(1.2,1.2)
\pspolygon[fillstyle=solid,fillcolor=pink](0,0)(1.2,0)(1.2,0.3)(0,0.3)(0,0)\rput(0.6,0.15){$_3$}
\rput(0,0.9){\pspolygon[fillstyle=solid,fillcolor=pink](0,0)(1.2,0)(1.2,0.3)(0,0.3)(0,0)\rput(0.6,0.15){$_3$}}
\pspolygon[fillstyle=solid,fillcolor=lightlightblue](0,0.3)(1.2,0.3)(1.2,0.9)(0,0.9)(0,0.3)
\psline[linecolor=blue,linewidth=1.5pt]{-}(0.2,0.3)(0.2,0.9)
\psline[linecolor=blue,linewidth=1.5pt]{-}(0.6,0.3)(0.6,0.9)
\psline[linecolor=blue,linewidth=1.5pt]{-}(1.0,0.3)(1.0,0.9)
\end{pspicture}
\ + \omega \
\begin{pspicture}[shift=-0.5](0,0.00)(1.2,1.2)
\pspolygon[fillstyle=solid,fillcolor=pink](0,0)(1.2,0)(1.2,0.3)(0,0.3)(0,0)\rput(0.6,0.15){$_3$}
\rput(0,0.9){\pspolygon[fillstyle=solid,fillcolor=pink](0,0)(1.2,0)(1.2,0.3)(0,0.3)(0,0)\rput(0.6,0.15){$_3$}}
\pspolygon[fillstyle=solid,fillcolor=lightlightblue](0,0.3)(1.2,0.3)(1.2,0.9)(0,0.9)(0,0.3)
\rput(0,0.6){
\psbezier[linecolor=blue,linewidth=1.5pt]{-}(0.2,0.3)(0.2,0.1)(-0.02,0)(-0.05,-0.05)
\rput(1.2,0){\psbezier[linecolor=blue,linewidth=1.5pt]{-}(-0.2,-0.3)(-0.2,-0.1)(0.02,0)(0.05,0.05)}
\multiput(0.4,0)(0.4,0){2}{\psbezier[linecolor=blue,linewidth=1.5pt]{-}(-0.2,-0.30)(-0.2,-0.0)(0.2,0.0)(0.2,0.30)}
\psframe[fillstyle=solid,linecolor=white,linewidth=0pt](-0.1,-0.4)(-0.02,0.4)
\psframe[fillstyle=solid,linecolor=white,linewidth=0pt](1.22,-0.4)(1.3,0.4)
}
\end{pspicture}
\ + \omega^2 \
\begin{pspicture}[shift=-0.5](0,0.00)(1.2,1.2)
\pspolygon[fillstyle=solid,fillcolor=pink](0,0)(1.2,0)(1.2,0.3)(0,0.3)(0,0)\rput(0.6,0.15){$_3$}
\rput(0,0.9){\pspolygon[fillstyle=solid,fillcolor=pink](0,0)(1.2,0)(1.2,0.3)(0,0.3)(0,0)\rput(0.6,0.15){$_3$}}
\pspolygon[fillstyle=solid,fillcolor=lightlightblue](0,0.3)(1.2,0.3)(1.2,0.9)(0,0.9)(0,0.3)
\psbezier[linecolor=blue,linewidth=1.5pt]{-}(0.2,0.3)(0.2,0.6)(1.0,0.6)(1.0,0.9)
\psbezier[linecolor=blue,linewidth=1.5pt]{-}(0.6,0.3)(0.6,0.48)(1.0,0.6)(1.21,0.67)
\psbezier[linecolor=blue,linewidth=1.5pt]{-}(1.0,0.3)(1.0,0.40)(1.1,0.48)(1.21,0.52)
\psbezier[linecolor=blue,linewidth=1.5pt]{-}(0.6,0.9)(0.6,0.72)(0.2,0.6)(-0.01,0.53)
\psbezier[linecolor=blue,linewidth=1.5pt]{-}(0.2,0.9)(0.2,0.8)(0.1,0.72)(-0.01,0.68)
\end{pspicture}
\ - \frac1{\omega^2+\omega^{-2}+\beta} \ 
\begin{pspicture}[shift=-0.5](0,0.00)(1.2,1.2)
\pspolygon[fillstyle=solid,fillcolor=pink](0,0)(1.2,0)(1.2,0.3)(0,0.3)(0,0)\rput(0.6,0.15){$_3$}
\rput(0,0.9){\pspolygon[fillstyle=solid,fillcolor=pink](0,0)(1.2,0)(1.2,0.3)(0,0.3)(0,0)\rput(0.6,0.15){$_3$}}
\pspolygon[fillstyle=solid,fillcolor=lightlightblue](0,0.3)(1.2,0.3)(1.2,0.9)(0,0.9)(0,0.3)
\psarc[linecolor=blue,linewidth=1.5pt]{-}(0,0.3){0.2}{0}{90}
\psarc[linecolor=blue,linewidth=1.5pt]{-}(1.2,0.3){0.2}{90}{180}
\psarc[linecolor=blue,linewidth=1.5pt]{-}(0,0.9){0.2}{-90}{0}
\psarc[linecolor=blue,linewidth=1.5pt]{-}(1.2,0.9){0.2}{180}{270}
\psline[linecolor=blue,linewidth=1.5pt]{-}(0.6,0.3)(0.6,0.9)
\end{pspicture}
\ \ \right],
\ee
with $\omega = \gamma^{1/3} \eE^{2\pi \iI r/3}$ and $r \in \{0,1,2\}$.\medskip 

In general, the projectors are given by
\begin{subequations}\label{eq:WJ.uaTL}
\begin{alignat}{3}
&\qatl_n(\beta,\gamma):
\qquad &&Q_{n,r} = \
\begin{pspicture}[shift=-0.75](0,-0.2)(1.4,1.5)
\pspolygon[fillstyle=solid,fillcolor=pink](0,0)(1.4,0)(1.4,0.3)(0,0.3)(0,0)\rput(0.7,0.15){$_{n}$}
\rput(0,1){\pspolygon[fillstyle=solid,fillcolor=pink](0,0)(1.4,0)(1.4,0.3)(0,0.3)(0,0)\rput(0.7,0.15){$_{n}$}}
\pspolygon[fillstyle=solid,fillcolor=lightlightblue](0,0.3)(1.4,0.3)(1.4,1)(0,1)(0,0.3)\rput(0.7,0.65){$\Pi_{n,r}$}
\end{pspicture}\ + 
 \sum_{k=1}^{\frac{n-1}2} \sum_{\ell = 0,\frac12,\dots}^{\frac{n-2k-1}2}\Gamma_{k,\ell} \
\begin{pspicture}[shift=-1.05](-0.2,0)(2.2,2.3)
\pspolygon[fillstyle=solid,fillcolor=pink](0,0)(2.2,0)(2.2,0.3)(0,0.3)(0,0)\rput(1.1,0.15){$_{n}$}
\rput(0,2){\pspolygon[fillstyle=solid,fillcolor=pink](0,0)(2.2,0)(2.2,0.3)(0,0.3)(0,0)\rput(1.1,0.15){$_{n}$}}
\pspolygon[fillstyle=solid,fillcolor=lightlightblue](0,0.3)(2.2,0.3)(2.2,0.8)(0,0.8)(0,0.3)
\pspolygon[fillstyle=solid,fillcolor=lightlightblue](0,2)(2.2,2)(2.2,1.5)(0,1.5)(0,2)
\pspolygon[fillstyle=solid,fillcolor=lightlightblue](0.4,0.8)(1.8,0.8)(1.8,1.5)(0.4,1.5)(0.4,0.8)
\psarc[linecolor=blue,linewidth=1.5pt]{-}(2.2,0.3){0.2}{90}{180}
\psarc[linecolor=blue,linewidth=1.5pt]{-}(0,0.3){0.2}{0}{90}
\psarc[linecolor=blue,linewidth=1.5pt]{-}(2.2,2){0.2}{180}{270}
\psarc[linecolor=blue,linewidth=1.5pt]{-}(0,2){0.2}{270}{0}
\psline[linecolor=blue,linewidth=1.5pt]{-}(1.1,0.3)(1.1,0.8)
\psline[linecolor=blue,linewidth=1.5pt]{-}(1.1,1.5)(1.1,2)
\rput(-0.15,0.5){$_k$}\rput(2.35,0.5){$_k$}
\rput(-0.15,1.8){$_k$}\rput(2.35,1.8){$_k$}
\rput(1.55,0.55){$_{n-2k}$}
\rput(1.55,1.75){$_{n-2k}$}
\rput(1.1,1.15){$\Omega^{2\ell}$}
\end{pspicture}\ ,
\\[0.3cm]
&\qatla_n(\beta,\alpha):
\qquad
&&Q_{n,r} = \
\begin{pspicture}[shift=-0.75](0,-0.2)(1.4,1.5)
\pspolygon[fillstyle=solid,fillcolor=pink](0,0)(1.4,0)(1.4,0.3)(0,0.3)(0,0)\rput(0.7,0.15){$_{n}$}
\rput(0,1){\pspolygon[fillstyle=solid,fillcolor=pink](0,0)(1.4,0)(1.4,0.3)(0,0.3)(0,0)\rput(0.7,0.15){$_{n}$}}
\pspolygon[fillstyle=solid,fillcolor=lightlightblue](0,0.3)(1.4,0.3)(1.4,1)(0,1)(0,0.3)\rput(0.7,0.65){$\Pi_{n,r}$}
\end{pspicture}\ + 
 \sum_{k=1}^{\frac{n-2}2} \sum_{\ell = 0,\frac12,\dots}^{\frac{n-2k-1}2}\Gamma_{k,\ell} \
\begin{pspicture}[shift=-1.05](-0.2,0)(2.2,2.3)
\pspolygon[fillstyle=solid,fillcolor=pink](0,0)(2.2,0)(2.2,0.3)(0,0.3)(0,0)\rput(1.1,0.15){$_{n}$}
\rput(0,2){\pspolygon[fillstyle=solid,fillcolor=pink](0,0)(2.2,0)(2.2,0.3)(0,0.3)(0,0)\rput(1.1,0.15){$_{n}$}}
\pspolygon[fillstyle=solid,fillcolor=lightlightblue](0,0.3)(2.2,0.3)(2.2,0.8)(0,0.8)(0,0.3)
\pspolygon[fillstyle=solid,fillcolor=lightlightblue](0,2)(2.2,2)(2.2,1.5)(0,1.5)(0,2)
\pspolygon[fillstyle=solid,fillcolor=lightlightblue](0.4,0.8)(1.8,0.8)(1.8,1.5)(0.4,1.5)(0.4,0.8)
\psarc[linecolor=blue,linewidth=1.5pt]{-}(2.2,0.3){0.2}{90}{180}
\psarc[linecolor=blue,linewidth=1.5pt]{-}(0,0.3){0.2}{0}{90}
\psarc[linecolor=blue,linewidth=1.5pt]{-}(2.2,2){0.2}{180}{270}
\psarc[linecolor=blue,linewidth=1.5pt]{-}(0,2){0.2}{270}{0}
\psline[linecolor=blue,linewidth=1.5pt]{-}(1.1,0.3)(1.1,0.8)
\psline[linecolor=blue,linewidth=1.5pt]{-}(1.1,1.5)(1.1,2)
\rput(-0.15,0.5){$_k$}\rput(2.35,0.5){$_k$}
\rput(-0.15,1.8){$_k$}\rput(2.35,1.8){$_k$}
\rput(1.55,0.55){$_{n-2k}$}
\rput(1.55,1.75){$_{n-2k}$}
\rput(1.1,1.15){$\Omega^{2\ell}$}
\end{pspicture}
\ \ + \Gamma_{n/2,0} \ \
\begin{pspicture}[shift=-0.7](-0.2,0)(1.6,1.6)
\pspolygon[fillstyle=solid,fillcolor=pink](0,0)(1.4,0)(1.4,0.3)(0,0.3)(0,0)\rput(0.7,0.15){$_{n}$}
\rput(0,1.3){\pspolygon[fillstyle=solid,fillcolor=pink](0,0)(1.4,0)(1.4,0.3)(0,0.3)(0,0)\rput(0.7,0.15){$_{n}$}}
\pspolygon[fillstyle=solid,fillcolor=lightlightblue](0,0.3)(1.4,0.3)(1.4,1.3)(0,1.3)(0,0.3)
\psarc[linecolor=blue,linewidth=1.5pt]{-}(1.4,0.3){0.2}{90}{180}
\psarc[linecolor=blue,linewidth=1.5pt]{-}(0,0.3){0.2}{0}{90}
\psarc[linecolor=blue,linewidth=1.5pt]{-}(1.4,1.3){0.2}{180}{270}
\psarc[linecolor=blue,linewidth=1.5pt]{-}(0,1.3){0.2}{270}{0}
\rput(-0.15,0.5){$_{\frac n2}$}\rput(1.55,0.5){$_{\frac n2}$}
\rput(-0.15,1.1){$_{\frac n2}$}\rput(1.55,1.1){$_{\frac n2}$}
\end{pspicture}\ ,
\\[0.3cm]
&\qatlb_n(\beta,\gamma):
\qquad
&&Q_{n,r} = \
\begin{pspicture}[shift=-0.75](0,-0.2)(1.4,1.5)
\pspolygon[fillstyle=solid,fillcolor=pink](0,0)(1.4,0)(1.4,0.3)(0,0.3)(0,0)\rput(0.7,0.15){$_{n}$}
\rput(0,1){\pspolygon[fillstyle=solid,fillcolor=pink](0,0)(1.4,0)(1.4,0.3)(0,0.3)(0,0)\rput(0.7,0.15){$_{n}$}}
\pspolygon[fillstyle=solid,fillcolor=lightlightblue](0,0.3)(1.4,0.3)(1.4,1)(0,1)(0,0.3)\rput(0.7,0.65){$\Pi_{n,r}$}
\end{pspicture}\ + 
 \sum_{k=1}^{\frac{n-2}2} \sum_{\ell = 0,\frac12,\dots}^{\frac{n-2k-1}2}\Gamma_{k,\ell} \
\begin{pspicture}[shift=-1.05](-0.2,0)(2.4,2.3)
\pspolygon[fillstyle=solid,fillcolor=pink](0,0)(2.2,0)(2.2,0.3)(0,0.3)(0,0)\rput(1.1,0.15){$_{n}$}
\rput(0,2){\pspolygon[fillstyle=solid,fillcolor=pink](0,0)(2.2,0)(2.2,0.3)(0,0.3)(0,0)\rput(1.1,0.15){$_{n}$}}
\pspolygon[fillstyle=solid,fillcolor=lightlightblue](0,0.3)(2.2,0.3)(2.2,0.8)(0,0.8)(0,0.3)
\pspolygon[fillstyle=solid,fillcolor=lightlightblue](0,2)(2.2,2)(2.2,1.5)(0,1.5)(0,2)
\pspolygon[fillstyle=solid,fillcolor=lightlightblue](0.4,0.8)(1.8,0.8)(1.8,1.5)(0.4,1.5)(0.4,0.8)
\psarc[linecolor=blue,linewidth=1.5pt]{-}(2.2,0.3){0.2}{90}{180}
\psarc[linecolor=blue,linewidth=1.5pt]{-}(0,0.3){0.2}{0}{90}
\psarc[linecolor=blue,linewidth=1.5pt]{-}(2.2,2){0.2}{180}{270}
\psarc[linecolor=blue,linewidth=1.5pt]{-}(0,2){0.2}{270}{0}
\psline[linecolor=blue,linewidth=1.5pt]{-}(1.1,0.3)(1.1,0.8)
\psline[linecolor=blue,linewidth=1.5pt]{-}(1.1,1.5)(1.1,2)
\rput(-0.15,0.5){$_k$}\rput(2.35,0.5){$_k$}
\rput(-0.15,1.8){$_k$}\rput(2.35,1.8){$_k$}
\rput(1.55,0.55){$_{n-2k}$}
\rput(1.55,1.75){$_{n-2k}$}
\rput(1.1,1.15){$\Omega^{2\ell}$}
\end{pspicture}\ .
\end{alignat}
\end{subequations}
The values of the constants $\Gamma_{k,\ell}$ are given in \cref{sec:Gammas}.

%
\section{Calculation of the constants \texorpdfstring{$\boldsymbol{\Gamma_{k,\ell}}$}{Gamma kl}}\label{sec:Gammas}
%

In this section, we compute the constants $\Gamma_{k,\ell}$ that arise in the definition of the projectors $Q_n$ and $Q_{n,r}$. We obtain, in \cref{Prop.const.rel} below, a 
 linear recurrence relation satisfied by the constants $\Gamma_{k,\ell}$. As explained in \cref{Sec:Fourier,sec:kernel}, this gives a way to compute all the constants. Then in \cref{Sec:Conj.Values}, we use this method to compute the values of $\Gamma_{k,\ell}$ for the lowest values of $k$, which then enables us to give closed-form formulas for all $k$. The closed-form formulas are proved in \cref{Sec:proof.closed.form}. 

\subsection{Linear relations satisfied by the constants}
\label{sec:linear.relations}

Let us start by defining the following three diagrammatic objects:
\be
\label{eq:ZYX}
Z_{k,\ell} = \ 
\begin{pspicture}[shift=-1.05](-0.2,0)(2.4,2.3)
\pspolygon[fillstyle=solid,fillcolor=pink](0,0)(2.2,0)(2.2,0.3)(0,0.3)(0,0)\rput(1.1,0.15){$_{n}$}
\rput(0,2){\pspolygon[fillstyle=solid,fillcolor=pink](0,0)(2.2,0)(2.2,0.3)(0,0.3)(0,0)\rput(1.1,0.15){$_{n}$}}
\pspolygon[fillstyle=solid,fillcolor=lightlightblue](0,0.3)(2.2,0.3)(2.2,0.8)(0,0.8)(0,0.3)
\pspolygon[fillstyle=solid,fillcolor=lightlightblue](0,2)(2.2,2)(2.2,1.5)(0,1.5)(0,2)
\pspolygon[fillstyle=solid,fillcolor=lightlightblue](0.4,0.8)(1.8,0.8)(1.8,1.5)(0.4,1.5)(0.4,0.8)
\psarc[linecolor=blue,linewidth=1.5pt]{-}(2.2,0.3){0.2}{90}{180}
\psarc[linecolor=blue,linewidth=1.5pt]{-}(0,0.3){0.2}{0}{90}
\psarc[linecolor=blue,linewidth=1.5pt]{-}(2.2,2){0.2}{180}{270}
\psarc[linecolor=blue,linewidth=1.5pt]{-}(0,2){0.2}{270}{0}
\psline[linecolor=blue,linewidth=1.5pt]{-}(1.1,0.3)(1.1,0.8)
\psline[linecolor=blue,linewidth=1.5pt]{-}(1.1,1.5)(1.1,2)
\rput(-0.15,0.5){$_k$}\rput(2.35,0.5){$_k$}
\rput(-0.15,1.8){$_k$}\rput(2.35,1.8){$_k$}
\rput(1.55,0.55){$_{n-2k}$}
\rput(1.55,1.75){$_{n-2k}$}
\rput(1.1,1.15){$\Omega^{2\ell}$}
\end{pspicture}
\ \ ,
\qquad
Y_{k,\ell} = \ 
\begin{pspicture}[shift=-2.05](-0.5,-1)(2.7,2.3)
\pspolygon[fillstyle=solid,fillcolor=lightlightblue](0,-1)(2.2,-1)(2.2,0.8)(0,0.8)(0,-1)
\rput(0,-0.3){\pspolygon[fillstyle=solid,fillcolor=pink](0.4,0)(2.2,0)(2.2,0.3)(0.4,0.3)(0.4,0)\rput(1.3,0.15){$_{n-1}$}}
\rput(0,2){\pspolygon[fillstyle=solid,fillcolor=pink](0,0)(2.2,0)(2.2,0.3)(0,0.3)(0,0)\rput(1.1,0.15){$_{n}$}}
\pspolygon[fillstyle=solid,fillcolor=lightlightblue](0,2)(2.2,2)(2.2,1.5)(0,1.5)(0,2)
\pspolygon[fillstyle=solid,fillcolor=lightlightblue](0.4,0.8)(1.8,0.8)(1.8,1.5)(0.4,1.5)(0.4,0.8)
\psarc[linecolor=blue,linewidth=1.5pt]{-}(2.2,0.3){0.2}{90}{180}
\psline[linecolor=blue,linewidth=1.5pt]{-}(0,0.5)(0.6,0.5)
\psarc[linecolor=blue,linewidth=1.5pt]{-}(0.6,0.3){0.2}{0}{90}
\psline[linecolor=blue,linewidth=1.5pt]{-}(0.8,0.3)(0.8,0)
\psarc[linecolor=blue,linewidth=1.5pt]{-}(2.2,2){0.2}{180}{270}
\psarc[linecolor=blue,linewidth=1.5pt]{-}(0,2){0.2}{270}{0}
\psline[linecolor=blue,linewidth=1.5pt]{-}(1.1,0)(1.1,0.8)
\psline[linecolor=blue,linewidth=1.5pt]{-}(1.1,1.5)(1.1,2)
\psline[linecolor=blue,linewidth=1.5pt]{-}(1.1,-1)(1.1,-0.3)
\psline[linecolor=blue,linewidth=1.5pt]{-}(2.0,0)(2.0,0.3)
\psarc[linecolor=blue,linewidth=1.5pt]{-}(0.4,0){0.2}{0}{180}
\psline[linecolor=blue,linewidth=1.5pt]{-}(0.2,0)(0.2,-0.3)
\psarc[linecolor=blue,linewidth=1.5pt]{-}(0,-0.3){0.2}{-90}{0}
\psarc[linecolor=blue,linewidth=1.5pt]{-}(2.2,-0.3){0.2}{180}{270}
\psarc[linecolor=blue,linewidth=1.5pt]{-}(0,-1.0){0.2}{0}{90}
\psarc[linecolor=blue,linewidth=1.5pt]{-}(2.2,-1.0){0.2}{90}{180}
\rput(-0.35,0.5){$_{k-1}$}\rput(2.55,0.5){$_{k-1}$}
\rput(-0.15,1.8){$_{k}$}\rput(2.35,1.8){$_{k}$}
\rput(1.55,0.55){$_{n-2k}$}
\rput(1.55,1.75){$_{n-2k}$}
\rput(1.48,-0.65){$_{n-2}$}
\rput(-0.15,-0.5){$_{1}$}\rput(2.35,-0.5){$_{1}$}
\rput(-0.15,-0.8){$_{1}$}\rput(2.35,-0.8){$_{1}$}
\rput(1.1,1.15){$\Omega^{2\ell-1}$}
\end{pspicture}
\ \ ,\qquad
X_{k,\ell} = \ 
\begin{pspicture}[shift=-1.15](-0.5,-0.1)(2.7,2.3)
\pspolygon[fillstyle=solid,fillcolor=lightlightblue](0,-0.1)(2.2,-0.1)(2.2,0.8)(0,0.8)(0,-0.1)
\pspolygon[fillstyle=solid,fillcolor=lightlightblue](0,2)(2.2,2)(2.2,1.5)(0,1.5)(0,2)
\pspolygon[fillstyle=solid,fillcolor=lightlightblue](0.4,0.8)(1.8,0.8)(1.8,1.5)(0.4,1.5)(0.4,0.8)
\rput(0,-0.1){\pspolygon[fillstyle=solid,fillcolor=pink](0.4,0)(1.8,0)(1.8,0.3)(0.4,0.3)(0.4,0)\rput(1.1,0.15){$_{n-2}$}}
\rput(0,2){\pspolygon[fillstyle=solid,fillcolor=pink](0,0)(2.2,0)(2.2,0.3)(0,0.3)(0,0)\rput(1.1,0.15){$_{n}$}}
\psarc[linecolor=blue,linewidth=1.5pt]{-}(1.8,0.2){0.2}{90}{180}
\psarc[linecolor=blue,linewidth=1.5pt]{-}(0.4,0.2){0.2}{0}{90}
\psarc[linecolor=blue,linewidth=1.5pt]{-}(2.2,-0.1){0.2}{90}{180}
\psarc[linecolor=blue,linewidth=1.5pt]{-}(0,-0.1){0.2}{0}{90}
\psarc[linecolor=blue,linewidth=1.5pt]{-}(2.2,2){0.2}{180}{270}
\psarc[linecolor=blue,linewidth=1.5pt]{-}(0,2){0.2}{270}{0}
\psline[linecolor=blue,linewidth=1.5pt]{-}(1.1,0.2)(1.1,0.8)
\psline[linecolor=blue,linewidth=1.5pt]{-}(1.1,1.5)(1.1,2)
\psline[linecolor=blue,linewidth=1.5pt]{-}(0,0.4)(0.4,0.4)
\psline[linecolor=blue,linewidth=1.5pt]{-}(1.8,0.4)(2.2,0.4)
\rput(-0.15,1.8){$_k$}\rput(2.35,1.8){$_k$}
\rput(-0.15,0.1){$_1$}\rput(2.35,0.1){$_1$}
\rput(-0.35,0.4){$_{k-1}$}\rput(2.55,0.4){$_{k-1}$}
\rput(1.55,0.6){$_{n-2k}$}
\rput(1.55,1.75){$_{n-2k}$}
\rput(1.1,1.15){$\Omega^{2\ell}$}
\end{pspicture} \ \ .
\ee
Here and below, we use the short-hand notation 
\be 
m_k = \frac{n-2k}{2}.
\ee 
The six families of projectors are defined as
\begin{subequations}
\label{eq:all.projectors}
\begin{alignat}{2}
\qptla_n(\beta,\alpha)&: \qquad Q_n = \sum_{k=0}^{\frac{n-2}2} \sum_{\ell = 0}^{m_k-1} \Gamma_{k,\ell} Z_{k,\ell} +  \Gamma_{n/2,0} Z_{n/2,0}\, ,
\label{eq:projectors.upa}
\\
\qptlb_n(\beta,\gamma)&: \qquad Q_n = \sum_{k=0}^{\frac{n-2}2} \sum_{\ell = 0}^{m_k-1} \Gamma_{k,\ell} Z_{k,\ell}\, , 
\label{eq:projectors.upb}
\\
\qptl_n(\beta,\gamma)&: \qquad Q_n = \sum_{k=0}^{\frac{n-1}2} \sum_{\ell = 0}^{2m_k-1} \Gamma_{k,\ell} Z_{k,\ell}\, ,\label{eq:projectors.up}
\\
\qatla_n(\beta,\alpha)&: \qquad Q_{n,r} = \sum_{k=0}^{\frac{n-2}2} \sum_{\ell = 0,\frac12,1, \dots}^{m_k-\frac12} \Gamma_{k,\ell} Z_{k,\ell} + \Gamma_{n/2,0} Z_{n/2,0}\,, \label{eq:projectors.uaa}
\\
\qatlb_n(\beta,\gamma)&: \qquad Q_{n,r} = \sum_{k=0}^{\frac{n-2}2} \sum_{\ell = 0,\frac12,1, \dots}^{m_k-\frac12} \Gamma_{k,\ell} Z_{k,\ell}\,,
\label{eq:projectors.uab}
\\
\qatl_n(\beta,\gamma)&: \qquad Q_{n,r} = \sum_{k=0}^{\frac{n-1}2} \sum_{\ell = 0,\frac12,1, \dots}^{m_k-\frac12} \Gamma_{k,\ell} Z_{k,\ell}\,.
\label{eq:projectors.ua}
\end{alignat}
\end{subequations}
The constants $\Gamma_{0,\ell}$ are given by
\be
\label{eq:initial.conditions}
\Gamma_{0,\ell} = \left\{\begin{array}{cc} 
\delta_{\ell,0} & \textrm{for }\qptla_n(\beta,\alpha), \qptlb_n(\beta,\gamma), \qptl_n(\beta,\gamma),
\\[0.2cm]
\tfrac1n\,\omega_{n,r}^{-2\ell} & \textrm{for }\qatla_n(\beta,\alpha), \qatlb_n(\beta,\gamma), \qatl_n(\beta,\gamma).
\end{array}\right.
\ee
We derive constraints satisfied by the constants $\Gamma_{k,\ell}$ using the relation $e_0 Q_n = 0$. We therefore compute the action of $e_0$ on $Z_{k,\ell}$, and write it as a linear combination of the diagrams $X_{k,\ell}$. To present the results in a uniform matter, we choose to view the projectors $Q_n$ of the uncoiled periodic algebras as elements of their corresponding uncoiled affine algebras. In this context, we have the equality 
\be
Z_{k,\ell+m_k} = \widehat\gamma\, Z_{k,\ell}, 
\qquad
\textrm{where}
\qquad
\widehat\gamma = \left\{\begin{array}{ll}
\gamma&\qptl_n(\beta,\gamma), \qptlb_n(\beta,\gamma),\qatl_n(\beta,\gamma), \qatlb_n(\beta,\gamma),\\[0.15cm]
1&\qptla_n(\beta,\alpha), \qatla_n(\beta,\alpha).
\end{array}\right.
\ee
 By choosing the convention $\Gamma_{k,\ell+m_k} = \wh\gamma^{-1}\,\Gamma_{k,\ell}$, we can, for instance, write
\be
\qptl_n(\beta,\gamma): \qquad Q_n = \sum_{k=0}^{\frac{n-1}2} \sum_{\ell = 0,\frac12,1,\dots}^{m_k-\frac12} \Gamma_{k,\ell} Z_{k,\ell}.
\ee
Likewise, the expressions \eqref{eq:projectors.uaa}, \eqref{eq:projectors.uab} and \eqref{eq:projectors.ua} for the projectors $Q_{n,r}$ can with this convention be rewritten with $\ell$ running over integers and extending up to $2m_k-1$. In the calculations below, we alternate between these two presentations when convenient.\medskip

Using the relations \eqref{Pm.props} satisfied by the projectors $P_n$, we find
\begin{subequations}
\begin{alignat}{3}
e_0 Z_{0,\ell} &= \tfrac{[n-2\ell-1]}{[n-1]} X_{1,\ell} + \tfrac{[2\ell]}{[n]} Y_{1,\ell},
\qquad &&0 \le \ell < m_0,
\\[0.15cm]
e_0 Z_{k,0} &= -\tfrac{[n]}{[n-1]} X_{k,0} + \tfrac{[n-k]}{[n]}Y_{k,1} + \tfrac{[k]}{[n]} Y_{k,0},
\\[0.15cm]
e_0 Z_{k,\ell} &= -\tfrac{[n]}{[n-1]} X_{k,\ell} + \tfrac{[n-k]}{[n]}Y_{k,\ell+1} + \tfrac{[k+2 \ell]}{[n]}Y_{k+1,\ell} + \tfrac{[k]}{[n]} Y_{k,\ell}, 
\qquad&&\tfrac12\le\ell<m_k,
\\[0.15cm]
e_0 Z_{\frac n2,0} &= -\tfrac{[n]}{[n-1]} X_{\frac n2,0}+ \alpha \tfrac{[\frac n2]}{[n]}Y_{\frac n2,\frac12}, &&(n \textrm{ even})
\\[0.15cm]
Y_{k,0} &= \tfrac{[n-k]}{[n-1]}X_{k,-1} + \tfrac{[k+1]}{[n-1]}X_{k+1,0} + \tfrac{[k]}{[n-1]}X_{k,0},
\\[0.15cm]
Y_{k,\frac12} &= \tfrac{[n-k]}{[n-1]}X_{k,-\frac12} + \tfrac{[k]}{[n-1]}X_{k,\frac12}, 
\\[0.15cm]
Y_{k,\ell} &= \tfrac{[n-k]}{[n-1]}X_{k,\ell-1} + \tfrac{[n-k-2\ell+1]}{[n-1]}X_{k+1,\ell-1} + \tfrac{[k]}{[n-1]}X_{k,\ell}, 
\qquad &&1\le\ell\le m_k,
\\[0.15cm]
Y_{k,m_k+\frac12} &= \tfrac{[n-k]}{[n-1]}X_{k,m_k-\frac12} + \tfrac{[k]}{[n-1]}X_{k,m_k+\frac12},
\\[0.15cm]
Y_{\frac n2,\frac12} &= \alpha\tfrac{[\frac n2]}{[n-1]} X_{\frac n2,0}, &&(n \textrm{ even})
\end{alignat}
\end{subequations}
where $1 \le k \le \lfloor\tfrac{n-1}2\rfloor$ and $\ell \in \frac12 \mathbb Z$. Putting these relations together, we find
\begin{alignat}{2}
e_0 Z_{k,\ell} &= f^1_k X_{k,\ell} + f^2_k (X_{k,\ell-1} + X_{k,\ell+1}) + \big[(1-\delta_{\ell,m_k-\frac12})f^{3\textrm{a}}_{k,\ell}+f^{3\textrm{b}}_{k,\ell}\big] X_{k+1,\ell} \nonumber
\\&+ (1-\delta_{\ell,0})\big[f^{4\textrm{a}}_{k,\ell}+(1-\delta_{\ell,\frac12})f^{4\textrm{b}}_{k,\ell}\big] X_{k+1,\ell-1} + (1-\delta_{\ell,0}-\delta_{\ell,\frac12}-\delta_{\ell,m_k-\frac12})f^5_{k,\ell} X_{k+2,\ell-1},
\label{eq:e0Z.all}
\end{alignat}
where 
\be
\label{eq:fikL}
\begin{array}{llll}
f^1_k = -\tfrac{[n-k][k][2n]}{[n-1][n]^2},\quad
&
f^{3\textrm{a}}_{k,\ell} = \tfrac{[n-k][n-k-2\ell-1]}{[n-1][n]},\quad
&
f^{4\textrm{a}}_{k,\ell} = \tfrac{[n-k-1][k+2\ell]}{[n-1][n]},
\\[0.25cm]
f^2_k = \tfrac{[n-k][k]}{[n-1][n]},
&
f^{3\textrm{b}}_{k,\ell} = \tfrac{[k+2\ell][k+1]}{[n-1][n]},\quad
&
f^{4\textrm{b}}_{k,\ell} = \tfrac{[n-k-2\ell+1][k]}{[n-1][n]},\quad
\end{array}
\quad
f^5_{k,\ell} = \tfrac{[n-k-2\ell][k+2\ell]}{[n-1][n]}.
\ee
The relation \eqref{eq:e0Z.all} is valid for $1 \le k \le \lfloor\tfrac{n-1}2\rfloor$ and $\ell \in \frac12 \mathbb Z$, for the algebras $\qptl_n(\beta,\gamma)$, $\qptlb_n(\beta,\gamma)$, $\qatl_n(\beta,\gamma)$ and $\qatlb_n(\beta,\gamma)$. For the largest values of $k$, they hold using the convention $X_{k,\ell} = 0$ for $k>\frac n2$. 
For $\qptla_n(\beta,\alpha)$ and $\qatla_n(\beta,\alpha)$, this relation is valid for $1 \le k \le \tfrac {n-4}2$. The extra relations for $k = \tfrac {n-2}2$ and $k=\tfrac {n}2$ are
\begin{subequations}
\begin{alignat}{2}
e_0 Z_{\frac {n-2}2,0} &= (f^1_{\frac {n-2}2} + 2 f^2_{\frac {n-2}2}) X_{\frac {n-2}2,0} + \frac{[2][\frac n2]^2}{[n][n-1]}X_{\frac n2,0}\,,
\\[0.2cm]
e_0 Z_{\frac {n-2}2,\frac12} &= (f^1_{\frac {n-2}2} + 2 f^2_{\frac {n-2}2}) X_{\frac {n-2}2,\frac12} + \frac{\alpha[\frac n2]^2}{[n][n-1]}X_{\frac n2,0}\,,
\\[0.2cm]
e_0 Z_{\frac n2,0} &= \frac{\alpha^2[\frac n2]^2-[n]^2}{[n][n-1]}\, X_{\frac n2,0}\,.
\end{alignat}
\end{subequations}

To compute $e_0$ on the projectors $Q_n$ or $Q_{n,r}$, we apply \eqref{eq:e0Z.all} to each term in the sums \eqref{eq:all.projectors}. We then shift the indices $k$ and $\ell$ of the sums in such a way that each term becomes proportional to $X_{k,\ell}$. These objects can all be shown to be linearly independent, so setting $e_0Q_n=0$ or $e_0Q_{n,r}=0$ is done by setting the coefficients of each $X_{k,\ell}$ to zero. This results in the following proposition.

\begin{Proposition}\label{Prop.const.rel}
For $k=1,2,\dots, \lfloor\tfrac{n-1}{2}\rfloor
$ and $\ell=0,\tfrac{1}{2},1,\dots, m_k-\tfrac{1}{2}$, the constants $\Gamma_{k,\ell}$ satisfy the linear conditions
 \begin{alignat}{3}
&0=f^1_k \Gamma_{k,\ell} + f^2_k(\Gamma_{k,\ell-1}+\Gamma_{k,\ell+1}) + f^3_{k-1,\ell} \Gamma_{k-1,\ell} + f^4_{k-1,\ell+1} \big(1+\tfrac{\gamma^2}{[2]}\delta_{k,\frac{n-1}2}\big)
\Gamma_{k-1,\ell+1} + f^5_{k-2,\ell+1} \Gamma_{k-2,\ell + 1}\nonumber\\[0.2cm]
\label{eq:Gamma.constraint}
&+\delta_{\ell,0}(f^3_{k-1,m_k} \Gamma_{k-1,-1}+f^5_{k-2,m_k+1} \Gamma_{k-2,-1}) + \delta_{\ell,\frac12} f^{3\textrm{b}}_{k-1,m_k+\frac12}\Gamma_{k-1,-\frac12} + \delta_{\ell,m_k-\frac12} f^{4\textrm{a}}_{k-1,\frac12}\Gamma_{k-1,m_k+\frac32},
\end{alignat}
where $f^3_{k,\ell} = f^{3\textrm{a}}_{k,\ell} + f^{3\textrm{b}}_{k,\ell}$ and $f^4_{k,\ell} = f^{4\textrm{a}}_{k,\ell} + f^{4\textrm{b}}_{k,\ell}$. For $\ell \in \{ 0, \frac12, m_k-1,m_k-\frac12\}$, some of the constants $\Gamma_{k,\ell}$ that arise in \eqref{eq:Gamma.constraint} have $\ell\notin[0,m_k-\frac12]$, and in this case, one uses the relation $\Gamma_{k,\ell + m_k} = \widehat\gamma^{-1}\, \Gamma_{k,\ell}$. Moreover, for $k = 1$, the above relation holds with the initial condition \eqref{eq:initial.conditions} and $\Gamma_{-1,\ell} = 0$.\medskip

Finally, for the algebras $\qptla_n(\beta,\alpha)$ and $\qatla_n(\beta,\alpha)$, there is an extra linear relation between the constants:
\be
\label{eq:extra.Gamma.constraint}
\frac{\alpha^2[\frac n2]^2-[n]^2}{[n][n-1]}\, \Gamma_{\frac n2,0} + \frac{[\frac n2]^2}{[n][n-1]} \big([2]\,\Gamma_{\frac{n-2}2,0} + \alpha\, \Gamma_{\frac{n-2}2,\frac12}\big) + f^5_{\frac{n-4}2,1} \Gamma_{\frac{n-4}2,1} = 0.
\ee
\end{Proposition}

The solutions to this recursive system of equations for the different uncoiled algebras are given in \cref{Prop.const.uptla.uptlb,Prop.const.uptl,Prop.const.qatl.qatla.qatlb}. Their proofs are given in \cref{Sec:proof.closed.form}. The following two sections present the path we took to obtain these formulas. 

\subsection{Fourier series and convolution sums}\label{Sec:Fourier}

The set of equations \eqref{eq:Gamma.constraint} is a triangular system with respect to the variable $k$ which can be solved explicitly using Fourier series. The argument works differently for $n$ even and odd. For $n$ even, $m_k$ is an integer and thus the set of linear relations with the periodicity property $\Gamma_{k,\ell + m_k} = \widehat\gamma^{-1}\, \Gamma_{k,\ell}$ really consists of two disjoint sets of relations. We thus define two sets of Fourier series:
\begin{subequations}
\begin{alignat}{3}
\widehat \Gamma^{\textrm{a}}_{k,s} &= \sum_{\ell = 0,1,\dots}^{m_k - 1} y_{k,s}^\ell \Gamma_{k,\ell},
\qquad
&&\Gamma_{k,\ell} =\, \frac1{m_k} \sum_{s = 0,1,\dots}^{m_k-1} y_{k,s}^{-\ell} \widehat \Gamma_{k,s}^{\textrm{a}},
\qquad &&\ell \in \mathbb Z,
\\
\widehat \Gamma^{\textrm{b}}_{k,s} &= \sum_{\ell = \frac12,\frac32,\dots}^{m_k - \frac12} y_{k,s}^{\ell-\frac12}\, \Gamma_{k,\ell},
\qquad
&&\Gamma_{k,\ell} =\, \frac1{m_k} \sum_{s = 0,1,\dots
}^{m_k-1} y_{k,s}^{-\ell+\frac12}\, \widehat \Gamma_{k,s}^{\textrm{b}},
\qquad &&\ell \in \mathbb Z+\tfrac12,
\end{alignat}
\end{subequations}
where $y_{k,s} = \widehat\gamma^{1/m_k} \eE^{2\pi \iI s/m_k}$, and 
$s\in \{0,1,\dots, m_k-1\}$ in the leftmost equations. Applying the Fourier series to \eqref{eq:Gamma.constraint}, we find
\begin{subequations}
\begin{alignat}{2}
- f^2_k\, \widehat\Gamma_{k,s}^{\textrm{a}} (y_{k,s} + y_{k,s}^{-1}-q^n-q^{-n}) = &\sum_{\ell = 0,1,\dots}^{m_k-1} y_{k,s}^\ell \Big( f^3_{k-1,\ell} \Gamma_{k-1,\ell} + f^4_{k-1,\ell+1} \Gamma_{k-1,\ell+1} + f^5_{k-2,\ell+1} \Gamma_{k-2,\ell+1}\Big) \nonumber\\
& + f^3_{k-1,m_k} \Gamma_{k-1,-1}+ f^5_{k-2,m_k+1} \Gamma_{k-2,-1} 
\nonumber\\
 &\hspace{-4.7cm}=\sum_{\ell = 0,1,\dots}^{m_k} y_{k,s}^\ell  f^3_{k-1,\ell} \Gamma_{k-1,\ell} + \sum_{\ell = 0,1,\dots}^{m_k-1} y_{k,s}^\ell f^4_{k-1,\ell+1} \Gamma_{k-1,\ell+1} + \sum_{\ell = 0,1,\dots}^{m_k} y_{k,s}^\ell f^5_{k-2,\ell+1} \Gamma_{k-2,\ell+1}\,,
 \label{eq:Gamma.hat.a}
 \\
- f^2_k\, \widehat\Gamma_{k,s}^{\textrm{b}} (y_{k,s} + y_{k,s}^{-1}-q^n-q^{-n}) = &\sum_{\ell = \frac12,\frac32,\dots}^{m_k-\frac12} y_{k,s}^{\ell-\frac12}\, \Big( f^3_{k-1,\ell} \Gamma_{k-1,\ell} + f^4_{k-1,\ell+1} \Gamma_{k-1,\ell+1} + f^5_{k-2,\ell+1} \Gamma_{k-2,\ell+1}\Big) \nonumber\\
& + f^{3\textrm{b}}_{k-1,m_k+\frac12} \Gamma_{k-1,-\frac12}+ y_{k,s}^{m_k-1} f^{4\textrm{a}}_{k-1,\frac12} \Gamma_{k-1,m_k+\frac32}\,. 
\label{eq:Gamma.hat.b}
\end{alignat}
\end{subequations}
For $\qptla_n(\beta,\alpha)$, \eqref{eq:Gamma.hat.b} is automatically satisfied because in this case, each $\Gamma_{k,\ell}$ vanishes for $\ell \in \mathbb Z + \frac12$. Applying the inverse transform, we find
\begin{subequations}
\label{eq:Gamma.kL.Fourier}
\begin{alignat}{2}
\Gamma_{k,\ell} &= -\frac{1}{f_k^2}\bigg[
\sum_{\ell'=0,1,\dots}^{m_k} J_k(\ell'-\ell) f^3_{k-1,\ell'}\Gamma_{k-1,\ell'}
+\sum_{\ell'=0,1,\dots}^{m_k-1} J_k(\ell'-\ell) f^4_{k-1,\ell'+1}\Gamma_{k-1,\ell'+1}
\label{eq:Gamma.a}
\\&\hspace{7cm}+\sum_{\ell'=0,1,\dots}^{m_k} J_k(\ell'-\ell) f^5_{k-2,\ell'+1}\Gamma_{k-2,\ell'+1} \bigg], \qquad \ell \in \mathbb Z,
\nonumber
\end{alignat}
\begin{alignat}{2}
\Gamma_{k,\ell} &= -\frac{1}{f_k^2}\bigg[
\sum_{\ell'=\frac12,\frac32,\dots}^{m_k-\frac12} J_k(\ell'-\ell) \big(f^3_{k-1,\ell'}\Gamma_{k-1,\ell'}
+ f^4_{k-1,\ell'+1}\Gamma_{k-1,\ell'+1} + f^5_{k-2,\ell'+1}\Gamma_{k-2,\ell'+1}\big)
\label{eq:Gamma.b} 
\\ &\hspace{2cm}+ J_k(\tfrac12-\ell) f^{3\textrm{b}}_{k-1,m_k+\frac12} \Gamma_{k-1,-\tfrac12} + J_k(m_k-\tfrac12-\ell) f^{4\textrm{a}}_{k-1,\frac12} \Gamma_{k-1,m_k+\frac32} 
\bigg], \qquad \ell \in \mathbb Z+\tfrac12,
\nonumber
\end{alignat}
\end{subequations}
where
\be
J_k(\ell) = \frac1{m_k} \sum_{s=0}^{m_k-1} \frac{y_{k,s}^\ell}{y_{k,s}+y_{k,s}^{-1}-q^n-q^{-n}}.
\ee
The result is thus expressed in terms of a convolution sum, with $J_k(\ell)$ playing the role of the kernel.\medskip

For $n$ odd, $m_k$ is a half-integer and thus the relations \eqref{eq:Gamma.constraint} do not divide into two disjoint subsets. In this case, we define a single Fourier series:
\be
\label{eq:Fourier.n.odd}
\widehat \Gamma_{k,s} = \sum_{\ell = 0}^{2m_k - 1} z_{k,s}^\ell \Gamma_{k,\ell},
\qquad
\Gamma_{k,\ell} = \frac1{2m_k} \sum_{s = 0}^{2m_k-1} z_{k,s}^{-\ell} \widehat \Gamma_{k,s},
\qquad
z_{k,s} =\gamma^{1/m_k} \eE^{\pi \iI s/m_k}.
\ee
We then find
\begin{alignat}{2}
-&f^2_k\,\widehat\Gamma_{k,s} (z_{k,s} + z_{k,s}^{-1}-q^n-q^{-n}) = \sum_{\ell = 0}^{m_k+\frac12} z_{k,s}^\ell f^3_{k-1,\ell} \Gamma_{k-1,\ell} + \sum_{\ell = 0}^{m_k-\frac12} z_{k,s}^\ell f^4_{k-1,\ell+1} \Gamma_{k-1,\ell+1} 
\nonumber\\&
+\sum_{\ell = 0}^{m_k-\frac12} z_{k,s}^\ell  f^5_{k-2,\ell+1} \Gamma_{k-2,\ell+1}
+ \sum_{\ell = m_k-\frac12}^{2m_k} z_{k,s}^\ell f^3_{k-1,\ell-m_k} \Gamma_{k-1,\ell+1} + \sum_{\ell = m_k+\frac12}^{2m_k-1} z_{k,s}^\ell f^4_{k-1,\ell+1-m_k} \Gamma_{k-1,\ell+2} 
\nonumber\\&
+ \sum_{\ell = m_k+\frac12}^{2m_k} z_{k,s}^\ell f^5_{k-2,\ell+1-m_k} \Gamma_{k-2,\ell+3}
+ z_{k,s}^{m_k-\frac12}f^{4\textrm{a}}_{k-1,\frac12}\Gamma_{k-1,m_k+\frac32} + z_{k,s}^{m_k+\frac12}f^{3\textrm{b}}_{k-1,m_k+\frac12} \Gamma_{k-1,m_k+\frac12}.
\end{alignat}
Applying the inverse transform yields
\begingroup
\allowdisplaybreaks
\begin{alignat}{2}
\Gamma_{k,\ell} &= -\frac1{f^2_k}
\bigg[\tilde J_k(m_k-\tfrac12 - \ell)f^{4\textrm{a}}_{k-1,\frac12}\Gamma_{k-1,m_k+\frac32} + \tilde J_k(m_k+\tfrac12 - \ell)f^{3\textrm{b}}_{k-1,m_k+\frac12}\Gamma_{k-1,m_k+\frac12} 
\nonumber\\ &
+ \sum_{\ell'=0}^{m_k-\frac12} \tilde J_k(\ell'-\ell) f^3_{k-1,\ell'} \Gamma_{k-1,\ell'}
+ \sum_{\ell'=0}^{m_k-\frac12} \tilde J_k(\ell'-\ell) f^4_{k-1,\ell'+1} \Gamma_{k-1,\ell'+1}
+ \sum_{\ell'=0}^{m_{k}-\frac12} \tilde J_k(\ell'-\ell) f^5_{k-2,\ell'+1} \Gamma_{k-2,\ell'+1} 
\nonumber\\ &
+ \sum_{\ell'=m_k+\frac12}^{2m_k} \tilde J_k(\ell'-\ell) f^3_{k-1,\ell'-m_k} \Gamma_{k-1,\ell'+1} 
+ \sum_{\ell'=m_k+\frac12}^{2m_k-1} \tilde J_k(\ell'-\ell) f^4_{k-1,\ell'-m_k+1}\Gamma_{k-1,\ell+2}
\nonumber\\ &
+ \sum_{\ell'=m_k+\frac12}^{2m_k} \tilde J_k(\ell'-\ell) f^5_{k-2,\ell'-m_k+1}\Gamma_{k-2,\ell+3}
\bigg],
\label{eq:Gamma.odd}
\end{alignat}
\endgroup
where
\be
\label{eq:J.tilde}
\tilde J_k(\ell) = \frac1{2m_k} \sum_{s=0}^{2m_k-1} \frac{z_{k,s}^\ell}{z_{k,s}+z_{k,s}^{-1}-q^n-q^{-n}}
=J_k(\ell) \bigg|_{\substack{\hspace{-0.5cm}\widehat\gamma\, \to\, \gamma^2 \\[0.05cm] m_k\,\to\,2 m_k}}.
\ee

\subsection{Evaluation of the kernel}\label{sec:kernel}

In this section, we prove the following lemma for the evaluation of the kernel function.
\begin{Lemma} The function $J_k(\ell)$ evaluates to
\be
\label{eq:Jk.explicit}
J_k(\ell) = -\frac1{q^n-q^{-n}} \sum_{\sigma = \pm 1} \frac{\sigma q^{\sigma n\ell}}{\widehat\gamma^{-1} q^{\sigma nm_k} -1}, \qquad 0 \le \ell \le m_k.
\ee
\end{Lemma}
\proof Let us define the function
\be
g(y) = \frac1{2 \pi \iI} \frac{y^{\ell}}{(y-q^n)(y-q^{-n})}\frac{1}{\widehat\gamma^{-1}y^{m_k}-1}.
\ee
This function has poles at $y = q^{\pm n}$ and at each $y = y_{k,s}$. For $\ell \ge 0$, it has no pole at $y = 0$. This function is also such that
\be
2 \pi \iI \,\textrm{Res}(g(y),y = y_{k,s}) = \frac1{m_k} \frac{y_{k,s}^\ell}{y_{k,s}+y_{k,s}^{-1}-q^n-q^{-n}},
\ee
which is precisely the summand for $J_k(\ell)$. As a result, we have
\be
J_k(\ell) = \oint_C \dd y \, g(y)
\ee
where $C$ is a contour in the complex plane that encircles all the poles $y_{k,s}$ in the counter-clockwise direction, but does not encircle the poles at $y = q^{\pm n}$. Let us now define the modified contour $C'$ that encircles all the poles counter-clockwise, and the corresponding function
\be
J'_k(\ell) = \oint_{C'} \dd y \, g(y).
\ee
Clearly,
\be
\label{eq:J.prime}
J'_k(\ell) = J_k(\ell) + 2 \pi \iI \, \textrm{Res}(g(y),y = q^n)  + 2 \pi \iI \, \textrm{Res}(g(y),y = q^{-n}). 
\ee
Moreover, $J'_k(\ell)$ can be computed in a second way, by changing variables to $y = z^{-1}$:
\be
J'_k(\ell) = \oint_{0} \frac{\dd z}{z^2} \, g(z^{-1}) = \oint_{0} \dd z \frac1{2 \pi \iI} \frac{z^{m_k-\ell}}{(1-zq^n)(1-zq^{-n})}\frac{1}{\widehat\gamma^{-1}-z^{m_k}}.
\ee
Here, the closed contour encircles the origin in the counter-clockwise direction, and does not encircle any of the poles at $z = q^{\pm n}$ and $z = y_{k,s}^{-1}$. For $\ell \le m_k$, the integrand has no pole at $z=0$, thus implying that $J'_k(\ell) = 0$. The residues at $y = q^{\pm n}$ in \eqref{eq:J.prime} are easily computed, giving us the expression \eqref{eq:Jk.explicit} for $J_k(\ell)$ and ending the proof.
\eproof

The expressions \eqref{eq:Gamma.kL.Fourier} for $\Gamma_{k,\ell}$ also require that we know $J_k(\ell)$ for $\ell<0$. For $\ell \notin \{0, 1, \dots, m_k\}$, $J_k(\ell)$ can be computed using the relation $y_{k,s}^{m_k} =\widehat\gamma$. With this idea, we find
\be
J_k(-\ell) = -\frac{1}{q^n-q^{-n}}\sum_{\sigma = \pm 1}\frac{\sigma q^{\sigma n\ell}}{\widehat\gamma\, q^{\sigma nm_k} -1}, \qquad 0 \le \ell \le m_k.
\ee
Finally, using the last equality in \eqref{eq:J.tilde}, we also find
\be
\tilde J_k(\ell) = -\frac1{q^n-q^{-n}} \sum_{\sigma = \pm 1} \frac{\sigma q^{\sigma n\ell}}{\gamma^{-2} q^{2\sigma nm_k} -1}, 
\qquad
\tilde J_k(-\ell) = -\frac{1}{q^n-q^{-n}}\sum_{\sigma = \pm 1}\frac{\sigma q^{\sigma n\ell}}{\gamma^2\, q^{2\sigma nm_k} -1},
\ee
valid for $\ell=0,1,\dots,2m_k$.

\subsection{Values of the constants}\label{Sec:Conj.Values}

In this section, we write down the values of the constants for the six uncoiled algebras. The relations \eqref{eq:Gamma.kL.Fourier} and \eqref{eq:Gamma.odd}, along with the kernel functions computed in \cref{sec:kernel}, uniquely fix $\Gamma_{1,\ell}$ from the known values of $\Gamma_{-1,\ell}$ and $\Gamma_{0,\ell}$. This then allows us to solve for $\Gamma_{2,\ell}$, and so on. The derivation at every step is tedious yet straightforward, as it only involves evaluating a number of geometric series. Here, we first give the expressions for $\Gamma_{k,\ell}$ for the smallest values of $k$. This allowed us to conjecture the general closed-forms for arbitrary values of $k$, which we then were able to prove in \cref{Sec:proof.closed.form} by showing that they satisfy the recurrence relation of \cref{Prop.const.rel}. 
\medskip

Before proceeding, let us recall that the $q$-factorials and $q$-binomials are given by
\be
[\kappa]! = \prod_{j=1}^{\kappa}[j],
\qquad
\left[\begin{matrix} \kappa\\ \tau\end{matrix}\right] = \frac{[\kappa]!}{[\tau]![\kappa-\tau]!}.
\ee

\paragraph{Constants for $\boldsymbol{\qptla_n(\beta, \alpha)}$ and $\boldsymbol{\qptlb_n(\beta,\gamma)}$}
Using the relations \eqref{eq:Gamma.a}, we find after some simplifications
\begin{subequations}
\begin{alignat}{2}
&\hspace{-1cm}\Gamma_{1,\ell} = \frac1{q-q^{-1}} \sum_{\sigma = \pm 1}\frac{\sigma q^{\sigma n\ell}}{\widehat\gamma\, q^{\sigma n m_1}-1},\label{eq:Gamma.1L}
\\[0.15cm]
&\hspace{-1cm}\Gamma_{2,\ell} = \frac1{(q-q^{-1})^3 [2]} \sum_{\sigma = \pm 1}\sigma q^{\sigma n\ell}\bigg(\frac{1}{\widehat\gamma\, q^{\sigma n m_2}-1}-\frac{[2m_2-2\ell]+q^{\sigma n}[2\ell]}{[n-2](\widehat\gamma\, q^{\sigma n m_1}-1)}\bigg),
\\[0.15cm]
&\hspace{-1cm}\Gamma_{3,\ell} = \frac1{(q-q^{-1})^5[2]^2[3]} \sum_{\sigma = \pm 1}\sigma q^{\sigma n\ell}\bigg(\frac{1}{\widehat\gamma\, q^{\sigma n m_3}-1} - [2]\frac{[2m_3-2\ell]+q^{\sigma n}[2 \ell]}{[n-3](\widehat\gamma\, q^{\sigma n m_2}-1)}
\\[0.15cm] &\hspace{3cm}
+\frac{[2m_3-2\ell][2m_3-2\ell+1]+q^{\sigma n}[2][2\ell][2m_3-2\ell]+q^{2\sigma n}[2\ell][2\ell+1]}{[n-3][n-2](\widehat\gamma\, q^{\sigma n m_1}-1)}\bigg).
\nonumber
\end{alignat}
\end{subequations}

\begin{Proposition}\label{Prop.const.uptla.uptlb}
For $\qptla_n(\beta,\alpha)$ and $\qptlb_n(\beta,\gamma)$, the general formula for the constants is
\begin{alignat}{2}
\Gamma_{k,\ell} =  \frac1{(q-q^{-1})^{2k-1}[k]![k-1]!} \sum_{\sigma = \pm 1}\sum_{\kappa = 1}^{k}\sum_{\tau = 0}^{k-\kappa} &\frac{(-1)^{k+\kappa} \sigma q^{\sigma n(\ell+\tau)}}{\widehat\gamma\, q^{\sigma n m_{\kappa}}-1} 
\frac{\left[\begin{matrix} k-1\\ \kappa-1\end{matrix}\right]}
{\left[\begin{matrix} n-\kappa-1\\ n-k-1\end{matrix}\right]}
\nonumber\\[0.15cm]&\times
\left[\begin{matrix} 2 m_k -2\ell + k - \kappa - \tau - 1\\ k-\kappa-\tau\end{matrix}\right]
\left[\begin{matrix} 2\ell + \tau - 1\\ \tau\end{matrix}\right],
\label{eq:Gamma.final}
\end{alignat}
valid for $1 \le k \le \frac{n-2}2$ and $\ell=0,1,\dots,m_k-1$.
For $\qptla_n(\beta,\alpha)$, the extra coefficient is given by
\be
\label{eq:Gamma.n/2.0}
\Gamma_{n/2,0} = -\frac{1}{(q-q^{-1})^{n-2}[\frac{n-2}2]!^2}\frac{[n][\frac n2]}{\alpha^2[\frac n2]^2-[n]^2}.
\ee
\end{Proposition}
Here and below, we use the convention
\be
{\left[\begin{matrix} 2\ell+\tau-1\\ \tau\end{matrix}\right]} \xrightarrow{\ell\to0} \delta_{\tau,0}
\ee
which is needed in \eqref{eq:Gamma.final} for $\ell = 0$.\medskip

Finally, let us note that the constants satisfy the reflection symmetry
\be
\label{eq:crossing.sym}
\Gamma_{k,\ell}\big|_{\widehat\gamma \to \widehat\gamma^{-1}} = \widehat\gamma\, \Gamma_{k,m_k-\ell}\, .
\ee

\paragraph{Constants for $\boldsymbol{\qptl_n(\beta,\gamma)}$}

In this case, the constants for $k=1$ and $k=2$ with $\ell = 0, \frac12, \dots, m_k-\frac12$ are
\begin{subequations}
\begin{alignat}{2}
&\Gamma_{1,\ell} = \frac1{q-q^{-1}} \sum_{\sigma = \pm 1}\frac{\sigma q^{\sigma n\ell}\mu_{\ell,1,\sigma}}{\gamma^2 q^{2\sigma n m_1}-1},
\\[0.15cm]
&\Gamma_{2,\ell} =\frac1{(q-q^{-1})^3 [2]} \sum_{\sigma = \pm 1}\sigma q^{\sigma n\ell}\bigg(\frac{\mu_{\ell,2,\sigma}}{\gamma^2\, q^{2\sigma n m_2}-1}-\mu_{\ell,1,\sigma}\frac{[2m_2-2\ell]+q^{\sigma n}[2\ell]}{[n-2](\gamma^2\, q^{2\sigma n m_1}-1)}\bigg) 
\end{alignat}
\end{subequations}
where 
\be
\mu_{\ell,\kappa,\sigma} = 
\left\{\begin{array}{ll}
1 & \ell \in \mathbb Z\,,\\[0.1cm]
\gamma \, q^{\sigma n m_\kappa} & \ell \in \mathbb Z + \frac12\,.
\end{array}\right.
\ee
The equivalent formulas with $\ell \ge m_k+\frac12$ are obtained using the relation $\Gamma_{k,\ell+m_k} = \gamma^{-1}\,\Gamma_{k,\ell}$.

\begin{Proposition}\label{Prop.const.uptl}
For $\qptl_n(\beta,\gamma)$, the general formula for the constants is
\begin{alignat}{2}
\Gamma_{k,\ell} 
= 
\frac1{(q-q^{-1})^{2k-1}[k]![k-1]!} 
\sum_{\sigma = \pm 1}\sum_{\kappa = 1}^{k}\sum_{\tau = 0}^{k-\kappa}& \frac{(-1)^\kappa \sigma q^{\sigma n(\ell+\tau)}\mu_{\ell,\kappa,\sigma}}{\gamma^2\, q^{2\sigma n m_{k-\kappa}}-1} 
\frac{\left[\begin{matrix} k-1\\ \kappa-1\end{matrix}\right]}
{\left[\begin{matrix} n-\kappa-1\\ n-k-1\end{matrix}\right]}
\nonumber\\&
\times\left[\begin{matrix} 2 m_k -2\ell + k - \kappa - \tau - 1\\ k-\kappa-\tau\end{matrix}\right]
\left[\begin{matrix} 2\ell + \tau - 1\\ \tau\end{matrix}\right]
\label{eq:Gamma.final.odd}
\end{alignat}
for $1 \le k \le \frac{n-1}2$ and $\ell=0,\frac12,\dots,m_k-\frac12$.
\end{Proposition}
This formula satisfies the reflection symmetry
\be
\Gamma_{k,\ell}\big|_{\gamma \to \gamma^{-1}} = \gamma^2\, \Gamma_{k,2m_k-\ell}\,.
\ee

\paragraph{Constants for $\boldsymbol{\qatla_n(\beta, \alpha)}$, $\boldsymbol{\qatlb_n(\beta,\gamma)}$ and $\boldsymbol{\qatl_n(\beta,\gamma)}$.}

Even if the argument with Fourier series is different for $n$ even and odd, it turns out that the solution takes the same form for the two parities of $n$. Indeed, for $k=1$ and $k=2$, we find after some simplifications
\begin{subequations}
\begin{alignat}{2}
\Gamma_{1,\ell} &= \frac{\omega^{-2\ell}}{n(q-q^{-1})}\sum_{\sigma = \pm 1} \frac{\sigma q^{2 \sigma \ell}}{\omega^2 q^{2 \sigma m_1}-1},\\
\Gamma_{2,\ell} &= \frac{\omega^{-2\ell}}{n(q-q^{-1})^3[2]}\sum_{\sigma = \pm 1} \sigma \bigg(\frac{q^{4 \sigma \ell}}{\omega^2 q^{2 \sigma m_2}-1}-q^{2 \sigma \ell}\frac{[2m_2-2\ell] + q^{\sigma n} [2\ell]}{[n-2](\omega^2 q^{2\sigma m_1}-1)}\bigg).
\end{alignat}
\end{subequations}

\begin{Proposition}\label{Prop.const.qatl.qatla.qatlb}
For $\qatla_n(\beta, \alpha)$, $\qatlb_n(\beta,\gamma)$ and $\qatl_n(\beta,\gamma)$, the general formula for the constants is
\begin{alignat}{2}
\Gamma_{k,\ell} 
= 
\frac{\omega^{-2\ell}}{n(q-q^{-1})^{2k-1}[k]![k-1]!} 
\sum_{\sigma = \pm 1}\sum_{\kappa = 1}^{k}\sum_{\tau = 0}^{k-\kappa} &\frac{(-1)^{k+\kappa} \sigma q^{\sigma (2\ell \kappa+n \tau)}}{\omega^2 q^{2\sigma m_{\kappa}}-1} 
\frac{\left[\begin{matrix} k-1\\ \kappa-1\end{matrix}\right]}
{\left[\begin{matrix} n-\kappa-1\\ n-k-1\end{matrix}\right]}
\nonumber\\[0.15cm]&\times
\left[\begin{matrix} 2 m_k -2\ell + k - \kappa - \tau - 1\\ k-\kappa-\tau\end{matrix}\right]
\left[\begin{matrix} 2\ell + \tau - 1\\ \tau\end{matrix}\right]
\label{eq:Gamma.final.affine}
\end{alignat}
for $1 \le k \le \lfloor\frac{n-1}2\rfloor$ and $\ell=0,\frac12,\dots,m_k-\frac12$.
For $\qatla_n(\beta,\alpha)$,
the extra coefficient is given by
\be
\label{eq:Gamma.n/2.0.affine}
\Gamma_{n/2,0} = \left\{\begin{array}{cc}
\displaystyle-\frac12 \frac1{(q-q^{-1})^{n-2} [\frac{n-2}2]!^2}\frac{[\frac n2]}{\alpha[\frac n2]-[n]}
& r = 0,
\\[0.5cm]
\displaystyle\frac12 \frac1{(q-q^{-1})^{n-2} [\frac{n-2}2]!^2}\frac{[\frac n2]}{\alpha[\frac n2]+[n]}
& r = \frac n2,
\\[0.5cm]
0 & \textrm{otherwise.}
\end{array}
\right.
\ee
\end{Proposition}
These constants satisfy the reflection symmetry
\be
\label{eq:reflection.symmetry}
\Gamma_{k,\ell}\big|_{\omega \to \omega^{-1}} = \widehat\gamma\, \Gamma_{k,m_k-\ell}\,.
\ee
Here, mapping $\omega\to\omega^{-1}$ with $\omega = \omega_{n,r}$ is equivalent  to simultaneously replacing $\widehat\gamma \to \widehat\gamma^{-1}$ and $r \to n-r$. This result is non-trivial and is proven in \cref{Sec:proof.closed.form}.

%
\section{Markov traces}\label{sec:Markov}
%

In this section, we discuss the Markov traces for the affine Temperley--Lieb algebras and their uncoiled quotients. We obtain simple formulas for their evaluation on the Wenzl--Jones projectors for the algebras $\qatl_n(\beta,\gamma = 1)$ for $n$ odd and $\qatla_n(\beta,\alpha)$ for $n$ even. 

\subsection[Definition of the Markov trace for $\atl_n(\beta)$]{Definition of the Markov trace for $\boldsymbol{\atl_n(\beta)}$}

Let us denote the greatest common divisor of two integers $i$ and $j$ by $i\wedge j$, and use the convention $i \wedge 0 = 0 \wedge i = i$. The \emph{Markov trace} $\mathcal F_{\boldsymbol \alpha}: \atl_n(\beta) \to \mathbb C$ is defined on the connectivities $c$ as
\be
\label{eq:def.F}
\mathcal F_{\boldsymbol \alpha}(c) = \beta^{n_\beta(c)} \prod_{i \wedge j=1} \alpha_{i,j}^{n_{i,j}(c)}.
\ee
To compute the numbers $n_\beta(c)$ and $n_{i,j}(c)$, we first connect the loop segments attached to the top and bottom segments of $c$, thus embedding this connectivity on a torus (see \eqref{eq:Tr.Cycl.Dessin} below for a visualisation). In the resulting diagram, these loop segments form closed loops. Some of them can be deformed into a point and are thus contractible. Their number for $c$ is $n_\beta(c)$. Some loops may wind around the torus and are thus non-contractible. These loops can wind around the torus $i$ in the horizontal direction and $j$ times in the vertical direction, with $i$ and $j$ coprime numbers (namely $i \wedge j = 1$). By convention, we choose $i \in \mathbb Z$ and $j \in \mathbb Z_{\ge 0}$. Following a curve moving upwards in the rectangular box delimiting $c$, we assign it a positive value of $i$ if it winds around the torus horizontally by moving towards the right, and a negative integer if it winds by moving towards the left. For example, embedding the connectivities $\Omega$ and $\Omega^{-1}$ of $\atl_{n=1}(\beta)$ on a torus, we obtain one non-contractible loop with the windings $(i,j) = (1,1)$ and $(-1,1)$, respectively. Thus the product in \eqref{eq:def.F} runs over all such coprime pairs $(i,j)$. The number $n_{i,j}(c)$ counts the number of non-contractible loops of $c$ with the winding $(i,j)$. For any given $c$, at most one $n_{i,j}(c)$ can be non-zero. Thus in \eqref{eq:def.F}, each loop with the winding $(i,j)$ is assigned the weight $\alpha_{i,j}$. Finally, we note that $\boldsymbol \alpha$ in bold refers to the set of weights $\alpha_{i,j}$.\medskip

The Markov trace satisfies the cyclicity property
\be
\mathcal F_{\boldsymbol \alpha}(c_1 c_2) = \mathcal F_{\boldsymbol \alpha}(c_2 c_1).\label{eq:Tr.Cycl}
\ee
This property can easily be proven by a diagrammatic manipulation of the trace, sliding one diagram along the cylinder:
\be\label{eq:Tr.Cycl.Dessin}
\begin{pspicture}[shift=-2.4](-0.25,0)(4.25,5)
\psframe[framearc=0.5](0,0)(4,5)
\psframe[framearc=0.5](1.25,1)(2.75,4)
\psframe[linecolor=blue,linewidth=1.5pt,framearc=0.5](.25,.25)(3.75,4.75)
\psframe[linecolor=blue,linewidth=1.5pt,framearc=0.5](1,.75)(3,4.25)
\pspolygon[fillstyle=solid,fillcolor=lightlightblue,linewidth=0pt](0,1.5)(0,2.5)(1.25,2.5)(1.25,1.5)
\rput(0.65,3.75){$\ldots$}
\rput (0.625,2){$c_1$} 
\pspolygon[fillstyle=solid,fillcolor=lightlightblue,linewidth=0pt](0,2.5)(0,3.5)(1.25,3.5)(1.25,2.5)
\rput (0.625,3){$c_2$} 
\end{pspicture} =
\begin{pspicture}[shift=-2.4](-0.25,0)(4.25,5)
\psframe[framearc=0.5](0,0)(4,5)
\psframe[framearc=0.5](1.25,1)(2.75,4)
\psframe[linecolor=blue,linewidth=1.5pt,framearc=0.5](.25,.25)(3.75,4.75)
\psframe[linecolor=blue,linewidth=1.5pt,framearc=0.5](1,.75)(3,4.25)
\pspolygon[fillstyle=solid,fillcolor=lightlightblue,linewidth=0pt](0,2)(0,3)(1.25,3)(1.25,2)
\rput (0.625,2.5){$c_1$} 
\pspolygon[fillstyle=solid,fillcolor=lightlightblue,linewidth=0pt](2.75,2)(2.75,3)(4,3)(4,2)
\rput{180}(3.375,2.5){$c_2$}
\rput(0.65,3.2){$\ldots$}
\end{pspicture} =
\begin{pspicture}[shift=-2.4](-0.25,0)(4.25,5)
\psframe[framearc=0.5](0,0)(4,5)
\psframe[framearc=0.5](1.25,1)(2.75,4)
\psframe[linecolor=blue,linewidth=1.5pt,framearc=0.5](.25,.25)(3.75,4.75)
\psframe[linecolor=blue,linewidth=1.5pt,framearc=0.5](1,.75)(3,4.25)
\pspolygon[fillstyle=solid,fillcolor=lightlightblue,linewidth=0pt](0,1.5)(0,2.5)(1.25,2.5)(1.25,1.5)
\rput (0.625,2){$c_2$} 
\pspolygon[fillstyle=solid,fillcolor=lightlightblue,linewidth=0pt](0,2.5)(0,3.5)(1.25,3.5)(1.25,2.5)
\rput (0.625,3){$c_1$} 
\rput(0.65,3.75){$\ldots$}
\end{pspicture}.
\ee

We denote as $\mathcal F_\alpha$ the \textit{homogeneous Markov trace}, where each $\alpha_{i,j}$ is set to the same value $\alpha$.

\subsection{Relation to standard modules}

As shown in \cite{RJ07,MDPR13}, the Markov trace can be expressed as
\be
\label{eq:F.Cdj}
\mathcal F_{\boldsymbol \alpha}(c) = \left\{\begin{array}{ll}
\displaystyle\sum_{d=1,3,...,n} \sum_{j=-\infty}^\infty 2\, T_{d \wedge j} (\tfrac12 \alpha_{\frac j{d \wedge j},\frac{d}{d \wedge j}})\, C_{d,j} & n \textrm{ odd},
\\
\displaystyle\textrm{tr}_{\raisebox{-0.05cm}{\tiny$\repW_{n,0,z}$}} (c) \big|_{\alpha \to \alpha_{1,0}} + \sum_{d=2,4,...,n} \sum_{j=-\infty}^\infty 2\, T_{d \wedge j} (\tfrac12 \alpha_{\frac j{d \wedge j},\frac{d}{d \wedge j}})\, C_{d,j} & n \textrm{ even},
\end{array}\right.
\ee
where $T_n(x)$ is the $n$-th Chebyshev polynomial of the first kind: $T_n(\cos \theta) = \cos n \theta$. The coefficients~$C_{d,j}$ are defined from the Laurent series expansion of the matrix trace of $c$ over the standard module $\repW_{n,d,z}$:
\be
\label{eq:trace.Cdj}
\textrm{tr}_{\raisebox{-0.05cm}{\tiny$\repW_{n,d,z}$}} (c) = \sum_{j = -\infty}^\infty  z^{-j} C_{d,j}\, .
\ee
Any given connectivity $c$ of the algebra has finitely many non-zero $C_{d,j}$.\medskip

Following the arguments presented in \cite{MDKP24} for the periodic dilute Temperley--Lieb algebra, we now show that the homogeneous Markov trace can be rewritten explicitly in terms of traces over standard modules. For $d \ge 1$, we write
\be
\sum_{j = -\infty}^{\infty} T_{d \wedge j}(\tfrac \alpha 2)\, C_{d,j} 
= \sum_{x = 0}^{d-1} \sum_{y = -\infty}^\infty T_{d \wedge (x+yd)}(\tfrac \alpha 2) C_{d,x+yd} 
= \sum_{x = 0}^{d-1} T_{d \wedge x}(\tfrac \alpha 2) M_{d,x},
\ee
where 
\be
M_{d,x} = \sum_{y =-\infty}^\infty  C_{d,x+yd}\,.
\ee

Let us now consider the trace of the transfer matrix over the standard modules with the twist set to $z = \eE^{2 \iI \pi m/d}$ for some $m \in \{0,1, \dots, d-1\}$. For $d \ge 1$, we have
\be
\textrm{tr}_{\raisebox{-0.05cm}{\tiny$\repW_{n,d, \exp(2 \iI \pi m/d)}$}}(c) = \sum_{x = 0}^{d-1}\sum_{y=-\infty}^\infty \eE^{-2 \iI \pi m (x+yd)/d} C_{d,x+yd} = \sum_{x=0}^{d-1} \eE^{-2 \iI \pi mx/d} M_{d,x}.
\ee
We now show that
\be
2 \sum_{x=0}^{d-1} T_{d \wedge x}(\tfrac \alpha 2) M_{d,x} = \sum_{m=0}^{d-1} \Theta_{d,m} \,\textrm{tr}_{\raisebox{-0.05cm}{\tiny$\repW_{n,d,\exp(2 \iI \pi m/d)}$}}(c),
\ee
where 
\be
\Theta_{d,m} = \frac 2d \sum_{j=0}^{d-1} \eE^{2 \iI \pi j m/d}\, T_{d \wedge j}(\tfrac \alpha 2)\,.
\ee
This stems from a direct computation:
\begin{alignat}{2}
\sum_{m=0}^{d-1}
 \Theta_{d,m}\, \textrm{tr}_{\raisebox{-0.05cm}{\tiny$\repW_{n,d,\exp(2 \iI \pi m/d)}$}}(c) 
&= \frac 2d \sum_{x=0}^{d-1} \sum_{j=0}^{d-1} T_{d \wedge j}(\tfrac \alpha 2) M_{d,x} \sum_{m=0}^{d-1} \eE^{2 \iI \pi m(j-x)/d}
\nonumber\\&= 2\sum_{x=0}^{d-1} \sum_{j=0}^{d-1} T_{d \wedge j}(\tfrac \alpha 2) M_{d,x} \delta_{j,x}
\nonumber\\&= 2\sum_{x=0}^{d-1}T_{d \wedge x}(\tfrac \alpha 2) M_{d,x}\, .
\end{alignat}
As a result, we have
\be
\label{eq:F.trW}
\mathcal F_{\alpha}(c) = \left\{\begin{array}{ll}
\displaystyle\sum_{d=1,3,...,n}\sum_{m=0}^{d-1} \Theta_{d,m} \,\textrm{tr}_{\raisebox{-0.05cm}{\tiny$\repW_{n,d,\exp(2 \iI \pi m/d)}$}}(c)
& n \textrm{ odd},
\\
\displaystyle\textrm{tr}_{\raisebox{-0.05cm}{\tiny$\repW_{n,0,z}$}} (c) + \sum_{d=2,4,...,n} 
\sum_{m=0}^{d-1} \Theta_{d,m} \,\textrm{tr}_{\raisebox{-0.05cm}{\tiny$\repW_{n,d,\exp(2 \iI \pi m/d)}$}}(c) 
& n \textrm{ even}. 
\end{array}\right. 
\ee
As shown in \cite{MDKP24}, the functions $\Theta_{d,m}$ are closely related to certain functions $\Lambda(M,N)$ introduced in Section 3 of \cite{dFSZ87} in the computation of torus partition functions of the critical $O(n)$ model using Coulomb Gas arguments.

\subsection{Markov traces for the uncoiled algebras and their projectors}
The standard modules arising in \eqref{eq:F.trW} are precisely those of the algebras $\qatl_n(\beta,\gamma = 1)$ for $n$ odd and $\qatla_n(\beta,\alpha)$ for $n$ even; see~\eqref{eq:Std.Mod.uaTL.odd} and~\eqref{eq:Std.Mod.uaTLa.even}. From the known evaluations of the projectors $Q_{n,r}$ on the standard modules, we directly obtain their homogeneous Markov traces. Summing over~$r$, we then readily obtain the homogeneous Markov traces of the projectors $Q_n$ of the corresponding uncoiled periodic algebras. 
 
\begin{Proposition}\label{prop:MarkovTrace.FQnr.et.FQn}
The homogeneous Markov trace for the uncoiled affine and periodic Temperley--Lieb algebras $\qatla_n(\beta,\alpha)$, $\qptla_n(\beta,\alpha)$ for $n$ even and $\qatl_n(\beta,\gamma = 1)$, $\qptl_n(\beta,\gamma=1)$ for $n$ odd are given by 
\be
\label{eq:FQnr.et.T}
\mathcal F_\alpha (Q_{n,r}) = \Theta_{n,r} = \frac 2n \sum_{j=0}^{n-1} \eE^{2 \iI \pi j r/n} T_{n \wedge j}(\tfrac \alpha 2) \qquad \textrm{for } \qatla_n(\beta,\alpha) \textrm{ and }  \qatl_n(\beta,\gamma = 1),
\ee
and
\be
\label{eq:FQn.et.T}
\mathcal F_\alpha(Q_n) = \sum_{r=0}^{n-1} \mathcal F_\alpha(Q_{n,r}) = 2\, T_n(\tfrac \alpha 2)
\qquad \textrm{for } \qptla_n(\beta,\alpha) \textrm{ and }  \qptl_n(\beta,\gamma = 1).
\ee
\end{Proposition}

Thus $\mathcal F_\alpha$ appears to be a good Markov trace for these algebras. This indeed makes sense:  the homogeneous Markov trace, where all loops are assigned the same weight independently of their winding around the torus, is a good trace for the uncoiled algebras for which the windings are uncoiled with unit factors.\medskip

This, of course, begs the question of whether there are also such traces for the algebras $\qatlb_n(\beta,\gamma)$ and $\qatl_n(\beta,\gamma \neq 1)$, and likewise for the corresponding periodic uncoiled algebras. Clearly,  $\mathcal F_\alpha$ is not a good candidate, as it gives inconsistent results. To illustrate why, let us compute $\mathcal F_{\boldsymbol \alpha}(\Omega^n)$ in the two above-mentioned uncoiled affine algebras:
\be
\alpha_{1,1}^n = \mathcal F_{\boldsymbol \alpha}(\Omega^n) = \gamma\, \mathcal F_{\boldsymbol \alpha}(\id) = \gamma\, \alpha_{1,0}^n.
\ee
Setting $\alpha_{1,1} = \alpha_{1,0} = \alpha$ yields a contradiction, except if $\gamma = 1$ or $\alpha = 0$. Are there other specialisations of the weights $\alpha_{i,j}$ that lead to consistent Markov traces, for the other values of $\gamma$? A trivial solution is to set $\alpha_{i,j} = 0$ for all $i,j$. The following proposition addresses this question in the case $\alpha_{1,0} \neq 0$.

\begin{Proposition}
Let $\alpha_{1,0}\neq 0$. Then there is no specialisation of the loop weights $\alpha_{i,j}$ that make $\mathcal F_{\boldsymbol\alpha}$ a consistent Markov trace for $\qatlb_n(\beta,\gamma)$ for $\gamma \neq 1$, and likewise for $\qatl_n(\beta,\gamma)$ for $\gamma^2 \neq 1$.
\end{Proposition} 
\proof
We have
\be
\alpha_{1,1}^{n-2k} = \mathcal F_{\boldsymbol \alpha}(e_0 e_2 \dots e_{2k-2}\Omega^n e_1 e_3 \dots e_{2k-1}) = \gamma\,\mathcal F_{\boldsymbol \alpha}(e_0 e_2 \dots e_{2k-2} e_1 e_3 \dots e_{2k-1}) = \gamma\,\alpha_{1,0}^{n-2k},
\ee
for $k \in \{ 0,1, \dots, \lfloor\tfrac n2\rfloor\}$, and therefore
\be
\label{eq:ratios.of.alphas}
\Big(\frac{\alpha_{1,1}}{\alpha_{1,0}}\Big)^{n-2k} = \gamma, \qquad k = 0,1,2, \dots, \lfloor\tfrac n2\rfloor.
\ee
For $n$ even, let us consider a pair $(k,k')$ satisfying $k+k' = \frac n2$. Then we have
\be
\gamma^2 = \Big(\frac{\alpha_{1,1}}{\alpha_{1,0}}\Big)^{n-2k} \Big(\frac{\alpha_{1,1}}{\alpha_{1,0}}\Big)^{n-2k'} = \Big(\frac{\alpha_{1,1}}{\alpha_{1,0}}\Big)^{n} = \gamma,
\ee
where at the last line we used \eqref{eq:ratios.of.alphas} for $k = 0$. Recalling that $\gamma \in \mathbb C^\times$, this identity holds only for $\gamma = 1$. For $n$ odd, we instead choose a triple $(k,k',k'')$ satisfying $k+k'+k'' = n$. Then
\be
\gamma^3 = \Big(\frac{\alpha_{1,1}}{\alpha_{1,0}}\Big)^{n-2k} \Big(\frac{\alpha_{1,1}}{\alpha_{1,0}}\Big)^{n-2k'} \Big(\frac{\alpha_{1,1}}{\alpha_{1,0}}\Big)^{n-2k''} = \Big(\frac{\alpha_{1,1}}{\alpha_{1,0}}\Big)^{n} =  \gamma
\ee
which holds only for $\gamma \in \{+1,-1\}$.
\eproof

We conclude that, if there exist Markov traces with $\alpha_{1,0} \neq 0$ that are consistent for the uncoiled algebras $\qatlb_n(\beta,\gamma \neq 1)$ and $\qatl_n(\beta,\gamma \neq \pm 1)$, they cannot be of the form \eqref{eq:def.F}.

%
\section{Conclusion}\label{sec:conclusion}
%

In this paper, we defined and studied quotients of the affine and periodic Temperley--Lieb algebras: the uncoiled affine and periodic Temperley--Lieb algebras. We described them in terms of sandwich diagrams, computed their dimensions, and showed that they are affine cellular, except in one special case where skew sandwich cellularity is needed. We also constructed Wenzl--Jones projectors for these algebras in terms of the sandwich diagrams and the similar projector of the ordinary Temperley--Lieb algebra. The construction involves certain coefficients $\Gamma_{k,\ell}$, for which we derived linear recurrence relations; see \cref{Prop.const.rel}. The closed forms for these constants are then given in \cref{Prop.const.uptla.uptlb,Prop.const.uptl,Prop.const.qatl.qatla.qatlb}. Finally, we computed the homogeneous Markov trace of the Wenzl--Jones projectors and obtained simple formulas in terms of Chebyshev polynomials of the first kind for the uncoiled affine algebras $\qatla_n(\beta,\alpha)$ and $\qatl_n(\beta,\gamma \neq 1)$, and for their periodic counterparts.\medskip 

This work marks a first step in the investigation of these algebras and leaves open a number of questions. A first step is to investigate the representation theory of the uncoiled algebras, as a way to approach the same investigations for the affine and periodic Temperley--Lieb algebras. To illustrate, let us focus on the affine case for $n$ odd. To be precise, we wish to study a functor between the categories of finite-dimensional modules of the affine Temperley--Lieb algebra and of the uncoiled affine Temperley--Lieb algebras
\be
\label{eq:Functor.uncoiled.affine}
\mathcal{F}: \ \operatorname{Mod}\atl_n(\beta) \ \longrightarrow\
\bigcup_{\gamma\in\field{C}} \operatorname{Mod} \qatl_n({\beta,\gamma}).
\ee
The functor $\mathcal{F}$ decomposes into functors $\mathcal{F}_{\gamma}$ projecting solely onto $\operatorname{Mod} \qatl_n(\beta,\gamma)$. Each functor is full. The full functor $\mathcal{F}$ sends a standard module $\repW_{n,d,z}$ of $\atl_n(\beta)$ to a standard module $\repW_{n,d,z}$ of an uncoiled algebra $\qatl_n(\beta,\gamma)$. The same analysis could, of course, be done for $\ptl_n(\beta)$ using its uncoiled algebras $\qptl_n(\beta,\gamma)$. The properties of $\mathcal F$ will hopefully shed some light on the representation theory of $\atl_n(\beta)$.\medskip

We then note that there exists a larger family of uncoiled algebras labelled by a positive integer~$k$, with the uncoiled algebras studied in this paper corresponding to $k=1$. For instance, for $\atl_n(\beta)$ with $n$ odd, this is achieved by replacing the unwinding relation~\eqref{eq:unwinding.ptl.atl.odd} by the relation 
\begin{equation}\label{eq:unwinding.atl.odd.k}
\Omega^{kn} = \gamma \id, \qquad k \in \mathbb Z_{\ge 1}.
\end{equation}
This means that the curves may wind along the periodic boundary condition
multiple times, and the quotient relation only equates a connectivity to a multiple of the identity if all its curves wind $k$ times. We therefore define the \emph{$k$-uncoiled affine Temperley--Lieb algebras} $\qatl_{n,k}(\beta,\gamma)$ as
\be
\label{eq:odd.quotient.k}
\qatl_{n,k}(\beta,\gamma) = \big\langle \Omega, \Omega^{-1}, e_0, e_1,\dots, e_{n-1} \big\rangle \Big/ \big\{\eqref{eq:PTL.relations},\eqref{eq:ATL.relations},\eqref{eq:unwinding.atl.odd.k}\big\}.
\ee
Their dimension increases with $k$ and is easily computed:
\be
\dim \mathsf{uaTL}_{n,k}(\beta,\gamma) = \sum_{d=1,3,\dots, n} k\,d\,|\mathsf{B}_{n,d}|^2 = k \dim \mathsf{uaTL}_{n,1}(\beta,\gamma).
\ee
One could then investigate the representation theory of $\mathsf{uaTL}_{n,k}(\beta,\gamma)$ for each $k$. Studying the limiting process $k\to \infty$ could be interesting, as it may give us some insight to understand the representation theory of the affine Temperley--Lieb algebra, for which the classification of the indecomposable modules is an open problem.\medskip

A second direction is to explore the characteristic $p$ case. The work done here always assumed that the variables $\beta$, $\alpha$, $\gamma$ and $z$ are complex numbers. However, most of the definitions and results can be extended to positive characteristics and thus fit in the framework of the recent work of Spencer~\cite{MS21,S21}. We note, however, that the expressions in \cref{Prop.const.uptla.uptlb,Prop.const.uptl,Prop.const.qatl.qatla.qatlb} for the constants $\Gamma_{k,\ell}$ are written in terms of $q$ instead of $\beta$. To study the similar objects for the modular uncoiled Temperley--Lieb algebra, one would need to rewrite these constants purely in terms of $\beta$.\medskip

We note that finding closed-form formulas for the Wenzl--Jones projectors as sums of connectivities is a hard problem, even for the Temperley--Lieb case. In this case, recursive formulas for the coefficients in these sums were found long ago~\cite{W88,FK97} and closed-form formulas were conjectured by Ocneanu~\cite{O02}. However, their proofs are known only in a few limited cases~\cite{R07}. These coefficients have deep relations with Kazhdan--Luzstig theory. Indeed, it was recently shown that they can be expressed as sums of Kazhdan--Lusztig polynomials~\cite{B24}, for which no closed formulas are 
known. Our results directly allow us to obtain closed formulas for the coefficients of the projectors for the uncoiled algebras in the basis of connectivities, as sums of the coefficients conjectured by Ocneanu for the Temperley--Lieb case. 
Moreover, it would be interesting to see if our methods with contour integrals can be adapted to also give a complete proof of Ocneanu's formulas for the coefficients of the Wenzl--Jones projectors of $\tl_n(\beta)$.\medskip

After the completion of this project, we discovered that it is possible to construct a proper cellular basis for the uncoiled affine algebras following the arguments from~Propositions~3.4.2 and~3.5.6 of \cite{G98}, with minor adaptations for the cases with $n$ even. In short, the sandwich diagram basis is replaced by certain linear combinations of diagrams. This is not a big change for the purposes of this paper; in particular, the cell modules are the same. In this work, we presented the sandwich basis, as it is the simplest and because it also works with the uncoiled periodic algebras. For further investigations of the representation theory of the uncoiled algebras, it could be useful to use
the proper cellular bases. 

\subsection*{Acknowledgments}

We thank T.~Pinet for suggesting the idea of the $k$-uncoiled algebras as well as for useful discussions with ALR at the Les Houches School of Physics {\it Applications of Hecke and Related Algebras: Representations, Integrability and Physics} in February--March 2023. We also thank Y.~Ikhlef and J.L.~Jacobsen for useful discussions, D.~Tubbenhauer for precision on their work on sandwich cellular algebras, and R.~Spencer for collaboration on preliminary computations. We further thank the referee for their careful reading of the manuscript and for noticing subtleties on the presentations of the uncoiled periodic algebras and their sandwich cellularity, which led us to realize the importance of the notion of skew sandwich cellularity. 

\subsection*{Declarations}
{\bf Ethical Approval.} Not applicable.

\noindent{\bf Competing interests.} The authors declare no conflicts of interest associated with this work.

\noindent{\bf Authors' contributions.} ALR and AMD took an active role in all parts of the research. ALR and AMD wrote and reviewed the manuscript together.

\noindent{\bf Funding.} ALR and AMD were supported by the EOS Research Project [project number 30889451]. Furthermore, ALR benefitted from doctoral and postdoctoral scholarships of the Fonds de recherche du Qu\'ebec -- Nature et technologies [grant numbers 270527 and 326641], and was funded by the Deutsche Forschungsgemeinschaft (DFG, German Research Foundation) under Germany's Excellence Strategy - GZ 2047/1, Projekt-ID 390685813. 

\noindent{\bf Availability of data and materials.} Data sharing is not applicable to this article as no datasets were generated or analysed
during the study.


\appendix

\section{Central elements for \texorpdfstring{$\boldsymbol{\ptl_n(\beta)}$}{pTLn(beta)}}\label{sec:non.isomorphic}
 
In this section, we construct central elements in $\ptl_n(\beta)$. They can be used to study the possible homomorphisms between standard modules over this algebra. Let us first recall the definition of known central elements in $\atl_n(\beta)$; see~\cite{MDSA13}. The {\em braid transfer matrices} are defined as 
\begin{subequations}\label{eq:Fb.def}
\begin{alignat}{2}
\Fb &= \ 
\psset{unit=.95cm}
\begin{pspicture}[shift=-1.1](-0.2,-0.7)(5.2,1.1)
\facegrid{(-0,0)}{(5,1)}
\psline[linewidth=1.5pt,linecolor=blue]{-}(-0.2,0.5)(3,0.5)
\psline[linewidth=1.5pt,linecolor=blue]{-}(0.5,0)(0.5,0.35)
\psline[linewidth=1.5pt,linecolor=blue]{-}(0.5,0.65)(0.5,1)
\psline[linewidth=1.5pt,linecolor=blue]{-}(1.5,0)(1.5,0.35)
\psline[linewidth=1.5pt,linecolor=blue]{-}(1.5,0.65)(1.5,1)
\psline[linewidth=1.5pt,linecolor=blue]{-}(2.5,0)(2.5,0.35)
\psline[linewidth=1.5pt,linecolor=blue]{-}(2.5,0.65)(2.5,1)
\rput(3.53,.51){$\dots$}
\psline[linewidth=1.5pt,linecolor=blue]{-}(4,0.5)(5.2,0.5)
\psline[linewidth=1.5pt,linecolor=blue]{-}(4.5,0)(4.5,0.35)
\psline[linewidth=1.5pt,linecolor=blue]{-}(4.5,0.65)(4.5,1)
\psline[linewidth=1pt]{<->}(0,-0.15)(5,-0.15)
\rput(2.5,-0.35){$_n$}
\end{pspicture}
\ , \qquad \quad&&
\psset{unit=.95cm}
\begin{pspicture}[shift=-.4](1,1)
\facegrid{(0,0)}{(1,1)}
\psline[linewidth=1.5pt,linecolor=blue]{-}(0,0.5)(1,0.5)
\psline[linewidth=1.5pt,linecolor=blue]{-}(0.5,0)(0.5,0.35)
\psline[linewidth=1.5pt,linecolor=blue]{-}(0.5,0.65)(0.5,1)
\end{pspicture}
\ = q^{1/2} \
\begin{pspicture}[shift=-.4](1,1)
\facegrid{(0,0)}{(1,1)}
\rput(0,0){\loopa}
\end{pspicture}
\ + q^{-1/2} \ 
\begin{pspicture}[shift=-.4](1,1)
\facegrid{(0,0)}{(1,1)}
\rput(0,0){\loopb}
\end{pspicture}\ ,
\\
\Fbb&=\ 
\psset{unit=.95cm}
\begin{pspicture}[shift=-1.1](-0.2,-0.7)(5.2,1.1)
\facegrid{(-0,0)}{(5,1)}
\rput(3.53,.51){$\dots$}
\psline[linewidth=1.5pt,linecolor=blue]{-}(0,0.5)(-0.2,0.5)
\psline[linewidth=1.5pt,linecolor=blue]{-}(5,0.5)(5.2,0.5)
\psline[linewidth=1.5pt,linecolor=blue]{-}(0.5,0)(0.5,1)
\psline[linewidth=1.5pt,linecolor=blue]{-}(1.5,0)(1.5,1)
\psline[linewidth=1.5pt,linecolor=blue]{-}(2.5,0)(2.5,1)
\psline[linewidth=1.5pt,linecolor=blue]{-}(4.5,0)(4.5,1)
\psline[linewidth=1.5pt,linecolor=blue]{-}(0,0.5)(0.35,0.5)
\psline[linewidth=1.5pt,linecolor=blue]{-}(0.65,0.5)(1.35,0.5)
\psline[linewidth=1.5pt,linecolor=blue]{-}(1.65,0.5)(2.35,0.5)
\psline[linewidth=1.5pt,linecolor=blue]{-}(2.65,0.5)(3,0.5)
\psline[linewidth=1.5pt,linecolor=blue]{-}(4,0.5)(4.35,0.5)
\psline[linewidth=1.5pt,linecolor=blue]{-}(4.65,0.5)(5,0.5)
\psline[linewidth=1pt]{<->}(0,-0.15)(5,-0.15)
\rput(2.5,-0.35){$_n$}
\end{pspicture}
\ , \qquad \quad &&
\psset{unit=.95cm}
\begin{pspicture}[shift=-.4](1,1)
\facegrid{(0,0)}{(1,1)}
\psline[linewidth=1.5pt,linecolor=blue]{-}(0.5,0)(0.5,1)
\psline[linewidth=1.5pt,linecolor=blue]{-}(0,0.5)(0.35,0.5)
\psline[linewidth=1.5pt,linecolor=blue]{-}(0.65,0.5)(1,0.5)
\end{pspicture}
\ = q^{-1/2} \
\begin{pspicture}[shift=-.4](1,1)
\facegrid{(0,0)}{(1,1)}
\rput(0,0){\loopa}
\end{pspicture}
\ + q^{1/2} \ 
\begin{pspicture}[shift=-.4](1,1)
\facegrid{(0,0)}{(1,1)}
\rput(0,0){\loopb}
\end{pspicture}\ ,
\end{alignat}
\end{subequations}
where we recall that $\beta = -q-q^{-1}$. These two elements, along with $\Omega^n$ and $\Omega^{-n}$,  are known to be central in $\atl_n(\beta)$. On the standard modules, they evaluate to multiples of the identity:
\be
\label{eq:FO.eigenvalues}
\Fb - \bigg(z q^{d/2}+\frac1{z q^{d/2}}\bigg)\id \Big|_{\repW_{n,d,z}} = 0,
\qquad
\Fbb - \bigg(\frac{z}{q^{d/2}} + \frac{q^{d/2}}{z}\bigg)\id \Big|_{\repW_{n,d,z}} = 0,
\qquad
\Omega^{\pm n} - z^{\pm d} \id \Big|_{\repW_{n,d,z}} = 0.
\ee

We proceed to construct central elements in $\ptl_n(\beta)$. Let us first note that
\be
\Fb = q^{n/2}\Omega + q^{-n/2}\Omega^{-1} + \dots,
\qquad 
\Fbb = q^{-n/2}\Omega + q^{n/2}\Omega^{-1} + \dots,
\ee
where the missing terms are sums of odd connectivities which involve at least one generator $e_j$ in their decompositions into words. Clearly, $\Fb$ and $\Fbb$ are not elements of $\ptl_n(\beta)$. We can nonetheless construct a central element in $\ptl_n(\beta)$ by noting that
\be
\Fb^2 = q^{n}\Omega^2 + q^{-n}\Omega^{-2} + \dots,
\qquad
\Fbb^2 = q^{-n}\Omega^2 + q^{n}\Omega^{-2}  + \dots,
\qquad
\Fb\Fbb = \Omega^2 + \Omega^{-2} + \dots,
\ee
where the extra terms are sums of even connectivities in $\ptl_n(\beta)$. Then
\be
\Gb = \Fb^2 + \Fbb^2 - (q^n+q^{-n}) \Fb \Fbb
\ee
is a central element in $\ptl_n(\beta)$. One can also define other central elements by taking linear combinations of higher even powers of $\Fb$ and $\Fbb$.\medskip

It is also possible to construct central elements of $\ptl_n(\beta)$ that combine $\Omega^n$, $\Omega^{-n}$ and either $\Fb$ or $\Fbb$. To construct these elements, we consider polynomials in $\Fb$ of higher degrees and write
\be
\begin{array}{ll}
\Fb^m = (q^{n/2}\Omega + q^{-n/2}\Omega^{-1})^m + \dots, 
&\qquad 
2\,T_m(\tfrac12\Fb) = q^{nm/2}\Omega^m + q^{-nm/2}\Omega^{-m} + \dots,
\\[0.15cm]
\Fbb^m = (q^{-n/2}\Omega + q^{n/2}\Omega^{-1})^m + \dots, 
&\qquad
2\,T_m(\tfrac12\Fbb) = q^{-nm/2}\Omega^m + q^{nm/2}\Omega^{-m} + \dots,
\end{array}
\ee
where $m \in \mathbb N$ and $T_m(x)$ is the $m$-th Chebyshev polynomial of the first kind. The other terms are again words that have at least one generator $e_j$. For $m$ even, they are elements of $\ptl_n(\beta)$. Then setting $m = 2nk$ with $k$ a non-negative integer or half-integer, we define
\be
\label{eq:Hk}
\begin{array}{l}
\Hb_k =  2\,T_{2nk}(\tfrac12\Fb) - q^{n^2k}\Omega^{2nk} - q^{-n^2k}\Omega^{-2nk},\\[0.15cm]
\Hbb_k =  2\,T_{2nk}(\tfrac12\Fbb) - q^{-n^2k}\Omega^{2nk} - q^{n^2k}\Omega^{-2nk},
\end{array}
 \qquad
k \in
\left\{\begin{array}{rl}
\mathbb N & n \textrm{ odd},
\\[0.1cm]
\frac12\mathbb N & n \textrm{ even}.
\end{array}\right.
\ee
These choices of $k$ ensure that the elements $\Hb_k$ and $\Hbb_k$ are always even, as they then involve only even connectivities with at least one generator $e_j$. They are therefore central elements in $\ptl_n(\beta)$.

%
\section{Binomial identities and dimensions}\label{app:binom.identities}
%

In this appendix, we show that 
\be\label{eq:D(n).result}
D(n) = \sum^*_{1 \le d \le n} d \binom{n}{\frac{n-d}2}^2 =
\left\{\begin{array}{cl}
\displaystyle n\binom{n-1}{\tfrac {n-1}2}^2 & n \textrm{ odd},\\[0.5cm]
\displaystyle n\binom{n-1}{\tfrac n2}^2 & n \textrm{ even},
\end{array}\right.
\ee
where $\overset{*}{\sum}$ indicates that the sum runs over values with $d \equiv n \textrm{ mod } 2$. We write
\begin{alignat}{2}
D(n) &=  n  \sum^*_{1 \le d \le n} \binom{n}{\frac{n-d}2}^2 - \sum^*_{1 \le d \le n} (n-d) \binom{n}{\frac{n-d}2}^2 = n \sum^*_{1 \le d \le n} \binom{n}{\frac{n-d}2}^2 - 2n \sum^*_{1 \le d \le n-2}  \binom{n}{\frac{n-d}2}\binom{n-1}{\frac{n-d-2}2}
\nonumber\\
&
= n \sum^*_{1 \le d \le n} \binom{n}{\frac{n-d}2}^2 - 2n \sum^*_{1 \le d \le n-2}  \binom{n-1}{\frac{n-d-2}2}^2 - 2n \sum^*_{1 \le d \le n-2}  \binom{n-1}{\frac{n-d-2}2}\binom{n-1}{\frac{n-d}2}.
\label{eq:D(n).intermediate} \end{alignat}
Using the Vandermonde identity
\be
\bigg[\sum_{r=0}^n \binom{n}{r} x^r\bigg]\bigg[\sum_{s=0}^n \binom{n}{s} x^{n-s}\bigg] = \sum_{t=0}^{2n} \binom{2n}{t} x^t,
\ee
we find that these three sums evaluate to
\begin{subequations}
\begin{alignat}{2}
\sum^*_{1 \le d \le n} \binom{n}{\frac{n-d}2}^2 &= \frac12 \binom{2n}{n}-\frac12 \binom{n}{\frac n2}^2 \delta_{n\,\textrm{even}},
\\
\sum^*_{1 \le d \le n-2} \binom{n-1}{\frac{n-d-2}2}^2 &= \frac12 \binom{2n-2}{n-1}-\frac12 \binom{n-1}{\frac {n-1}2}^2 \delta_{n\,\textrm{odd}},
\\
\sum^*_{1 \le d \le n-2} \binom{n-1}{\frac{n-d-2}2}\binom{n-1}{\frac{n-d}2} &= \frac12 \binom{2n-2}{n}-\frac12 \binom{n-1}{\frac {n}2}^2 \delta_{n\,\textrm{even}}.
\end{alignat}
\end{subequations}
Thus $D(u)$ is written in terms of a finite combination of binomial coefficients, which is easily seen to simplify to \eqref{eq:D(n).result}.

%
\section{Proof of the formulas for the constants}\label{Sec:proof.closed.form}
%

In this section, we prove the formulas for the constants $\Gamma_{k,\ell}$ given in \cref{Sec:Conj.Values} by rewriting them as sums of contour integrals and showing that they satisfy the recursion relations obtained in \cref{sec:linear.relations}.

\paragraph{Relation between the periodic and affine constants.}

Let us first note that we have the identity
\begin{alignat}{2}
\frac 1n \sum_{r=0}^{n-1} \frac{\omega_{n,r}^{-2\ell}}{\omega_{n,r}^2 q^2 -1} 
&= - \frac{q^{-n+2\ell}}2 \bigg(\frac1{q^{-n}-\widehat\gamma} + \frac{(-1)^{n-2\ell}}{(-1)^n q^{-n}-\widehat\gamma} \bigg) \qquad \ell \in \{0,\tfrac12,1, \dots, \tfrac{n-1}2\}
\nonumber\\
& = \left\{\begin{array}{cll}
\displaystyle \frac{q^{2\ell}}{\widehat\gamma q^n - 1}& n \textrm{ even}, &\ell \in \{0,1,\dots, \frac{n-1}2\},
\\[0.3cm]
0& n \textrm{ even}, &\ell \in \{\frac12,\frac32,\dots, \frac{n-2}2\},
\\[0.3cm]
\displaystyle \frac{q^{2\ell}}{\widehat\gamma^2 q^{2n} - 1}& n \textrm{ odd}, &\ell \in \{0,1,\dots, \frac{n-2}2\},
\\[0.3cm]
\displaystyle \frac{\widehat\gamma\, q^{n+2\ell}}{\widehat\gamma^2 q^{2n} - 1}& n \textrm{ odd}, &\ell \in \{\frac12,\frac32,\dots, \frac{n-1}2\},
\end{array}\right.
\end{alignat}
where we recall that $\omega_{n,r}= \widehat\gamma^{1/n} \eE^{2 \pi \iI r/n}$. To show this, we use the same strategy as in \cref{sec:kernel}, namely we write
\be
\frac 1n \sum_{r=0}^{n-1} \frac{\omega_{n,r}^{-2\ell}}{\omega_{n,r}^2 q^2 -1}  = \oint \frac{\dd \omega}{2 \pi \iI} \frac{\omega^{n-2\ell-1}}{(\omega q-1)(\omega q+1)} \frac1{\omega^n - \widehat\gamma},
\ee
where the contour encircles all the poles at $\omega = \omega_{n,r}$ for $r \in \{0,1,\dots, n-1\}$, but not the two other poles of the integrand. For $0 \le \ell \le \frac{n-1}2$, there are no poles at the origin nor at infinity, and the integral can be evaluated from the residues at $\omega = q^{-1}$ and $\omega = -q^{-1}$, which then simplifies to the above expression.\medskip

We use this identity with $q \to q^{2 \sigma m_\kappa}$ to compute $\sum_{r=0}^{n-1} \Gamma_{k,\ell}$, where the constants $\Gamma_{k,\ell}$ are those for the uncoiled affine  algebras given in \eqref{eq:Gamma.final.affine} and evaluated at $\omega = \omega_{n,r}$. This yields precisely the formulas \eqref{eq:Gamma.final} and \eqref{eq:Gamma.final.odd} for the constants $\Gamma_{k,\ell}$ for the uncoiled periodic algebras, for $n$ even and odd respectively, as it should be according to \eqref{eq:Qn=sumQnr}. It is easy to check from \eqref{eq:Gamma.n/2.0} and \eqref{eq:Gamma.n/2.0.affine} that this also holds for $\Gamma_{n/2,0}$ for the algebras $\qatla_n(\beta,\alpha)$ and $\qptla_n(\beta,\alpha)$. We conclude that \cref{Prop.const.qatl.qatla.qatlb} implies \cref{Prop.const.uptla.uptlb,Prop.const.uptl}. We thus focus on proving \cref{Prop.const.qatl.qatla.qatlb} in the following.

\paragraph{Contour integral formulas.}
Our strategy to prove \cref{Prop.const.qatl.qatla.qatlb} is to rewrite the expressions~\eqref{eq:Gamma.final.affine}  for $\Gamma_{k,\ell}$ as contour integrals, and subsequently show that they satisfy the recursion relations of \cref{Prop.const.rel}. As a first step, we derive expressions for $\Gamma_{k,\ell}$ as sums of contour integrals. The $q$-binomial satisfies the following known relations
\begin{subequations}
\label{eq:qbin.rel}
\begin{alignat}{2}
&\prod_{j=\frac{1-n}2}^{\frac{n-1}2} (1+z q^{2j}) = \sum_{j=0}^{n}\left[\begin{matrix} n\\ j\end{matrix}\right] z^j, 
\quad
&&\oint \frac{\dd z}{2\pi \iI} \frac1{z^{m+1}} \prod_{j=\frac{1-n}2}^{\frac{n-1}2} (1+z q^{2j}) = \left[\begin{matrix} n\\ m\end{matrix}\right],
\label{eq:qbin.rel.a}\\\label{eq:qbin.rel.b}
&\prod_{j=\frac{1-n}2}^{\frac{n-1}2} \frac1{1-z q^{2j}} = 
\sum_{j=0}^{\infty}\left[\begin{matrix} n+j-1\\ j\end{matrix}\right] z^j, 
\quad
&&\oint \frac{\dd z}{2\pi \iI} \frac1{z^{m+1}} \prod_{j=\frac{1-n}2}^{\frac{n-1}2} \frac1{1-z q^{2j}} = \left[\begin{matrix} n+m-1\\ m\end{matrix}\right].
\end{alignat}
\end{subequations}
In these formulas, the contour integrals are over closed contours around the pole at the origin. Moreover, the products over $j$ use unit steps and thus run over integers and half-integers for $n$ odd and even, respectively.\medskip

The sums over $\kappa$ and $\tau$ in \eqref{eq:Gamma.final.affine} run over the values $\{1,2,\dots, k\}$ and $\{0,1,\dots,k-\kappa\}$, respectively. It is, however, convenient to extend these ranges to larger intervals using ${\left[\begin{smallmatrix} n\\ m\end{smallmatrix}\right]} = 0$ for $n\ge 0$ with $m<0$ or $m>n$. We therefore omit the bounds of the sums, understanding that they run over $\mathbb Z$. Using~\eqref{eq:qbin.rel.b}, we find 
\begin{alignat}{2}
\sum_{\tau} & q^{\sigma n \tau}
\left[\begin{matrix} 2m_k-2\ell+k-\kappa-\tau-1\\ k-\kappa-\tau\end{matrix}\right]
\left[\begin{matrix} 2\ell+\tau-1\\ \tau\end{matrix}\right] 
\nonumber\\&= 
\sum_{\tau} q^{\sigma n \tau}
\left[\begin{matrix} 2\ell+\tau-1\\ \tau\end{matrix}\right]
\oint \frac{\dd z}{2\pi \iI} \frac1{z^{k-\kappa-\tau+1}} \prod_{j=-(m_k-\ell-\frac12)}^{m_k-\ell-\frac12} \frac1{1-z q^{2j}}
\nonumber\\&= 
\oint \frac{\dd z}{2\pi \iI} \frac{1}{z^{k-\kappa+1}}\prod_{j=-(m_k-\ell-\frac12)}^{m_k-\ell-\frac12} \frac1{1-z q^{2j}} \prod_{j=-(\ell-\frac12)}^{\ell-\frac12} \frac1{1-z q^{\sigma n+2j}}\, .
\label{eq:tau.sum}
\end{alignat}
This yields
\begin{alignat}{2}
\label{eq:Gamma.Lambda.1}
\Gamma_{k,\ell} = \frac{\omega^{-2\ell}}{n(q-q^{-1})^{2k-1}[k]![k-1]!} \sum_{\sigma = \pm 1} &\sum_{\kappa}  \frac{(-1)^{k+\kappa} \sigma q^{2 \sigma \ell \kappa}}{\omega^2 q^{2\sigma m_\kappa}-1} 
\frac{\left[\begin{smallmatrix} k-1\\ \kappa-1\end{smallmatrix}\right]}
{\left[\begin{smallmatrix} n-\kappa-1\\ n-k-1\end{smallmatrix}\right]}
\\[0.2cm] \nonumber
&\times \oint \frac{\dd z}{2\pi \iI}\,
\frac{1}{z^{k-\kappa+1}}
\prod_{j=-(m_k-\ell-\frac12)}^{m_k-\ell-\frac12} \frac1{1-z q^{2j}} \prod_{j=-(\ell-\frac12)}^{\ell-\frac12} \frac1{1-z q^{\sigma n+2j}}\,.
\end{alignat}
A different expression for the left-hand side of \eqref{eq:tau.sum} is obtained by first replacing $\left[\begin{smallmatrix} 2\ell + \tau - 1\\ \tau\end{smallmatrix}\right]$ by a contour integral, and then evaluating the resulting sum over $\kappa$ containing the other $q$-binomial coefficient:
\begin{alignat}{2}
\label{eq:Gamma.Lambda.2}
\Gamma_{k,\ell} = \frac{\omega^{-2\ell}}{n(q-q^{-1})^{2k-1}[k]![k-1]!} \sum_{\sigma = \pm 1} &\sum_{\kappa}  \frac{(-1)^{k+\kappa} \sigma q^{\sigma(2 \ell \kappa+nk-n\kappa)}}{\omega^2 q^{2\sigma m_\kappa}-1} 
\frac{\left[\begin{smallmatrix} k-1\\ \kappa-1\end{smallmatrix}\right]}
{\left[\begin{smallmatrix} n-\kappa-1\\ n-k-1\end{smallmatrix}\right]}
\\[0.2cm] \nonumber
&\times \oint \frac{\dd z}{2\pi \iI}\,
\frac{1}{z^{k-\kappa+1}}
\prod_{j=-(\ell-\frac12)}^{\ell-\frac12} \frac1{1-z q^{2j}} \prod_{j=-(m_k-\ell-\frac12)}^{m_k-\ell-\frac12} \frac1{1-z q^{-\sigma n+2j}}\,.
\end{alignat}
Both of these formulas are valid for $0 \le k \le \lfloor\frac{n-1}2\rfloor$ and $\ell \in \{0,\frac12,1,\dots, m_k-\frac12\}$. For $\ell = 0$, they are valid provided that we use the convention
\be
\prod_{j=-(\ell-\frac12)}^{\ell-\frac12}  f(j) \xrightarrow{\ell\to0} 1\,.
\ee
Going back to \eqref{eq:Gamma.Lambda.1}, we use the partial fraction decomposition
\be
\label{eq:partial.fraction}
\frac1{y^s} \frac1{y^2 \alpha^2 -1} = \frac{\alpha^s}{y^2 \alpha^2 -1} \big(\id_{s \,\in\, 2 \mathbb Z}+\alpha y\, \id_{s \,\in\, 2 \mathbb Z+1}\big) - \sum_{\rho = 1}^s \frac{\alpha^{s-\rho}}{y^\rho}\,\id_{s-\rho \,\in\, 2 \mathbb Z}\,, \qquad s \in \mathbb Z_{\ge 0}\,.
\ee
We use this relation with $y = \omega$, $s = 2\ell$ and $\alpha = q^{\sigma m_\kappa}$ to re-express \eqref{eq:Gamma.Lambda.1} as a sum of two contributions
\begin{alignat}{2}
\label{eq:Gamma.Lambda.3}
\Gamma_{k,\ell} &= \sum_{\sigma = \pm 1} \sum_{\kappa}
\frac{1}{\omega^2 q^{2 \sigma m_{k}}-1}
\oint \frac{\dd z}{2\pi \iI}\, \Lambda_{k,\ell,\kappa,\sigma}
+ \sum_{\sigma = \pm 1} \sum_{\rho=1}^{2\ell}
\frac{1}{\omega^\rho}
\oint \frac{\dd z}{2\pi \iI}\, \Xi_{k,\ell,\rho,\sigma}\,,
\end{alignat}
where
\begin{subequations}
\begin{alignat}{2}
\Lambda_{k,\ell,\kappa,\sigma} &= \big(\id_{\ell \,\in\, \mathbb Z} + \omega q^{\sigma m_\kappa} \id_{\ell \,\in\, \mathbb Z+\frac12}\big)\,\frac{(-1)^{k+\kappa} \sigma q^{\sigma n \ell}}{n(q-q^{-1})^{2k-1}[k]![k-1]!} 
\frac{\left[\begin{smallmatrix} k-1\\ \kappa-1\end{smallmatrix}\right]}
{\left[\begin{smallmatrix} n-\kappa-1\\ n-k-1\end{smallmatrix}\right]}
\nonumber\\[0.15cm]
& \times\frac{1}{z^{k-\kappa+1}}
\prod_{j=-(m_k-\ell-\frac12)}^{m_k-\ell-\frac12} \frac1{1-z q^{2j}} \prod_{j=-(\ell-\frac12)}^{\ell-\frac12} \frac1{1-z q^{\sigma n+2j}}\,,
\\[0.15cm]
\Xi_{k,\ell,\rho,\sigma} & = \id_{2\ell + \rho\, \in \, 2\mathbb Z}\,\frac{(-1)^{k} \sigma q^{\sigma (n \ell-n \rho/2+\rho)}}{n(q-q^{-1})^{2k-1}[k]![k-1]!} 
\frac1{\left[\begin{smallmatrix} n-2\\ k-1\end{smallmatrix}\right]}
\\[0.15cm]
& \times
\frac1{z^k} \prod_{j=-\frac{n-3}2}^{\frac{n-3}2} (1-z q^{\sigma \rho+2j})
\prod_{j=-(m_k-\ell-\frac12)}^{m_k-\ell-\frac12} \frac1{1-z q^{2j}} \prod_{j=-(\ell-\frac12)}^{\ell-\frac12} \frac1{1-z q^{\sigma n+2j}}\,, \qquad 1 \le \rho \le 2\ell.
\nonumber 
\end{alignat}
\end{subequations}
In the second contribution, the sum over $\kappa$ was evaluated using \eqref{eq:qbin.rel.a}.\medskip

The expressions \eqref{eq:Gamma.Lambda.1}, \eqref{eq:Gamma.Lambda.2} and \eqref{eq:Gamma.Lambda.3} will all be useful below. In particular, we will show that the constants $\Gamma_{k,\ell}$ satisfy the recursion relations \eqref{eq:Gamma.constraint} by showing that certain closed contour integrals around $z=0$ vanish. These integrals usually vanish for one of the three following reasons:
\begin{enumerate}[(i)] 
\item The contour integral is of the form $\oint \dd z f(z)$ and vanishes because $f(z)$ vanishes identically. \label{item:reasoni}
\item The contour integral is of the form $\oint \dd z f(z)$ and vanishes because $f(z)$ is regular at $z=0$. We write this as $f(z) \simeq 0$.\label{item:reasonii}
\item The contour integral is the sum of two terms $\sum_{\sigma = \pm 1} \mathcal I_\sigma$, where $\mathcal I_\sigma = \oint \dd z f_\sigma(z)$. The functions $f_\sigma(z)$ have poles at the origin and at infinity only, and satisfy 
\be
\label{eq:f(z)}
f_\sigma(z^{-1}) = - z^2 f_{-\sigma}(z), \qquad \sigma \in \{+1,-1\}.
\ee 
It then follows that $\mathcal I_\sigma = -\mathcal I_{-\sigma}$, so that the sum of the two contributions indeed vanishes. \label{item:reasoniii}
\end{enumerate}

\paragraph{Reflection symmetry.} 
We now show that $\Gamma_{k,\ell}$ satisfies the reflection symmetry \eqref{eq:reflection.symmetry}, which we rewrite as 
\be
\label{eq:reflection.symmetry.v2}
\Gamma_{k,m_k-\ell} \big|_{\omega \to \omega^{-1}} = \widehat\gamma\,\Gamma_{k,\ell}\,.
\ee 
Using \eqref{eq:Gamma.Lambda.1}, we find 
\begin{alignat}{2}
\Gamma_{k,m_k-\ell} \big|_{\omega \to \omega^{-1}} &= \frac{\widehat\gamma\, \omega^{-2\ell-2k+2}}{n(q-q^{-1})^{2k-1}[k]![k-1]!} \sum_{\sigma = \pm 1} \sum_{\kappa}  \frac{(-1)^{k+\kappa} \sigma q^{-2 \sigma(m_k-\ell) \kappa}}{\omega^2 q^{2\sigma m_\kappa}-1} 
\frac{\left[\begin{smallmatrix} k-1\\ \kappa-1\end{smallmatrix}\right]}
{\left[\begin{smallmatrix} n-\kappa-1\\ n-k-1\end{smallmatrix}\right]}
\nonumber\\[0.2cm] 
&\times \oint \frac{\dd z}{2\pi \iI}\,
\frac{1}{z^{k-\kappa+1}}
\prod_{j=-(\ell-\frac12)}^{\ell-\frac12} \frac1{1-z q^{2j}} \prod_{j=-(m_k-\ell-\frac12)}^{m_k-\ell-\frac12} \frac1{1-z q^{-\sigma n+2j}},
\end{alignat}
where we used $\omega^n = \widehat\gamma$ and changed the summation variable $\sigma \to -\sigma$. We now use \eqref{eq:partial.fraction} with $y = \omega^2$, $s = k-1$ and $\alpha = q^{2 \sigma m_\kappa}$ to rewrite $\frac{\omega^{-2k+2}}{\omega^2 q^{2\sigma m_\kappa}-1}$ in terms of two contributions. The first contribution yields exactly $\widehat\gamma$ times the formula \eqref{eq:Gamma.Lambda.2} for $\Gamma_{k,\ell}$. For the second contribution, we reorganise the ratio of $q$-binomials, allowing us to evaluate the sum over $\kappa$ using \eqref{eq:qbin.rel.a}. After simplifications, the result reads 
\begin{alignat}{2}
\Gamma_{k,m_k-\ell} \big|_{\omega \to \omega^{-1}} & =  \widehat\gamma\, \Gamma_{k,\ell} + \widehat\gamma \sum_{\sigma = \pm 1} \sum_{\rho = 1}^{k-1} \frac1{\omega^{2(\rho+\ell)}} \frac{(-1)^k \sigma q^{\sigma (nk-n\rho-n+2\ell+2\rho)}}{n(q-q^{-1})^{2k-1}[k]![k-1]!}
\frac1{\left[\begin{smallmatrix} n-2\\ k-1\end{smallmatrix}\right]}
\\[0.15cm]\nonumber&\times
\oint \frac{\dd z}{2\pi \iI}\, \frac1{z^k}
\prod_{j=-\frac{n-3}2}^{\frac{n-3}2}(1-z q^{\sigma(-n+2\ell+2\rho)+2j})
\prod_{j=-(\ell-\frac12)}^{\ell-\frac12} \frac1{1-z q^{2j}} \prod_{j=-(m_k-\ell-\frac12)}^{m_k-\ell-\frac12} \frac1{1-z q^{-\sigma n+2j}}\,.
\end{alignat}
One can check that the integrand only has poles at the origin and at infinity, for all $\rho \in \{1,2, \dots, k-1\}$ and $\sigma \in \{+1,-1\}$. Changing variables to $z \to z^{-1}$, we find that this sum of contour integrals vanishes due to reason~(\ref{item:reasoniii}) explained above. This ends the proof of \eqref{eq:reflection.symmetry.v2}.

\paragraph{The recursion relations.}

The proof of \cref{Prop.const.qatl.qatla.qatlb} is inductive over increasing values of $k$. It is not hard to see that the expression \eqref{eq:Gamma.final} satisfies the initial conditions, namely $\Gamma_{-1,\ell} = 0$ and $\Gamma_{0,\ell} =\frac 1n\, \omega^{-2\ell}$. We now show that $\Gamma_{k,\ell}$ satisfies the recursive relation \eqref{eq:Gamma.constraint} for $1 \le k \le \frac{n-2}2$ and $\ell \in \{0, \frac 12, 1, \dots, m_k-\frac 12\}$ using the contour integral expressions derived above.\medskip 

The recursion relations for $\ell \in \{0,\frac12, m_k-1,m_k-\frac12\}$ must be treated with extra care. Indeed, the relation \eqref{eq:Gamma.constraint} specialised to $\ell = 0$ or $\ell = \frac12$ involves $\Gamma_{k,-1}$ or $\Gamma_{k,-1/2}$, which are not covered by the above integral equations and must instead be replaced by $\widehat\gamma\,\Gamma_{k,m_k-1}$ and $\widehat\gamma\,\Gamma_{k,m_k-1/2}$, respectively. Similarly, the relation \eqref{eq:Gamma.constraint} specialised to $\ell = m_k-1$ or $\ell = m_k-\frac12$ involves $\Gamma_{k,m_k}$ and $\Gamma_{k,m_k+1/2}$, which must be replaced by $\widehat\gamma^{-1}\Gamma_{k,0}$ and $\widehat\gamma^{-1}\Gamma_{k,\frac12}$, respectively. Using the reflection symmetry \eqref{eq:reflection.symmetry.v2}, we find that the two relations for $\ell = m_k-1$ or $\ell = m_k-\frac12$ are identical to the relations for $\ell = 1$ and $\ell = \frac12$, respectively. Thus only the cases $\ell = 0$ and $\ell = \frac12$ must be proved separately.\medskip

Let us therefore first investigate the cases $1 \le \ell \le m_k-\frac32$. After inserting the integral expressions~\eqref{eq:Gamma.Lambda.3} for $\Gamma_{k,\ell}$ in \eqref{eq:Gamma.constraint}, we find that the identity to be shown becomes 
\begin{alignat}{3} 
0 &= \sum_{\sigma = \pm 1} \sum_\kappa \frac1{\omega^2 q^{2 \sigma m_\kappa}-1} 
\oint \frac{\dd z}{2\pi \iI}\, \Big(f^1_k \Lambda_{k,\ell,\kappa,\sigma} + f^2_k(\Lambda_{k,\ell-1,\kappa,\sigma}+\Lambda_{k,\ell+1,\kappa,\sigma}) + f^3_{k-1,\ell} \Lambda_{k-1,\ell,\kappa,\sigma} 
\nonumber\\&
\hspace{5.4cm}+ f^4_{k-1,\ell+1} \Lambda_{k-1,\ell+1,\kappa,\sigma} + f^5_{k-2,\ell+1} \Lambda_{k-2,\ell+1,\kappa,\sigma}\Big)
\nonumber\\&
+\sum_{\sigma = \pm 1} \sum_{\rho=1}^{2\ell+2} \frac1{\omega^\rho}
\oint \frac{\dd z}{2\pi \iI}\, \Big(f^1_k \Xi_{k,\ell,\rho,\sigma} + f^2_k(\Xi_{k,\ell-1,\rho,\sigma}+\Xi_{k,\ell+1,\rho,\sigma}) + f^3_{k-1,\ell} \Xi_{k-1,\ell,\rho,\sigma} 
\nonumber\\&
\hspace{5.4cm}+ f^4_{k-1,\ell+1} \Xi_{k-1,\ell+1,\rho,\sigma} + f^5_{k-2,\ell+1} \Xi_{k-2,\ell+1,\rho,\sigma}\Big)\,,
 \label{eq:Lambda.constraint}
\end{alignat}
where for convenience we use the convention $\Xi_{k,\ell,\rho,\sigma} = 0$ for $\rho>2\ell$. It turns out that 
\be
f^1_k \Lambda_{k,\ell,\kappa,\sigma} + f^2_k(\Lambda_{k,\ell-1,\kappa,\sigma}+\Lambda_{k,\ell+1,\kappa,\sigma}) + f^3_{k-1,\ell} \Lambda_{k-1,\ell,\kappa,\sigma} 
+ f^4_{k-1,\ell+1} \Lambda_{k-1,\ell+1,\kappa,\sigma} + f^5_{k-2,\ell+1} \Lambda_{k-2,\ell+1,\kappa,\sigma} \simeq 0
\ee
for  $1 \le \kappa \le k$.
The left-hand side, in fact, vanishes identically for $1 \le \kappa \le k-2$. This can be shown directly by factoring the common factors and simplifying the resulting expression.
The corresponding integrals thus vanish due to reason~(\ref{item:reasoni}). For $\kappa = k-1$ and $\kappa = k$, the content of this parenthesis does not vanish, but the corresponding integrals vanish for reason~(\ref{item:reasonii}). The poles at the origin have degrees one or two in these cases, so the corresponding residues can be computed simply, yielding
\begin{subequations}
\begin{alignat}{2}
\Omega_{k,\ell,\sigma} &= \oint \frac{\dd z}{2\pi \iI}\, \Lambda_{k,\ell,k,\sigma} = \frac{(\id_{\ell \,\in\, \mathbb Z} + \omega q^{\sigma m_k} \id_{\ell \,\in\, \mathbb Z+1/2})\sigma q^{\sigma n \ell}}{n(q-q^{-1})^{2k-1}[k]![k-1]!}\,,
\\[0.15cm]
\Pi_{k,\ell,\sigma} &= \oint \frac{\dd z}{2\pi \iI}\, \Lambda_{k,\ell,k-1,\sigma} = -\frac{(\id_{\ell \,\in\, \mathbb Z} + \omega q^{\sigma m_{k-1}} \id_{\ell \,\in\, \mathbb Z+1/2})\sigma q^{\sigma n \ell}}{n(q-q^{-1})^{2k-1}[k]![k-2]![n-k]}([2\ell]q^{\sigma n}+[m_k-2\ell])\,.
\end{alignat}
\end{subequations}
These satisfy the relations
\begin{subequations}
\begin{alignat}{2}
0&=f^1_k \Omega_{k,\ell,\sigma} + f^2_k(\Omega_{k,\ell-1,\sigma}+\Omega_{k,\ell+1,\sigma})\,,
\\[0.15cm]
0&=f^1_k \Pi_{k,\ell,\sigma} + f^2_k(\Pi_{k,\ell-1,\sigma}+\Xi_{k,\ell+1,\sigma}) + f^3_{k-1,\ell}\Omega_{k-1,\ell,\sigma} + f^4_{k-1,\ell+1}\Omega_{k-1,\ell+1,\sigma}\,,
\end{alignat}
\end{subequations}
so the corresponding integrals for $\kappa = k$ and $\kappa = k-1$ indeed vanish.\medskip 

We similarly find that the content of the second parenthesis in \eqref{eq:Lambda.constraint} vanishes for $\rho \in \{1,2, \dots, 2\ell-1\} \cup \{2\ell +1\}$, so again the corresponding integrals vanish for reason~(\ref{item:reasoni}). For $\rho = 2\ell$ and $\rho = 2\ell + 2$, we instead find the identities
\begin{subequations}
\begin{alignat}{2}
&f^1_k \Xi_{k,\ell,2\ell,\sigma} + f^2_k(\Xi_{k,\ell-1,2\ell,\sigma}+\Xi_{k,\ell+1,2\ell,\sigma}) + f^3_{k-1,\ell} \Xi_{k-1,\ell,2\ell,\sigma} + f^4_{k-1,\ell+1} \Xi_{k-1,\ell+1,2\ell,\sigma} + f^5_{k-2,\ell+1} \Xi_{k-2,\ell+1,2\ell,\sigma}
\nonumber\\&= 
\frac{(-1)^{k+1} \sigma q^{\sigma(2\ell-n)}}{n (q-q^{-1})^{2k-1}} \frac{[n-k]!}{[n]![k-1]!} \frac{1}{z^k} \prod_{j=-\frac{n-3}2}^{\frac{n-3}2} (1-z q^{2 \sigma \ell+2j})
\prod_{j=-(m_k-\ell+\frac12)}^{m_k-\ell+\frac12} \frac1{1-z q^{2j}} \prod_{j=-(\ell-\frac32)}^{\ell-\frac32} \frac1{1-z q^{\sigma n+2j}}
\end{alignat} 
and
\begin{alignat}{2}
&f^1_k \Xi_{k,\ell,2\ell+2,\sigma} + f^2_k(\Xi_{k,\ell-1,2\ell+2,\sigma}+\Xi_{k,\ell+1,2\ell+2,\sigma}) + f^3_{k-1,\ell} \Xi_{k-1,\ell,2\ell+2,\sigma} + f^4_{k-1,\ell+1} \Xi_{k-1,\ell+1,2\ell+2,\sigma} 
\nonumber\\[0.25cm]&
+ f^5_{k-2,\ell+1} \Xi_{k-2,\ell+1,2\ell+2,\sigma} = \frac{(-1)^k \sigma q^{\sigma(2\ell+2)}}{n (q-q^{-1})^{2k-1}} \frac{[n-k]!}{[n]![k-1]!} \frac{(1-z q^{-\sigma n + 2 \ell + 1})(1-z q^{-\sigma n - 2 \ell - 1})}{z^k} 
\\&\hspace{4.25cm}\nonumber
\times\prod_{j=-\frac{n-3}2}^{\frac{n-3}2} (1-z q^{2\sigma (\ell+1)+2j})
\prod_{j=-(m_k-\ell+\frac12)}^{m_k-\ell+\frac12} \frac1{1-z q^{2j}} \prod_{j=-(\ell-\frac12)}^{\ell-\frac12} \frac1{1-z q^{\sigma n+2j}}\,.
\end{alignat} 
\end{subequations}
The two functions on the right-hand sides can be shown to satisfy the relation \eqref{eq:f(z)} and to have no pole except at the origin, so the resulting sum of integrals over $\sigma \in \{+1,-1\}$ vanishes due to reason~(iii).
\medskip

For the special case $\ell = \frac12$, we use the reflection symmetry and find that, for $1 \le k \le \lfloor\frac{n-3}2\rfloor$, the identity to be proven is\footnote{The special case with $n$ even and $(k,\ell) = (\frac{n-2}2,\frac12)$ is treated separately at the end of the section.}
\begin{alignat}{2}
\label{eq:Gamma.constraint.L=1/2}
0 &= f^1_k \Gamma_{k,\frac12} + f^2_k \Big(\Gamma_{k,\frac12}\big|_{\omega \to \omega^{-1}}+\Gamma_{k,\frac32}\Big)+f^3_{k-1,\frac12} \Gamma_{k-1,\frac12} + f^4_{k-1,\frac32} \Gamma_{k-1,\frac32} + f^5_{k-2,\frac32}\Gamma_{k-2,\frac32} 
\nonumber\\[0.15cm]
&+ f^{3\textrm{b}}_{k-1, m_{k-1}-\frac12} \Big(\Gamma_{k-1,\frac12}\big|_{\omega \to \omega^{-1}}\Big).
\end{alignat}
Because
\be
\Gamma_{k,\frac12}\big|_{\omega \to \omega^{-1}} = \sum_{\sigma = \pm 1} \sum_\kappa \frac1{\omega^2 q^{2 \sigma m_\kappa}-1} \oint \frac{\dd z}{2\pi \iI}\, \big(\!-\!q^{2 \sigma m_\kappa} \Lambda_{k,\frac12,\kappa,-\sigma}\big) + \omega \sum_{\sigma = \pm 1} \oint \frac{\dd z}{2\pi \iI}\, \Xi_{k,\frac12,1,\sigma}\,,
\ee
we must show that
\begin{alignat}{3} 
&0 = \sum_{\sigma = \pm 1} \sum_\kappa \frac1{\omega^2 q^{2 \sigma m_\kappa}-1} 
\oint \frac{\dd z}{2\pi \iI}\, \Big(
f^1_k \Lambda_{k,\frac12,\kappa,\sigma} 
- q^{2 \sigma m_\kappa} f^2_k \Lambda_{k,\frac12,\kappa,-\sigma} 
+ f^2_k \Lambda_{k,\frac32,\kappa,\sigma} 
+ f^3_{k-1,\frac12} \Lambda_{k-1,\frac12,\kappa,\sigma} 
\nonumber\\&
\hspace{5.4cm}
+ f^4_{k-1,\frac32} \Lambda_{k-1,\frac32,\kappa,\sigma} 
+ f^5_{k-2,\frac32} \Lambda_{k-2,\frac32,\kappa,\sigma} 
- q^{2 \sigma m_\kappa} f^{3\textrm{b}}_{k-1,m_k+\frac12}\Lambda_{k-1,\frac12,\kappa,-\sigma}\Big)
\nonumber\\&
+\sum_{\sigma = \pm 1} \sum_{\rho=1,3} \frac1{\omega^\rho}
\oint \frac{\dd z}{2\pi \iI}\, \Big(
f^1_k \Xi_{k,\frac12,\rho,\sigma} + f^2_k\Xi_{k,\frac32,\rho,\sigma} + f^3_{k-1,\frac12} \Xi_{k-1,\frac12,\rho,\sigma} 
+ f^4_{k-1,\frac32} \Xi_{k-1,\frac32,\rho,\sigma} + f^5_{k-2,\frac32} \Xi_{k-2,\frac32,\rho,\sigma}\Big)
\nonumber\\&
+ \omega \sum_{\sigma = \pm 1} \oint \frac{\dd z}{2\pi \iI}\, \Big(f^2_k \Xi_{k,\frac12,1,\sigma} + f^{3\textrm{b}}_{k-1,m_k+\frac12} \Xi_{k-1,\frac12,1,\sigma}\Big)\,.
 \label{eq:Lambda.constraint.L=1/2}
\end{alignat}
We find that the content of the first parenthesis vanishes for $1 \le \kappa \le k-2$, so the corresponding integrals vanish for reason~(\ref{item:reasoni}). For $\kappa = k-1$ and $\kappa = k$, these integrals instead vanish for reason~(\ref{item:reasonii}). This follows from
\begin{subequations}
\begin{alignat}{2}
0&=f^1_k \Omega_{k,\frac12,\sigma} - q^{2 \sigma m_k} f^2_k \Omega_{k,\frac12,-\sigma}+f^2_k\Omega_{k,\frac32,\sigma}\,,
\\[0.15cm]
0&=f^1_k \Pi_{k,\frac12,\sigma} - q^{2 \sigma m_{k-1}} f^2_k \Pi_{k,\frac12,-\sigma}+ f^2_k \Xi_{k,\frac32,\sigma} + f^3_{k-1,\frac12}\Omega_{k-1,\frac12,\sigma} + f^4_{k-1,\frac32}\Omega_{k-1,\frac32,\sigma} 
\nonumber\\&
- q^{2 \sigma m_{k-1}} f^{3\textrm{b}}_{k-1,m_k+\frac12} \Omega_{k-1,\frac12,-\sigma} \,.
\end{alignat}
\end{subequations}
For the second and third contributions in \eqref{eq:Lambda.constraint.L=1/2}, we find the identities
\begingroup
\allowdisplaybreaks
\begin{subequations}
\begin{alignat}{2}
f^1_k & \Xi_{k,\frac12,1,\sigma} + f^2_k\Xi_{k,\frac32,1,\sigma} + f^3_{k-1,\frac12} \Xi_{k-1,\frac12,1,\sigma} 
+ f^4_{k-1,\frac32} \Xi_{k-1,\frac32,1,\sigma} + f^5_{k-2,\frac32} \Xi_{k-2,\frac32,1,\sigma} \\\nonumber&= \frac{(-1)^{k+1} \sigma\, q^{\sigma(1-n)}}{n(q-q^{-1})^{2k-1}}\frac{[n-k]!}{[n]![k-1]!}\frac{1-z q^{\sigma n}}{z^k} \prod_{j=-\frac{n-3}2}^{\frac{n-3}2} (1-z q^{\sigma + 2j}) \prod_{j=-m_k}^{m_k}\frac1{1-z q^{2j}}\,,
\\
f^1_k & \Xi_{k,\frac12,3,\sigma} + f^2_k\Xi_{k,\frac32,3,\sigma} + f^3_{k-1,\frac12} \Xi_{k-1,\frac12,3,\sigma} 
+ f^4_{k-1,\frac32} \Xi_{k-1,\frac32,3,\sigma} + f^5_{k-2,\frac32} \Xi_{k-2,\frac32,3,\sigma} 
\\\nonumber&= \frac{(-1)^{k} \sigma\, q^{3\sigma}}{n(q-q^{-1})^{2k-1}}\frac{[n-k]!}{[n]![k-1]!}\frac{(1-z q^{-\sigma n-2})(1-z q^{-\sigma n+2})}{z^k(1-z q^{\sigma n})} \prod_{j=-\frac{n-3}2}^{\frac{n-3}2} (1-z q^{3\sigma + 2j}) \prod_{j=-m_k}^{m_k}\frac1{1-z q^{2j}}\,,
\\
f^2_k & \Xi_{k,\frac12,1,\sigma} + f^{3\textrm{b}}_{k-1,m_k+\frac12}\Xi_{k-1,\frac12,1,\sigma} = \frac{(-1)^{k} \sigma\, q^{\sigma}}{n(q-q^{-1})^{2k-1}}\frac{[n-k]!}{[n]![k-1]!}\frac{1-z q^{-\sigma n}}{z^k}
\\\nonumber& \hspace{5cm}\times 
\prod_{j=-\frac{n-3}2}^{\frac{n-3}2} (1-z q^{\sigma + 2j}) \prod_{j=-m_k}^{m_k}\frac1{1-z q^{2j}}\,.
\end{alignat}
\end{subequations}
\endgroup
All the corresponding integrals vanish because of reason~(\ref{item:reasoniii}). \medskip

For the special cases where $\ell = 0$ and $0 \le k \le \lfloor \frac{n-2}2\rfloor$, we need to prove the identity\footnote{The special case with $n$ odd and $(k,\ell) = (\frac{n-1}2,0)$ is treated separately at the end of the section.}
\begin{alignat}{2}
\label{eq:Gamma.constraint.L=0}
0 &= f^1_k \Gamma_{k,0} + f^2_k \Big(\Gamma_{k,1}\big|_{\omega \to \omega^{-1}}+\Gamma_{k,1}\Big)+f^3_{k-1,0} \Gamma_{k-1,0} + f^4_{k-1,1} \Gamma_{k-1,1} + f^5_{k-2,1}\Gamma_{k-2,1} 
\nonumber\\[0.15cm]
&+ f^3_{k-1, m_{k-1}-1} \Big(\Gamma_{k-1,1}\big|_{\omega \to \omega^{-1}}\Big) + f^5_{k-2, m_{k-2}-1} \Big(\Gamma_{k-2,1}\big|_{\omega \to \omega^{-1}}\Big)\, ,
\end{alignat}
where
\be
\Gamma_{k,1}\big|_{\omega \to \omega^{-1}} = -\sum_{\sigma = \pm 1} \sum_\kappa \bigg(\frac1{\omega^2 q^{2 \sigma m_\kappa}-1}+1\bigg) \oint \frac{\dd z}{2\pi \iI}\, \Lambda_{k,1,\kappa,-\sigma} + \omega^2 \sum_{\sigma = \pm 1} \oint \frac{\dd z}{2\pi \iI}\, \Xi_{k,1,2,\sigma}\,.
\ee
As a result, we must show that
\begin{alignat}{3} 
&0 = \sum_{\sigma = \pm 1} \sum_\kappa \frac1{\omega^2 q^{2 \sigma m_\kappa}-1} 
\oint \frac{\dd z}{2\pi \iI}\, \Big(
f^1_k \Lambda_{k,0,\kappa,\sigma} 
- f^2_k \Lambda_{k,1,\kappa,-\sigma} 
+ f^2_k \Lambda_{k,1,\kappa,\sigma} 
+ f^3_{k-1,0} \Lambda_{k-1,0,\kappa,\sigma} 
\nonumber\\&
\hspace{2.4cm}
+ f^4_{k-1,1} \Lambda_{k-1,1,\kappa,\sigma} 
+ f^5_{k-2,1} \Lambda_{k-2,1,\kappa,\sigma} 
- f^3_{k-1,m_{k-1}-1}\Lambda_{k-1,1,\kappa,-\sigma}
- f^5_{k-2,m_{k-2}-1}\Lambda_{k-2,1,\kappa,-\sigma}\Big)
\nonumber\\&
- \sum_{\sigma = \pm 1} \sum_\kappa \oint \frac{\dd z}{2\pi \iI} 
\Big(
f^2_k \Lambda_{k,1,\kappa,\sigma} 
+ f^3_{k-1,m_{\kappa-1}+1} \Lambda_{k-1,1,\kappa,\sigma} 
+ f^5_{k-2,m_{k-2}-1} \Lambda_{k-2,1,\kappa,\sigma} 
\Big)
\nonumber\\&
+ \omega^{-2}\sum_{\sigma = \pm 1} 
\oint \frac{\dd z}{2\pi \iI}\,
\Big(
f^2_k \Xi_{k,1,2,\sigma} + f^4_{k-1,1} \Xi_{k-1,2,\rho,\sigma} + f^5_{k-2,1} \Xi_{k-2,1,2,\sigma}\Big)
\nonumber\\&
+ \omega^2 \sum_{\sigma = \pm 1} \oint \frac{\dd z}{2\pi \iI} 
\Big(f^2_k \Xi_{k,1,2,\sigma} 
+ f^3_{k-1,m_{k-1}-1} \Xi_{k-1,1,2,\sigma}
+f^5_{k-2,m_{k-2}-1} \Xi_{k-2,1,2,\sigma}\Big).
 \label{eq:Lambda.constraint.L=0}
\end{alignat}
One can show that the content of the first parenthesis vanishes identically for $1\le \kappa \le k-2$. For $\kappa = k-1$ and $k$, the corresponding integrals vanish because
\begin{subequations}
\begin{alignat}{2}
0&=f^1_k \Omega_{k,0,\sigma} - f^2_k \Omega_{k,1,-\sigma}+f^2_k\Omega_{k,1,\sigma}\,,
\\[0.15cm]
0&=f^1_k \Pi_{k,0,\sigma} - f^2_k \Pi_{k,1,-\sigma}+ f^2_k \Xi_{k,1,\sigma} + f^3_{k-1,0}\Omega_{k-1,0,\sigma} + f^4_{k-1,1}\Omega_{k-1,1,\sigma} 
-f^3_{k-1,m_{k-1}-1}\Omega_{k-1,1,-\sigma}  \,.
\end{alignat}
\end{subequations}
The other integrals vanish because of the identities
\begingroup
\allowdisplaybreaks
\begin{subequations}
\begin{alignat}{2}
&f^2_k \Lambda_{k,1,\kappa,\sigma} + f^3_{k-1,m_{k-1}-1} \Lambda_{k-1,1,\kappa,\sigma} + f^5_{k-2,m_{k-2}-1} \Lambda_{k-2,1,\kappa, \sigma} 
\\[0.15cm]\nonumber
& \simeq \frac{(-1)^{k+\kappa} \sigma q^{\sigma n}}{n(q-q^{-1})^{2k-1}} \frac{[n-k]!}{[n]![k-1]!}
\left[\begin{matrix} n-2\\ \kappa-1\end{matrix}\right] \frac{(1-z q^{-\sigma n-1})(1-z q^{-\sigma n+1})}{z^{k-\kappa+1}}
\prod_{j=-(m_k+\frac12)}^{m_k+\frac12} \frac1{1-z q^{2j}}\,,
\\[0.3cm]
&f^2_k \Xi_{k,1,2,\sigma} + f^4_{k-1,1} \Lambda_{k-1,1,2,\sigma} + f^5_{k-2,1} \Lambda_{k-2,1,2, \sigma} 
\nonumber\\[0.15cm]& 
= f^2_k \Xi_{k,1,2,\sigma} + f^3_{k-1,m_{k-1}-1} \Lambda_{k-1,1,2,\sigma} + f^5_{k-2,m_{k-2}-1} \Lambda_{k-2,1,2, \sigma} 
\\[0.15cm]\nonumber
&= \frac{(-1)^{k} \sigma q^{2 \sigma}}{n(q-q^{-1})^{2k-1}} \frac{[n-k]!}{[n]![k-1]!}
\frac{(1-z q^{-\sigma n-1})(1-z q^{-\sigma n+1})}{z^k}
\prod_{j=-\frac{n-3}2}^{\frac{n-3}2} (1-z q^{2 \sigma + 2 j})
\prod_{j=-(m_k+\frac12)}^{m_k+\frac12} \frac1{1-z q^{2j}}\,.
\end{alignat}
\end{subequations}
\endgroup
One can then show that the second, third and fourth contributions in \eqref{eq:Lambda.constraint.L=0} vanish due to reason~(\ref{item:reasoniii}). We note that, to show this for the second contribution, one must evaluate the sum over $\kappa$ using \eqref{eq:qbin.rel.a}.

\paragraph{Extremal cases.} 
The only remaining relations to prove in the recursive system are the special cases with $(k,\ell)$ specialised to $(\frac{n-2}2,\frac12)$ and $(\frac{n-1}2,0)$ in \eqref{eq:Gamma.constraint}, as well as \eqref{eq:extra.Gamma.constraint} for the algebra $\qatla_n(\beta, \alpha)$. Let us start with the latter case. As a first step, we keep $\widehat \gamma$ generic and evaluate $\Gamma_{(n-2)/2,0}$, $\Gamma_{(n-2)/2,1/2}$ and $\Gamma_{(n-4)/2,1}$. In their corresponding integral expressions \eqref{eq:Gamma.Lambda.1}, there are either two or four simple poles away from the origin. The contour integrals can then be computed explicitly using the same arguments as those used in \cref{sec:kernel}, namely one checks that there is no pole at infinity and computes the integral from the sum of the residues of the other poles. After simplifications, the results read
\begin{subequations}
\label{eq:Gamma.v1}
\begin{alignat}{2}
\label{eq:Gamma.n-2/2,0}
\Gamma_{\frac{n-2}2,0} &= \frac{(-1)^{n/2+1}}{n (q-q^{-1})^{n-3}} \sum_{\sigma = \pm 1} \sum_{\kappa=1}^{(n-2)/2}\frac{(-1)^\kappa \sigma}{\omega^2 q^{2\sigma m_\kappa}-1} \frac{[\frac n2][\frac n2-\kappa]}{[\kappa-1]![n-\kappa-1]!}\,,
\\[0.15cm]\label{eq:Gamma.n-2/2,1/2}
\Gamma_{\frac{n-2}2, \frac12} &= \frac{(-1)^{n/2+1}\omega^{-1}}{n (q-q^{-1})^{n-2}} \sum_{\sigma = \pm 1} \sum_{\kappa=1}^{(n-2)/2}\frac{(-1)^\kappa q^{-\sigma m_\kappa}}{\omega^2 q^{2\sigma m_\kappa}-1} \frac{q^{\sigma n m_\kappa}-1}{[\kappa-1]![n-\kappa-1]!}\,,
\\[0.15cm]\label{eq:Gamma.n-4/2,1}
\Gamma_{\frac{n-4}2,1} &=  \frac{(-1)^{n/2}\omega^{-2}}{n (q-q^{-1})^{n-3}}\sum_{\sigma = \pm 1} \sum_{\kappa=1}^{(n-4)/2} \frac{\sigma (-1)^\kappa q^{\sigma(n/2-2)m_\kappa}}{\omega^2 q^{2\sigma m_\kappa}-1}  \frac{[\tfrac n2-1][(\tfrac n2+1)m_\kappa]-[\tfrac n2+1][(\tfrac n2-1)m_\kappa]}{[\kappa-1]![n-\kappa-1]!}\,. 
\end{alignat}
\end{subequations}
These can equivalently be written as
\begin{subequations}
\label{eq:Gamma.v2}
\begin{alignat}{2}
\label{eq:Gamma.n-2/2,0v2}
\Gamma_{\frac{n-2}2,0} &= \frac{(-1)^{n/2+1}}{n}\, [\tfrac n2]\, \omega^{n-2} \prod_{\kappa=1}^{\frac{n-2}2} \prod_{\sigma = \pm 1} \frac1{\omega^2 q^{2 \sigma m_\kappa}-1}\,,
\\[0.15cm]\label{eq:Gamma.n-2/2,1/2v2}
\Gamma_{\frac{n-2}2,\frac 12}
& = \frac{(-1)^{n/2+1}\omega^{-1}}{n}\,\frac{\omega^n -1}{\omega^2 - 1} \prod_{\kappa=1}^{\frac{n-2}2} \prod_{\sigma = \pm 1} \frac1{\omega^2 q^{2 \sigma m_\kappa}-1}\,,
\\[0.15cm]\label{eq:Gamma.n-4/2,1v2}
\Gamma_{\frac{n-4}2,1} &=  \frac{(-1)^{n/2+1}\omega^{-2}}{n} \frac{[\frac n2+1](\omega^{n} -\omega^2)-[\frac n2-1](\omega^{n+2} -1)}{\omega^2 - 1} \prod_{\kappa=1}^{\frac{n-2}2} \prod_{\sigma = \pm 1} \frac1{\omega^2 q^{2 \sigma m_\kappa}-1}\,.
\end{alignat}
\end{subequations}
To show this, we note that \eqref{eq:Gamma.v1} and \eqref{eq:Gamma.v2} give us expressions  for $\Gamma_{(n-2)/2,0}$, $\omega\,\Gamma_{(n-2)/2,1/2}$ and $\omega^2\,\Gamma_{(n-4)/2,1}$ that we can re-express as rational functions of the form $\mathcal N(\omega^2)/\mathcal D(\omega^2)$ in the variable~$\omega^2$, with the same denominator in the three cases:
\be
\mathcal D(\omega^2)=\prod_{\kappa=1}^{\frac{n-2}2}\prod_{\sigma = \pm 1} (\omega^2 q^{2 \sigma m_\kappa}-1)\,.
\ee 
Then $\mathcal N(\omega^2)$ is a polynomial in $\omega^2$, of degree strictly less than $n-2$, that is different for the three cases. To show that the two expressions for a given constant are equal, we show that the numerators are equal for $n-2$ values of $\omega^2$, namely those where $\mathcal D(\omega^2)$ vanishes: $\omega^2 = q^{-2\sigma m_\kappa}$ for $\kappa \in \{1,2,\dots, \frac{n-2}2\}$ and $\sigma \in \{+1,-1\}$. Equivalently, we must show that $\lim_{\omega^2 \to q^{-2\sigma m_\kappa}} (\omega^2q^{2\sigma m_\kappa}-1) \Gamma_{(n-2)/2,0}$ is identical for the two expressions (and likewise for the two other constants). This verification is straightforward.\medskip

Using the above formulas, we find
\be
[2] \Gamma_{\frac{n-2}2,0} + \Gamma_{\frac{n-4}2,1} = \frac{(-1)^{n/2+1}}{n} \frac{\omega^n-1}{\omega^2 - 1} \big( [\tfrac{n+2}2] - \omega^{-2} [\tfrac{n-2}2]\big) \prod_{\kappa=1}^{\frac{n-2}2} \prod_{\sigma = \pm 1} \frac1{\omega^2 q^{2 \sigma m_\kappa}-1}\,.
\ee
Setting $\omega = \eE^{2 \pi \iI r/n}$ with $r\in \{0,1, \dots, n\}$, we obtain
\begin{subequations}
\begin{alignat}{2}
&[2] \Gamma_{\frac{n-2}2,0} + \Gamma_{\frac{n-4}2,1} = \frac{[n]}{2(q-q^{-1})^{n-2}[\frac n2]![\frac{n-2}2]!} 
\times \left\{\begin{array}{cl} 
1 & r \in \{0,\frac n2\},\\[0.15cm]
0 & \textrm{otherwise},
\end{array}\right.
\\[0.2cm]
&\Gamma_{\frac{n-2}2,\frac12} = \frac{1}{2(q-q^{-1})^{n-2}[\frac{n-2}2]!^2} 
\times \left\{\begin{array}{cl} 
1 & r = 0,\\[0.15cm]
-1 & r = \frac n2,\\[0.15cm]
0 & \textrm{otherwise}.
\end{array}\right.
\end{alignat}
\end{subequations}
With these expressions, it is straightforward to check that \eqref{eq:extra.Gamma.constraint} holds with the formula \eqref{eq:Gamma.n/2.0.affine} for $\Gamma_{n/2,0}$, ending the proof in this case.\medskip

The two remaining cases are tackled using similar arguments. For $(k,\ell) = (\frac{n-2}2,\frac12)$ with $n$ even, we note that the two last terms in \eqref{eq:Gamma.constraint} are  non-zero. With the same arguments as above, we find
\begin{subequations}
\begin{alignat}{2}
\Gamma_{\frac{n-4}2,\frac12} &=  \frac{(-1)^{n/2+1}\omega^{-1}}{n} \frac{1-
\frac{[\frac n2][\frac{n+2}2]}{[2]}\,\omega^{n-2}+[\frac{n-2}2][\frac{n+2}2]\,\omega^n- \frac{[\frac{n-2}2][\frac n2]}{[2]}\,\omega^{n+2}}{\omega^2 - 1} \prod_{\kappa=1}^{\frac{n-2}2} \prod_{\sigma = \pm 1} \frac1{\omega^2 q^{2 \sigma m_\kappa}-1}\,,
\\[0.15cm]
\Gamma_{\frac{n-4}2,\frac32} &=  \frac{(-1)^{n/2+1}\omega^{-3}}{n} \frac{\frac{[\frac{n-2}2][\frac n2]}{[2]}- [\frac{n-2}2][\frac{n+2}2]\,\omega^2 + \frac{[\frac n2][\frac{n+2}2]}{[2]}\,\omega^{4}- \omega^{n+2}}{\omega^2 - 1} \prod_{\kappa=1}^{\frac{n-2}2} \prod_{\sigma = \pm 1} \frac1{\omega^2 q^{2 \sigma m_\kappa}-1}\,,
\\[0.15cm]
\Gamma_{\frac{n-6}2,\frac32} &=  \frac{(-1)^{n/2+1}\omega^{-3}}{n} \frac{\frac{[\frac{n-4}2][\frac{n-2}2]}{[2]}(\omega^{n+4}-1)-[\frac{n-4}2][\frac{n+4}2](\omega^{n+2}-\omega^2)+\frac{[\frac{n+2}2][\frac{n+4}2]}{[2]}(\omega^n-\omega^4)}{\omega^2 - 1} 
\\\nonumber&\times\prod_{\kappa=1}^{\frac{n-2}2} \prod_{\sigma = \pm 1} \frac1{\omega^2 q^{2 \sigma m_\kappa}-1}\,,
\end{alignat}
\end{subequations}
satisfying
\begin{alignat}{2}
0 &= \bigg(\gamma + \gamma^{-1}-\frac{[2n]}{[n]}\bigg)\Gamma_{\frac{n-2}2,\frac12} 
+ \bigg(\gamma^{-1}+\frac{[\frac{n-2}2]}{[\frac{n+2}2]}+\frac{[\frac{n}2][\frac{n+4}2]}{[\frac{n-2}2][\frac{n+2}2]} \bigg) \Gamma_{\frac{n-4}2,\frac12} 
\nonumber\\&
+ \bigg(\gamma + \frac{[\frac{n+2}2]}{[\frac{n-2}2]}+\frac{[\frac{n-4}2][\frac{n}2]}{[\frac{n-2}2][\frac{n+2}2]} \bigg) \Gamma_{\frac{n-4}2,\frac32}
+ \frac{[\frac n2]^2}{[\frac{n-2}2][\frac{n+2}2]} \Gamma_{\frac{n-6}2,\frac32}\,,
\end{alignat}
which is equivalent to \eqref{eq:Gamma.constraint} for $(k,\ell) = (\frac{n-2}2,\frac12)$.\medskip

Finally, for $(k,\ell) = (\frac{n-1}2,0)$ with $n$ odd, we find the expressions
\begingroup
\allowdisplaybreaks
\begin{subequations}
\begin{alignat}{2} 
\Gamma_{\frac{n-1}2,0} &= \frac{(-1)^{(n-1)/2}}{n}\, \omega^{n-1} \prod_{\kappa=1}^{\frac{n-1}2} \prod_{\sigma = \pm 1} \frac1{\omega^2 q^{2 \sigma m_\kappa}-1}\,,
\\[0.15cm]
\Gamma_{\frac{n-3}2,0} &= \frac{(-1)^{(n+1)/2}}{n}\, \frac{[\frac{n-1}2][\frac{n+1}2]}{[2]}\,\omega^{n-3} \prod_{\kappa=1}^{\frac{n-3}2} \prod_{\sigma = \pm 1} \frac1{\omega^2 q^{2 \sigma m_\kappa}-1}\,,
\\[0.15cm]
\Gamma_{\frac{n-3}2,\frac12} &= \frac{(-1)^{(n+1)/2} \omega^{-1}}{n}\, \big(1-[\tfrac{n+1}2]\,\omega^{n-1}+[\tfrac{n-1}2]\,\omega^{n+1}\big) \prod_{\kappa=1}^{\frac{n-1}2} \prod_{\sigma = \pm 1} \frac1{\omega^2 q^{2 \sigma m_\kappa}-1}\,,
\\[0.15cm]
\Gamma_{\frac{n-3}2,1} &= \frac{(-1)^{(n+1)/2} \omega^{-2}}{n}\, \big([\tfrac{n-1}2]-[\tfrac{n+1}2]\,\omega^{2}+\omega^{n+1}\big) \prod_{\kappa=1}^{\frac{n-1}2} \prod_{\sigma = \pm 1} \frac1{\omega^2 q^{2 \sigma m_\kappa}-1}\,,
\\[0.15cm]
\Gamma_{\frac{n-5}2,1} &= \frac{(-1)^{(n-1)/2} \omega^{-2}}{n}\, \bigg([\tfrac{n-3}2]-[\tfrac{n+3}2]\,\omega^{2}+\frac{[\tfrac{n+1}2][\tfrac{n+3}2]}{[2]}\,\omega^{n-1}-[\tfrac{n-3}2][\tfrac{n+3}2]\,\omega^{n+1}+\frac{[\tfrac{n-3}2][\tfrac{n-1}2]}{[2]}\,\omega^{n+3}\bigg) 
\nonumber\\&\times\prod_{\kappa=1}^{\frac{n-1}2} \prod_{\sigma = \pm 1} \frac1{\omega^2 q^{2 \sigma m_\kappa}-1}\,,
\\[0.15cm]
\Gamma_{\frac{n-5}2,\frac32} &= \frac{(-1)^{(n-1)/2} \omega^{-3}}{n}\, \bigg(\frac{[\tfrac{n-3}2][\tfrac{n-1}2]}{[2]}-[\tfrac{n-3}2][\tfrac{n+3}2]\,\omega^{2}+\frac{[\tfrac{n+1}2][\tfrac{n+3}2]}{[2]}\,\omega^{4}-[\tfrac{n+3}2]\,\omega^{n+1}+[\tfrac{n-3}2]\,\omega^{n+3}\bigg) 
\nonumber\\&\times\prod_{\kappa=1}^{\frac{n-1}2} \prod_{\sigma = \pm 1} \frac1{\omega^2 q^{2 \sigma m_\kappa}-1}\,,
\end{alignat}
\end{subequations}
\endgroup
satisfying
\begin{alignat}{2}
0 &= \bigg(\gamma^2 + \gamma^{-2}-\frac{[2n]}{[n]}\bigg)\Gamma_{\frac{n-1}2,0} 
+ \bigg( \frac{[\frac{n+3}2]}{[\frac{n-1}2]}+\frac{[\frac{n-3}2]}{[\frac{n+1}2]} \bigg) \Gamma_{\frac{n-3}2,0} 
+ \big(\gamma [2] + \gamma^{-1} \big) \Gamma_{\frac{n-3}2,\frac12}
\nonumber\\&
+ \big(\gamma^2 + [2]) \Gamma_{\frac{n-3}2,1}
+ \Gamma_{\frac{n-5}2,1}
+ \gamma\, \Gamma_{\frac{n-5}2,\frac32}\,,
\end{alignat}
which is equivalent to \eqref{eq:Gamma.constraint} for $(k,\ell) = (\frac{n-1}2,0)$. This ends the proof.

\addcontentsline{toc}{section}{References}

\end{document}